# Towards a taxonomy of atlases and of morphisms between them

Shmuel (Seymour J.) Metz

June 25, 2019

**Abstract**

Manifolds and fiber bundles, while superficially different, have strong parallels; in particular, they are both defined in terms of equivalence classes of atlases or in terms of maximal atlases, with the atlases treated as mere adjuncts. This paper presents a unified view of atlases for manifolds and fiber bundles as mathematical entities in their own right. It defines some convenient notation, defines categories of, e.g., atlases and defines functors among them.

[Metz, 2018] introduced the ideas presented here, but many of the details are not needed there. This paper fleshes out the concepts in more detail than would be relevant there.

# Part I
# Introduction

Historically, manifolds were formalized in terms of atlases based on Euclidean spaces and fiber bundles were formalized in terms of atlases based on product spaces, using either equivalence classes of atlases or maximal atlases. The concept of pseudogroups allowed unifying manifolds and manifolds with boundary. The definitions of fiber bundles and manifolds have strong parallels, and can be unified in a similar fashion. The standard presentations treat these atlases as secondary adjuncts.

[Metz, 2018] introduced the concept of an m-atlas and of morphisms among them only as preliminaries to unifying manfiolds and fiber bundles as special cases of local coordiinate spaces. This paper treats m-atlases as objects of interest in their own right, and includes concepts not relevant to [Metz, 2018].

This paper defines an approach using categories and commutative diagrams that is designed for easy exposition at the possible expense of abstractness and generality. In particular, I have chosen to simplify proofs by stating stricter auumptions than necessary and have assumed the Axiom of Choice (AOC).

---

I wish to gratefully thank Walter Hoffman (z"l), Milton Parnes, Dr. Stanley H. Levy, the Mathematics department of Wayne State University, the Mathematics department of the State University of New York at Buffalo and others who guided my education.

Although this paper incidentally defines partial equivalents to manifolds and fiber bundles using model spaces and model atlases, it does not explicitly reflect the role of the group in fiber bundles.

This paper defines functors among categories of atlases, categories of model spaces, categories of manifolds and categories of fiber bundles; it constructs more machinery than is customary in order to facilitate the presentation of those categories and functors.

*Remark* 0.1. The definitions of morphisms given here are nonstandard, but there are also definitions of classical morphisms that are much closer to the standard ones.

Parts II to V present nomenclature and give basic results.

Part VI defines m-atlases, m-atlas morphisms and categories of them; lemma 11.23 proves that the defined categories are indeed categories.

part VII (Equivalence of manifolds) on page 111 and part VIII (Equivalence of fiber bundles) on page 134 show the equivalence of manifolds and fiber bundles with m-atlases by explicitly exhibiting functors to and from the m-atlases.

*Remark* 0.2. The unconventional definitions of manifold and fiber bundle are intended to make their relationship to local coordinate spaces more natural. See [Metz, 2018] for details.

Most of the lemmata, theorems and corollaries in this paper should be substantially identical to results that are familiar to the reader. What is novel is the perspective and the material directly related to local coordinate spaces. The presentation assumes only a basic knowledge of Category Theory, such as may be found in the first chapter of [Mac Lane, 1998] or [Adámek, Herrlich, Strecker, 1990].

Many of the proofs are obvious and may be omitted from other versions of this paper.

# 1 New concepts and notation

This paper introduces a significant number of new concepts and some modifications of the definitions for some conventional concepts. It also introduces some notation of lesser importance. The following are the most important.

1. Nearly commutative diagram (NCD), NCD at a point, locally NCD and special cases with related nomenclature. These are cases where a diagram can be modified to make it commutative.

2. Model space and related concepts. A model space[1] is a topological space with a category of permissible open sets and mappings satisfying specified conditions. It is similar to a pseudo-group, with some important differences.

3. Model topology and M-paracompactness

4. An m-atlas is a generalization of atlases in manifolds.

---

[1]The phrase has been used before, but with a different meaning.

5. An m-atlas near morphism (m-atlas morphism) is a generalization of a $C^k$ map between differentiable manifolds. It differs in using both a function between the manifolds and a function between the coordinate spaces.

6. Linear space and related concepts A linear space is a subspace of a Banach space or of a Fréchet space that is suitable for defining coordinate patches in a manifold.

7. Trivial $C^k$ linear model space and related concepts. A trivial $C^k$ linear model space is a maximal model space with morphisms $C^k$ maps between open subsets of a linear space.

8. $G$-$\rho$ bundle atlases[2] and related concepts

# Part II
# Conventions

An arrow with an Equal-Tilde ($A \xrightarrow[\phi]{\tilde{=}} B$) represents an isomorphism. One with a hook ($A \xhookrightarrow{i} B$) represents an inclusion map. One with a tail ($A \xrightarrowtail{i} B$) represents a monomorphism. One with a double head ($A \xtwoheadrightarrow{\pi} B$) represents a surjection.

When a superscript ends in $-1$, e.g., $f^{i-1}$, it is to be taken as function inverse rather than subtraction of 1.

All diagrams shown are commutative; none are exact.

A definition of a base term and several related terms may specify modifiers in parenthese in the base definition and then give the definitions for each modifier; usually the modifiers will be restrictions and the base definition will not be repeated. E.g., "$f$ is a (semi-strict, strict) prestructure morphism of $P^1$ to $P^2$ iff" is followed by the definition of prestructure morphism, then by the restrictions for strict and semi-strict prestructure morphisms.

When a definition defines a base propositional function and variant propositional functions with a qualifier given as a superscript or underset, then the form base$^{(qualifiers)}$ will refer to either base$^{qualifier}$ or base and the form $\underset{(qualifiers)}{base}$ will refer to either $\underset{qualifier}{base}$ or base, e.g., isAtl$_{\underset{semi-strict}{Ar}}^{(classic,near)}$ refers to either isAtl$_{\underset{semi-strict}{Ar}}^{classic}$, isAtl$_{\underset{semi-strict}{Ar}}^{near}$ or isAtl$_{\underset{semi-strict}{Ar}}$.

The form base$^{qualifier(,qualifiers)}$ will refer to either base$^{qualifier}$ or to the variant with one or more of the additional qualifiers in the superscript, e.g., isAtl$_{\underset{semi-strict}{Ar}}^{classic(,near)}$ refers to either isAtl$_{\underset{semi-strict}{Ar}}^{classic}$ or isAtl$_{\underset{semi-strict}{Ar}}^{classic,near}$.

---

[2]Similar to coordinate bundles

Alternatively, a definition may specify a numbered list of alternatives, and subsequently specify additional numbered lists with items corresponding to those in the first list, e.g., $\boldsymbol{f}$ is also a

1. full
2. semi-maximal
3. maximal
4. full semi-maximal
5. full maximal

$E^1$-$E^2$ $C^k$ near morphism of $\boldsymbol{A}^1$ to $\boldsymbol{A}^2$ in the coordinate spaces $C^1$, $C^2$, abbreviated as

1. abbreviation for full
2. abbreviation for semi-maximal
3. abbreviation for maximal
4. abbreviation for full semi-maximal
5. abbreviation for full maximal

iff

1. definition for full
2. definition for semi-maximal
3. definition for maximal
4. definition for full semi-maximal
5. definition for full maximal

A Corollary, Lemma or Theorem that applies to related terms defined with the above convention will specify the modifiers in parentheses to indicate that it applies to the base term and to the varianr terms, e.g., "a (semi-strict, strict) $\boldsymbol{E}^i$-$\boldsymbol{E}^{i+1}$ m-atlas morphism" If it applies only to more restrictive terms them it will specify the first relevant modifier followed by the remaining relevant modifiers in parentheses, e.g., "If $\mathcal{S}^i \stackrel{\textbf{full--cat}}{\subseteq} \mathcal{S}'^i$ and $\boldsymbol{f}^i$ is a semi-strict (strict) prestructure morphism".

Blackboard bold upper case will denote specific sets, e.g., the Naturals.

Bold lower case italic letters will refer to sets, sequences and tuples of functions, e.g., $\boldsymbol{f} \stackrel{\text{def}}{=} (f_1, f_2)$.

Bold lower case Latin letters will refer to sequence valued functions of sequences and tuple valued functions of tuples, e.g., **range** yields the sequence of ranges of a sequence of functions.

Bold upper case italic letters will refer to sequences or tuples, e.g., $\boldsymbol{A} = (x, y, z)$, to sets of them, to sets of topological spaces or to sets of open sets.

Bold upper case script letters will refer to sequences of categories, e.g., $\boldsymbol{\mathcal{A}} \stackrel{\text{def}}{=} (\mathcal{A}_\alpha, \alpha \in A)$.

Fraktur will refer to topologies and to topology-valued functions, e.g., $\mathfrak{Top}$.

Functions have a range, domain and relation, not just a relation. Unless otherwise stated, they are assumed or proven to be continuous.

Groups are assumed to be topological groups. The ambiguous notation $x^{-1}$ will be used when it is obvious from context what the group operation $\star$ and the group identity $\mathbf{1}_G$ are.

Lower case Greek letters other than $\pi, \rho, \sigma, \phi$ and $\psi$ will refer to ordinals, possibly transfinite, and to formal labels. A letter with a Greek superscript and a letter with a Latin or numeric superscript always refer to distinct variables.

Lower case $\pi$ will refer to a projection operator

Lower case $\rho$ will refer to a continuous effective group action, i.e., a continuous representation of a group in a homeomorphism group.

Lower case $\sigma$ will refer to a sequence of ordinals, referred to as a signature.

Lower case $\phi$ will refer to a coordinate function.

Lower case italic and Latin letters will refer to

1. elements of a set or sequence

2. formal labels
   When the letters $u$ and $v$, with or without a superscript, refer to an element of a set then they may also be used as formal labels for sets associated with that element.

3. functions

4. morphisms of a category

5. natural numbers

6. objects of a category

7. ordinal numbers

Upper case Greek letters other than $\Sigma$ may refer to

1. ordinal used as the limit of a sequence of consecutive ordinals, e.g., $x_\alpha, \alpha \leq A$

2. ordinal used as the order type of a sequence of consecutive ordinals, e.g., $x_\alpha, \alpha < A$

Upper case $\Sigma$ will refer to a sequence of signatures

Upper case Latin letters will refer to

1. Natural numbers

2. Topological spaces

3. Open sets

4. Formal labels for elements of a sequence or tuple of functions, e.g., $f_E$ might be $f_0\colon E_1 \longrightarrow E_2$.

Upper case Script Latin letters will refer to categories and functors.
Upright Latin letters will be used for long names. In particular,

1. Ar, Hom, Mor and Ob are the conventional set-valued functions on categories:

    Ar($\mathscr{C}$)  All morphisms of $\mathscr{C}$

    Hom$_\mathscr{C}(a,b)$  All morphisms of $\mathscr{C}$ from $a$ to $b$

    Mor($\mathscr{C}$)  All morphisms of $\mathscr{C}$ Synonymous with Ar($\mathscr{C}$).

    Ob($\mathscr{C}$)  All objects of $\mathscr{C}$

    They will never be applied to large categories in this paper.

2. **Set** and **Top** are the conventional categories of sets and topological spaces.

3. The letters E, G, X and Y will be used for formal labels.

The term $C^k$ includes $C^\infty$ (smooth) and $C^\omega$ (analytic).
This paper uses the term morphism in preference to arrow, but uses the conventional Ar.
This paper uses the terms empty, null and void interchangeably.
The term sequence without an explicit reference to $\mathbb{N}$ will refer to a general ordinal sequence, possibly transfinite.
Sequence numbering, unlike tuple numbering, starts at 0 and the exposition assumes a von Neumann definition of ordinals, so that $\alpha \in \beta \equiv \alpha \prec \beta$.
Except where explicitly stated otherwise, all categories mentioned are small categories with underlying sets, but the morphisms will often not be set functions between the objects and there will not always be a forgetful functor to **Set** or **Top**. By abuse of language no distinction will be made between a category $\mathscr{A}$ of topological spaces and the concrete category $(\mathscr{A}, \mathscr{U})$ over **Top**. Similarly, no distinction will be made among the object $U \in \mathrm{Ob}(\mathscr{A})$, the topological space $\mathscr{U}(U)$ and the underlying set.
The notation $G^V$ will refer only to the set of continuous functions from $V$ to $G$, never to the set of all functions from $V$ to $G$.
When defining a category, the Ordered pair $(\boldsymbol{O}, \boldsymbol{M})$ refers to the smallest concrete category over **Set** or **Top** whose objects are the elements in $\boldsymbol{O}$, whose morphisms include all of the elements of $\boldsymbol{M}$ and whose morphisms from $o^1 \in \boldsymbol{O}$ to $o^2 \in \boldsymbol{O}$ are functions $f\colon o^1 \longrightarrow o^2$ and whose composition is function composition.
When defining a category, the Ordered triple $(\boldsymbol{O}, \boldsymbol{M}, C)$ refers to the small category whose objects are in $\boldsymbol{O}$, whose morphisms are in $\boldsymbol{M}$, whose Hom is

$$\mathrm{Hom}_{(\boldsymbol{O},\boldsymbol{M},C)}(o_1 \in \boldsymbol{O}, o_2 \in \boldsymbol{O}) \stackrel{\mathrm{def}}{=} \{(\boldsymbol{f}, o_1, o_2) \in \boldsymbol{M}\} \qquad (1.1)$$

and whose composition is C.

By abuse of language I may write "$\mathcal{S}$" for $\mathrm{Ob}(\mathcal{S})$, "$A \in \mathcal{A}$" for $A \in \mathrm{Ob}(\mathcal{A})$, "$A \subset \mathcal{A}$" for $A \subset \mathrm{Ob}(\mathcal{A})$, "$A \in \mathcal{A} \subset B \in \mathcal{B}$" for "the underlying set of $A$ is contained in the underlying set of $B$ and the inclusion $i\colon x \in A \hookrightarrow x \in B$ is a morphism" and "$f\colon A \longrightarrow B$" for $f \in \mathrm{Hom}_{\mathcal{C}}(A, B)$, where $\mathcal{C}$ is understood by context.

By abuse of language I shall use the same nomenclature for sequences and tuples, and will use parentheses around a single expression both for grouping and for a tuple with a single element; the intent should be clear from context.

By abuse of language I will omit parenthese around the operands of Functors when they can be assumed by context.

By abuse of language I shall use the $\times$ and $\bigtimes$ symbols for both Cartesian products of sets and Cartesian products of functions on those sets.

By abuse of language, and assuming AOC, I shall refer to some sets as ordinal sequences, e.g., "$(C_\alpha, \alpha \in A)$" for "$\{C_\alpha \mid \alpha \in A\}$", in cases where the order is irrelevant.

By abuse of language, I may omit universal quantifiers in cases where the intent is clear.

In some cases I define a notion similar to a conventional notion and also need to refer to the conventional notion. In those cases I prefix a letter or phrase to the term, e.g., m-paracompact versus paracompact.

# Part III
# General notions

This part of the paper describes nomenclature used throughout the paper. In some cases this reflects new nomenclature or notions, in others it simply makes a choice from among the various conventions in the literature.

**Definition 1.1** (Operations on categories). If $\mathcal{C}$ is a category then

1. $x \stackrel{\mathrm{Ob}}{\in} \mathcal{C} \equiv x \in \mathrm{Ob}(\mathcal{C})$

2. $x \stackrel{\mathrm{Ar}}{\in} \mathcal{C} \equiv x \in \mathrm{Ar}(\mathcal{C})$

3. If $\mathcal{C}$ is a concrete category over **Set** and $o \stackrel{\mathrm{Ob}}{\in} \mathcal{C}$, then $\mathbf{Set}_{\mathcal{C}}(o)$ is the underlying sete of $o$. By abuse of language we will write $\mathbf{Set}(o)$ when the category is understood from context.

4. If $\mathcal{C}$ is a concrete category over **Top** and $o \stackrel{\mathrm{Ob}}{\in} \mathcal{C}$, then $\mathrm{Top}_{\mathcal{C}}(o)$ is the underlying topological space of $o$. By abuse of language we will write $\mathrm{Top}(o)$ when the category is understood from context.

If $\mathcal{S}$ and $\mathcal{T}$ are categories then $\mathcal{S} \overset{\text{cat}}{\subseteq} \mathcal{T}$ iff S is a subcategory of $\mathcal{T}$ and $\mathcal{S} \overset{\text{full-cat}}{\subseteq} \mathcal{T}$ iff S is a full subcategory of $\mathcal{T}$.

If $\mathcal{S}$ and $\mathcal{T}$ are categories then the category union of $\mathcal{S}$ and $\mathcal{T}$, abbreviated $\mathcal{S} \overset{\text{cat}}{\cup} \mathcal{T}$, is the category whose objects are in $\mathcal{S}$ or in $\mathcal{T}$ and whose morphisms are in $\mathcal{S}$ or in $\mathcal{T}$.

**Definition 1.2** (Identity). $\text{Id}_S$ is the identity function on the space $S$,
$\text{Id}_o$ is the identity morphism for the object $o$[3],
$\text{Id}_{U,V}$, for $U \subseteq V$, is the inclusion map[4]. $\text{Id}_V^U$ is synonymous with $\text{Id}_{U,V}$.
$\text{Id}_{\mathscr{C}}$ is the identity functor on the category $\mathscr{C}$.
$\text{ID}_{\mathbf{S}^i}$, $i = 1, 2$, is the sequence of identity functions for the elements of the sequence $\mathbf{S}^i \overset{\text{def}}{=} (S^1_\alpha, \alpha \prec A)$. Let $\mathbf{S}^1 \overset{()}{\subseteq} \mathbf{S}^2$. Then $\text{ID}_{\mathbf{S}^1, \mathbf{S}^2}$ is the sequence of inclusion maps $(\text{Id}_{S^1_\alpha, S^2_\alpha})$, $\alpha \prec A$ for the elements of the sequences $\mathbf{S}^i$.

The subscript may be omitted when it is clear from context.

**Definition 1.3** (Images). $f[U] \overset{\text{def}}{=} \{f(x) | x \in U\}$ is the image of $U$ under $f$ and $f^{-1}[V] \overset{\text{def}}{=} \{x | f(x) \in V\}$ is the inverse image of $V$ under $f$.

*Remark* 1.4. This notation, adopted from [Kelley, 1955], avoids the ambiguity in the traditional $f(U)$ and $f^{-1}(V)$.

**Definition 1.5** (Projections). $\pi_\alpha$ is the projection function that maps a sequence into element $\alpha$ of the sequence. $\pi_i$ is also the projection function that maps a tuple into element $i$ of the tuple.

**Definition 1.6** (Topological category). A topological category is a small subcategory of **Top** or its concrete category over **Set**.

$\mathcal{T}$ is a full topological category iff it is a topological category and whenever $U^i, V^i \overset{\text{Ob}}{\in} \mathcal{T}$, $i = 1, 2$, $V^i \subseteq U^i$, $f \colon U^1 \longrightarrow U^2 \overset{\text{Ar}}{\in} \mathcal{T}$ and $f[V^1] \subseteq V^2$ then $f \colon V^1 \longrightarrow V^2 \overset{\text{Ar}}{\in} \mathcal{T}$.

**Lemma 1.7** (Inclusions in topological categories are morphisms). *Let $\mathcal{T}$ be a full topological category, $S^i \overset{\text{Ob}}{\in} \mathcal{T}$, $i = 1, 2$, and $S^1 \subseteq S^2$. Then $\text{Id}_{S^2}^{S^1}$ is a morphism of $\mathcal{T}$*

*Proof.* $\text{Id}_{S^2} \overset{\text{Ar}}{\in} \mathcal{T}$, $S^1 \subseteq S^2$ by hypothesis and $S^1 \subseteq S^2$, so $\text{Id}_{S^2}^{S^1} \overset{\text{Ar}}{\in} \mathcal{T}$ by definition 1.6. □

**Definition 1.8** (Local morphisms). Let $\mathcal{S}^i$, $i = 1, 2$, be a full topological category and $S^i \overset{\text{Ob}}{\in} \mathcal{S}^i$, A continuous function $f \colon S^1 \longrightarrow S^2$ is a local $\mathcal{S}^1$-$\mathcal{S}^2$ morphism of $S^1$ to $S^2$ iff $\mathcal{S}^1 \overset{\text{full-cat}}{\subseteq} \mathcal{S}^2$ and for every $u \in S^1$ there is an open neighborhood $U_u$ for $u$ and an open neighborhood $V_u$ for $v \overset{\text{def}}{=} f(u)$ such that $f[U_u] \subseteq V_u$ and $f \colon U_u \longrightarrow V_u$ is a morphism of $\mathcal{S}^2$.

---

[3]The object is often expressed as a tuple, e.g., $\text{Id}_{(A,B)}$ is the identity morphism for the object $(A, B)$
[4]$U$ and $V$ need not have the same topology.

*Remark* 1.9. The condition $f\colon U_u \longrightarrow V_u \overset{\text{Ar}}{\in} \mathcal{S}^2$ ensures that $U_u \overset{\text{Ob}}{\in} \mathcal{S}^1$ and $V_u \overset{\text{Ob}}{\in} \mathcal{S}^2$

Let $\mathcal{T}$ be a full topological category and $S^i \overset{\text{Ob}}{\in} \mathcal{T}, i = 1, 2$. A continuous function $f\colon S^1 \longrightarrow S^2$ is a local $\mathcal{T}$ morphism of $S^1$ to $S^2$ iff it is a local $\mathcal{T}$-$\mathcal{T}$ morphism of $S^1$ to $S^2$.

**Lemma 1.10** (Local morphisms). *Let $\mathcal{T}^i$, $i \in [1,3]$, be a full topological category, $S^i \overset{\text{Ob}}{\in} \mathcal{T}^i$ and $\mathcal{T}^i \overset{\text{full-cat}}{\subseteq} \mathcal{T}^{i+1}$.*

*If $f^i\colon S^i \longrightarrow S^{i+1} \overset{\text{Ar}}{\in} \mathcal{T}^{i+1}$ then $f^i$ is a local $\mathcal{T}^i$-$\mathcal{T}^{i+1}$ morphism of $S^i$ to $S^{i+1}$.*

*Proof.* Let $u \in S^i$ and $v \overset{\text{def}}{=} f^i(u) \in S^{i+1}$. $S^i$ is an open for $u$, $S^{i+1}$ is an open neighborhood of $v$ and $f^i\colon S^i \longrightarrow S^{i+1} \overset{\text{Ar}}{\in} \mathcal{T}^{i+1}$ by hypothesis. □

*If each $f^i\colon S^i \longrightarrow S^{i+1}$, is a local $\mathcal{T}^i$-$\mathcal{T}^{i+1}$ morphism of $S^i$ to $S^{i+1}$ then $f^2 \circ f^1\colon S^1 \longrightarrow S^3$ is a local $\mathcal{T}^1$-$\mathcal{T}^3$ morphism of $S^1$ to $S^3$.*

*Proof.* Since $\mathcal{T}^1 \overset{\text{full-cat}}{\subseteq} \mathcal{T}^2$ and $\mathcal{T}^2 \overset{\text{full-cat}}{\subseteq} \mathcal{T}^3$, $\mathcal{T}^1 \overset{\text{full-cat}}{\subseteq} \mathcal{T}^3$. Let $u \in S^1$, $v \overset{\text{def}}{=} f^1(u)$ and $w \overset{\text{def}}{=} f^2(v)$. There exist an open neighborhood $U_u$ for $u$, open neighborhoods $V_u$, $V'_v$ for $v$ and an open neighborhood $W_v$ of $w$ such that $f^1[U_u] \subseteq V_u$, $f^1\colon U_u \longrightarrow V_u$ is a morphism of $\mathcal{T}^2$, $f^2[V'_v] \subseteq W_v$ and $f^2\colon V'_v \longrightarrow W_v$ is a morphism of $\mathcal{T}^3$. Then $\hat{V}_u \overset{\text{def}}{=} V_u \cap V'_v \neq \emptyset$, $\hat{V}_u$ is an open neighborhood of $v$ and $\hat{U}_u \overset{\text{def}}{=} f_1^{i-1}[\hat{V}_u]$ is an open neighborhood of $u$. $f^1\colon \hat{U}_u \longrightarrow \hat{V}_u$ and $f^2\colon \hat{V}_u \longrightarrow W_v$ are morphisms of $\mathcal{T}^3$ by definition 1.6 (Topological category) on page 8 and thus $f^2 \circ f^1\colon \hat{U}_u \longrightarrow W_v$ is a morphism of $\mathcal{T}^3$. □

**Corollary 1.11** (Local morphisms). *Let $\mathcal{T}^i$, $i = 1, 2$, be a full topological category, $\mathcal{T}^i \overset{\text{full-cat}}{\subseteq} \mathcal{T}^{i+1}$, $S^i \overset{\text{Ob}}{\in} \mathcal{T}^i$ and $S^1 \subseteq S^2$. Then $\text{Id}_{S^2}^{S^1}$ is a local $\mathcal{T}^1$-$\mathcal{T}^2$ morphism of $S^1$ to $S^2$ and $\text{Id}_{S^i}$ is a local $\mathcal{T}$ morphism of $S^i$ to $S^i$.*

*Proof.* $S^1 \overset{\text{Ob}}{\in} \mathcal{T}^2$ because $S^1 \overset{\text{Ob}}{\in} \mathcal{T}^1$ and $\mathcal{T}^1 \overset{\text{cat}}{\subseteq} \mathcal{T}^2$, $S^2 \overset{\text{Ob}}{\in} \mathcal{T}^2$ by hypothesis and $S^1 \subseteq S^2$ by hypothesis, so $\text{Id}_{S^2}^{S^1} \overset{\text{Ar}}{\in} \mathcal{T}^2$ by lemma 1.7.

$\text{Id}_{S^i} \overset{\text{def}}{=} \text{Id}_{S^i,S^i}$. □

**Definition 1.12** (Sequence functions). Let $\boldsymbol{S} \overset{\text{def}}{=} (s_\alpha, \alpha \prec A)$ be a sequence of functions. Then

$$\textbf{domain}(\boldsymbol{S}) \overset{\text{def}}{=} (\text{domain}(s_\alpha), \alpha \prec A) \tag{1.2}$$

$$\textbf{range}(\boldsymbol{S}) \overset{\text{def}}{=} (\text{range}(s_\alpha), \alpha \prec A) \tag{1.3}$$

Let $\boldsymbol{T} \overset{\text{def}}{=} (t_\alpha, \alpha \prec A)$ be a sequence of functions with $\textbf{range}(\boldsymbol{S}) = \textbf{domain}(\boldsymbol{T})$. Then their composition is the sequence $\boldsymbol{T} \overset{()}{\circ} \boldsymbol{S} \overset{\text{def}}{=} (t_\alpha \circ s_\alpha, \alpha \prec A)$,

Let $\boldsymbol{S} \stackrel{\text{def}}{=} (s_\gamma, \gamma \leq \Gamma)$, then these functions extract information about the sequence:

$$\text{head}(\boldsymbol{S}, \Omega) \stackrel{\text{def}}{=} (s_\gamma, \gamma \prec \Omega) \tag{1.4}$$

$$\text{head}(\boldsymbol{S}) \stackrel{\text{def}}{=} \text{head}(\boldsymbol{S}, \Gamma) \tag{1.5}$$

$$\text{length0}(\boldsymbol{S}) \stackrel{\text{def}}{=} \Gamma \tag{1.6}$$

$$\text{tail}(\boldsymbol{S}) \stackrel{\text{def}}{=} s_\Gamma \tag{1.7}$$

Let $\boldsymbol{S} \stackrel{\text{def}}{=} (s_\gamma, \gamma \prec \Gamma)$, then

$$\text{length}(\boldsymbol{S}) \stackrel{\text{def}}{=} \Gamma \tag{1.8}$$

*Remark* 1.13. If length0($\boldsymbol{S}$) is defined then length($\boldsymbol{S}$) = length0($\boldsymbol{S}$) + 1. length0($\boldsymbol{S}$) is the ordinal type of head($\boldsymbol{S}$), not the ordinal type of $\boldsymbol{S}$.

Let $\boldsymbol{\mathcal{S}} \stackrel{\text{def}}{=} (\mathcal{S}_\alpha, \alpha \prec A)$ and $\boldsymbol{\mathcal{T}} \stackrel{\text{def}}{=} (\mathcal{T}_\alpha, \alpha \prec A)$ be sequences of categories. Then $\boldsymbol{\mathcal{S}}$ is a subcategory sequence of $\boldsymbol{\mathcal{T}}$, abbreviated $\boldsymbol{\mathcal{S}} \stackrel{\text{cat}}{\subseteq} \boldsymbol{\mathcal{T}}$, iff every category in $\boldsymbol{\mathcal{S}}$ is a subcategory of the corresponding category in $\boldsymbol{\mathcal{T}}$, i.e., $\left(\forall_{\alpha \prec A}\right) \mathcal{S}_\alpha \stackrel{\text{cat}}{\subseteq} \mathcal{T}_\alpha$, and $\boldsymbol{\mathcal{S}}$ is a full subcategory sequence of $\boldsymbol{\mathcal{T}}$, abbreviated $\boldsymbol{\mathcal{S}} \stackrel{\text{full-cat}}{\subseteq} \boldsymbol{\mathcal{T}}$, iff every category in $\boldsymbol{\mathcal{S}}$ is a full subcategory of the corresponding category in $\boldsymbol{\mathcal{T}}$, i.e., $\left(\forall_{\alpha \prec A}\right) \mathcal{S}_\alpha \stackrel{\text{full-cat}}{\subseteq} \mathcal{T}_\alpha$.

The category sequence union of $\boldsymbol{\mathcal{S}}$ and $\boldsymbol{\mathcal{T}}$, abbreviated $\boldsymbol{\mathcal{S}} \stackrel{\text{cat}}{\cup} \boldsymbol{\mathcal{T}}$, is the sequence of category unions of corresponding categories in $\boldsymbol{\mathcal{S}}$ and $\boldsymbol{\mathcal{T}}$, i.e., $(\mathcal{S}_\alpha \stackrel{\text{cat}}{\cup} \mathcal{T}_\alpha)$.

**Lemma 1.14** (Sequence functions). *Let $\boldsymbol{f}^i \stackrel{\text{def}}{=} (f_\alpha^i, \alpha \prec A)$, $i \in [1, 3]$, be sequences of functions with $\textbf{domain}(\boldsymbol{f}^2) = \textbf{range}(\boldsymbol{f}^1)$ and $\textbf{domain}(\boldsymbol{f}^3) = \textbf{range}(\boldsymbol{f}^2)$. Then $(\boldsymbol{f}^3 \stackrel{()}{\circ} \boldsymbol{f}^2) \stackrel{()}{\circ} \boldsymbol{f}^1 = \boldsymbol{f}^3 \stackrel{()}{\circ} (\boldsymbol{f}^2 \stackrel{()}{\circ} \boldsymbol{f}^1)$.*

*Proof.*

$$\begin{aligned}(\boldsymbol{f}^3 \stackrel{()}{\circ} \boldsymbol{f}^2) \stackrel{()}{\circ} \boldsymbol{f}^1 &= ((f_\alpha^3 \circ f_\alpha^2) \circ f_\alpha^1, \alpha \prec A) \\ &= (f_\alpha^3 \circ (f_\alpha^2 \circ f_\alpha^1), \alpha \prec A) \\ &= \boldsymbol{f}^3 \stackrel{()}{\circ} (\boldsymbol{f}^2 \stackrel{()}{\circ} \boldsymbol{f}^1)\end{aligned}$$

$\square$

Let $\boldsymbol{f} \stackrel{\text{def}}{=} (f_\alpha, \alpha \prec A)$ be a sequence of functions, $\boldsymbol{D} = \textbf{domain}(\boldsymbol{f})$ and $\boldsymbol{R} = \textbf{range}(\boldsymbol{f})$. Then $ID_{\boldsymbol{R}}$ is a left $\stackrel{()}{\circ}$ identity for $\boldsymbol{f}$ and $ID_{\boldsymbol{D}}$ is a right $\stackrel{()}{\circ}$ identity for $\boldsymbol{f}$.

*Proof.*

$$\mathrm{ID}_R \overset{()}{\circ} f = (\mathrm{Id}_{\mathrm{range}(f_\alpha)} \circ f_\alpha, \alpha \prec A)$$
$$= (f_\alpha, \alpha \prec A)$$
$$= f$$

$$f \overset{()}{\circ} \mathrm{ID}_D = (f_\alpha \circ \mathrm{Id}_{\mathrm{domain}(f_\alpha)}, \alpha \prec A)$$
$$= (f_\alpha, \alpha \prec A)$$
$$= f$$

□

**Definition 1.15** (Tuple functions). *Let $S \overset{\mathrm{def}}{=} (s_n, n \in [1, N])$ be a tuple of functions. Then*

$$\mathbf{domain}(S) \overset{\mathrm{def}}{=} (\mathrm{domain}(s_n), n \in [1, N]) \quad (1.9)$$

$$\mathbf{range}(S) \overset{\mathrm{def}}{=} (\mathrm{range}(s_n), n \in [1, N]) \quad (1.10)$$

Let $T \overset{\mathrm{def}}{=} (t_n, n \in [1, N])$ be a tuple of functions with $\mathbf{range}(S) = \mathbf{domain}(T)$, Then their composition is the tuple $T \overset{()}{\circ} S \overset{\mathrm{def}}{=} (t_n \circ s_n, n \in [1, N])$

Let $S \overset{\mathrm{def}}{=} (s_m, m \in [1, M])$ and $T \overset{\mathrm{def}}{=} (t_n, n \in [1, N])$ be tuples. Then the following are tuple functioms

$$\mathrm{head}(S, I) \overset{\mathrm{def}}{=} (s_m, m \in [1, I]) \quad (1.11)$$

$$\mathrm{head}(S) \overset{\mathrm{def}}{=} \mathrm{head}(S, M-1) \quad (1.12)$$

$$\mathrm{tail}(T, I) \overset{\mathrm{def}}{=} (t_n, n \in [I, N]) \quad (1.13)$$

$$\mathrm{tail}(T) \overset{\mathrm{def}}{=} t_N \quad (1.14)$$

$$\mathrm{join}(S, T) \overset{\mathrm{def}}{=} (s_1, \ldots, s_M, t_1, \ldots, t_N) \quad (1.15)$$

**Lemma 1.16** (Tuple functions). *Let $f \overset{\mathrm{def}}{=} (f_n, n \in [1, N])$ be a tuple of functions, $D \overset{\mathrm{def}}{=} \mathbf{domain}(f)$ and $R \overset{\mathrm{def}}{=} \mathbf{range}(f)$. Then $\mathrm{ID}_R$ is a left $\overset{()}{\circ}$ identity for $f$ and $\mathrm{ID}_D$ is a right $\overset{()}{\circ}$ identity for $f$.*

*Proof.*

$$\mathrm{ID}_R \overset{()}{\circ} f = (\mathrm{Id}_{\mathrm{range}(f_\alpha)} \circ f_n, n \in [1, N])$$
$$= (f_n, n \in [1, N])$$
$$= f$$

$$f \overset{()}{\circ} \mathrm{ID}_{\boldsymbol{D}} = (f_n \circ \mathrm{Id}_{\mathrm{domain}(f_n)}, n \in [1, N])$$
$$= (f_n, n \in [1, N])$$
$$= \boldsymbol{f}$$

□

**Lemma 1.17** (Tuple composition for unlabeled morphisms). *Let $\boldsymbol{f}^i \overset{\mathrm{def}}{=} (f_n^i, n \in [1, N])$, $i \in [1, 3]$, be tuples of functions with $\mathrm{domain}(\boldsymbol{f}^2) = \mathrm{range}(\boldsymbol{f}^1)$ and $\mathrm{domain}(\boldsymbol{f}^3) = \mathrm{range}(\boldsymbol{f}^2)$. Then $(\boldsymbol{f}^3 \overset{()}{\circ} \boldsymbol{f}^2) \overset{()}{\circ} \boldsymbol{f}^1 = \boldsymbol{f}^3 \overset{()}{\circ} (\boldsymbol{f}^2 \overset{()}{\circ} \boldsymbol{f}^1)$.*

*Proof.*

$$(\boldsymbol{f}^3 \overset{()}{\circ} \boldsymbol{f}^2) \overset{()}{\circ} \boldsymbol{f}^1 = ((f_n^3 \circ f_n^2) \circ f_n^1, n \in [1, N])$$
$$= (f_n^3 \circ (f_n^2 \circ f_n^1), n \in [1, N])$$
$$= \boldsymbol{f}^3 \overset{()}{\circ} (\boldsymbol{f}^2 \overset{()}{\circ} \boldsymbol{f}^1)$$

□

**Definition 1.18** (Tuple composition for labeled morphisms). Let $\boldsymbol{M}^i \overset{\mathrm{def}}{=} (\boldsymbol{f}^i, o_1^i, o_2^i)$, $i = 1, 2$, be tuples such that each $\boldsymbol{f}^i$ is a sequence of functions or each $\boldsymbol{f}^i$ is a tuple of functions, $\mathrm{range}(\boldsymbol{f}^1) = \mathrm{domain}(\boldsymbol{f}^2)$ and $0_2^1 = o_1^2$. Then

$$\boldsymbol{M}^2 \overset{A}{\circ} \boldsymbol{M}^1 \overset{\mathrm{def}}{=} (\boldsymbol{f}^2 \overset{()}{\circ} \boldsymbol{f}^1, o_1^1, o_2^2) \tag{1.16}$$

**Lemma 1.19** (Tuple composition for labeled morphisms). *Let $\boldsymbol{M}^i \overset{\mathrm{def}}{=} (\boldsymbol{f}^i, o_1^i, o_2^i)$, $i \in [1, 3]$, be a tuple such that $\boldsymbol{f}^i$ is a sequence or tuple of functions, $\mathrm{range}(\boldsymbol{f}^i) = \mathrm{domain}(\boldsymbol{f}^{i+1})$ and $o_2^i = o_1^{i+1}$, $i = 1, 2$. Then*

$$\boldsymbol{M}^3 \overset{A}{\circ} (\boldsymbol{M}^2 \overset{A}{\circ} \boldsymbol{M}^1) = (\boldsymbol{M}^3 \overset{A}{\circ} \boldsymbol{M}^2) \overset{A}{\circ} \boldsymbol{M}^1 \tag{1.17}$$

*Proof.* From definition 1.18 (Tuple composition for labeled morphisms), lemma 1.14 (Sequence functions) on page 10 and lemma 1.16 (Tuple functions), we have

$$\boldsymbol{M}^3 \overset{A}{\circ} (\boldsymbol{M}^2 \overset{A}{\circ} \boldsymbol{M}^1) = \boldsymbol{M}^3 \overset{A}{\circ} (\boldsymbol{f}^2 \overset{()}{\circ} \boldsymbol{f}^1, o_1^1, o_2^2)$$
$$= (\boldsymbol{f}^3 \overset{()}{\circ} \boldsymbol{f}^2 \overset{()}{\circ} \boldsymbol{f}^1, o_1^1, o_2^3)$$
$$= (\boldsymbol{f}^3 \overset{()}{\circ} \boldsymbol{f}^2, o_1^2, o_2^3) \overset{A}{\circ} \boldsymbol{M}^1$$
$$= (\boldsymbol{M}^3 \overset{A}{\circ} \boldsymbol{M}^2) \overset{A}{\circ} \boldsymbol{M}^1$$

□

Let $D^i \stackrel{\text{def}}{=} \textbf{domain}(f^i)$ and $R^i \stackrel{\text{def}}{=} \textbf{range}(f^i)$. Then

$$(\text{ID}_{R^i}, o_2^i, o_2^i) \stackrel{A}{\circ} M^i = M^i \tag{1.18}$$

$$M^i \stackrel{A}{\circ} (\text{ID}_{D^i}, o_1^i, o_1^i) = M^i \tag{1.19}$$

*Proof.*

$$(\text{ID}_{R^i}, o_2^i, o_2^i) \stackrel{A}{\circ} M^i = (\text{ID}_{R^i} \stackrel{()}{\circ} f^i, o_1^i, o_2^i)$$
$$= (f^i, o_1^i, o_2^i)$$
$$= M^i$$

$$M^i \stackrel{A}{\circ} (\text{ID}_{D^i}, o_1^i, o_1^i) = (f^i \stackrel{()}{\circ} \text{ID}_{D^i}, o_1^i, o_2^i)$$
$$= (f^i, o_1^i, o_2^i)$$
$$= M^i$$

$\square$

**Definition 1.20** (Cartesian product of sequence)**.** Let $S^i \stackrel{\text{def}}{=} (S_\alpha^i, \alpha \prec A)$, $i = 1, 2$, be a sequence and $f \stackrel{\text{def}}{=} (f_\alpha \colon S_\alpha^1 \longrightarrow S_\alpha^2, \alpha \prec A)$ be a sequence of functions, then $\bigtimes S^i \stackrel{\text{def}}{=} \bigtimes_{\alpha \prec A} S_\alpha^i$ is the generalized Cartesian product of the sequence $S^i$ and $\bigtimes f \colon S^1 \longrightarrow S^2 \stackrel{\text{def}}{=} \bigtimes_{\alpha \prec A} f_\alpha$ is the generalized Cartesian product of the function sequence $f$.

**Definition 1.21** (Topology functions)**.** Let $S$ and $T$ be a topological spaces and $Y$ a subset of $S$. Then

1. $S \stackrel{\text{subsp}}{\subseteq} T$ iff $S$ is a subspace of $T$. $S \stackrel{\text{open-subsp}}{\subseteq} T$ iff $S$ is an open subspace of $T$.

2. $\mathfrak{Top}(S)$ is the topology of $S$.

3. $\mathfrak{Top}(Y, S) \stackrel{\text{def}}{=} \{U \cap Y | U \in \mathfrak{Top}(S)\}$ is the relative topology of $Y$.

4. $\text{Top}(Y, S) \stackrel{\text{def}}{=} (Y, \mathfrak{Top}(Y, S))$ is $Y$ with the relative topology.

5. $S_{\text{op}} \stackrel{\text{def}}{=} \{(U, \mathfrak{Top}(U, S)) | U \in \mathfrak{Top}(S) \setminus \{\emptyset\}\}$ is the set of all non-null open subspaces of $S$.

Let $\mathbf{S}$ be a set of topological spaces. Then $\mathbf{S}_{\text{op}} \stackrel{\text{def}}{=} \bigcup_{S \in \mathbf{S}} S_{\text{op}}$ is the set of open subspaces in $\mathbf{S}$.

Let $S$ and $T'$ be spaces, $T \subseteq T'$ be a subspace and $f \colon S \longrightarrow T$ a function. Then $f \colon S \longrightarrow T' \stackrel{\text{def}}{=} \text{Id}_{T,T'} \circ f$ is $f$ considered as a function from $S$ to $T'$.

Let $S'$ and $T'$ be spaces, $S \subseteq S', T \subseteq T'$ be subspaces and $f' \colon S' \longrightarrow T'$ a function such that $f'[S] \subseteq T$. Then $f' \colon S \longrightarrow T$, also written $f'\restriction_{S,T}$, is $f'\restriction_S$ considered as a function from $S$ to $T$.

Let $\boldsymbol{S}^i \stackrel{\text{def}}{=} (S_\alpha^i, \alpha \preceq A)$, $i = 1, 2$, be a sequence of spaces, $\boldsymbol{S}^1 \stackrel{()}{\subseteq} \boldsymbol{S}^2$ and $f^2 \colon \operatorname{head}(\boldsymbol{S}^2) \longrightarrow \operatorname{tail}(\boldsymbol{S}^2)$ a function. $f^2\restriction_{\operatorname{head}(\boldsymbol{S}^1)} \stackrel{\text{def}}{=} f^2\restriction_{\bigtimes \operatorname{head}(\boldsymbol{S}^1)}$. If $f^2\restriction_{\operatorname{head}(\boldsymbol{S}^1)}[\bigtimes \operatorname{head}(\boldsymbol{S}^1)] \subseteq \operatorname{tail}(\boldsymbol{S}^1)$ then $f^2 \colon \operatorname{head}(\boldsymbol{S}^1) \longrightarrow \operatorname{tail}(\boldsymbol{S}^1)$, also written $f^2\restriction_{\operatorname{head}(\boldsymbol{S}^1),\operatorname{tail}(\boldsymbol{S}^1)}$, is $f^2\restriction_{\operatorname{head}(\boldsymbol{S}^1)}$ considered as a function from $\bigtimes \operatorname{head}(\boldsymbol{S}^1)$ to $\operatorname{tail}(\boldsymbol{S}^1)$.

**Definition 1.22** (Sequence inclusion). Let $\boldsymbol{S} \stackrel{\text{def}}{=} (S_\alpha, \alpha \prec A)$ and $\boldsymbol{T} \stackrel{\text{def}}{=} (T_\alpha, \alpha \prec A)$ be sequences. $\boldsymbol{S} \stackrel{()}{\in} \boldsymbol{T}$ iff $\left(\forall_{\alpha \prec A}\right) S_\alpha \in T_\alpha$ or $\left(\forall_{\alpha \prec A}\right) S_\alpha \stackrel{\text{Ob}}{\in} T_\alpha$.

**Lemma 1.23** (Sequence inclusion). *Let $\boldsymbol{S} \stackrel{\text{def}}{=} (S_\alpha, \alpha \prec A)$ be a sequence, $\boldsymbol{\mathcal{T}}^i \stackrel{\text{def}}{=} (\mathcal{T}_\alpha^i, \alpha \prec A)$, $i = 1, 2$, a sequence of categories, $\boldsymbol{\mathcal{T}}^1 \stackrel{\text{cat}}{\subseteq} \boldsymbol{\mathcal{T}}^2$ and $\boldsymbol{S} \stackrel{()}{\in} \boldsymbol{\mathcal{T}}^1$ Then $\boldsymbol{S} \stackrel{()}{\in} \boldsymbol{\mathcal{T}}^2$.*

*Proof.* If $\left(\forall_{\alpha \prec A}\right) S_\alpha \stackrel{\text{Ob}}{\in} \mathcal{T}_\alpha^1$ then $\left(\forall_{\alpha \prec A}\right) S_\alpha \stackrel{\text{Ob}}{\in} \mathcal{T}_\alpha^2$. □

# Part IV
# Nearly commutative diagrams

The notion of commutative diagrams is very useful in, e.g., Algebraic Topology. Often one encounters commutative diagrams in which two outgoing terminal nodes can be connected by a bridging function such that the resulting diagram is still commutative. This paper uses the term nearly commutative to describe a restricted class of such diagrams.

Let $\mathscr{C}$ be a full topological category and $D$ a tree with two branches, whose nodes are topological spaces $U_i$ and $V^j$ and whose links are continuous functions $f_i \colon U_i \longrightarrow U_{i+1}$ and $f'_j \colon U_j \longrightarrow U_{j+1}$ between the sets:

$$D \stackrel{\text{def}}{=} \{f_0 \colon U_0 = V_0 \longrightarrow U_1, \ldots, f_{m-1} \colon U_{m-1} \longrightarrow U_m, \qquad (1.20)$$
$$f'_0 \colon U_0 = V_0 \longrightarrow V_1, \ldots, f'_{m-1} \colon V_{m-1} \longrightarrow V_n\}$$

with $U_0 = V_0$, $U_m \stackrel{\text{Ob}}{\in} \mathscr{C}$ and $V_n \stackrel{\text{Ob}}{\in} \mathscr{C}$, as shown in fig. 1.

**Definition 1.24** (Nearly commutative diagrams in category $\mathscr{C}$). $D$ is left (right) nearly commutative in category $\mathscr{C}$ iff the two final nodes are in $\mathscr{C}$ and there is a morphism $\hat{f} \colon U_m \longrightarrow V_n \stackrel{\text{Ar}}{\in} \mathscr{C}$ ($\hat{f}' \colon V_n \longrightarrow U_m \stackrel{\text{Ar}}{\in} \mathscr{C}$) making the graph a commutative diagram, as shown in fig. 2. $D$ is nearly commutative in category $\mathscr{C}$ iff it is eiher left nearly commutative or right nearly commutative.

$D$ is strongly nearly commutative in category $\mathscr{C}$ iff the morphism is an isomorphism of $\mathscr{C}$,

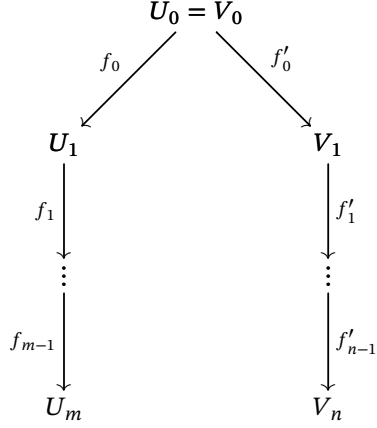

Figure 1: Uncompleted nearly commutative diagram

**Definition 1.25** (Nearly commutative diagrams in category $\mathscr{C}$ at a point)**.** Let $D$ be as above and $x$ be an element of the initial node. $D$ is (left,right,strongly) nearly commutative in $\mathscr{C}$ at $x$ iff there are subobjects of the nodes such that the tree $D'$ formed by replacing the nodes is (left,right,strongly) nearly commutative in $\mathscr{C}$, and $x$ is in the new initial node, as shown in fig. 3 (Local nearly commutative diagram) on page 16 below:

$$x \in U'_0 = V'_0 \tag{1.21}$$
$$U'_i \subseteq U_i \tag{1.22}$$
$$V'_j \subseteq V_j \tag{1.23}$$
$$f_i \!\upharpoonright_{U'_i} : U'_i \longrightarrow U'_{i+1} \tag{1.24}$$
$$f'_j \!\upharpoonright_{V'_j} : V'_j \longrightarrow V'_{j+1} \tag{1.25}$$

$$D' \stackrel{\text{def}}{=} \{f_0 \colon U'_0 = V'_0 \longrightarrow U'_1, \dots, f_{m-1} \colon U'_{m-1} \longrightarrow U'_m, \\ f'_0 \colon U'_0 = V'_0 \longrightarrow V'_1, \dots, f'_{m-1} \colon V'_{m-1} \longrightarrow V'_n\} \tag{1.26}$$

**Definition 1.26** (Locally nearly commutative diagrams in category $\mathscr{C}$)**.** Let $D$ be as above. $D$ is (left,right,strongly) locally nearly commutative in $\mathscr{C}$ iff it is (left,right,strongly) nearly commutative in $\mathscr{C}$ at $x$ for every $x$ in the initial node.

*Remark* 1.27. It will often be clear from context what the relevant categories are. This paper may use the term "nearly commutative" without explicitly identifying the categories in which the modes are found.

**Lemma 1.28** (Locally nearly commutative diagrams in category $\mathscr{C}$)**.** *Let $\mathscr{C}$ and $D$ be as above. If $U_0 = \emptyset$ then $D$ is (left, right, strongly) locally nearly commutative in $\mathscr{C}$.*

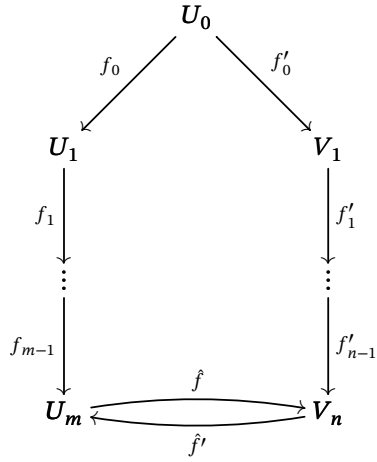

Figure 2: Completed nearly commutative diagram

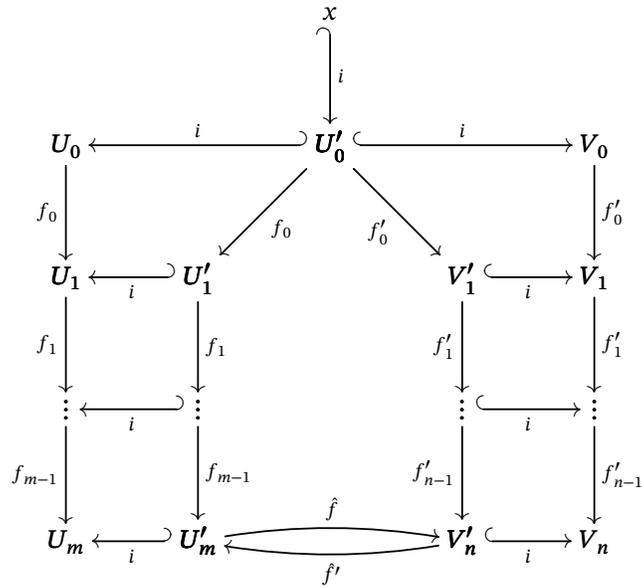

Figure 3: Local nearly commutative diagram

*Proof.* $D$ is vacuously locally nearly commutative at every $x \in U_0$ since there is no such $x$. □

*Let $D$ be locally nearly commutative in $\mathscr{C}$ and let $\hat{U}_0 = \hat{V}_0 \subseteq U_0$ and $\hat{D}$ be $D$ with $U_0 = V_0$ replaced by $\hat{U}_0 = \hat{V}_0$. Then $\hat{D}$ is locally nearly commutative in $\mathscr{C}$.*

*Proof.* If $x \in \hat{U}_0$ then $x \in U_0$ and hence $D$ is locally nearly commutative in $\mathscr{C}$ at $x$. Replacing $U_0'$ with $U_0' \cap \hat{U}_0$ in the definition shows that $\hat{D}$ is locally nearly commutative in $\mathscr{C}$ at $x$. □

**Lemma 1.29** (Nearly commutative diagrams in category $\mathscr{C}'$). *Let $\mathscr{C} \stackrel{\text{cat}}{\subseteq} \mathscr{C}'$ and $D$ be as in diagram* (1.20) *above.*

*If $D$ is (left, right, strongly) nearly commutative in $\mathscr{C}$ then $D$ is (left, right, strongly) commutative in $\mathscr{C}'$.*

*Proof.* If $U_m \stackrel{\text{Ob}}{\in} \mathscr{C}$ then $U_m \stackrel{\text{Ob}}{\in} \mathscr{C}'$. If $V_n \stackrel{\text{Ob}}{\in} \mathscr{C}$ then $V_n \stackrel{\text{Ob}}{\in} \mathscr{C}'$. If $\hat{f} \colon U_m \rightarrowtail\!\!\!\stackrel{\tilde{=}}{\twoheadrightarrow} V_n \stackrel{\text{Ar}}{\in} \mathscr{C}$ then $\hat{f} \stackrel{\text{Ar}}{\in} \mathscr{C}'$ □

**Corollary 1.30** (Locally nearly commutative diagrams in category $\mathscr{C}'$ at a point). *Let $\mathscr{C} \stackrel{\text{cat}}{\subseteq} \mathscr{C}'$ and $D$ be as in diagram* (1.20) *above.*

*If $D$ is locally nearly commutative in $\mathscr{C}$ at $x$ then $D$ is locally nearly commutative in $\mathscr{C}'$ at $x$.*

*Proof.* Let $U_i', V_j'$ and diagram (1.26) be as in definition 1.25 (Nearly commutative diagrams in category $\mathscr{C}$ at a point) on page 15, as shown in diagram fig. 3 (Local nearly commutative diagram). $D'$ is nearly commutative in category $\mathscr{C}$, and hence is nearly commutative in category $\mathscr{C}$. □

**Corollary 1.31** (Locally nearly commutative diagrams in category $\mathscr{C}'$). *Let $\mathscr{C} \stackrel{\text{cat}}{\subseteq} \mathscr{C}'$ and $D$ be as in diagram* (1.20) *above.*

*If $D$ is is locally nearly commutative in $\mathscr{C}$ then $D$ is locally nearly commutative in $\mathscr{C}'$.*

*Proof.* Let $x \in U_0 = V_0$. $D$ is is locally nearly commutative in $\mathscr{C}$ at $x$, hence locally nearly commutative in $\mathscr{C}'$ at $x$. □

# Part V
# Model spaces and allied notions

Let $S$ be a topological space. We need to formalize the notions of an open cover by sets that are "well behaved" in some sense, e.g., convex, sufficiently small, and of "well behaved" functions among those sets, e.g., preserving fibers, smooth. We do this by associating a category of acceptable sets and functions.

*Remark* 1.32. Using pseudo-groups, as in [Kobayashi, 1996, p. 1], would not allow restricting model neighborhoods to, e.g., convex sets.

## 2 Model spaces

**Definition 2.1** (Model spaces)**.** Let $S$ be a topological space and $\mathcal{S}$ a small category whose objects are open subsets of $S$ and whose morphisms are continuous functions. $\boldsymbol{S} \stackrel{\text{def}}{=} (S, \mathcal{S})$ is a model space for $S$ iff

1. $\text{Ob}(\mathcal{S})$ is an open cover for $S$. Note that it need not be a basis for $S$.

2. $\text{Ob}(\mathcal{S})$ is closed under finite intersections.

3. The morphisms of $\mathcal{S}$ are continuous functions in $S$.

4. If $f\colon A \longrightarrow B$ is a morphism, $A' \stackrel{\text{Ob}}{\in} \mathcal{S} \subseteq A \stackrel{\text{Ob}}{\in} \mathcal{S}$, $B' \stackrel{\text{Ob}}{\in} \mathcal{S} \subseteq B \stackrel{\text{Ob}}{\in} \mathcal{S}$ and $f[A'] \subseteq B'$ then $f\restriction_{A'}\colon A' \longrightarrow B'$ is a morphism.

5. If $A' \stackrel{\text{Ob}}{\in} \mathcal{S} \subseteq A \stackrel{\text{Ob}}{\in} \mathcal{S}$ then the inclusion map $\text{Id}_{A',A}\colon A' \hookrightarrow A$ is a morphism.

    *Remark* 2.2. This is actually a consequence of item 4, but it is convenient to give it here.

6. Restricted sheaf condition: informally, consistent morphisms can be glued together. Whenever

    (a) $U_\alpha$ and $V_\alpha$, $\alpha \prec A$, are objects of $\mathcal{S}$.

    (b) $f_\alpha\colon U_\alpha \longrightarrow V_\alpha$ are morphisms of $\mathcal{S}$.

    (c) $U \stackrel{\text{def}}{=} \bigcup_{\alpha \prec A} U_\alpha$ is an object of $\mathcal{S}$.

    (d) $V \stackrel{\text{def}}{=} \bigcup_{\alpha \prec A} V_\alpha$ is an object of $\mathcal{S}$.

    (e) $f\colon U \longrightarrow V$ is a continuous function and for every $\alpha \prec A$, $f$ agrees with $f_\alpha$ on $U_\alpha$

    then $f$ is a morphism of $\mathcal{S}$.

*Remark* 2.3. A model space is similar to a pseudo-group, but has some important differences:

1. Open subsets of objects of $\mathcal{S}$ need not be objects of $\mathcal{S}$.

2. The morphisms of $\mathcal{S}$ need not be homeomorphisms or even 1-1, thus they need not have inverses.

3. The morphisms of $\mathcal{S}$ need not satisfy the sheaf condition, but only the restricted sheaf condition in item 6.

$\boldsymbol{S}$ is fine grained iff every open subset of an object of $\mathcal{S}$ is an object of $\mathcal{S}$. Note that it is not sufficient that $\mathrm{Ob}(\mathcal{S})$ be a basis for $S$.

$\mathrm{Top}(\boldsymbol{S}) \stackrel{\text{def}}{=} S = \pi_1(\boldsymbol{S})$ is the topological space of $\boldsymbol{S}$.

Let $\boldsymbol{S}^i$, $i = 1, 2$, be a model space. $\mathrm{Cat}(\boldsymbol{S}^i) \stackrel{\text{def}}{=} \mathcal{S}^i = \pi_2(\boldsymbol{S}^i)$ is the category of model space $\boldsymbol{S}^i$. By abuse of language we write $\boldsymbol{S}^1 \stackrel{\text{subsp}}{\subseteq} \boldsymbol{S}^2$ ($\boldsymbol{S}^1 \stackrel{\text{open-subsp}}{\subseteq} \boldsymbol{S}^2$) for $\mathrm{Top}(\boldsymbol{S}^1) \stackrel{\text{subsp}}{\subseteq} \mathrm{Top}(\boldsymbol{S}^2)$ ($\mathrm{Top}(\boldsymbol{S}^1) \stackrel{\text{open-subsp}}{\subseteq} \mathrm{Top}(\boldsymbol{S}^2)$).

Let $U \stackrel{\mathrm{Ob}}{\in} \mathcal{S}$. Then $\mathrm{Top}(U, \boldsymbol{S}) \stackrel{\text{def}}{=} \mathrm{Top}(U, S)$ is $U$ with the relative topology.

$\mathcal{T}\!op(\boldsymbol{S}) \stackrel{\text{def}}{=} \left( \left\{ \mathrm{Top}(U, S) \,\middle|\, U \stackrel{\mathrm{Ob}}{\in} \mathcal{S} \right\}, \mathrm{Ar}(\mathcal{S}) \right)$ is the topological category of $\boldsymbol{S}$.

By abuse of language we write $U \subseteq \boldsymbol{S}$ for $U \stackrel{\mathrm{Ob}}{\in} \mathcal{S}$, $f$ is a morphism of $\boldsymbol{S}$ for $f$ is a morphism of $\mathcal{S}$ and $f$ is an isomorphism of $\boldsymbol{S}$ for $f$ is an isomorphism of $\mathcal{S}$.

Let $\boldsymbol{C}$ be a set of model spaces. $\boldsymbol{C}$ is fine grained iff every $(C, \mathcal{C}) \in \boldsymbol{C}$ is fine grained.

**Lemma 2.4** (The topological category of a model space is a full topological category)**.** *Let $\boldsymbol{M} \stackrel{\text{def}}{=} (S, \mathcal{S})$ be a model space for $S$. Then $\mathcal{T}\!op(\boldsymbol{M})$ is a full topological category.*

*Proof.* $\mathcal{T}\!op(\boldsymbol{M})$ is a small subcategory of **Top** by construction. $\mathcal{T}\!op(\boldsymbol{M})$ is a full topological category by item 4 of definition 2.1. □

**Definition 2.5** (Model neighborhoods)**.** Let $\boldsymbol{S} \stackrel{\text{def}}{=} (S, \mathcal{S})$ be a model space for $S$. Then the objects of $\mathcal{S}$ are model neighborhoods of $\boldsymbol{S}$. If $u \in U \stackrel{\mathrm{Ob}}{\in} \mathcal{S}$ then $U$ is a model neighborhood of $\boldsymbol{S}$ for $u$. If $U_i \stackrel{\mathrm{Ob}}{\in} \mathcal{S}$, $i = 1, 2$, and $U_1 \subseteq U_2$ then $U_1$ is a model subneighborhood of $U_2$; if $u^1 \in U_1$ then $U_1$ is also a model subneighborhood for $u^1$.

If $\emptyset \stackrel{\mathrm{Ob}}{\in} \mathcal{S}$ them $\emptyset$ is a degenerate model neighborhood of $\boldsymbol{S}$; any nonvoid model neighborhood of $\boldsymbol{S}$ is a nonvoid model neighborhood of $\boldsymbol{S}$.

**Definition 2.6** (Model subspaces)**.** $\boldsymbol{S} \stackrel{\text{def}}{=} (S, \mathcal{S})$ is a model subspace of $\boldsymbol{T} \stackrel{\text{def}}{=} (T, \mathcal{T})$, abbreviated $\boldsymbol{S} \stackrel{\text{mod}}{\subseteq} \boldsymbol{T}$, iff

1. $\boldsymbol{S}$ and $\boldsymbol{T}$ are model spaces

2. $S$ is a subspace of $T$

3. $\mathcal{S} \stackrel{\text{full-cat}}{\subseteq} \mathcal{T}$

4. every intersection of $S$ with an object of $\mathcal{T}$ is an object of $\mathcal{S}$

$$\left( \forall_{U \stackrel{\mathrm{Ob}}{\in} \mathcal{T}} \right) U \cap S \stackrel{\mathrm{Ob}}{\in} \mathcal{S} \tag{2.1}$$

5. Every object of $\mathcal{T}$ contained in $S$ is an object of $\mathcal{S}$

$$\left(\forall_{U \stackrel{\text{Ob}}{\in} \mathcal{T} : U \subseteq S}\right) U \stackrel{\text{Ob}}{\in} \mathcal{S} \tag{2.2}$$

*Remark* 2.7. This is actually a consequence of item 4, but it is convenient to state it here.

$\boldsymbol{S}$ is also a full model subspace of $\boldsymbol{T}$, abbreviated $\boldsymbol{S} \stackrel{\text{full-mod}}{\subseteq} \boldsymbol{T}$, iff $S = T$.

$\boldsymbol{S}$ is also a strict model subspace of $\boldsymbol{T}$, abbreviated $\boldsymbol{S} \stackrel{\text{strict-mod}}{\subseteq} \boldsymbol{T}$, iff $\text{Id}_{S,T}$ is open and continuous.

Let $\boldsymbol{S} \stackrel{\text{def}}{=} (S, \mathcal{S})$ be a model space and $T$ a model neighborhood of $\boldsymbol{S}$. Then $\text{Mod}(T, \boldsymbol{S})$, the relative model space of $T$, is $(T, \mathcal{T})$, where $\mathcal{T}$ is the full subcategory of $\mathcal{S}$ containing all model subneighborhoods of $T$.

By abuse of language we write $T$ for $\text{Mod}(T, \boldsymbol{S})$ when $\boldsymbol{S}$ is understood by context.

**Lemma 2.8** (Model subspaces). $\stackrel{\text{mod}}{\subseteq}$ and $\stackrel{\text{strict-mod}}{\subseteq}$ are transitive.

*Proof.* Let $\boldsymbol{S}^i \stackrel{\text{def}}{=} (S^i, \mathcal{S}^i)$ $i = 1, 2, 3$, be model spaces. If $\boldsymbol{S}^i \stackrel{\text{mod}}{\subseteq} \boldsymbol{S}^{i+1}$, $i = 1, 2$, then

1. $\boldsymbol{S}^1$ and $\boldsymbol{S}^3$ are model spaces

2. Since $S^i$ is a subspace of $S^{i+1}$, $i = 1, 2$, then $S^1$ is a subspace of $S^3$.

3. Since $\mathcal{S}^i \stackrel{\text{full-cat}}{\subseteq} \mathcal{S}^{i+1}$, $i = 1, 2$, then $\mathcal{S}^1 \stackrel{\text{full-cat}}{\subseteq} \mathcal{S}^3$.

4. Since every intersection of $S^i$ with an object of $\mathcal{S}^{i+1}$ is an object of $\mathcal{S}^i$, $i = 1, 2$, every intersection of $S^1$ with an object of $\mathcal{S}^3$ is an object of $\mathcal{S}^1$.

5. Let $U \stackrel{\text{Ob}}{\in} \mathcal{S}^3$ and $U \subseteq S^1$. $U \cap S^1 = U$ is an object of $\mathcal{S}^1$.

6. If $\boldsymbol{S}^i \stackrel{\text{strict-mod}}{\subseteq} \boldsymbol{S}^{i+1}$, $i = 1, 2$, then $\text{Id}_{\boldsymbol{S}^i,\boldsymbol{S}^{i+1}}$ is open and continuous, hence $\text{Id}_{\boldsymbol{S}^1,\boldsymbol{S}^3} = \text{Id}_{\boldsymbol{S}^2,\boldsymbol{S}^3} \circ \text{Id}_{\boldsymbol{S}^1,\boldsymbol{S}^2}$ is open and continuous and thus $\boldsymbol{S}^1 \stackrel{\text{strict-mod}}{\subseteq} \boldsymbol{S}^3$. □

Let $\boldsymbol{S} \stackrel{\text{def}}{=} (S.\mathcal{S})$ be a model space and $T \subseteq \boldsymbol{S}$. Then $\text{Mod}(T, \boldsymbol{S}) \stackrel{\text{strict-mod}}{\subseteq} \boldsymbol{S}$

*Proof.* Let $\boldsymbol{T} = \left((T.\mathfrak{T}), \mathcal{T} \, bigr\right) \stackrel{\text{def}}{=} \text{Mod}(T, \boldsymbol{S})$.

1. $\boldsymbol{S}$ and $\boldsymbol{T}$ are model spaces.

2. $T$ is a subspace of $S$.

3. $\mathcal{T} \stackrel{\text{strictcat}}{\subseteq} \mathcal{S}$ by construction.

4. Let $U \stackrel{\text{Ob}}{\in} \boldsymbol{S}$ be contained in $T$. Then $U$ is a model neighborhood of $\text{Mod}(T, \boldsymbol{S})$ by construction.

$\square$

*Let $\boldsymbol{S}^i \stackrel{\text{def}}{=} (S^i, \mathcal{S}^i)$, $i = 1, 2$, and $\boldsymbol{S}^1 \stackrel{\text{mod}}{\subseteq} \boldsymbol{S}^2$. Then $S^1$ is an open subspace of $S^2$.*

*Proof.* The model neighborhoods of $\boldsymbol{S}^1$ are model neighborhoods of $\boldsymbol{S}^2$, hence open in $S^2$. Their union, which is open in $S^2$, is $S^1$ by item 1 of definition 2.6 (Model subspaces) on page 19. $\square$

*If $\text{Ob}(\mathcal{S}^1)$ is a basis for $S^1$ then $\text{Id}_{S^2}^{S^1}$ is open.*

*Proof.* Let $U$ be open in $S^1$. $U$ is a union of model neighborhoods of $\boldsymbol{S}^1$, each of which is also a model neighborhood of $\boldsymbol{S}^2$ by item 3 of definition 2.6 (Model subspaces) on page 19, and hence open in $S^2$. Then their union $U$ is open in $S^2$. $\square$

*If $\text{Ob}(\mathcal{S}^2)$ is a basis for $S^2$ then $\text{Id}_{S^2}^{S^1}$ is continuous.*

*Proof.* Let $V$ be open in $S^2$. $V \cap S^1$ is open in $S^2$, and by hypothesis a union of model neighborhoods of $\boldsymbol{S}^2$, each of which is also a model neighborhood of $\boldsymbol{S}^1$ by item 5 of definition 2.6, and hence open in $S^1$. Then their union $U$ is open in $S^1$. $\square$

**Corollary 2.9.** *Let $\boldsymbol{S}^i \stackrel{\text{def}}{=} (S^i, \mathcal{S}^i)$, $i = 1, 2$, be a model space, $\boldsymbol{S}^1 \stackrel{\text{mod}}{\subseteq} \boldsymbol{S}^2$ and $\text{Ob}(\mathcal{S}^i)$ be a basis for $S^i$. Then $\boldsymbol{S}^1 \stackrel{\text{strict-mod}}{\subseteq} \boldsymbol{S}^2$.*

*Proof.* $\text{Id}_{S^2}^{S^1}$ is open and continuous, the conditions of definition 2.6 (Model subspaces) above. $\square$

# 3 M-nearly commutative diagrams

Let $\boldsymbol{M} \stackrel{\text{def}}{=} (M, \mathcal{M})$ be a model space and $D$ a tree with two branches, whose nodes are topological spaces $U_i$ and $V_j$ and whose links are continuous functions $f_i \colon U_i \longrightarrow U_{i+1}$ and $f'_j \colon V_j \longrightarrow V_{j+1}$ between the spaces:

$$D \stackrel{\text{def}}{=} \{f_0 \colon U_0 = V_0 \longrightarrow U_1, \ldots, f_{m-1} \colon U_{m-1} \longrightarrow U_m,$$
$$f'_0 \colon U_0 = V_0 \longrightarrow V_1, \ldots, f'_{m-1} \colon V_{m-1} \longrightarrow V_n\}$$

with $U_0 = V_0$, $U_m \stackrel{\text{Ob}}{\in} \mathcal{M}$ and $V_n \stackrel{\text{Ob}}{\in} \mathcal{M}$, as shown in fig. 1 (Uncompleted nearly commutative diagram) on page 15.

**Definition 3.1** (M-nearly commutative diagrams). *$D$ is (left,right,strongly) M-nearly commutative in model space $\boldsymbol{M}$ iff $D$ is (left,right,strongly) nearly commutative in category $\mathcal{M}$.*

**Definition 3.2** (M-nearly commutative diagrams at a point)**.** Let $M$ and $D$ be as above and $x$ be an element of the initial node. $D$ is (left,right,strongly) M-nearly commutative in $M$ at $x$ iff $D$ is (left,right,strongly) nearly commutative in category $\mathcal{M}$ at $x$.

**Definition 3.3** (M-locally nearly commutative diagrams)**.** Let $\mathcal{C}$ and $D$ be as above. $D$ is (left,right,strongly) M-locally nearly commutative in $M$ iff $D$ is (left,right,strongly) nearly commutative in category $\mathcal{M}$ at every point.

# 4 Trivial model spaces

Informally, a trivial model space of a specific type is one that does not restrict the potential objects and morphisms of its type.

**Definition 4.1** (Trivial model spaces and categories)**.** Let $S$ be a topological space and $\mathcal{S}$ the category of all continuous functions between open sets of $S$.

$S_{\text{triv}} \overset{\text{def}}{=} (S, \mathcal{S})$ is the trivial model space of $S$ and $S_{\text{triv}}$ is a trivial model space.

Let $\boldsymbol{S}$ be a set of topological spaces.

$\boldsymbol{S}_{\text{triv}} \overset{\text{def}}{=} \left\{ S'_{\text{triv}} \,\middle|\, S' \in \boldsymbol{S} \right\}$ is the set of all trivial model spaces in $\boldsymbol{S}$.

$\boldsymbol{S}_{\mathcal{T}riv}$, the category of all continuous functions between elements of $\boldsymbol{S}_{\text{triv}}$, is the trivial model category of $\boldsymbol{S}$.

$\boldsymbol{S}_{\text{op−triv}}$, the category of open trivial model spaces in $\boldsymbol{S}$, is the category whose objects are $\left\{ U_{\text{triv}} \,\middle|\, U \in \boldsymbol{S}_{\text{op}} \right\}$, the trivial model spaces of non-null open sets of spaces in $\boldsymbol{S}$, and whose morphisms are all the continuous functions among them.

Let $\boldsymbol{S} \overset{\text{def}}{=} (S, \mathcal{S})$ be a model space.

$\boldsymbol{S}_{\text{submod−triv}} \overset{\text{def}}{=} \left\{ M \,\middle|\, M \overset{\text{mod}}{\subseteq} \boldsymbol{S} \right\}$ is the set of all subspaces of $\boldsymbol{S}$.

$\boldsymbol{M}_{\text{submod−}\mathcal{T}riv}$, the category of all model functions between elements of $\boldsymbol{S}_{\text{submod−triv}}$, is the trivial subspace model category of $\boldsymbol{S}$.

$\boldsymbol{S}_{\text{full−submod−triv}} \overset{\text{def}}{=} \left\{ M \,\middle|\, M \overset{\text{full−mod}}{\subseteq} \boldsymbol{S} \right\}$ is the set of all full subspaces of $\boldsymbol{S}$.

$\boldsymbol{M}_{\text{full−submod−}\mathcal{T}riv}$, the category of all model functions between elements of $\boldsymbol{S}_{\text{full−submod−triv}}$, is the trivial full subspace model category of $\boldsymbol{S}$.

Let $\boldsymbol{S}$ be a set of model spaces.

$\boldsymbol{S}_{\text{submod−triv}} \overset{\text{def}}{=} \left\{ M \,\middle|\, \left(\exists_{M' \in \boldsymbol{S}}\right) M \overset{\text{mod}}{\subseteq} M' \right\}$ is the set of all subspaces in $\boldsymbol{S}$.

$\boldsymbol{S}_{\text{submod−}\mathcal{T}riv}$, the category of all model functions between elements of $\boldsymbol{S}_{\text{submod−triv}}$, is the trivial subspace model category of $\boldsymbol{S}$.

$$\mathbf{S}_{\textbf{full–submod–triv}} \stackrel{\text{def}}{=} \left\{ M \,\middle|\, \left(\exists_{M' \in S}\right) M \stackrel{\text{full–mod}}{\subseteq} M' \right\} \text{ is the set of all full subspaces}$$

in $\mathbf{S}$.

$\mathbf{S}_{\textbf{submod–}\mathcal{T}riv}$, the category of all model functions between elements of $\mathbf{S}_{\textbf{full–submod–triv}}$, is the trivial full subspace model category of $\mathbf{S}$.

**Lemma 4.2** (Trivial model spaces). *Let $S^i$, $i = 1, 2$, be a topological space.*

1. $S^i_{\text{triv}}$ *is a model space.*

   *Proof.* Let $\mathcal{S}^i \stackrel{\text{def}}{=} \text{Cat}(S^i_{\text{triv}})$, $f \colon A \longrightarrow B$ be a morphism of $\mathcal{S}^i$, $A' \stackrel{\text{Ob}}{\in} \mathcal{S}^i \subseteq A$ and $B' \stackrel{\text{Ob}}{\in} \mathcal{S}^i \subseteq B$. Then the definitions of continuity and of $S^i_{\text{triv}}$ imply each of the following.

   (a) $\text{Ob}(\mathcal{S}^i) = \mathfrak{Top}(S^i)$ and thus is an open cover of $S^i$.
   
   (b) $\text{Ob}(\mathcal{S}^i) = \mathfrak{Top}(S^i)$ and thus is closed under finite intersections.
   
   (c) All morphisms in $\text{Ar}(\mathcal{S}^i)$ are continuous.
   
   (d) Since $f \colon A \longrightarrow B$ is a morphism of $\mathcal{S}^i$, $f$ is continuous.
   
   (e) Since $A' \stackrel{\text{Ob}}{\in} \mathcal{S}^i$ then $A'$ is open.
   
   (f) Since $B' \stackrel{\text{Ob}}{\in} \mathcal{S}^i$ then $B'$ is open.
   
   (g) Since $f \colon A \longrightarrow B$ is continuous then $f\restriction_{A'} \colon A' \longrightarrow B$ is continuous.
   
   (h) Since $f[A'] \subseteq B'$ then $f\restriction_{A'} \colon A' \longrightarrow B'$ is well defined and continuous, hence a morphism.
   
   (i) If $A' \stackrel{\text{Ob}}{\in} \mathcal{S}^i \subseteq A \stackrel{\text{Ob}}{\in} \mathcal{S}^i$ then
   
      i. $A$ and $A'$ are open.
   
      ii. The inclusion map $\text{Id}_{A'} \colon A' \hookrightarrow A$ is continuous.
   
      iii. The inclusion map $i \colon A' \hookrightarrow A$ is a morphism of $\mathcal{S}^i$ by the definition of $S^i_{\text{triv}}$.

   $\square$

2. *Let $S^1$ be an open subspace of $S^2$. Then $S^1_{\text{triv}} \stackrel{\text{strict–mod}}{\subseteq} S^2_{\text{triv}}$.*

   *Proof.* Let $\mathcal{S}^i \stackrel{\text{def}}{=} \text{Cat}(S^i_{\text{triv}})$, $i = 1, 2$.

   $S^1_{\text{triv}}$ and $S^2_{\text{triv}}$ satisfy the three conditions of definition 2.6 (Model subspaces) on page 19:

**Top($S^1_{\text{triv}}$) is a subspace of Top($S^2_{\text{triv}}$):**

Top($S^i_{\text{triv}}$) = $S^i$ by definition 4.1 (Trivial model spaces and categories) on page 22. $S^1$ is a subspace of $S^2$ by hypothesis.

$\mathcal{S}^1 \overset{\text{full-cat}}{\subseteq} \mathcal{S}^2$:

Let $f\colon U \longrightarrow V$ be a morphism of $S^1_{\text{triv}}$. $U$ and $V$ are objects of $S^1_{\text{triv}}$, hence open in $S^1$, and open in $S^2$ because $S^1$ is an open subspace of $S^2$ by hypothesis.

Let $W$ be an open set in $S^2$. $V \cap W$ is open in $S^2$, and open in $S^1$ because $S^1$ is an open subspace of $S^2$ by hypothesis. $f^{-1}[V] = f^{-1}[V \cap W]$ is open in $S^1$ and thus open in $S^2$.

**Every object of $S^2_{\text{triv}}$ contained in $S^1$ is an object of $S^1_{\text{triv}}$:**

$\text{Id}^{S^1}_{S^2}$ is continuous by hypothesis.

$\text{Id}^{S^1}_{S^2}$ **is open and continuous:** The topology of $\boldsymbol{S}^i$ is the same as the topology of $S^i_{\text{triv}}$, and $\text{Id}^{S^1}_{S^2}$ is open and continuous by hypothesis.

□

3. *Let $\boldsymbol{S}^i \overset{\text{def}}{=} (S^i, \mathcal{S}^i), i = 1, 2$, be a model space and $\boldsymbol{S}^1 \overset{\text{strict-mod}}{\subseteq} \boldsymbol{S}^2$. Then $\boldsymbol{S}^1_{\text{triv}} \overset{\text{strict-mod}}{\subseteq} \boldsymbol{S}^2_{\text{triv}}$.*

*Proof.* Let $\mathcal{S}'^i \overset{\text{def}}{=} \pi_2 \boldsymbol{S}^i_{\text{triv}}, i = 1, 2$.

$\boldsymbol{S}^1_{\text{triv}}$ and $\boldsymbol{S}^2_{\text{triv}}$ satisfy the three conditions of definition 2.6 (Model subspaces) on page 19:

**Top($S^1_{\text{triv}}$) is a subspace of Top($S^2_{\text{triv}}$):**

By hypothesis, $S^1 \subseteq S^2$ and $\text{Id}^{S^1}_{S^2}$ is open.

$\mathcal{S}'^1 \overset{\text{full-cat}}{\subseteq} \mathcal{S}'^2$:

Let $f\colon U \longrightarrow V$ be a morphism of $S^1_{\text{triv}}$. $U$ and $V$ are objects of $S^1_{\text{triv}}$, hence open in $S^1$, and open in $S^2$ because $\text{Id}^{S^1}_{S^2}$ is open by hypothesis.

Let $W$ be an open set in $S^2$. $V \cap W$ is open in $S^2$, and open in $S^1$ because $\text{Id}^{S^1}_{S^2}$ is continuous by hypothesis. $f^{-1}[V] = f^{-1}[V \cap W]$ is open in $S^1$ and thus open in $S^2$.

**Every object of $S^2_{\text{triv}}$ contained in $S^1$ is an object of $S^1_{\text{triv}}$:**

$\text{Id}^{S^1}_{S^2}$ is continuous by hypothesis.

$\mathrm{Id}_{S^2}^{S^1}$ **is open and continuous:** The topology of $\boldsymbol{S}^i$ is the same as the topology of $\underset{\mathbf{triv}}{S^i}$, and $\mathrm{Id}_{S^2}^{S^1}$ is open and continuous by hypothesis.

$\square$

4. Let $S^i$, $i = 1, 2$ be a topological space and $S^1$ be an open subspace of $S^2$. Then $\underset{\mathbf{triv}}{S^1} \overset{\text{strict-mod}}{\subseteq} \underset{\mathbf{triv}}{S^2}$.

    *Proof.* By definition 4.1 (Trivial model spaces and categories) on page 22, $\underset{\mathbf{triv}}{S^i} = (S^i, \mathcal{S}^i)$, where $\mathcal{S}^i$ is the category of all continuous functions between open sets of $S^i$,

    $\underset{\mathbf{triv}}{S^1} \overset{\text{mod}}{\subseteq} \underset{\mathbf{triv}}{S^2}$:

    (a) $S^1$ is a subspace of $S^2$ by hypothesis.
    (b) Every object of $\mathcal{S}^1$ is open in $S^1$, thus open in $S^2$ and thus an object of $\mathcal{S}^2$. If $f\colon U \longrightarrow V \overset{\text{Ar}}{\in} \mathcal{S}^1$ then $U$ and $V$ are open in $S^1$ and $f$ is continuous in $S^1$, thus $U$ and $V$ are open in $S^2$ and $f$ is continuous in $S^2$ and $f\colon U \longrightarrow V \overset{\text{Ar}}{\in} \mathcal{S}^2$. If $f\colon U \longrightarrow V$ is a morphism of $\mathcal{S}^2$, $U \overset{\text{Ob}}{\in} \mathcal{S}^1$ and $V \overset{\text{Ob}}{\in} \mathcal{S}^1$, then $U$ and $V$ are open in $S^1$ and $f$ is continuous in $S^2$, thus continuous in $S^1$, and $f\colon U \longrightarrow V \overset{\text{Ar}}{\in} \mathcal{S}^1$.
    (c) If $U \overset{\text{Ob}}{\in} \mathcal{S}^2$ then $U$ is open in $S^2$, hence $U \cap S^2$ is open in $S^2$, hence $U \cap S^2$ is open in $S^2$ and $U \overset{\text{Ob}}{\in} \mathcal{S}^2$.
    (d) If $U \overset{\text{Ob}}{\in} \mathcal{S}^2$ and $U \subseteq S^1$ then U is open is $S^2$, hence open is $S^1$, hence $U \overset{\text{Ob}}{\in} \mathcal{S}^1$.

    $\underset{\mathbf{triv}}{S^1} \overset{\text{strict-mod}}{\subseteq} \underset{\mathbf{triv}}{S^2}$: $\mathrm{Id}_{S^1,S^2}$ is open and continuous because $S^1$ is an open subspace of $S^2$.

$\square$

5. Let $\boldsymbol{S} \overset{\text{def}}{=} (S, \mathcal{S})$ be a model space. $\underset{\mathbf{full-submod-triv}}{\boldsymbol{S}} \overset{\text{fullcat}}{\subseteq} \underset{\text{submod}-\mathcal{T}riv}{\boldsymbol{M}}$.

    *Proof.* $\left(\forall \underset{\boldsymbol{M} \subseteq \boldsymbol{S}}{\text{full-mod}}\right) \boldsymbol{M} \overset{\text{mod}}{\subseteq} \boldsymbol{S}$, so $\mathrm{Ob}(\underset{\mathbf{full-submod-triv}}{\boldsymbol{S}}) \subseteq \mathrm{Ob}(\underset{\text{submod}-\mathcal{T}riv}{\boldsymbol{M}})$. It remains to show that for a pair of objects in both categories, the Hom sets are the same. Let $\boldsymbol{M}^i \overset{\text{full-mod}}{\subseteq} \boldsymbol{S}$, $i = 1.2$.

(a) Let $f\colon \boldsymbol{M}^1 \longrightarrow \boldsymbol{M}^2 \overset{\text{Ar}}{\in} \boldsymbol{S}_{\text{full--submod--triv}}$. Then $f$ is a model function and thus $f\colon \boldsymbol{M}^1 \longrightarrow \boldsymbol{M}^2 \overset{\text{Ar}}{\in} \boldsymbol{S}_{\text{submod--triv}}$.

(b) Let $f\colon \boldsymbol{M}^1 \longrightarrow \boldsymbol{M}^2 \overset{\text{Ar}}{\in} \boldsymbol{S}_{\text{submod--triv}}$. Then $f$ is a model function and thus $f\colon \boldsymbol{M}^1 \longrightarrow \boldsymbol{M}^2 \overset{\text{Ar}}{\in} \boldsymbol{S}_{\text{full--submod--triv}}$.

$\square$

6. Let $f\colon S^1 \longrightarrow S^2$ be continuous. Then $f$ is a model function from $S^1_{\text{triv}}$ to $S^2_{\text{triv}}$.

   *Proof.* Let $U^i$ be a model neighborhood of $S^i_{\text{triv}}$, $i = 1, 2$.

   (a) Since $U^2$ is a model neighborhood of $S^2_{\text{triv}}$, $U^2$ is open, $f^{-1}[U^2]$ is open and thus a model neighborhood of $S^1_{\text{triv}}$.

   (b) $f[U^1] \subseteq S^2$. Since $S^2$ is open, it is a model neighborhood of $S^2_{\text{triv}}$.

$\square$

**Corollary 4.3** (Trivial model spaces)**.** *Let $\boldsymbol{S}^i \overset{\text{def}}{=} (S^i, \mathcal{S}^i)$, $i = 1, 2$, be a model space, $\boldsymbol{S}^1 \overset{\text{mod}}{\subseteq} \boldsymbol{S}^2$ and $\mathrm{Ob}(\mathcal{S}^i)$ be a basis for $S^i$. Then $S^1_{\text{triv}} \overset{\text{strict--mod}}{\subseteq} S^2_{\text{triv}}$.*

*Proof.* $\mathrm{Ob}(\mathcal{S}^i)$ is a basis for $S^1$, $i = 1, 2$, by hypothesis. $\boldsymbol{S}^1 \overset{\text{strict--mod}}{\subseteq} \boldsymbol{S}^2$ by corollary 2.9 () on page 21. $\square$

# 5 Minimal model spaces

Even if a set of open sets fails one of items 1 and 2 or a set of functions among them fails one of items 3 to 6 in the definition of a model space, there is a minimal model space containing them.

**Definition 5.1** (Minimal model spaces)**.** Let $S$ be a topological space, $\boldsymbol{O}$ a set of open sets in $S$, $\boldsymbol{f}$ a set of continuous functions between elements of $\boldsymbol{O}$ and $\mathcal{S}$ the smallest concrete category over **Top** having all sets in $\boldsymbol{O}$ as objects, having all functions in $\boldsymbol{f}$ as morphisms and satisfying items 3 to 6 of definition 2.1 (Model spaces) on page 18 .
Then
$$\operatorname*{Mod}_{\min}(S, \boldsymbol{O}, \boldsymbol{f}) \overset{\text{def}}{=} \left(\mathrm{Top}(\bigcup \boldsymbol{O}, S), \mathcal{S}\right) \tag{5.1}$$
is the minimal model space of $S$ with neighborhoods $\boldsymbol{O}$ and neighborhood mappings $\boldsymbol{f}$.

*Remark* 5.2. The trivial model space $S_{\text{triv}}$ is a special case.

Let $\boldsymbol{S} \stackrel{\text{def}}{=} (S, \mathcal{S})$ be a model space, $\boldsymbol{O}$ a set of model neighborhoods of $\boldsymbol{S}$ and $\boldsymbol{f}$ a set of morphisms between elements of $\boldsymbol{O}$.

Then

$$\underset{\min}{\operatorname{Mod}}(\boldsymbol{S}, \boldsymbol{O}, \boldsymbol{f}) \stackrel{\text{def}}{=} \underset{\min}{\operatorname{Mod}}(S, \boldsymbol{O}, \boldsymbol{f}) \tag{5.2}$$

is the minimal submodel space of $\boldsymbol{S}$ with neighborhoods $\boldsymbol{O}$ and neighborhood mappings $\boldsymbol{f}$.

*Remark* 5.3. while every model neighborhood of $\underset{\min}{\operatorname{Mod}}(\boldsymbol{S}, \boldsymbol{O}, \boldsymbol{f})$ is a model neighborhood of $\underset{\min}{\operatorname{Mod}}(S, \boldsymbol{O}, \boldsymbol{f})$, the converse need not be true.

**Lemma 5.4** (Minimal model spaces are model spaces). *Let $S$ be a topological space, $\boldsymbol{O}$ a set of open sets in $S$ and $\boldsymbol{f}$ a set of continuous functions between elements of $\boldsymbol{O}$. Then $(C, \mathscr{C}) \stackrel{\text{def}}{=} \underset{\min}{\operatorname{Mod}}(S, \boldsymbol{O}, \boldsymbol{f})$ is a model space.*

*Proof.*

1. Finite intersections of open sets are open and $C = \bigcup_{U \in \boldsymbol{O}} U$ by construction
2. $\operatorname{Ob}(\mathscr{C})$ is closed under finite intersections by construction
3. Compositions of continuous functions, inclusion maps and restrictions of continuous functions are continuous
4. Restrictions of morphisms are morphisms by construction
5. Inclusion maps are morphisms by construction
6. The restricted sheaf condition holds by construction

$\square$

Let $\boldsymbol{S} \stackrel{\text{def}}{=} (S, \mathcal{S})$ be a model space, $\boldsymbol{O}$ a set of model neighborhoods of $\boldsymbol{S}$ and $\boldsymbol{f}$ a set of morphisms between elements of $\boldsymbol{O}$. Then $\underset{\min}{\operatorname{Mod}}(\boldsymbol{S}, \boldsymbol{O}, \boldsymbol{f})$ is a model space.

*Proof.* $\underset{\min}{\operatorname{Mod}}(\boldsymbol{S}, \boldsymbol{O}, \boldsymbol{f}) = \underset{\min}{\operatorname{Mod}}(S, \boldsymbol{O}, \boldsymbol{f})$ and $\underset{\min}{\operatorname{Mod}}(S, \boldsymbol{O}, \boldsymbol{f})$ is a model space. $\square$

# 6 M-paracompact model spaces

Paracompactness is an important property for topological spaces because of partitions of unity. There is an analogous property for model spaces.

**Definition 6.1** (Model topology). Let $\boldsymbol{M} \stackrel{\text{def}}{=} (S, \mathcal{S})$ be a model space for $S$. Then the model topology $\mathfrak{M}^*$ for $M$ is the topology generated by $\operatorname{Ob}(\mathcal{S})$.

*Remark* 6.2. $\mathfrak{M}^*$ is not guarantied to be T0 even if $S$ is T4. However, $\mathfrak{M}^*$ may be normal or regular even if $S$ is not.

**Definition 6.3** (m-paracompactness). A model space $\boldsymbol{M} \stackrel{\text{def}}{=} (S, \mathcal{S})$ is m-paracompact iff $\mathfrak{M}^*$ is regular and every cover of $S$ by model neighborhoods has a locally finite refinement by model neighborhoods.

*Remark* 6.4. This is a stronger condition than merely requiring $\mathfrak{M}^*$ to be paracompact.

**Theorem 6.5** (m-paracompactness and paracompactness). *If $\boldsymbol{M} = (S, \mathcal{S})$ is m-paracompact and the model neighborhoods form a basis for $S$ then $S$ is paracompact.*

*Proof.* Let $R$ be an open cover of $S$. Since the model neighborhoods form a basis, every set in $R$ is a union of model neighborhoods and thus there is a refinement $R_1$ by model neighborhoods. Since by hypothesis $\boldsymbol{M}$ is m-paracompact, $R_1$ has a locally finite refinement $R_2$ by model neighborhoods. Since model neighborhoods are open sets, $R_2$ is a locally finite refinement of $R$ in the conventional sense. □

# 7 Model functions and model categories

It is convenient to have a notion of mappings between model spaces that are well behaved in some sense, e.g., fiber preserving; using that notion it is then possible to group model spaces into categories.

## 7.1 Model functions

**Definition 7.1** (Model functions). Let $\boldsymbol{S}^i \stackrel{\text{def}}{=} (S^i, \mathcal{S}^i)$, $i = 1, 2$, be a model space and $f \colon S^1 \longrightarrow S^2$ be a continuous function. $f$ is a model function of $\boldsymbol{S}^1$ to $\boldsymbol{S}^2$ iff the inverse images of model neighborhoods are model neighborhoods and the images of model neighborhoods are contained in model neighborhoods.

$$\left(\forall_{V \stackrel{\text{Ob}}{\in} \mathcal{S}^2}\right) f^{-1}[V] \stackrel{\text{Ob}}{\in} \mathcal{S}^1 \tag{7.1}$$

$$\left(\forall_{U \stackrel{\text{Ob}}{\in} \mathcal{S}^1}\right)\left(\exists_{V \stackrel{\text{Ob}}{\in} \mathcal{S}^2}\right) f[U] \subseteq V \tag{7.2}$$

$f$ is also a constrained model function of $\boldsymbol{S}^1$ to $\boldsymbol{S}^2$ iff $f[S^1]$ is contained in a model neighborhood of $\boldsymbol{S}^2$.

$$\left(\exists_{V \stackrel{\text{Ob}}{\in} \mathcal{S}^2}\right) f[S^1] \subseteq V \tag{7.3}$$

*Remark* 7.2. When $f \colon \boldsymbol{S}^1 \longrightarrow \boldsymbol{S}^2$ is constrained, $S^1$ is a model neighborhood of $\boldsymbol{S}^1$.

By abuse of language we write $f\colon \boldsymbol{S}^1 \longrightarrow \boldsymbol{S}^2$ both for $f$ considered as a model function and for $f$ considered as a continuous function.

Let $\boldsymbol{S}^1 \stackrel{\text{def}}{=} (S^1, \mathcal{S}^1)$ and $\boldsymbol{S}'^2 \stackrel{\text{def}}{=} (S'^2, \mathcal{S}'^2)$ be model spaces, $\boldsymbol{S}^2 \stackrel{\text{def}}{=} (S^2, \mathcal{S}^2) \stackrel{\text{mod}}{\subseteq} \boldsymbol{S}'^2$ be a model subspace and $f\colon S^1 \longrightarrow S^2$ be a model function of $\boldsymbol{S}^1$ to $\boldsymbol{S}^2$. Then $f\colon \boldsymbol{S}^1 \longrightarrow \boldsymbol{S}'^2 \stackrel{\text{def}}{=} \mathrm{Id}_{\boldsymbol{S}^2, \boldsymbol{S}'^2} \circ f$ is $f$ considered as a model function from $\boldsymbol{S}^1$ to $\boldsymbol{S}'^2$.

Let $\boldsymbol{S}'^i \stackrel{\text{def}}{=} (S'^i, \mathcal{S}'^i)$, $i = 1, 2$, be a model space, $\boldsymbol{S}^i \stackrel{\text{def}}{=} (S^i, \mathcal{S}^i) \stackrel{\text{mod}}{\subseteq} \boldsymbol{S}'^i$ be a model subspace and $f'$ be a model function of $\boldsymbol{S}'^1$ to $\boldsymbol{S}'^2$ such that $f'[S^1] \subseteq S^2$. Then $f\colon \boldsymbol{S}^1 \longrightarrow \boldsymbol{S}^2$, also written $f'\!\restriction_{\boldsymbol{S}^1, \boldsymbol{S}^2}$, is $f'\!\restriction_{S^1}$ considered as a model function of $\boldsymbol{S}^1$ to $\boldsymbol{S}^2$.

**Lemma 7.3** (Model functions). *Let $\boldsymbol{S}^i \stackrel{\text{def}}{=} (S^i, \mathcal{S}^i)$ and $\boldsymbol{S}'^i \stackrel{\text{def}}{=} (S'^i, \mathcal{S}'^i)$, $i = 1, 2$, be model spaces, $\boldsymbol{S}^i \stackrel{\text{mod}}{\subseteq} \boldsymbol{S}'^i$ and $f'\colon S'^1 \longrightarrow S'^2$ be a model function of $\boldsymbol{S}'^1$ to $\boldsymbol{S}'^2$.*

$f \stackrel{\text{def}}{=} f'\!\restriction_{\boldsymbol{S}^1, \boldsymbol{S}^2}$ *is a model function of $\boldsymbol{S}^1$ to $\boldsymbol{S}^2$.*

*Proof.* $f'$ is a model function by hypothesis. $f'$ is continuous by definition 7.1 above, so $f \stackrel{\text{def}}{=} f'\!\restriction_{\boldsymbol{S}^1, \boldsymbol{S}^2}$ is continuous.

Let $U$ be a model neighborhood of $\boldsymbol{S}^1$. $\boldsymbol{S}^1 \stackrel{\text{mod}}{\subseteq} \boldsymbol{S}'^1$ by hypothesis, so $U$ is a model neighborhood of $\boldsymbol{S}'^1$. $f'$ is a model function of $\boldsymbol{S}'^1$ to $\boldsymbol{S}'^2$ by hypothesis. □

**Lemma 7.4** (Model functions and trivial model spaces). *Let $S^i$, $i = 1, 2$, be a topological space and $f\colon S^1 \longrightarrow S^2$ be a continuous function. Then $f\colon \underset{\text{triv}}{S^1} \longrightarrow \underset{\text{triv}}{S^2}$ is a model function.*

*Proof.* Let $U^i$ be a model neighborhood of $\underset{\text{triv}}{S^i}$, $i = 1, 2$.

$f^{-1}[U^2]$ is open, hence a model neighborhood of $\underset{\text{triv}}{S^1}$.

$S^2$ is a model neighborhood of $\underset{\text{triv}}{S^2}$ and $f[U^1] \subseteq U^2$. □

Let $\boldsymbol{S}^i \stackrel{\text{def}}{=} (S^i, \mathcal{S}^i)$, $i = 1, 2$, be a model space and $f\colon \boldsymbol{S}^1 \longrightarrow \boldsymbol{S}^2$ be a model function. Then $f\colon \underset{\text{triv}}{S^1} \longrightarrow \underset{\text{triv}}{S^2}$ is a model function.

*Proof.* $f\colon S^1 \longrightarrow S^2$ is continuous. □

Let $\boldsymbol{S}^i \stackrel{\text{def}}{=} (S^i, \mathcal{S}^i)$, $i = 1, 2$, be a model space and $f\colon \boldsymbol{S}^1 \longrightarrow \boldsymbol{S}^2$ be a model function. Then $f\colon \underset{\text{triv}}{S^1} \twoheadrightarrow \boldsymbol{S}^2$ is a model function iff $f\colon \boldsymbol{S}^1 \twoheadrightarrow \boldsymbol{S}^2$ is a constrained model function.

*Proof.* let $U^2$ be a model neighborhood of $\boldsymbol{S}^2$. $f^{-1}[U^2]$ is open, hence a model neighborhood of $\underset{\text{triv}}{S^1}$.

$S^1$ is a model neighborhood of $\underset{\text{triv}}{S^1}$. If $f\colon \underset{\text{triv}}{S^1} \longrightarrow \underset{\text{triv}}{S^2}$ is a model function then $f[S^1]$ is contained in a model neighbothood of $\boldsymbol{S}^2$ and thus $f$ is constrained.

If $f$ is constrained, let $U^1$ be a model neighborhood of $S^1_{\text{triv}}$. Then $f[U^1] \subseteq f[S^1]$ is contained in a model neighbothood of $\boldsymbol{S}^2$. □

**Lemma 7.5** (Composition of model functions). *Let $\boldsymbol{S}_i \stackrel{\text{def}}{=} (S_i, \mathcal{S}_i)$, $i \in [1,3]$, be a model space and $f_i \colon \boldsymbol{S}_i \longrightarrow \boldsymbol{S}_{i+1}$, $i = 1, 2$, be a model function.*
*$f_2 \circ f_1$ is a model function.*

*Proof.* Let $U_i \stackrel{\text{Ob}}{\in} \mathcal{S}_i$, $i = 1, 2, 3$.

Since $f_2$ is a model function, $f_2^{-1}[U_3] \stackrel{\text{Ob}}{\in} \mathcal{S}_2$. Since $f_1$ is a model function, $(f_2 \circ f_1)^{-1}[U_3] = f_1^{-1}[f_2^{-1}[U_3]] \stackrel{\text{Ob}}{\in} \mathcal{S}_1$. Since $f_1$ is a model functions, $f_1[U_1]$ is contained in a model neighborhood $V_2$. Since $f_2$ is a model functions, $f_2[V_2]$ is contained in a model neighborhood $V_3$. Then $(f_2 \circ f_1)[U_1] \subseteq V_3$. □

*If each $f_i$ is constrained then $f_2 \circ f_1$ is constrained.*

*Proof.* $f_1[S^1] \subseteq S^2$ and $f_2[S^2] \subseteq S^3$, hence $f_2 \circ f_1[S^1] \subseteq S^3$. □

**Definition 7.6** (Model homeomorphisms). Let $\boldsymbol{S}^i \stackrel{\text{def}}{=} (S^i, \mathcal{S}^i)$, $i = 1, 2$, be a model space and $f \colon S^1 \longrightarrow S^2$ be a model function. $f$ is a model homeomorphism iff it is also invertible and its inverse is a model function.

## 7.2 Model categories

**Definition 7.7** (Model categories). A category $\mathcal{M}$ is a model category iff

1. the objects of $\mathcal{M}$ are model spaces

2. the morphisms of $\mathcal{M}$ are model functions.

3. composition is functional composition.

4. If $\boldsymbol{S} \stackrel{\text{def}}{=} (S, \mathcal{S}) \stackrel{\text{Ob}}{\in} \mathcal{M}$ and $S^i \stackrel{\text{Ob}}{\in} \mathcal{S}$, $i = 1, 2$, then $\boldsymbol{S}^i \stackrel{\text{def}}{=} \text{Mod}(S^i, \boldsymbol{S}) \stackrel{\text{Ob}}{\in} \mathcal{M}$. If further $S^1 \subseteq S^2$ then the inclusion map $\text{Id}^{S^1}_{S^2} \colon \boldsymbol{S}^1 \rightarrowtail \boldsymbol{S}^2$ is a morphism of $\mathcal{S}$.

5. If $\boldsymbol{S}^i \stackrel{\text{def}}{=} (S^i, \mathcal{S}^i) \stackrel{\text{Ob}}{\in} \mathcal{M}$, $U^i \stackrel{\text{Ob}}{\in} \mathcal{S}^i$, $i = 1, 2$, $f \colon \boldsymbol{S}^1 \longrightarrow \boldsymbol{S}^2$ is a morphism of $\mathcal{M}$ and $f[U^1] \subseteq U^2$ then $f \colon \text{Mod}(U^1, \boldsymbol{S}^1) \longrightarrow \text{Mod}(U^2, \boldsymbol{S}^2)$ is a morphism of $\mathcal{M}$.

6. If $\boldsymbol{S}^i \stackrel{\text{def}}{=} (S^i, \mathcal{S}^i)$, $i = 1, 2$, is a subspace of $\boldsymbol{S} \stackrel{\text{def}}{=} (S, \mathcal{S}) \stackrel{\text{Ob}}{\in} \mathcal{M}$ and $f \colon S^1 \longrightarrow S^2$ is a morphism of $\mathcal{S}$ then $f \colon \boldsymbol{S}^1 \longrightarrow \boldsymbol{S}^2$ is a morphism of $\mathcal{M}$.

A model category $\mathcal{M}$ satisfies the restricted sheaf condition iff whenever

1. $U$ and $V$ are objects of $\mathcal{S}$.

2. $U_\alpha \stackrel{\text{full-mod}}{\subseteq} U$ and $V_\alpha \stackrel{\text{full-mod}}{\subseteq} V$, $\alpha \prec A$, are objects of $\mathcal{M}$.

3. $f_\alpha \colon U_\alpha \longrightarrow V_\alpha$ are morphisms of $\mathcal{M}$.

4. $f \colon U \longrightarrow V$ is a model function and for every $\alpha < A$, $f$ agrees with $f_\alpha$ on $U_\alpha$

then $f$ is a morphism of $\mathcal{M}$.

**Lemma 7.8** (Model categories). *Let $\mathcal{M}$ be a model category, $\boldsymbol{S}^i \stackrel{\text{def}}{=} (S^i, \mathcal{S}^i) \stackrel{\text{Ob}}{\in} \mathcal{M}$. $i = 1, 2$, $U^i \stackrel{\text{Ob}}{\in} \mathcal{S}^i$, $\boldsymbol{U}^i \stackrel{\text{def}}{=} \mathrm{Mod}(U^i, \boldsymbol{S}^i)$ and $f \colon \boldsymbol{U}^1 \longrightarrow \boldsymbol{U}^2 \stackrel{\text{Ar}}{\in} \mathcal{M}$.*

$f \colon \boldsymbol{U}^1 \longrightarrow \boldsymbol{S}^2 \stackrel{\text{Ar}}{\in} \mathcal{M}$.

*Proof.* Since $\mathrm{Id}_{\boldsymbol{S}^2}^{\boldsymbol{U}^2} \stackrel{\text{Ar}}{\in} \mathcal{M}$ by item 4 of definition 7.7 above and $f \colon \boldsymbol{U}^1 \longrightarrow \boldsymbol{U}^2 \stackrel{\text{Ar}}{\in} \mathcal{M}$ by hypothesis, $f \colon \boldsymbol{U}^1 \longrightarrow \boldsymbol{S}^2 = \mathrm{Id}_{\boldsymbol{S}^2}^{\boldsymbol{U}^2} \circ f \colon \boldsymbol{U}^1 \longrightarrow \boldsymbol{U}^2 \stackrel{\text{Ar}}{\in} \mathcal{M}$ □

**Definition 7.9** (Trivial model categories). Let $\boldsymbol{M}$ be a set of model spaces. Then $\mathrm{Mod}(\boldsymbol{M})^{\mathrm{triv}}$ is the category of all model functions between model subspaces of model spaces in $\boldsymbol{M}$

## 7.3 (Local) m-morphisms

**Definition 7.10** (Local m-morphisms). Let $\mathcal{M}^i$ be a model category, $i = 1, 2$, and $\boldsymbol{S}^i \stackrel{\text{def}}{=} (S^i, \mathcal{S}^i) \stackrel{\text{Ob}}{\in} \mathcal{M}^i$.

Let $U^i \subseteq \boldsymbol{S}^i$, $i = 1, 2$. be open. A function $f \colon U^1 \longrightarrow U^2$ is a local $\boldsymbol{S}^1$-$\boldsymbol{S}^2$ m-morphism of $U^1$ to $U^2$ iff $\boldsymbol{S}^1 \stackrel{\text{mod}}{\subseteq} \boldsymbol{S}^2$ and for every $u \in U^1$ there is a model neighborhood $U_u \subseteq U^1$ for $u$ and a model neighborhood $V_u \subseteq U^2$ for $v \stackrel{\text{def}}{=} f(u)$ such that $f[U_u] \subseteq V_u$ and $f \colon U_u \longrightarrow V_u$ is a morphism of $\boldsymbol{S}^2$. It is also a strict local $\boldsymbol{S}^1$-$\boldsymbol{S}^2$ m-morphism of $U^1$ to $U^2$ iff $\boldsymbol{S}^1 \stackrel{\text{strict-mod}}{\subseteq} \boldsymbol{S}^2$.

Let $U^1 \subseteq \boldsymbol{S}^i$ and $U^2 \subseteq \boldsymbol{S}^i$, $i = 1, 2$, be open. A function $f \colon U^1 \longrightarrow U^2$ is a (strict) local $\boldsymbol{S}^i$ m-morphism of $U^1$ to $U^2$ iff it is a (strict) local $\boldsymbol{S}^i$-$\boldsymbol{S}^i$ m-morphism of $U^1$ to $U^2$.

A model function $f \colon \boldsymbol{S}^1 \longrightarrow \boldsymbol{S}^2$ is a (strict) local m-morphism of $\boldsymbol{S}^1$ to $\boldsymbol{S}^2$ iff it is a (strict) local $\boldsymbol{S}^1$-$\boldsymbol{S}^2$ m-morphism of $S^1$ to $S^2$.[5]

A model function $f \colon \boldsymbol{S}^i \longrightarrow \boldsymbol{S}^i$ is a (strict) local m-morphism of $\boldsymbol{S}^i$ iff it is a (strict) local m-morphism of $\boldsymbol{S}^i$ to $\boldsymbol{S}^i$.

Let $U^i \stackrel{\text{Ob}}{\in} \mathcal{S}^i$, $i = 1, 2$. A function $f \colon U^1 \longrightarrow U^2$ is a local $\boldsymbol{S}^1$-$\boldsymbol{S}^2$-$\mathcal{M}^1$-$\mathcal{M}^2$ m-morphism of $U^1$ to $U^2$ iff $\mathcal{M}^1 \stackrel{\text{full-cat}}{\subseteq} \mathcal{M}^2$ and for every $u \in U^1$ there is a model neighborhood $U_u \subseteq U^1$ for $u$ and a model neighborhood $V_u \subseteq U^2$ for $v \stackrel{\text{def}}{=} f(u)$ such that $f[U_u] \subseteq V_u$ and $f \colon \mathrm{Mod}(U_u, \boldsymbol{S}^1) \longrightarrow \mathrm{Mod}(V_u, \boldsymbol{S}^2)$ is a morphism of $\mathcal{M}^2$.

---

[5]This definition does not require that $S^i$ be a model neighborhod of $S^i$.

*Remark* 7.11. We do not need a notion of a strict local $\boldsymbol{S}^1$-$\boldsymbol{S}^2$-$\mathcal{M}^1$-$\mathcal{M}^2$ m-morphism of $U^1$ to $U^2$ since the definition of a local $\boldsymbol{S}^1$-$\boldsymbol{S}^2$-$\mathcal{M}^1$-$\mathcal{M}^2$ m-morphism of $U^1$ to $U^2$ already includes the requirement that $\mathcal{M}^1 \stackrel{\text{full-cat}}{\subseteq} \mathcal{M}^2$.

Let $U^1 \subseteq \boldsymbol{S}^i$ and $U^2 \subseteq \boldsymbol{S}^i$, be open. A model function $f\colon U^1 \longrightarrow U^2$ is a local $\boldsymbol{S}^i$-$\mathcal{M}^i$ m-morphism of $U^1$ to $U^2$ iff it is a local $\boldsymbol{S}^i$-$\boldsymbol{S}^i$-$\mathcal{M}^i$-$\mathcal{M}^i$ m-morphism of $U^1$ to $U^2$.

A model function $f\colon \boldsymbol{S}^1 \longrightarrow \boldsymbol{S}^2$ is a local $\mathcal{M}^1$-$\mathcal{M}^2$ m-morphism of $\boldsymbol{S}^1$ to $\boldsymbol{S}^2$ iff it is a local $\boldsymbol{S}^1$-$\boldsymbol{S}^2$-$\mathcal{M}^1$-$\mathcal{M}^2$ m-morphism of $S^1$ to $S^2$.[5]

A model function $f\colon \boldsymbol{S}^i \longrightarrow \boldsymbol{S}^i$ is a local $\mathcal{M}^i$ m-morphism iff it is a local $\mathcal{M}^i$-$\mathcal{M}^i$ m-morphism of $\boldsymbol{S}^i$ to $\boldsymbol{S}^i$.

The phrases "locally a(n) (strict) ... m-morphism" and "a (strict) local ... m-morphism" are equivalent.

**Definition 7.12** (M-morphisms). Let $\mathcal{M}^i$ be a model category, $i = 1, 2$, and $\boldsymbol{S}^i \stackrel{\text{def}}{=} (S^i, \mathcal{S}^i) \stackrel{\text{Ob}}{\in} \mathcal{M}^i$.

Let $U^i \subseteq \boldsymbol{S}^i$, $i = 1, 2$. be open. A model function $f\colon U^1 \longrightarrow U^2$ is an $\boldsymbol{S}^1$-$\boldsymbol{S}^2$ m-morphism of $U^1$ to $U^2$ iff $\boldsymbol{S}^1 \stackrel{\text{mod}}{\subseteq} \boldsymbol{S}^2$ and $f$ is a morphism of $\boldsymbol{S}^2$. It is also a strict $\boldsymbol{S}^1$-$\boldsymbol{S}^2$ m-morphism of $U^1$ to $U^2$ iff $\boldsymbol{S}^1 \stackrel{\text{strict-mod}}{\subseteq} \boldsymbol{S}^2$.

A model function $f\colon \boldsymbol{S}^1 \longrightarrow \boldsymbol{S}^2$ is a (strict) m-morphism of $\boldsymbol{S}^1$ to $\boldsymbol{S}^2$ iff it is a (strict) $\boldsymbol{S}^1$-$\boldsymbol{S}^2$ m-morphism of $S^1$ to $S^2$.[5]

A model function $f\colon \boldsymbol{S}^i \longrightarrow \boldsymbol{S}^i$ is a (strict) m-morphism of $\boldsymbol{S}^i$ iff it is a (strict) m-morphism of $\boldsymbol{S}^i$ to $\boldsymbol{S}^i$.

Let $U^i \subseteq \boldsymbol{S}^i$, $i = 1, 2$. be model neighborhoods. A model function $f\colon U^1 \longrightarrow U^2$ is an $\boldsymbol{S}^1$-$\boldsymbol{S}^2$-$\mathcal{M}^1$-$\mathcal{M}^2$ m-morphism of $U^1$ to $U^2$ iff $\mathcal{M}^1 \stackrel{\text{full-cat}}{\subseteq} \mathcal{M}^2$ and $f\colon \mathrm{Mod}(U^1, \boldsymbol{S}^1) \longrightarrow \mathrm{Mod}(U^2, \boldsymbol{S}^2)$ is a morphism of $\mathcal{M}^2$.

*Remark* 7.13. We do not need a notion of a strict $\boldsymbol{S}^1$-$\boldsymbol{S}^2$-$\mathcal{M}^1$-$\mathcal{M}^2$ m-morphism of $U^1$ to $U^2$ since the definition of a $\boldsymbol{S}^1$-$\boldsymbol{S}^2$-$\mathcal{M}^1$-$\mathcal{M}^2$ $\mathcal{M}^1$-$\mathcal{M}^2$ m-morphism of $U^1$ to $U^2$ already includes the requirement that $\mathcal{M}^1 \stackrel{\text{full-cat}}{\subseteq} \mathcal{M}^2$.

Let $U^1 \subseteq \boldsymbol{S}^i$ and $U^2 \subseteq \boldsymbol{S}^i$, $i = 1, 2$, be model neighborhoods. A model function $f\colon U^1 \longrightarrow U^2$ is an $\boldsymbol{S}^i$-$\mathcal{M}^i$ m-morphism of $U^1$ to $U^2$ iff it is an $\boldsymbol{S}^i$-$\boldsymbol{S}^i$-$\mathcal{M}^i$-$\mathcal{M}^i$ m-morphism of $U^1$ to $U^2$.

A model function $f\colon \boldsymbol{S}^1 \longrightarrow \boldsymbol{S}^2$ is an $\mathcal{M}^1$-$\mathcal{M}^2$ m-morphism of $\boldsymbol{S}^1$ to $\boldsymbol{S}^2$ iff $\mathcal{M}^1 \stackrel{\text{full-cat}}{\subseteq} \mathcal{M}^2$ and $f$ is a morphism of $\mathcal{M}^2$.[5]

*Remark* 7.14. This definition does not require that $S^i$ be a model neighborhod of $\boldsymbol{S}^i$.

A model function $f\colon \boldsymbol{S}^i \longrightarrow \boldsymbol{S}^i$ is an $\mathcal{M}^i$ m-morphism of $\boldsymbol{S}^i$ iff it is an $\mathcal{M}^i$-$\mathcal{M}^i$ m-morphism of $\boldsymbol{S}^i$ to $\boldsymbol{S}^i$.

**Lemma 7.15** ((Local) m-morphisms). *Let $\mathcal{M}'^i$, $\mathcal{M}^i \stackrel{\text{full-cat}}{\subseteq} \mathcal{M}'^i$, $i = 1, 2$, be model categories, $\boldsymbol{S}^i \stackrel{\text{def}}{=} (S^i, \mathcal{S}^i) \stackrel{\text{Ob}}{\in} \mathcal{M}^i$, $\boldsymbol{S}'^i \stackrel{\text{def}}{=} (S'^i, \mathcal{S}'^i) \stackrel{\text{Ob}}{\in} \mathcal{M}'^i$ and $\boldsymbol{S}^i \stackrel{\text{mod}}{\subseteq} \boldsymbol{S}'^i$.*

1. Let $U^i \subseteq S^i$, $i = 1, 2$, be open. If $f \colon U^1 \longrightarrow U^2$ is a (strict) $S^1$-$S^2$ m-morphism of $U^1$ to $U^2$ then $f \colon U^1 \longrightarrow U^2$ is a (strict) local $S^1$-$S^2$ m-morphism of $U^1$ to $U^2$.

   *Proof.*

   (a) Since $f$ is an $S^1$-$S^2$ m-morphism of $U^1$ to $U^2$ then $S^1 \overset{\text{mod}}{\subseteq} S^2$ and $f$ is a morphism of $S^2$,

   (b) $S^1$ and $S^2$ are model neighborhoods of $S^2$ and $S^1 = f^{-1}[S^2]$ is a model neighborhood of $S^1$.

   (c) Let $u^1 \in U^1$, $u^2 \overset{\text{def}}{=} f(u^1) \in U^2$, Then $U^1$ is a model neighborhood for $u^1$ and $U^2$ is a model neighborhood for $u^2$. $f[U^1] \subseteq U^2$ by hypothesis. $f \colon U^1 \longrightarrow U^2$ is a morphism of $S^2$ as shown above.

   (d) The condition $S^1 \overset{\text{strict-mod}}{\subseteq} S^2$ is the same for local m-morphisms and m-morphisms.

   □

2. Let $U^i \subseteq S^i$ be a model neighborhood of $S^i$, $i = 1, 2$. If $f \colon U^1 \longrightarrow U^2$ is a (strict) $S^1$-$S^2$ local m-morphism of $U^1$ to $U^2$ then $f \colon U^1 \longrightarrow U^2$ is a (strict) $S^1$-$S^2$ m-morphism of $U^1$ to $U^2$.

   *Proof.*

   (a) $S^1 \overset{\text{mod}}{\subseteq} S^2$ by hypothesis. $U^1$ is a model neighborhood of $S^1$ by hypothesis, thus a model neighborhood of $S^2$

   (b) For each $u \in U^1$, let $U_u$ be a model subneighborhood of $U^1$ for $u$ and $V_u$ be a model subneighborhood of $U^2$ for $v \overset{\text{def}}{=} f(u)$ such that $f[U_u] \subseteq V_u$ and $f \colon U_u \longrightarrow V_u$ is a morphism of $S^2$. Then $f \colon S^1 \longrightarrow S^2$ is a morphism of $S^2$ by item 6 of definition 2.1 (Model spaces) on page 18.

   (c) The condition $S^1 \overset{\text{strict-mod}}{\subseteq} S^2$ is the same for local m-morphisms and m-morphisms.

   □

3. Let $U^i \subseteq S^i$ be a model neighborhood, $i = 1, 2$, If $f \colon U^1 \longrightarrow U^2$ is a (strict) $S^1$-$S^2$-$\mathcal{M}^1$-$\mathcal{M}^2$ m-morphism of $U^1$ to $U^2$ then $f \colon U^1 \longrightarrow U^2$ is a (strict) local $S^1$-$S^2$-$\mathcal{M}^1$-$\mathcal{M}^2$ m-morphism of $U^1$ to $U^2$.

   *Proof.*

   (a) Each $U^i$ is a model neighborhood of $S^i$ by hypothesis

(b) Since $f\colon U^1 \longrightarrow U^2$ is an $S^1$-$S^2$-$\mathcal{M}^1$-$\mathcal{M}^2$ m-morphism of $U^1$ to $U^2$, $\mathcal{M}^1 \stackrel{\text{full-cat}}{\subseteq} \mathcal{M}^2$ and $f\colon \text{Mod}(U^1, S^1) \longrightarrow \text{Mod}(U^2, S^2)$ is a morphism of $\mathcal{M}^2$ by definition 7.12 (M-morphisms) on page 32.

(c) Let $u^1 \in U^1$, $u^2 \stackrel{\text{def}}{=} f(u^1) \in U^2$.

   i. Each $U^i$ is a model neighborhood of $S^i$ for $u^i$,
   ii. Each $U^i \subseteq U^i$.
   iii. $f[U^1] \subseteq U^2$.

   □

4. Let $\mathcal{M}^2$ satisfy the restricted sheaf condition and each $U^i \subseteq S^i$ be a model neighborhood, $i = 1, 2$. Every $f\colon U^1 \longrightarrow U^2$ that is a $S^1$-$S^2$-$\mathcal{M}^1$-$\mathcal{M}^2$ local m-morphism of $U^1$ to $U^2$ is a $S^1$-$S^2$-$\mathcal{M}^1$-$\mathcal{M}^2$ m-morphism of $U^1$ to $U^2$.

*Proof.*

(a) Since each $U^i$ is a model neighborhood of $S^i$ by hypothesis, each $\text{Mod}(U^i, S^i) \stackrel{\text{Ob}}{\in} \mathcal{M}^i$ by item 4 of definition 7.7

(b) $f\colon \text{Mod}(U^1, S^1) \longrightarrow \text{Mod}(U^2, S^2)$ is a model function by definition 7.10

(c) Each $S^i \stackrel{\text{Ob}}{\in} \mathcal{M}^i$ by hypothesis

(d) $\mathcal{M}^1 \stackrel{\text{cat}}{\subseteq} \mathcal{M}^2$ by hypothesis

(e) $\mathcal{M}^2$ satisfies the restricted sheaf condition by hypothesis

(f) For every $u^1 \in U^1$ there is a model neighborhood $U^1_{u^1}$ for $u^1$ and a model neighborhood $U^2_{u^1}$ for $u^2 \stackrel{\text{def}}{=} f(u^1)$ such that $f[U^1_{u^1}] \subseteq U^2_{u^1}$ and $f\colon \text{Mod}(U^1_{u^1}, S^i) \longrightarrow \text{Mod}(U^2_{u^1}, S^i)$ is a morphism of $\mathcal{M}^i$.

   i. Each $\text{Mod}(U^1_{u^1}, S^1) \stackrel{\text{full-mod}}{\subseteq} S^1$.
   ii. Each $f\colon \text{Mod}(U^1_{u^1}, S^1) \longrightarrow \text{Mod}(U^2_{u^1}, S^2)$ agrees with $f$ on $U^1_{u^1}$.
   iii. Since $\bigcup_{u^1 \in U^1} U^1_{u^1} = U^1$ and $\bigcup_{u^1 \in U^1} U^2_{u^1} = U^2$, then $f \stackrel{\text{Ar}}{\in} \mathcal{M}^2$ by definition 7.7 (Model categories) on page 30.

   □

5. If $f\colon S^1 \longrightarrow S^2$ is a (strict) local m-morphism of $S^1$ to $S^2$ and $f\colon S^1 \longrightarrow S'^2$ is continuous, then $f$ is a (strict) local m-morphism of $S^1$ to $S'^2$. If $f\colon S^1 \longrightarrow S^2$ is strict and $f\colon S^1 \longrightarrow S'^2$ is open then $f$ is a strict local m-morphism of $S^1$ to $S'^2$.

*Proof.* $f\colon S^1 \longrightarrow S'^2$ is continuous by hypothesis.

Let $u \in S^1$, $U_u$ be a model neighborhood of $\boldsymbol{S}^1$ for $u$ and $V_u$ be a model neighborhood of $\boldsymbol{S}^2$ for $v \stackrel{\text{def}}{=} f(u)$ such that $f[U_u] \subseteq V_u$ and $f\colon U_u \longrightarrow V_u$ is a morphism of $\boldsymbol{S}^2$.

Since $\boldsymbol{S}^2 \stackrel{\text{mod}}{\subseteq} \boldsymbol{S}'^2$ by hypothesis, $f\colon U_u \longrightarrow V_u$ is a morphism of $\boldsymbol{S}'^2$. □

6. If $f\colon S^1 \longrightarrow S^2$ is an m-morphism of $\boldsymbol{S}^1$ to $\boldsymbol{S}^2$, then $f$ is an m-morphism of $\boldsymbol{S}^1$ to $\boldsymbol{S}'^2$.

*Proof.* $\boldsymbol{S}^1 \stackrel{\text{mod}}{\subseteq} \boldsymbol{S}^2$ by definition 7.12 (M-morphisms) on page 32 and $\boldsymbol{S}^2 \stackrel{\text{mod}}{\subseteq} \boldsymbol{S}'^2$ by hypothesis, so $\boldsymbol{S}^1 \stackrel{\text{mod}}{\subseteq} \boldsymbol{S}'^2$.

$f$ is a morphism of $\boldsymbol{S}^1$ by definition 7.12, hence a morphism of $\boldsymbol{S}'^2$. □

7. If $f\colon S^1 \longrightarrow S^2$ is a local $\mathcal{M}^1$-$\mathcal{M}^2$ m-morphism of $\boldsymbol{S}^1$ to $\boldsymbol{S}^2$, then $f$ is a local $\mathcal{M}'^1$-$\mathcal{M}'^2$ m-morphism of $\boldsymbol{S}^1$ to $\boldsymbol{S}'^2$.

*Proof.* Let $u \in S^1$, $U_u$ be a model neighborhood of $\boldsymbol{S}^1$ for $u$ and $V_u$ be a model neighborhood of $\boldsymbol{S}^2$ for $v \stackrel{\text{def}}{=} f(u)$ such that $f[U_u] \subseteq V_u$ and $f\colon U_u \longrightarrow V_u$ is a morphism of $\mathcal{S}^2$.

Since $\mathcal{S}^2 \stackrel{\text{full-cat}}{\subseteq} \mathcal{S}'^2$ by hypothesis, $f\colon U_u \longrightarrow V_u$ is a morphism of $\mathcal{S}'^2$. □

8. If $f\colon S^1 \longrightarrow S^2$ is an $\mathcal{M}^1$-$\mathcal{M}^2$ m-morphism of $\boldsymbol{S}^1$ to $\boldsymbol{S}^2$, then $f$ is an $\mathcal{M}'^1$-$\mathcal{M}'^2$ m-morphism of $\boldsymbol{S}^1$ to $\boldsymbol{S}'^2$.

*Proof.* $\mathcal{S}^1 \stackrel{\text{full-cat}}{\subseteq} \mathcal{S}^2$ by definition 7.12 (M-morphisms) on page 32 and $\mathcal{S}^2 \stackrel{\text{full-cat}}{\subseteq} \mathcal{S}'^2$ by hypothesis, so $\mathcal{S}^1 \stackrel{\text{full-cat}}{\subseteq} \mathcal{S}'^2$.

$f$ is a morphism of $\mathcal{S}^1$ by definition 7.12, hence a morphism of $\mathcal{S}'^2$. □

9. If $\boldsymbol{S}^1 \stackrel{\text{strict-mod}}{\subseteq} \boldsymbol{S}^2$ and $f\colon S^1 \longrightarrow S^2$ is continuous, then $f$ is a strict m-morphism of $\boldsymbol{S}^1_{\text{triv}}$ to $\boldsymbol{S}^2_{\text{triv}}$.

*Proof.*

**$f$ is continuous:**
$S^1 \subseteq S^2$ by hypothesis. $S^1$ is open in $S^2$ by lemma 2.8 (Model subspaces) on page 20 and a model neighborhood of $\boldsymbol{S}^2_{\text{triv}}$ by definition 4.1 (Trivial model spaces and categories) on page 22. $S^2$ is a model neighborhood of $S^2_{\text{triv}}$ by definition 4.1.

Let $U$ be a model meighborhood of $S^1_{\text{triv}}$. $f[U] \subseteq S^2$ and $S^2$ is open in $S^2$, hence a model meighborhood of $S^2_{\text{triv}}$.

Let $V$ be a model neighborhood of $S^2_{\text{triv}}$. $V$ is open in $S^2$, $\text{Id}^{S^1}_{S^2}$ is continuous by hypothesis and $f$ is continuous by hypothesis and thus $f^{-1}[V]$ is open in $S^1$. $\text{Id}^{S^1}_{S^2}$ is open by hypothesis, so $f^{-1}[V]$ open in $S^2$

**$f$ maps model neighborhoods into model neighborhoods:**
$S^2$ is open in $S^2$, hence a model neighborhood of $S^2_{\text{triv}}$. Let $U$ be a model neighborhood of $S^1_{\text{triv}}$. Then $f[U] \subset S^2$.

**Inverse images of model neighborhoods are model neighborhoods:**
Let $V$ be a model neighborhood of $S^2_{\text{triv}}$. $V$ is open in $S^2$ and $f$ is continuous, so $f^{-1}[V]$ is open in $S^1$ and hence a model neighborhood of $S^1_{\text{triv}}$.

**$f$ is strict:** $\boldsymbol{S}^1 \overset{\text{strict-mod}}{\subseteq} \boldsymbol{S}^2$ by hypothesis.

$\square$

10. If $\boldsymbol{S}^1 \overset{\text{mod}}{\subseteq} \boldsymbol{S}^2$ then $\text{Id}^{S^1}_{S^2}$ is a local m-morphism of $\boldsymbol{S}^1$ to $\boldsymbol{S}^2$. If each $S^i$ is a model neighborhood of $\boldsymbol{S}^i$ then $\text{Id}^{S^1}_{S^2}$ is a morphism of $\mathcal{S}^2$.

*Proof.* Let $u \in S^1$ and $U_u$ a model neighborhood of $\boldsymbol{S}^1$ for $u$. Since $\boldsymbol{S}^1 \overset{\text{mod}}{\subseteq} \boldsymbol{S}^2$, $U_u$ is also a model neighborhood of $\boldsymbol{S}^2$ for $u$, and hence $\text{Id}_{U_u} \overset{\text{Ar}}{\in} \mathcal{S}^2$.

Each $S^i$ is a model neighborhood of $\boldsymbol{S}^i$ by hypothesis. Since $\boldsymbol{S}^1 \overset{\text{mod}}{\subseteq} \boldsymbol{S}^2$ by hypothesis. $S^1$ is a model neighborhood of $\boldsymbol{S}^2$. The result follows by item 5 of definition 2.1 (Model spaces) on page 18. $\square$

11. If $\mathcal{M}^1 \overset{\text{cat}}{\subseteq} \mathcal{M}^2$ and $\boldsymbol{S}^1 \overset{\text{mod}}{\subseteq} \boldsymbol{S}^2$ then $\text{Id}^{S^1}_{S^2}$ is a morphism of $\mathcal{M}^2$.

*Proof.* $\boldsymbol{S}^1 \overset{\text{Ob}}{\in} \mathcal{M}^2$ because $\mathcal{M}^1 \overset{\text{cat}}{\subseteq} \mathcal{M}^2$. $\text{Id}_{S^2} \overset{\text{Ar}}{\in} \mathcal{M}^2$. Then $\text{Id}^{S^1}_{S^2} \overset{\text{Ar}}{\in} \mathcal{M}^2$ by item 4 of definition 7.7 (Model categories) on page 30. $\square$

**Corollary 7.16** ((Local) m-morphisms). *Let $\mathcal{M}'^i$, $\mathcal{M}^i \overset{\text{full-cat}}{\subseteq} \mathcal{M}'^i$, $i = 1, 2$, be model categories, $\boldsymbol{S}^i \overset{\text{def}}{=} (S^i, \mathcal{S}^i) \overset{\text{Ob}}{\in} \mathcal{M}^i$, $\boldsymbol{S}'^i \overset{\text{def}}{=} (S'^i, \mathcal{S}'^i) \overset{\text{Ob}}{\in} \mathcal{M}'^i$ and $\boldsymbol{S}^i \overset{\text{mod}}{\subseteq} \boldsymbol{S}'^i$.*

1. If $f\colon S^1 \longrightarrow S^2$ is a (strict) m-morphism of $\boldsymbol{S}^1$ to $\boldsymbol{S}^2$ then $f\colon S^1 \longrightarrow S^2$ is a (strict) local m-morphism of $\boldsymbol{S}^1$ to $\boldsymbol{S}^2$.

   *Proof.* Since $f\colon S^1 \longrightarrow S^2$ is a (strict) $\boldsymbol{S}^1$-$\boldsymbol{S}^2$ m-morphism of $S^1$ to $S^2$ then $f\colon S^1 \longrightarrow S^2$ is a (strict) local $\boldsymbol{S}^1$-$\boldsymbol{S}^2$ m-morphism of $S^1$ to $S^2$. □

2. If $U^i$, $i = 1, 2$, is a model neighborhood of $\boldsymbol{S}^i$ then every function $f\colon U^1 \longrightarrow U^2$ that is a (strict) local $\boldsymbol{S}^i$ m-morphism of $U^1$ to $U^2$ is a (strict) $\boldsymbol{S}^i$ m-morphism of $U^1$ to $U^2$.

   *Proof.* Since $f\colon U^1 \longrightarrow U^2$ is a (strict) local $\boldsymbol{S}^i$-$\boldsymbol{S}^i$ m-morphism of $U^1$ to $U^2$ then it is a (strict) $\boldsymbol{S}^i$-$\boldsymbol{S}^i$ m-morphism of $U^1$ to $U^2$. □

3. If $f\colon S^1 \longrightarrow S^2$ is a (strict) local m-morphism of $\boldsymbol{S}^1$ to $\boldsymbol{S}^2$ and each $S^i$ is a model neighborhood of $\boldsymbol{S}^i$, then $f\colon S^1 \longrightarrow S^2$ is a (strict) m-morphism of $\boldsymbol{S}^1$ to $\boldsymbol{S}^2$.

   *Proof.* Since $f\colon S^1 \longrightarrow S^2$ is a (strict) $\boldsymbol{S}^1$-$\boldsymbol{S}^2$ local m-morphism of $S^1$ to $S^2$ then it is a (strict) $\boldsymbol{S}^1$-$\boldsymbol{S}^2$ m-morphism of $S^1$ to $S^2$. □

4. If $f\colon S^1 \longrightarrow S^2$ is a $\mathcal{M}^1$-$\mathcal{M}^2$ m-morphism of $\boldsymbol{S}^1$ to $\boldsymbol{S}^2$ then $f\colon S^1 \longrightarrow S^2$ is a local $\mathcal{M}^1$-$\mathcal{M}^2$ m-morphism of $\boldsymbol{S}^1$ to $\boldsymbol{S}^2$.

   *Proof.* Since $f$ is a $\boldsymbol{S}^1$-$\boldsymbol{S}^2$-$\mathcal{M}^1$-$\mathcal{M}^2$ m-morphism of $S^1$ to $S^2$ then $f$ is a local $\boldsymbol{S}^1$-$\boldsymbol{S}^2$-$\mathcal{M}^1$-$\mathcal{M}^2$ m-morphism of $S^1$ to $S^2$. □

5. If $\mathcal{M}^i$ satisfies the restricted sheaf condition and $U^j$ is a model neighborhood of $\boldsymbol{S}^i$, $j = 1, 2$, then every function $f\colon U^1 \longrightarrow U^2$ that is a local $\boldsymbol{S}^i$-$\mathcal{M}^i$ m-morphism of $U^1$ to $U^2$ is a $\boldsymbol{S}^i$-$\mathcal{M}^i$ m-morphism of $U^1$ to $U^2$.

   *Proof.* (a) $f\colon \operatorname{Mod}(U^1, \boldsymbol{S}^i) \longrightarrow \operatorname{Mod}(U^2, \boldsymbol{S}^i)$ is a model function by definition 7.10
   
   (b) $\boldsymbol{S}^i \stackrel{\text{Ob}}{\in} \mathcal{M}^i$ by hypothesis
   
   (c) Since each $U^j$ is a model neighborhood of $\boldsymbol{S}^i$ by hypothesis, each $\operatorname{Mod}(U^j, \boldsymbol{S}^i) \stackrel{\text{Ob}}{\in} \mathcal{M}^i$ by item 4 of definition 7.7
   
   (d) For every $u^1 \in U^1$ there is a model neighborhood $U_{u^1}$ for $u^1$ and a model neighborhood $U_{u^2}$ for $u^2 \stackrel{\text{def}}{=} f(u^1)$ such that $f[U_{u^1}] \subseteq U_{u^2}$ and $f\colon \operatorname{Mod}(U_{u^1}, \boldsymbol{S}^i) \longrightarrow \operatorname{Mod}(U_{u^2}, \boldsymbol{S}^i)$ is a morphism of $\mathcal{M}^i$. Each $\operatorname{Mod}(U_{u^1}, \boldsymbol{S}^i) \stackrel{\text{full-mod}}{\subseteq} \boldsymbol{S}^i$. Each $f\colon \operatorname{Mod}(U_{u^1}, \boldsymbol{S}^i) \longrightarrow \operatorname{Mod}(U_{u^2}, \boldsymbol{S}^i)$ agrees with $f$ on $U_{u^1}$. Since $\bigcup_{u^1 \in U^1} U_{u^1} = U^1$ and $\bigcup_{u^1 \in U^1} U^2 = U^2$, $f \stackrel{\text{Ar}}{\in} \mathcal{M}^i$ by definition 7.7. □

6. Let each $S^i$ be a model neighborhood of $\boldsymbol{S}^i$, Then every $f\colon S^1 \longrightarrow S^2$ that is a (strict) $\mathcal{M}^1$-$\mathcal{M}^2$ local m-morphism of $\boldsymbol{S}^1$ to $\boldsymbol{S}^2$ is a (strict) $\mathcal{M}^1$-$\mathcal{M}^2$ m-morphism of $\boldsymbol{S}^1$ to $\boldsymbol{S}^2$.

1. If $f\colon S^1 \longrightarrow S^2$ is a (strict) local m-morphism of $\boldsymbol{S}^1$ to $\boldsymbol{S}^2$ and $\boldsymbol{S}^2 \overset{\text{strict}-\text{mod}}{\subseteq} \boldsymbol{S}'^2$ then $f$ is a (strict) local m-morphism of $\boldsymbol{S}^1$ to $\boldsymbol{S}'^2$.

   *Proof.* $f\colon S^1 \longrightarrow S'^2$ is continuous by definition 7.10 (Local m-morphisms) on page 31. □

2. If $f\colon S^1 \longrightarrow S^2$ is a strict m-morphism of $\boldsymbol{S}^1$ to $\boldsymbol{S}^2$ then $f$ is a strict m-morphism of $\boldsymbol{S}^1_{\text{triv}}$ to $\boldsymbol{S}^2_{\text{triv}}$.

   *Proof.* $f$ is continuous by definition 7.10 (Local m-morphisms) on page 31. □

3. If $\boldsymbol{S}^i$ is a model neighborhood of $\boldsymbol{S}^i$ then $\operatorname{Id}_{\boldsymbol{S}^i}$ is a morphism of $\mathcal{S}^i$.

   *Proof.* $\boldsymbol{S}^i \overset{\text{mod}}{\subseteq} \boldsymbol{S}^i$. □

4. If $f\colon S^1 \longrightarrow S^2$ is a local m-morphism of $\boldsymbol{S}^1$ to $\boldsymbol{S}^2$ and $f^2\colon S^2 \longrightarrow S^3$ is an m-morphism of $\boldsymbol{S}^2$ to $\boldsymbol{S}^3$ them $f^2 \circ f$ is a local m-morphism of $\boldsymbol{S}^1$ to $\boldsymbol{S}^3$.

   *Proof.* Since $f^2\colon S^2 \longrightarrow S^3$ is an m-morphism of $\boldsymbol{S}^2$ to $\boldsymbol{S}^3$ then $f^2$ is a local m-morphism of $\boldsymbol{S}^2$ to $\boldsymbol{S}^3$. □

5. If $f\colon S^1 \longrightarrow S^2$ is an m-morphism of $\boldsymbol{S}^1$ to $\boldsymbol{S}^2$ and $f^2\colon S^2 \longrightarrow S^3$ is a local m-morphism of $\boldsymbol{S}^2$ to $\boldsymbol{S}^3$ them $f^2 \circ f$ is a local m-morphism of $\boldsymbol{S}^1$ to $\boldsymbol{S}^3$.

   *Proof.* Since $f\colon S^1 \longrightarrow S^2$ is an m-morphism of $\boldsymbol{S}^1$ to $\boldsymbol{S}^2$ then $f$ is a local m-morphism of $\boldsymbol{S}^2$ to $\boldsymbol{S}^2$. □

6. If $f\colon S^1 \longrightarrow S^2$ is a local $\mathcal{M}^1$-$\mathcal{M}^2$ m-morphism of $\boldsymbol{S}^1$ to $\boldsymbol{S}^2$ and $f^2\colon S^2 \longrightarrow S^3$ is an $\mathcal{M}^2$-$\mathcal{M}^3$ m-morphism of $\boldsymbol{S}^2$ to $\boldsymbol{S}^3$ them $f^2 \circ f$ is a local $\mathcal{M}^1$-$\mathcal{M}^3$ m-morphism of $\boldsymbol{S}^1$ to $\boldsymbol{S}^3$.

   *Proof.* Since $f^2\colon S^2 \longrightarrow S^3$ is an $\mathcal{M}^2$-$\mathcal{M}^3$ m-morphism of $\boldsymbol{S}^2$ to $\boldsymbol{S}^3$ then $f^2$ is a local $\mathcal{M}^2$-$\mathcal{M}^3$ m-morphism of $\boldsymbol{S}^2$ to $\boldsymbol{S}^3$. □

7. If $f\colon S^1 \longrightarrow S^2$ is an $\mathcal{M}^1$-$\mathcal{M}^2$ m-morphism of $\boldsymbol{S}^1$ to $\boldsymbol{S}^2$ and $f^2\colon S^2 \longrightarrow S^3$ is a local $\mathcal{M}^2$-$\mathcal{M}^3$ m-morphism of $\boldsymbol{S}^2$ to $\boldsymbol{S}^3$ them $f^2 \circ f$ is a local $\mathcal{M}^1$-$\mathcal{M}^3$ m-morphism of $\boldsymbol{S}^1$ to $\boldsymbol{S}^3$.

*Proof.* Since $f\colon S^1 \longrightarrow S^2$ is an $\mathcal{M}^1$-$\mathcal{M}^2$ m-morphism of $\boldsymbol{S}^1$ to $\boldsymbol{S}^2$ then $f$ is a local $\mathcal{M}^1$-$\mathcal{M}^2$ m-morphism of $\boldsymbol{S}^2$ to $\boldsymbol{S}^2$. □

**Lemma 7.17** (Composition of local m-morphisms). *Let $\mathcal{M}^i$, $i = 1, 2, 3$, be model categories, $\boldsymbol{S}^i \stackrel{\text{def}}{=} (S^i, \mathcal{S}^i) \stackrel{\text{Ob}}{\in} \mathcal{M}^i$ and $U^i \subseteq S^i$ be open.*

1. *If $f^i\colon U^i \longrightarrow U^{i+1}$, $i = 1, 2$, is a (strict) local $\boldsymbol{S}^i$-$\boldsymbol{S}^{i+1}$ m-morphism of $U^i$ to $U^{i+1}$ then $f^2 \circ f^1$ a (strict) local $\boldsymbol{S}^1$-$\boldsymbol{S}^3$ m-morphism of $U^i$ to $U^{i+1}$.*

   *Proof.* Since $\boldsymbol{S}^1 \stackrel{\text{mod}}{\subseteq} \boldsymbol{S}^2$ and $\boldsymbol{S}^2 \stackrel{\text{mod}}{\subseteq} \boldsymbol{S}^3$ then $\boldsymbol{S}^1 \stackrel{\text{mod}}{\subseteq} \boldsymbol{S}^3$. Let $u^1 \in U^1, u^2 \stackrel{\text{def}}{=} f^1(u^1)$ and $u^3 \stackrel{\text{def}}{=} f^2(u^2)$. Let $U_{u^i}$ and $V_{u^i}$ be model neighborhoods of $\boldsymbol{S}^i$ and $\boldsymbol{S}^{i+1}$ such that that $f^i[U_{u^i}] \subseteq V_{u^i}$ and $f^i\colon \text{Mod}(U_{u^i}, \boldsymbol{S}^i) \longrightarrow \text{Mod}(V_{u^i}, \boldsymbol{S}^{i+1})$ is a morphism of $\boldsymbol{S}^{i+1}$.

   Let $U'_{u^1} \stackrel{\text{def}}{=} U_{u^1} \cap (f^2 \circ f^1)^{-1}[V_{u^3}]$

   Then $f^2 \circ f^1\colon U'_{u^1} \longrightarrow V_{u^3}$ is a morphism of $\boldsymbol{S}^3$. □

2. *If $f^i\colon S^i \longrightarrow S^{i+1}$, $i = 1, 2$, is a (strict) local $\mathcal{M}^i$-$\mathcal{M}^{i+1}$ m-morphism of $\boldsymbol{S}^i$ to $\boldsymbol{S}^{i+1}$ then $f^2 \circ f^1$ a (strict) local $\mathcal{M}^i$-$\mathcal{M}^{i+2}$ m-morphism of $\boldsymbol{S}^1$ to $\boldsymbol{S}^3$.*

   *Proof.* Since $\mathcal{M}^1 \stackrel{\text{fullcat}}{\subseteq} \mathcal{M}^2$ and $\boldsymbol{M}^2 \stackrel{\text{fullcat}}{\subseteq} \mathcal{M}^3$ then $\boldsymbol{M}^1 \stackrel{\text{fullcat}}{\subseteq} \mathcal{M}^3$. Let $u^1 \in \boldsymbol{S}^1$, $u^2 \stackrel{\text{def}}{=} f^1(u)$ and $u^3 \stackrel{\text{def}}{=} f^2(u^2)$. Let $U_{u^i}$ and $V_{u^i}$ be model neighborhoods of $\boldsymbol{S}^i$ and $\boldsymbol{S}^{i+1}$ such that $f^i\colon U_{u^i} \longrightarrow V_{u^i}$ is a morphism of $\mathcal{M}^{i+1}$. Let $U'_{u^1} \stackrel{\text{def}}{=} U_{u^1} \cap (f^2 \circ f^2)^{-1}[V_{u^3}]$

   Then $f^2 \circ f^2\colon U'_{u^1} \longrightarrow V_{u^3}$ is a morphism of $\mathcal{M}^3$ by item 5 of definition 7.7 (Model categories) on page 30. □

*If each $f^i\colon \boldsymbol{S}^i \longrightarrow \boldsymbol{S}^{i+1}$, $i = 1, 2$, is a (strict) local m-morphism of $\boldsymbol{S}^i$ to $\boldsymbol{S}^{i+1}$, then $f^2 \circ f^1\colon S^1 \longrightarrow S^3$ is a (strict) local m-morphism of $\boldsymbol{S}^1\colon S^1 \longrightarrow S^3$. If each $f^i$ is strict then $f^2 \circ f^1$ is strict.*

*Proof.* Since $\boldsymbol{M}^1 \stackrel{\text{mod}}{\subseteq} \boldsymbol{M}^2$ and $\boldsymbol{M}^2 \stackrel{\text{mod}}{\subseteq} \boldsymbol{M}^3$, $\boldsymbol{M}^1 \stackrel{\text{mod}}{\subseteq} \boldsymbol{M}^3$. Let $u \in S^i$, $v \stackrel{\text{def}}{=} f^1(u)$ and $w \stackrel{\text{def}}{=} f^2(v)$. There exist a model neighborhood $U_u$ for $u$, model neighborhoods $V_u$, $V'_u$ for $v$ and a model neighborhood $W_v$ of $w$ such that $f_1[U_u] \subseteq V'_u, f_2[V_v] \subseteq W_v$, $f_1\colon U_u \longrightarrow V'_u$ is a morphism of $\boldsymbol{S}^2$ and $f_2\colon V_v \longrightarrow W_v$ is a morphism of $\mathcal{S}^3$. Then $\hat{V}_u \stackrel{\text{def}}{=} V_v \cap V'_u \neq \emptyset$, $\hat{V}_u$ is a model neighborhood of $v$ and $\hat{U}_u \stackrel{\text{def}}{=} f_1^{-1}[\hat{V}_u]$ is a model neighborhood for $u$. $f_1\colon \hat{U}_u \longrightarrow \hat{V}_u$ and $f_2\colon \hat{V}_u \longrightarrow W_v$ are morphisms of $\mathcal{S}^3$ by item 4 of definition 2.1 (Model spaces) on page 18 and thus $f_2 \circ f_1\colon \hat{U}_u \longrightarrow W_v$ is a morphism of $\boldsymbol{S}^3$.

If $\boldsymbol{M}^1 \stackrel{\text{strict-mod}}{\subseteq} \boldsymbol{M}^2$ and $\boldsymbol{M}^2 \stackrel{\text{strict-mod}}{\subseteq} \boldsymbol{M}^3$, then $\boldsymbol{M}^1 \stackrel{\text{strict-mod}}{\subseteq} \boldsymbol{M}^3$. □

If each $f^i\colon \boldsymbol{S}^i \longrightarrow \boldsymbol{S}^{i+1}$ is a (strict) local $\mathcal{M}^i$-$\mathcal{M}^{i+1}$ m-morphism of $\boldsymbol{S}^i$ to $\boldsymbol{S}^{i+1}$ then $f^2 \circ f^1\colon \boldsymbol{S}^1 \longrightarrow \boldsymbol{S}^3$ is a (strict) local $\mathcal{M}^1$-$\mathcal{M}^3$ m-morphism of $\boldsymbol{S}^1$ to $\boldsymbol{S}^3$.

*Proof.* Since $\mathcal{M}^1 \overset{\text{full-cat}}{\subseteq} \mathcal{M}^2$ and $\mathcal{M}^2 \overset{\text{full-cat}}{\subseteq} \mathcal{M}^3$, $\mathcal{M}^1 \overset{\text{full-cat}}{\subseteq} \mathcal{M}^3$. Let $u \in \boldsymbol{S}^1$, $v \overset{\text{def}}{=} f^1(u)$ and $w \overset{\text{def}}{=} f^2(v)$. There exist a model neighborhood $U_u$ for $u$, model neighborhoods $V_v, V'_u$ for $v$ and a model neighborhood $W_v$ of $w$ such that $f^1[U_u] \subseteq V'_u$, $f^1\colon U_u \longrightarrow V'_u$ is a morphism of $\mathcal{M}^2$, $f^2[V_v] \subseteq W_v$ and $f^2\colon V_v \longrightarrow W_v$ is a morphism of $\mathcal{M}^3$. Then $\hat{V}_u \overset{\text{def}}{=} V_v \cap V'_u \neq \emptyset$, $\hat{V}_u$ is a model neighborhood of $v$ and $\hat{U}_u \overset{\text{def}}{=} f_1^{i-1}[\hat{V}_u]$ is a model neighborhood for $u$. $f^1\colon \hat{U}_u \longrightarrow \hat{V}_u$ and $f^2\colon \hat{V}_u \longrightarrow W_v$ are morphisms of $\mathcal{M}^3$ by item 4 of definition 7.7 (Model categories) on page 30 and thus $f^2 \circ f^1\colon \hat{U}_u \longrightarrow W_v$ is a morphism of $\mathcal{M}^3$. □

**Corollary 7.18** ([Composition of local m-morphisms). *Let $\mathcal{M}^i$, $i = 1, 2, 3$, be model categories, $\boldsymbol{S}^i \overset{\text{def}}{=} (S^i, \mathcal{S}^i) \overset{\text{Ob}}{\in} \mathcal{M}^i$.*

1. *Let $U^j \subseteq S^i$, $j = 1, 2, 3$, be open. If $f^j\colon U^j \longrightarrow U^{j+1}$, $j = 1, 2$, is a (strict) local $\boldsymbol{S}^i$ m-morphism of $U^j$ to $U^{j+1}$ then $f^2 \circ f^1$ a (strict) local $\boldsymbol{S}^i$ m-morphism of $U^1$ to $U^3$.*

   *Proof.* Since each $f^j\colon U^j \longrightarrow U^{j+1}$ is a (strict) local $\boldsymbol{S}^i$-$\boldsymbol{S}^i$ m-morphism of $U^j$ to $U^{j+1}$ then $f^2 \circ f^1$ a (strict) local $\boldsymbol{S}^i$-$\boldsymbol{S}^i$ m-morphism of $U^j$ to $U^{j+1}$. □

2. *If $f^j\colon S^j \longrightarrow S^{j+1}$, $j = 1, 2$, is a (strict) local $\boldsymbol{S}^i$ m-morphism of $\boldsymbol{S}^j$ to $\boldsymbol{S}^{j+1}$ then $f^2 \circ f^1$ a (strict) local $\boldsymbol{S}^i$ m-morphism of $\boldsymbol{S}^1$ to $\boldsymbol{S}^3$.*

   *Proof.* Since each $f^j\colon S^j \longrightarrow S^{j+1}$ is a (strict) local $\boldsymbol{S}^i$-$\boldsymbol{S}^i$ m-morphism of $S^j$ to $S^{j+1}$ then $f^2 \circ f^1$ a (strict) local $\boldsymbol{S}^i$-$\boldsymbol{S}^i$ m-morphism of $S^j$ to $S^{j+1}$. □

**Lemma 7.19** (Composition of m-morphisms). *Let $\mathcal{M}^i$, $i = 1, 2, 3$, be model categories and $\boldsymbol{S}^i \overset{\text{def}}{=} (S^i, \mathcal{S}^i) \overset{\text{Ob}}{\in} \mathcal{M}^i$.*

1. *If $f^i\colon \boldsymbol{S}^i \longrightarrow \boldsymbol{S}^{i+1}$, $i = 1, 2$, is an m-morphism of $\boldsymbol{S}^i$ to $\boldsymbol{S}^{i+1}$ then $f^2 \circ f^1\colon \boldsymbol{S}^1 \longrightarrow \boldsymbol{S}^3$ is an m-morphism of $\boldsymbol{S}^1\colon \boldsymbol{S}^1 \longrightarrow \boldsymbol{S}^3$.*
   *If each $f^i$ is strict then $f^2 \circ f^1$ is strict.*

   *Proof.* Since each $f^i$ is a morphism of $\boldsymbol{S}^{i+1}$ by definition 7.12 (M-morphisms) on page 32 and $\boldsymbol{S}^1 \overset{\text{mod}}{\subseteq} \boldsymbol{S}^2 \overset{\text{mod}}{\subseteq} \boldsymbol{S}^3$, then each $f^i$ is a morphism of $\boldsymbol{S}^3$ by definition 7.12 (M-morphisms) on page 32, and thus $f^2 \circ f^1\colon \boldsymbol{S}^1 \longrightarrow \boldsymbol{S}^3$ is a morphism of $\boldsymbol{S}^3$.
   If $\boldsymbol{S}^1 \overset{\text{strict-mod}}{\subseteq} \boldsymbol{S}^2$ and $\boldsymbol{S}^2 \overset{\text{strict-mod}}{\subseteq} \boldsymbol{S}^3$, then $\boldsymbol{S}^1 \overset{\text{strict-mod}}{\subseteq} \boldsymbol{S}^3$. □

2. *If $f^i\colon S^i \longrightarrow S^{i+1}$, $i = 1, 2$, is an $\mathcal{M}^i$-$\mathcal{M}^{i+1}$ m-morphism of $S^i$ to $S^{i+1}$ then $f^2 \circ f^1$ an $\mathcal{M}^i$-$\mathcal{M}^{i+2}$ m-morphism of $S^1$ to $S^3$.*

   *Proof.* Since $\mathcal{M}^1 \overset{\text{fullcat}}{\subseteq} \mathcal{M}^2$ and $\mathbf{M}^2 \overset{\text{fullcat}}{\subseteq} \mathcal{M}^3$ then $\mathbf{M}^1 \overset{\text{fullcat}}{\subseteq} \mathcal{M}^3$. Since each $f^i\colon S^i \longrightarrow S^{i+1}$ is a morphism of $\mathcal{S}^3$, then so is $f^2 \circ f^1$. □

3. *If each $f^i\colon \mathbf{S}^i \longrightarrow \mathbf{S}^{i+1}$ is an $\mathcal{M}^i$-$\mathcal{M}^{i+1}$ m-morphism of $\mathbf{S}^i$ to $\mathbf{S}^{i+1}$ then $f^2 \circ f^1\colon \mathbf{S}^1 \longrightarrow \mathbf{S}^3$ is an $\mathcal{M}^1$-$\mathcal{M}^3$ m-morphism of $\mathbf{S}^1$ to $\mathbf{S}^3$.*

   *Proof.* Since $\mathcal{M}^1 \overset{\text{full–cat}}{\subseteq} \mathcal{M}^2$ and $\mathcal{M}^2 \overset{\text{full–cat}}{\subseteq} \mathcal{M}^3$, $\mathcal{M}^1 \overset{\text{full–cat}}{\subseteq} \mathcal{M}^3$. $f^i$ is a morphism of $\mathcal{S}^{i+1}$ by definition 7.12 (M-morphisms) on page 32. Since $\mathcal{M}^2 \overset{\text{fullcat}}{\subseteq} \mathcal{M}^3$, $f^1$ is a morphism of $\mathcal{S}^3$. and thus $f^2 \circ f^1\colon S^1 \longrightarrow S^3$ is a morphism of $\mathcal{S}^3$. □

# 8 Spaces and proper functions

Several of the following definitions involve spaces of a restrictive character and specific types of mappings among them. Except where otherwise qualified, the word *space* will have this restricted meaning.

**Definition 8.1** (Spaces). A space is a topological space, a model space or either with an additional associated structure.

**Definition 8.2** (Proper functions). Let $S^1$, $S^2$ be spaces and $f\colon S^1 \longrightarrow S^2$ a continuous function; $f$ need not preserve any associated algebraic structure[6]. $f$ is a proper[7] function iff

1. $S^1$ and $S^2$ are both Truthspace or both **Truthspace**, and $f(\text{True}) = \text{True}$.

2. $S^1$, $S^2$ are topological spaces other than Truthspace.

3. $S^1$ is a topological space other than Truthspace, $S^2$ is a model space other than **Truthspace** and the images of open sets are contained in model neighborhoods.

4. $S^1$ is a model space other than **Truthspace**, $S^2$ is a topological space other than Truthspace and the inverse images of open sets are model neighborhoods.

---

[6] However, in practice commutation relations will often enforce the preservation of algebraic structures.

[7] In the spirit of [Kelley, 1955, footnote, p. 112]

> This nomenclature is an excellent example of the time-honored custom of referring to a problem we cannot handle as abnormal, irregular, improper, degenerate, inadmissible, and otherwise undesirable.

5. $S^1$, $S^2$ are both model spaces other than **Truthspace** and $f$ is a model function.

**Definition 8.3** (Singleton categories and sequences of singleton categories). Let $S$ be a topological space. Then the singleton category of $S$, abbreviated $\underset{Sing}{S}$, is the category whose sole object is $S$ and whose sole morphism is the identity morphism of $S$ to itself.

Let $\boldsymbol{S}$ be a model space. Then the singleton category of $\boldsymbol{S}$, abbreviated $\underset{Sing}{\boldsymbol{S}}$, is the category whose sole object is $\boldsymbol{S}$ and whose sole morphism is the identity morphism of $\boldsymbol{S}$ to itself.

Let $\boldsymbol{S} \stackrel{\text{def}}{=} (S^\alpha, \alpha \prec A)$ be a sequence of spaces. Then the singleton category sequence of $\boldsymbol{S}$, abbreviated $\underset{Sing}{\boldsymbol{S}}$, is $(\underset{Sing}{S^\alpha}, \alpha \prec A)$.

Let $\boldsymbol{S}$ be a set of spaces. Then the singleton category of $\boldsymbol{S}$, abbreviated $\underset{Sing}{\boldsymbol{S}}$, is $\bigcup_{S \in \boldsymbol{S}} \underset{Sing}{S}$.

*Remark* 8.4. By abuse of language the notation $\underset{Sing}{\boldsymbol{S}}$ will be used to name sequences of categories constructed with this and similar functions.

# Part VI
# M-charts and m-atlases

The literature defines fiber bundles and manifolds using the language of charts, atlases and transition functions; it has multiple equivalent definitions. Some authors start with topological spaces and define atlases over them. Some start with abstract sets and define atlases over them, deriving the topology from the atlas. Some start with an indexed set of open patches in a coordinate space and transition functions among them satisfying a cocycle (compatibility) constraint, and then derive the total space as a quotient space of the disjoin union of the patches. Some authors use maximal atlases while others use equivalence classes of atlases.

This paper uses the first approach, explicitly making the relevant spaces topological spaces, but modifies the definitions of charts in order to make them fit more natuarally into the context of local coordinate spaces. It adds a prefix, e.g., $C^k$, M, in order to avoid confusion with the conventional definitions. An m-atlas based on a topological spaces is closer to the conventional definitions of an atlas for a manifold while an m-atlas based on a model space is suitable for defining both manifolds and fiber bundles. Although simple manifolds and fiber bundles could both be defined directly in terms of maximal m-atlases, this paper has a different perspective, and uses the m-atlases as part of the more general Local Coordinate Space (LCS), presented in [Metz, 2018].

Sections 9 to 10 define m-charts and m-atlases. Maximal m-atlases are essentially the same as manifolds.

Section 11 defines morphisms between m-atlasses in a fashion tailored to use in defining local coordinate spaces. It constructs categories of M-atlases and functors, and proves some basic results. It also defines a related concept of near morphisms, although the restriction needed to use them as morphisms in a category of m-atlases requires them to be morphisms.

Section 13 defines morphisms between m-atlasses in a fashion similar to that conventionally used for $C^k$ maps between differentiable manifolds, in order to provide context. It constructs categories of M-atlases and functors, and proves some basic results.

## 9 M-charts

**Definition 9.1** (M-charts). Let $\boldsymbol{E} \stackrel{\text{def}}{=} (E, \mathscr{E})$ and $\boldsymbol{C} \stackrel{\text{def}}{=} (C, \mathscr{C})$ be model spaces. An m-chart $(U, V, \phi)$ of $\boldsymbol{E}$ in the coordinate space $\mathscr{C}$ consists of

1. A nonvoid model neighborhood $U \stackrel{\text{Ob}}{\in} \mathscr{E}$, known as a coordinate patch

2. A model neighborhood $V \stackrel{\text{Ob}}{\in} \mathscr{C}$

3. A model homeomorphism $\phi\colon U \rightarrowtail\stackrel{\widetilde{=}}{\twoheadrightarrow} V$, known as a coordinate function

*Remark* 9.2. I consider it clearer to explicate the range, rather than the conventional usage of specifying only the domain and function or the minimalist usage of specifying only the function.

For evey $x \in U$, $(U, V, \phi)$ is an m-chart of $\boldsymbol{E}$ in the coordinate space $\mathscr{C}$ at $x$. $\boldsymbol{E}$ is the total model space or total space for the chart and $\boldsymbol{C}$ is the coordinate model space or coordinate space for the chart.

Let $E$ be a topological space and $\boldsymbol{C} \stackrel{\text{def}}{=} (C, \mathscr{C})$ be a model space. $(U, V, \phi)$ is an m-chart of $E$ in the coordinate space $\mathscr{C}$ iff it is an m-chart of $\underset{\text{triv}}{E}$ in the coordinate space $\mathscr{C}$.

For evey $x \in U$, $(U, V, \phi)$ is an m-chart of $E$ in the coordinate space $\mathscr{C}$ at $x$. $E$ is the total space for the chart and $\boldsymbol{C}$ is the coordinate model space or coordinate space for the chart.

**Lemma 9.3** (M-charts).

1. Let $\boldsymbol{E} \stackrel{\text{def}}{=} (E, \mathscr{E})$ and $\boldsymbol{C} \stackrel{\text{def}}{=} (C, \mathscr{C})$ be model spaces, and $(U, V, \phi)$ be an m-chart of $\boldsymbol{E}$ in the coordinate space $\boldsymbol{C}$ such that every open subset of $V$ is a model neighborhood of $\boldsymbol{C}$. Then $(U, V, \phi)$ is an m-chart of $E$ in the coordinate space $\boldsymbol{C}$.

   *Proof.* $U$, $V$ and $\phi$ satisfy the conditions of the definition:

   (a) $U$ is a model neighborhood of $\boldsymbol{E}$, hence open and a model neighborhood of $\underset{\text{triv}}{E}$.

(b) $V$ is a model neighborhood of $\boldsymbol{C}$,

(c) Since $\phi$ is a model homeomorphism from $\text{Mod}(U, \boldsymbol{E})$ to $\text{Mod}(V, \boldsymbol{C})$, it is a homeomorphism. If $V'$ is a model subneighborhood of $V$, then $\phi^{-1}[V']$ is open and hence a model neighborhood of $\underset{\text{triv}}{E}$. If $U' \subseteq U$ is a model neighborhood of $\underset{\text{triv}}{E}$, then it is open and $\phi[U'] \subseteq V$ is open. Every open subset of $V$ is a model neighborhood of $\boldsymbol{C}$ by hypothesis.

$\square$

2. Let $(U, V, \phi)$ be an m-chart of $E$ in the coordinate space $\boldsymbol{C}$. $(U, V, \phi)$ is an m-chart of $\boldsymbol{E}$ in the coordinate space $\boldsymbol{C}$ iff every open subset of $U$ is a model neighborhood of $\boldsymbol{E}$.

*Proof.* If $(U, V, \phi)$ is an m-chart of $E$ in the coordinate space $\boldsymbol{C}$ and also an m-chart of $\boldsymbol{E}$ in the coordinate space $\boldsymbol{C}$ then let $U' \subseteq U$ be open. Since $(U, V, \phi)$ is an m-chart of $E$ in the coordinate space $\boldsymbol{C}$, $\phi^{-1}$ is a model function from $\text{Mod}(V, \boldsymbol{C})$ to $\underset{\text{triv}}{E}$ and $\phi[U']$ is a model neighborhood of $\boldsymbol{C}$. Since $(U, V, \phi)$ is an m-chart of $\boldsymbol{E}$ in the coordinate space $\boldsymbol{C}$, $\phi$ is a model function from $\text{Mod}(U, \boldsymbol{E})$ to $\boldsymbol{C}$ and thus $\phi^{-1}[\phi[U']]$ is a model neighborhood of $\boldsymbol{E}$.

If $(U, V, \phi)$ is an m-chart of $E$ in the coordinate space $\boldsymbol{C}$ and every open $U' \subseteq U$ is a model neighborhood of $\boldsymbol{E}$ then

(a) Since $U \subseteq U$ is open, $U$ is a model neighborhood of $\boldsymbol{E}$ by hypothesis.

(b) $V$ is a model neighborhood of $\boldsymbol{C}$ by hypothesis.

(c) If $V' \subseteq V$ is a model neighborhood of $\boldsymbol{C}$ then $\phi^{-1}[V']$ is open and hence a model neighborhod of $\boldsymbol{E}$ by hypothesis.

$\phi^{-1}$ is a model function from $\text{Mod}(U, \underset{\text{triv}}{E})$ to $\text{Mod}(V, \boldsymbol{C})$. if $U' \subseteq U$ is a model neighborhood of $\boldsymbol{E}$ then $U'$ is open and hence a model neighborhood of $\underset{\text{triv}}{E}$. Then $\phi[U']$ is a model neighborhood of $\boldsymbol{C}$.

$\square$

**Corollary 9.4** (M-charts). *Let $\boldsymbol{E} \stackrel{\text{def}}{=} (E, \mathscr{E})$ and $\boldsymbol{C} \stackrel{\text{def}}{=} (C, \mathscr{C})$ be model spaces, $\boldsymbol{C}$ be fine grained and $(U, V, \phi)$ be an m-chart of $E$ in the coordinate space $\boldsymbol{C}$. Then $(U, V, \phi)$ is an m-chart of $\boldsymbol{E}$ in the coordinate space $\boldsymbol{C}$.*

*Proof.* Since $V$ is open and $\boldsymbol{C}$ is fine grained, every open subset of $V$ is a model neighborhood of $\boldsymbol{C}$. The result follows by item 1 of lemma 9.3 (M-charts) on page 43. $\square$

**Definition 9.5** (Subcharts). Let $(U, V, \phi)$ be an m-chart of $\boldsymbol{E} \stackrel{\text{def}}{=} (E, \mathscr{E})$ in the coordinate space $(C, \mathscr{C})$ and $U' \subseteq U$ a nonvoid model neighborhood of $(E, \mathscr{E})$. Then $(U', V', \phi') \stackrel{\text{def}}{=} (U', \phi[U'], \phi\!\restriction_{U', V'})$ is a subchart of $(U, V, \phi)$.

Let $(U, V, \phi)$ be an m-chart of $E$ in the coordinate space $(C, \mathscr{C})$ and $U' \subseteq U$ an open set of $E$. Then $(U', V', \phi') \stackrel{\text{def}}{=} (U', \phi[U'], \phi\!\restriction_{U',V'})$ is a subchart of $(U, V, \phi)$.

$(U', V', \phi')$ is a subchart of $(U, V, \phi)$ at $x$ iff $(U', V', \phi')$ is a subchart of $(U, V, \phi)$ and $x \in U'$.

**Lemma 9.6** (Subcharts). *Let $\boldsymbol{E} \stackrel{\text{def}}{=} (E, \mathscr{E})$ and $\boldsymbol{C} \stackrel{\text{def}}{=} (C, \mathscr{C})$ be model spaces.*
*Let $(U, V, \phi)$ be an m-chart of $\boldsymbol{E}$ in the coordinate space $\boldsymbol{C}$ and $(U', V', \phi')$ a subchart of $(U, V, \phi)$. Then $(U', V', \phi')$ is an m-chart of $\boldsymbol{E}$ in the coordinate space $\boldsymbol{C}$.*

*Proof.* $U'$, $V'$ and $\phi'$ satisfy the conditions of the definition of an m-chart:

1. $U'$ is a model neighborhood of $\boldsymbol{E}$ by the definition of subchart.

2. Since $U'$ is a model neighborhood and $\phi$ is a model homeomorphism, $V' = \phi[U']$ is a model neighborhood.

3. Since $\phi$ and $\phi^{-1}$ are model functions, so are $\phi' = \phi\!\restriction_{U',V'}$ and $\phi'^{-1} = \phi^{-1}\!\restriction_{V',U'}$. Thus $\phi'$ is a model homeomorphism.

$\square$

*Let $(U, V, \phi)$ be an m-chart of $E$ in the coordinate space $\boldsymbol{C}$ and $(U', V', \phi')$ a subchart of $(U, V, \phi)$. Then $(U', V', \phi')$ is an m-chart of $E$ in the coordinate space $\boldsymbol{C}$.*

*Proof.* $U'$, $V'$ and $\phi'$ satisfy the conditions of the definition of an m-chart:

1. $U'$ is open by the definition of subchart.

2. Since $\phi$ is a homeomorphism and $U'$ is open, $V' = \phi[U']$ is open.

3. Since $\phi$ is a homeomorphism that maps open sets into model neighborhoods, $\phi\!\restriction_{U',V'}$ is a homeomorphism that maps open sets into model neighborhoods.

$\square$

**Definition 9.7** (M-compatibility). Let $\boldsymbol{E} \stackrel{\text{def}}{=} (E, \mathscr{E})$ and $\boldsymbol{C} \stackrel{\text{def}}{=} (C, \mathscr{C})$ be model spaces and $(u^i, v^i, \phi^i)$, $i = 1, 2$, be an m-chart of $\boldsymbol{E}$ in the coordinate space $\boldsymbol{C}$. Then $(u^1, v^1, \phi^1)$ is compatible with $(u^2, v^2, \phi^2)$ iff either

1. $U^1$ and $U^2$ are disjoint

2. The transition function $t = \phi^2 \circ \phi^{1-1}\!\restriction_{\phi[U^1 \cap U^2]}$ is an isomorphism of $\mathscr{C}$

Let $E$ be a topological space, $\boldsymbol{C} \stackrel{\text{def}}{=} (C, \mathscr{C})$ be a model space and $(u^i, v^i, \phi^i)$, $i = 1, 2$, be an m-chart of $E$ in the coordinate space $\boldsymbol{C}$. Then $(u^1, v^1, \phi^1)$ is compatible with $(u^2, v^2, \phi^2)$ iff either

1. $U^1$ and $U^2$ are disjoint

2. The transition function $t = \phi^2 \circ \phi^{1-1}\!\restriction_{\phi[U^1 \cap U^2]}$ is an isomorphism of $\mathscr{C}$

**Lemma 9.8** (Symmetry of m-compatibility). *Let $(U, V, \phi)$ and $(U', V', \phi')$ be m-charts of $(E, \mathcal{E})$ in the coordinate space $(C, \mathcal{C})$. Then $(U, V, \phi)$ is m-compatible with $(U', V', \phi')$ in the coordinate space $(C, \mathcal{C})$ iff $(U', V', \phi')$ is m-compatible with $(U, V, \phi)$ in the coordinate space $(C, \mathcal{C})$.*

*Let $(U, V, \phi)$ and $(U', V', \phi')$ be m-charts of $E$ in the coordinate space $(C, \mathcal{C})$. Then $(U, V, \phi)$ is m-compatible with $(U', V', \phi')$ iff $(U', V', \phi')$ is m-compatible with $(U, V, \phi)$.*

*Proof.* The same proof applies to both cases. It suffices to prove the implication in only one direction.

1. $U \cap U' = U' \cap U$.

2. Since the transition function $t = \phi' \circ \phi^{-1} \restriction_{\phi[U \cap U']}$ is an isomorphism of $\mathcal{C}$, so is $t^{-1} = \phi \circ \phi'^{-1} \restriction_{\phi'[U \cap U']}$.

$\square$

**Lemma 9.9** (M-compatibility of subcharts). *Let $(U_i, V_i, \phi_i)$, $i = 1, 2$, be m-charts of $(E, \mathcal{E})$ in the coordinate space $(C, \mathcal{C})$, $(U'_i, V'_i, \phi'_i)$ be subcharts and $(U_1, V_1, \phi_1)$ be m-compatible with $(U_2, V_2, \phi_2)$. Then $(U'_1, V'_1, \phi'_1)$ is m-compatible with $(U'_2, V'_2, \phi'_2)$.*

*Proof.* Since subcharts are charts, $U'_i$ and $V'_i$ are model neighborhoods. If $U_1 \cap U_2 = \emptyset$ then $U'_1 \cap U'_2 = \emptyset$. If $U'_1 \cap U'_2 = \emptyset$ then $(U'_1, V'_1, \phi'_1)$ is m-compatible with $(U'_2, V'_2, \phi'_2)$. Otherwise, the transition function $t^1_2 \stackrel{\text{def}}{=} \phi_2 \circ \phi_1^{-1} \restriction_{\phi_1[U_1 \cap U_2]}$ is a model homeomorphism and hence $t^1_2 \restriction_{\phi_1[U'_1 \cap U'_2]} : \phi_1[U'_1 \cap U'_2] \xrightarrow{\tilde{=}} \phi_2[U'_1 \cap U'_2]$ is a model homeomorphism. $\square$

*Let $(U_i, V_i, \phi_i)$, $i = 1, 2$, be m-charts of $E$ in the coordinate space $(C, \mathcal{C})$ and $(U'_i, V'_i, \phi'_i)$ subcharts. Then $(U'_i, V'_i, \phi'_i)$ is m-compatible with $(U_i, V_i, \phi_i)$.*

*Proof.* Since subcharts are charts, $U'_i$ and $V'_i$ are open. If $U_1 \cap U_2 = \emptyset$ then $U'_1 \cap U'_2 = \emptyset$. If $U'_1 \cap U'_2 = \emptyset$ then $(U'_1, V'_1, \phi'_1)$ is m-compatible with $(U'_2, V'_2, \phi'_2)$. Otherwise, the transition function $t^1_2 \stackrel{\text{def}}{=} \phi_2 \circ \phi_1^{-1} \restriction_{\phi_1[U_1 \cap U_2]}$ is a model homeomorphism and hence $t^1_2 \restriction_{\phi_1[U'_1 \cap U'_2]} : \phi_1[U'_1 \cap U'_2] \xrightarrow{\tilde{=}} \phi_2[U'_1 \cap U'_2]$ is a model homeomorphism. $\square$

**Corollary 9.10** (M-compatibility with subcharts). *Let $(U, V, \phi)$ be an m-chart of $(E, \mathcal{E})$ in the coordinate space $(C, \mathcal{C})$ and $(U', V', \phi')$ a subchart. Then $(U', V', \phi')$ is m-compatible with $(U, V, \phi)$.*

*Proof.* $(U, V, \phi)$ is m-compatible with itself and is a subchart of itself, $\square$

*Let $(U, V, \phi)$ be an m-chart of $E$ in the coordinate space $(C, \mathcal{C})$ and $(U', V', \phi')$ a subchart. Then $(U', V', \phi')$ is m-compatible with $(U, V, \phi)$.*

*Proof.* $(U, V, \phi)$ is m-compatible with itself and is a subchart of itself, $\square$

**Definition 9.11** (Covering by m-charts)**.** Let $A$ be a set of charts of the topological space $E$ in the coordinate space $\boldsymbol{C} = (C, \mathscr{C})$. $A$ covers $E$ iff $\pi_1[A]$ covers $E$.

Let $A$ be a set of charts of the model space $\boldsymbol{E} = (E, \mathscr{E})$ in the coordinate space $\boldsymbol{C} = (C, \mathscr{C})$. $A$ covers $\boldsymbol{E}$ iff $\pi_1[A]$ covers $E$.

## 10 M-atlases

A set of charts can be atlases for different coordinate model spaces even if it is for the same total model space. In order to aggregate atlases into categories, there must be a way to distinguish them. Including the two[8] spaces in the definitions of the categories serves the purpose.

**Definition 10.1** (M-atlases)**.** Let $A$ be a set of mutually m-compatible m-charts of $\boldsymbol{E} = (E, \mathscr{E})$ in the coordinate space $\boldsymbol{C} = (C, \mathscr{C})$. Then $A$ is an m-atlas of $\boldsymbol{E}$ in the coordinate space $\boldsymbol{C}$, abbreviated $\text{isAtl}_{\text{Ob}}(A, \boldsymbol{E}, \boldsymbol{C})$, iff $A$ covers $E$. $A$ is a full atlas of $\boldsymbol{E}$ in the coordinate space $\boldsymbol{C}$, abbreviated $\text{isAtl}_{\text{Ob}}_{\text{full}}(A, E, \boldsymbol{C})$, iff $A$ covers $E$ and $\pi_2[A]$ covers $C$. The triple $(A, \boldsymbol{E}, \boldsymbol{C})$ refers to $A$ considered as an m-atlas of $\boldsymbol{E}$ in the coordinate space $\boldsymbol{C}$. $\boldsymbol{E}$ is the total model space or total space for the atlas and $\boldsymbol{C}$ is the coordinate model space or coordinate space for the atlas.

Let $E$ be a topological space, $\boldsymbol{C} = (C, \mathscr{C})$ a model space and $A$ be a set of mutually m-compatible m-charts of $E$ in the coordinate space $\boldsymbol{C}$. Then $A$ is an m-atlas of $E$ in the coordinate space $\boldsymbol{C}$, abbreviated as $\text{isAtl}_{\text{Ob}}(A, E, \boldsymbol{C})$, iff $A$ is an m-atlas of $\underset{\text{triv}}{E}$ in the coordinate space $\boldsymbol{C}$. $A$ is a full atlas of $E$ in the coordinate space $\boldsymbol{C}$, abbreviated $\text{isAtl}_{\text{Ob}}_{\text{full}}(A, E, \boldsymbol{C})$, iff $A$ is a full atlas of $\underset{\text{triv}}{E}$ in the coordinate space $\boldsymbol{C}$. The triple $(A, E, \boldsymbol{C})$ refers to $A$ considered as an m-atlas of $E$ in the coordinate space $\boldsymbol{C}$. $E$ is the total space for the atlas and $\boldsymbol{C}$ is the coordinate model space or coordinate space for the atlas.

By abuse of language we write $U \in A$ for $U \in \pi_1[A]$.

**Lemma 10.2** (M-atlases)**.** *Let $A$ be an m-atlas of $\boldsymbol{E} \stackrel{\text{def}}{=} (E, \mathscr{E})$ in the coordinate space $\boldsymbol{C} \stackrel{\text{def}}{=} (C, \mathscr{C})$. If $\boldsymbol{C}$ is fine grained then $A$ is an m-atlas of $E$ in the coordinate space $\boldsymbol{C}$.*

*Proof.* Let $(U, V, \phi) \in A$. $(U, V, \phi)$ is an m-chart of $E$ in the coordinate space $\boldsymbol{C}$ by corollary 9.4 (M-charts) on page 44. □

**Definition 10.3** (M-compatibility with m-atlases)**.** An m-chart $(U, V, \phi)$ of $(E, \mathscr{E})$ is m-compatible with an m-atlas $A$ iff it is m-compatible with every chart in $A$.

An m-chart $(U, V, \phi)$ of $E$ is m-compatible with an m-atlas $A$ iff it is m-compatible with every chart in $A$.

**Lemma 10.4** (M-compatibility of subcharts with atlases)**.** *Let $A$ be an m-atlas of $(E, \mathscr{E})$ in the coordinate space $(C, \mathscr{C})$ and $\boldsymbol{C}_1 = (U_1, V_1, \phi_1)$ an m-chart in $A$. Then any subchart of $\boldsymbol{C}_1$ is m-compatible with $A$.*

---

[8]The spaces are redundant, but convenient.

Let $A$ be an m-atlas of $E$ in the coordinate space $(C, \mathscr{C})$ and $C_1 = (U_1, V_1, \phi_1)$ an m-chart in $A$. Then any subchart of $C_1$ is m-compatible with $A$.

*Proof.* The same proof applies in both cases. Let $C' = (U', V', \phi')$ be a subchart of $C_1$ and $C_2 = (U_2, V_2, \phi_2)$ another chart in $A$.

1. If $U_1 \cap U_2 = \emptyset$, then $U' \cap U_2 = \emptyset$

2. If $U' \cap U_2 = \emptyset$ then $C'$ is m-compatible with $C_2$

3. Otherwise the transition function $t_2^1 \stackrel{\text{def}}{=} \phi_2 \circ \phi_1^{-1} \restriction_{\phi_1[U_1 \cap U_2]}$ is an isomorphism of $(C, \mathscr{C})$. Since $\phi_1[U_1 \cap U_2]$ and $V'$ are model neighborhoods, so is $\phi_1[U' \cap U_2]$ and thus $t_2^1 \restriction_{\phi_1[U' \cap U_2]}$ is an isomorphism of $(C, \mathscr{C})$.

□

**Lemma 10.5** (Extensions of m-atlases). *Let $E = (E, \mathscr{E})$ and $C = (C, \mathscr{C})$ be model spaces, $A$ an m-atlas of $E$ in the coordinate space $C$, and $(U, V, \phi)$, $(U', V', \phi')$ m-charts of $E$ in the coordinate space $C$ m-compatible with $A$ in the coordinate space $C$. Then $(U, V, \phi)$ is m-compatible with $(U', V', \phi')$ in the coordinate space $C$.*

*Proof.* If $U \cap U' = \emptyset$ then $(U, V, \phi)$ is m-compatible with $(U', V', \phi')$. Otherwise, $\phi$ is a model homeomorphism, $U \cap U'$ is a model neighborhood of $E$, $\phi[U \cap U']$ is a model neighborhood of $C$, $\phi'[U \cap U']$ is a model neighborhood of $C$ and $\phi' \circ \phi^{-1} \colon \phi[U \cap U'] \rightarrowtail\!\!\!\twoheadrightarrow \phi'[U \cap U']$ is a homeomorphism. It remains to show that $\phi' \circ \phi^{-1} \restriction_{\phi[U \cap U']}$ is a model homeomorphism. Let $(U_\alpha, V_\alpha, \phi_\alpha)$, $\alpha \prec A$, be charts in $A$ such that $U \cap U' \subseteq \bigcup_{\alpha \prec A} U_\alpha$ and $U \cap U' \cap U_\alpha \neq \emptyset$, $\alpha \prec A$. Since the charts are m-compatible with $(U_\alpha, V_\alpha, \phi_\alpha)$, $\phi' \circ \phi_\alpha^{-1} \restriction_{\phi_\alpha[U \cap U' \cap U_\alpha]}$ and $\phi_\alpha \circ \phi^{-1} \restriction_{\phi[U \cap U' \cap U_\alpha]}$ are model homeomorphisms and thus $\phi' \circ \phi^{-1} \restriction_{\phi[U \cap U' \cap U_\alpha]} = \phi' \circ \phi_\alpha^{-1} \circ \phi_\alpha \circ \phi^{-1} \restriction_{\phi[U \cap U' \cap U_\alpha]}$ is a model homeomorphism. Then by item 6 of definition 2.1 (Model spaces) on page 18, $\phi' \circ \phi^{-1}$ is a model homeomorphism. □

**Definition 10.6** (Maximal (semi-maximal) m-atlases). Let $E$ and $C$ be model spaces and $A$ an m-atlas of $E$ in the coordinate space $C$. $A$ is a maximal m-atlas of $E$ in the coordinate space $C$, abbreviated $\text{isAtl}_{\text{Ob}}^{\max}(A, E, C)$, iff $A$ cannot be extended by adding an additional m-compatible chart. $A$ is a semi-maximal m-atlas of $E$ in the coordinate space $C$, abbreviated $\text{isAtl}_{\text{Ob}}^{\text{S-max}}(A, E, C)$, iff whenever $(U, V, \phi) \in A$, $U' \subseteq U$, $V' \subseteq V$ and $V'' \stackrel{\text{Ob}}{\in} \mathscr{C}$ are model neighborhoods, $\phi[U'] = V'$ and $\phi' \colon V' \rightarrowtail\!\!\!\twoheadrightarrow V''$ is an isomorphism of $\mathscr{C}$ then $(U', V'', \phi' \circ \phi \restriction_{U'}) \in A$.

$$\text{isAtl}_{\text{Ob}}^{\text{max-full}}(A, E, C) \stackrel{\text{def}}{=} \text{isAtl}_{\text{Ob}}^{\text{full}}(A, E, C) \wedge \text{isAtl}_{\text{Ob}}^{\max}(A, E, C) \tag{10.1}$$

$$\text{isAtl}_{\text{Ob}}^{\text{S-max-full}}(A, E, C) \stackrel{\text{def}}{=} \text{isAtl}_{\text{Ob}}^{\text{full}}(A, E, C) \wedge \text{isAtl}_{\text{Ob}}^{\text{S-max}}(A, E, C) \tag{10.2}$$

Let $E$ be a topological space, $\boldsymbol{C} = (C, \mathscr{C})$ be a model space and $\boldsymbol{A}$ be an m-atlas of $E$ in the coordinate space $\boldsymbol{C}$. $\boldsymbol{A}$ is a maximal m-atlas of $E$ in the coordinate space $\boldsymbol{C}$, abbreviated $\text{isAtl}_{\text{Ob}}^{\max}(\boldsymbol{A}, E, \boldsymbol{C})$, iff $\boldsymbol{A}$ is an m-atlas that cannot be extended by adding an additional m-compatible chart. $\boldsymbol{A}$ is a semi-maximal m-atlas of $E$ in the coordinate space $\boldsymbol{C}$, abbreviated $\text{isAtl}_{\text{Ob}}^{\text{S-max}}(\boldsymbol{A}, E, \boldsymbol{C})$, iff whenever $(U, V, \phi) \in \boldsymbol{A}$, $U' \subseteq U$, $V' \subseteq V$ and $V'' \stackrel{\text{Ob}}{\in} \mathscr{C}$ are model neighborhoods, $\phi[U'] = V'$ and $\phi' \colon V' \stackrel{\tilde{=}}{\rightarrowtail\!\!\!\twoheadrightarrow} V''$ is an isomorphism of $\mathscr{C}$ then $(U', V'', \phi' \circ \phi\!\restriction_{U'}) \in \boldsymbol{A}$.

*Remark* 10.7. There is no sheaf condition[9]; the union of a set of coordinate patches in the maximal atlas whose coordinate functions match on the intersections need not be a coordinate patch in the atlas.

$$\text{isAtl}_{\text{Ob}}^{\text{max-full}}(\boldsymbol{A}, E, \boldsymbol{C}) \stackrel{\text{def}}{=} \text{isAtl}_{\text{Ob}}^{\text{full}}(\boldsymbol{A}, E, \boldsymbol{C}) \wedge \text{isAtl}_{\text{Ob}}^{\max}(\boldsymbol{A}, E, \boldsymbol{C}) \tag{10.3}$$

$$\text{isAtl}_{\text{Ob}}^{\text{S-max-full}}(\boldsymbol{A}, E, \boldsymbol{C}) \stackrel{\text{def}}{=} \text{isAtl}_{\text{Ob}}^{\text{full}}(\boldsymbol{A}, E, \boldsymbol{C}) \wedge \text{isAtl}_{\text{Ob}}^{\text{S-max}}(\boldsymbol{A}, E, \boldsymbol{C}) \tag{10.4}$$

**Lemma 10.8** (Maximal m-atlases are semi-maximal m-atlases). *Let $\boldsymbol{E} = (E, \mathscr{E})$ and $\boldsymbol{C} = (C, \mathscr{C})$ be model spaces and $\boldsymbol{A}$ a maximal m-atlas of $\boldsymbol{E}$ in the coordinate space $\boldsymbol{C}$. Then $\boldsymbol{A}$ is a semi-maximal m-atlas of $\boldsymbol{E}$ in the coordinate space $\boldsymbol{C}$.*

*Proof.* Let $(U, V, \phi) \in \boldsymbol{A}$, $U' \subseteq U$, $V' \subseteq V$ and $V'' \stackrel{\text{Ob}}{\in} \mathscr{C}$ be model neighborhoods, $\phi[U'] = V'$ and $\phi' \colon V' \stackrel{\tilde{=}}{\rightarrowtail\!\!\!\twoheadrightarrow} V''$ be an isomorphism of $\mathscr{C}$. $(U', V', \phi)$ is a subchart of $(U, V, \phi)$ and by lemma 10.4 (M-compatibility of subcharts with atlases) on page 47 is m-compatible with the charts of $\boldsymbol{A}$. Since $\phi'$ is a model homeomorphism, $(U', V'', \phi' \circ \phi)$ is m-compatible with the charts of $\boldsymbol{A}$. Since $\boldsymbol{A}$ is maximal, $(U', V'', \phi' \circ \phi)$ is a chart of $\boldsymbol{A}$. $\square$

Let $E$ be a topological space, $\boldsymbol{C} = (C, \mathscr{C})$ a model space and $\boldsymbol{A}$ a maximal m-atlas of $E$ in the coordinate space $\boldsymbol{C}$. Then $\boldsymbol{A}$ is a semi-maximal m-atlas of $E$ in the coordinate space $\boldsymbol{C}$.

*Proof.* Let $(U, V, \phi) \in \boldsymbol{A}$, $U' \subseteq U$ open, $V' \subseteq V$ and $V'' \stackrel{\text{Ob}}{\in} \mathscr{C}$ be model neighborhoods, $\phi[U'] = V'$ and $\phi' \colon V' \stackrel{\tilde{=}}{\rightarrowtail\!\!\!\twoheadrightarrow} V''$ be an isomorphism of $\mathscr{C}$. $(U', V', \phi)$ is a subchart of $(U, V, \phi)$ and by lemma 10.4 (M-compatibility of subcharts with atlases) on page 47 is m-compatible with the charts of $\boldsymbol{A}$. Since $\phi'$ is a model homeomorphism, $(U', V'', \phi' \circ \phi)$ is m-compatible with the charts of $\boldsymbol{A}$. Since $\boldsymbol{A}$ is maximal, $(U', V'', \phi' \circ \phi)$ is a chart of $\boldsymbol{A}$. $\square$

**Theorem 10.9** (Existence and uniqueness of maximal m-atlases). *Let $\boldsymbol{A}$ be an m-atlas of $\boldsymbol{E} \stackrel{\text{def}}{=} (E, \mathscr{E})$ in the coordinate space $\boldsymbol{C} \stackrel{\text{def}}{=} (C, \mathscr{C})$. Then there exists a unique maximal m-atlas $\text{Atlas}_{\max}(\boldsymbol{A}, \boldsymbol{E}, \boldsymbol{C})$ of $\boldsymbol{E}$ in the coordinate space $\boldsymbol{C}$ m-compatible with $\boldsymbol{A}$.*

---

[9] However, note item 6 (restricted sheaf condition) of definition 2.1 (Model spaces) on page 18.

*Proof.* Let $P$ be the set of all m-atlases of $E$ in the coordinate space $C$ containing $A$ and m-compatible in the coordinate space $C$ with all of the m-charts in $A$. Let $\underset{\max}{P}$ be a maximal chain of $A$. Then $A' = \bigcup \underset{\max}{P}$ is a maximal m-atlas of $E$ in the coordinate space $C$ m-compatible with $A$. Uniqueness follows from lemma 10.5 (Extensions of m-atlases) on page 48. $\square$

**Corollary 10.10** (Existence and uniqueness of maximal m-atlases). *Let $A$ be an m-atlas of $E$ in the coordinate space $C \overset{\text{def}}{=} (C, \mathscr{C})$. Then there exists a unique maximal m-atlas $\underset{\max}{\mathrm{Atlas}}(A, E, C)$ of $E$ in the coordinate space $C$ m-compatible with $A$.*

*Proof.* $\underset{\max}{\mathrm{Atlas}}(A, E, C) \overset{\text{def}}{=} \underset{\max}{\mathrm{Atlas}}(A, \underset{\text{triv}}{E}, C)$. $\square$

**Definition 10.11** (Sets of m-atlases). Let $E$ and $C$ be model spaces.

$$\underset{(\text{full,S-max,max,S-max-full,max-full})}{\mathscr{A}\mathscr{t}\mathscr{l}_{\mathrm{Ob}}}(E, C) \overset{\text{def}}{=} \left\{ (A, E, C) \;\middle|\; \underset{(\text{full,S-max,max,S-max-full,max-full})}{\mathrm{isAtl}_{\mathrm{Ob}}}(A, E, C) \right\} \quad (10.5)$$

Let $E$ and $C$ be sets of model spaces.

$$\underset{(\text{full,S-max,max,S-max-full,max-full})}{\mathscr{A}\mathscr{t}\mathscr{l}_{\mathrm{Ob}}}(E, C) \overset{\text{def}}{=} \left\{ (A, (E, \mathscr{E}), (C, \mathscr{C})) \;\middle|\; (E, \mathscr{E}) \in E \wedge (C, \mathscr{C}) \in C \wedge \underset{(\text{full,S-max,max,S-max-full,max-full})}{\mathrm{isAtl}_{\mathrm{Ob}}}(A, (E, \mathscr{E}), (C, \mathscr{C})) \right\} \quad (10.6)$$

Let $E$ be a set of topological spaces and $C$ a set of model spaces. Then

$$\underset{(\text{full,S-max,max,S-max-full,max-full})}{\mathscr{A}\mathscr{t}\mathscr{l}_{\mathrm{Ob}}}(E, C) \overset{\text{def}}{=} \left\{ (A, E, (C, \mathscr{C})) \;\middle|\; E \in E \wedge (C, \mathscr{C}) \in C \wedge \underset{(\text{full,S-max,max,S-max-full,max-full})}{\mathrm{isAtl}_{\mathrm{Ob}}}(A, E, (C, \mathscr{C})) \right\} \quad (10.7)$$

**Lemma 10.12** (Sets of m-atlases). *Let $E$ and $C$ be sets of model spaces.*

$$\underset{\text{full (S-max,max,S-max-full,max-full)}}{\mathscr{A}\mathscr{t}\mathscr{l}_{\mathrm{Ob}}}(E, C) \subseteq \mathscr{A}\mathscr{t}\mathscr{l}_{\mathrm{Ob}}(E, C) \quad (10.8)$$

$$\underset{\text{S-max-full (max-full)}}{\mathscr{A}\mathscr{t}\mathscr{l}_{\mathrm{Ob}}}(E, C) \subseteq \underset{\text{S-max (max)}}{\mathscr{A}\mathscr{t}\mathscr{l}_{\mathrm{Ob}}}(E, C) \quad (10.9)$$

*Proof.* The result follows from definition 10.1 (M-atlases) on page 47, definition 10.6 (Maximal (semi-maximal) m-atlases) on page 48 and definition 10.11 above. $\square$

$$\underset{\text{max-full}}{\mathcal{A}t\ell_{\text{Ob}}(E,C)} \subseteq \underset{\text{S-max-full}}{\mathcal{A}t\ell_{\text{Ob}}(E,C)} \subseteq \underset{\text{full}}{\mathcal{A}t\ell_{\text{Ob}}(E,C)} \tag{10.10}$$

*Proof.* The result follows from definition 10.6 (Maximal (semi-maximal) m-atlases) on page 48, lemma 10.8 (Maximal m-atlases are semi-maximal m-atlases) on page 49 and definition 10.11 above. □

**Corollary 10.13** (Sets of m-atlases). *Let $E$ be a set of topological spaces and $C$ a set of model spaces. Then*

$$\underset{\text{full (S-max,max,S-max-full,max-full)}}{\mathcal{A}t\ell_{\text{Ob}}(E,C)} \subseteq \mathcal{A}t\ell_{\text{Ob}}(E,C) \tag{10.11}$$

$$\underset{\text{S-max-full (max-full)}}{\mathcal{A}t\ell_{\text{Ob}}(E,C)} \subseteq \underset{\text{S-max (max)}}{\mathcal{A}t\ell_{\text{Ob}}(E,C)} \tag{10.12}$$

$$\underset{\text{max-full}}{\mathcal{A}t\ell_{\text{Ob}}(E,C)} \subseteq \underset{\text{S-max-full}}{\mathcal{A}t\ell_{\text{Ob}}(E,C)} \subseteq \underset{\text{full}}{\mathcal{A}t\ell_{\text{Ob}}(E,C)} \tag{10.13}$$

*Proof.* The result follows from definition 10.1 (M-atlases) on page 47, definition 10.6 (Maximal (semi-maximal) m-atlases) on page 48, lemma 10.8 (Maximal m-atlases are semi-maximal m-atlases) on page 49 and definition 10.11 above. □

## 11 M-atlas (near) morphisms and functors

This section introduces a taxonomy for morphisms between m-atlases, defines categories of m-atlases, defines functors amomg them and proves some basic reasults.

### 11.1 M-atlas (near) morphisms

This subsection introduces the notions of m-atlas near morphisms, m-atlas morphisms, semi-strict (near) m-atlas morphisms and strict (near) m-atlas morphisms. Only (semi-strict, strict) m-atlas morphisms are needed for the exposition of local coordinate spaces, as m-atlas near morphisms between maximal atlases are proven to be m-atlas morphisms.

#### 11.1.1 Definitions of m-atlas (near) morphisms

**Definition 11.1** (M-atlas near morphisms for model spaces). Let $\mathscr{E}^i, \mathscr{C}^i, i = 1, 2$, be model categories, $E^i \overset{\text{Ob}}{\in} \mathscr{E}^i$, $C^i \overset{\text{Ob}}{\in} \mathscr{C}^i$, $A^i$ be an m-atlas of $E^i$ in the coordinate space $C^i$ and $f \overset{\text{def}}{=} (f_0 \colon E^1 \longrightarrow E^2, f_1 \colon C^1 \longrightarrow C^2)$ a pair of model functions.

$f$ is an $E^1$-$E^2$ m-atlas near morphism of $A^1$ to $A^2$ in the coordinate spaces $C^1, C^2$, abbreviated as $\text{isAtl}^{\text{near}}_{\text{Ar}}(A^1, E^1, C^1, A^2, E^2, C^2, f_0, f_1)$ and an $\mathscr{E}^1$-$\mathscr{E}^2$ m-atlas near morphism of $A^1$ to $A^2$ in the coordinate model categories $\mathscr{C}^1, \mathscr{C}^2$, abbreviated as $\text{isAtl}^{\text{near}}_{\text{Ar}}(A^1, \mathscr{E}^1, \mathscr{C}^1, A^2, \mathscr{E}^2, \mathscr{C}^2, f_0, f_1)$ iff for any charts $(U^i, V^i, \phi^i \colon U^i \overset{\tilde{=}}{\rightarrowtail\!\!\!\twoheadrightarrow} V^i) \in$ $A^i$, $i = 1, 2$, with $I \stackrel{\text{def}}{=} U^1 \cap f_0^{-1}[U^2] \neq \emptyset$, diagram (11.1) below is strongly M-locally nearly commutative in $\boldsymbol{C}^2$, i.e., for any $u^1 \in I$ there are model neighborhoods $U'^1 \subseteq I$, $V'^1 \subseteq V^1$, $U'^2 \subseteq U^2$, $V'^2 \subseteq V^2$, $\hat{V}'^2 \subseteq C^2$ and an isomorphism $\hat{f}\colon V'^2 \rightarrowtail\!\!\!\xrightarrow{\tilde{=}}\!\!\!\twoheadrightarrow \hat{V}'^2$[10] such that eqs. (11.2) to (11.7) below hold, as shown in figs. 4 to 5.

$$D \stackrel{\text{def}}{=} \left\{ \begin{array}{l} f_0\colon I \stackrel{\text{def}}{=} U^1 \cap f_0^{-1}[U^2] \longrightarrow U^2, \phi^2\colon U^2 \rightarrowtail\!\!\!\xrightarrow{\tilde{=}}\!\!\!\twoheadrightarrow V^2, \\ \phi^1\colon I \rightarrowtail V^1, f_1\colon V^1 \longrightarrow \boldsymbol{C}^2 \end{array} \right\} \quad (11.1)$$

$$u^1 \in U'^1 \quad (11.2)$$

$$f_0[U'^1] \subseteq U'^2 \quad (11.3)$$

$$\phi^1[U'^1] = V'^1 \quad (11.4)$$

$$f_1[V'^1] \subseteq \hat{V}'^2 \quad (11.5)$$

$$\phi^2[U'^2] = V'^2 \quad (11.6)$$

$$\hat{f} \circ \phi^2 \circ f_0 = f_1 \circ \phi^1 \quad (11.7)$$

*Remark* 11.2. $(U'^i, V'^i, \phi^i\restriction_{U'^i, V'^i})$ need not be a chart of $\boldsymbol{A}^i$.

The triple $\left(\boldsymbol{f}, (\boldsymbol{A}^1, \boldsymbol{E}^1, \boldsymbol{C}^1), (\boldsymbol{A}^2, \boldsymbol{E}^2, \boldsymbol{C}^2)\right)$ will refer to $\boldsymbol{f}$ considered as an $\boldsymbol{E}^1$-$\boldsymbol{E}^2$ m-atlas near morphism of $\boldsymbol{A}^1$ to $\boldsymbol{A}^2$ in the coordinate spaces $\boldsymbol{C}^1, \boldsymbol{C}^2$.

**Definition 11.3** (M-atlas near morphisms for topological spaces). Let $\mathscr{C}^i$, $i = 1, 2$, be a model category, $\boldsymbol{C}^i \stackrel{\text{Ob}}{\in} \mathscr{C}^i$, $E^i$ be a topological space, $\boldsymbol{A}^i$ be an m-atlas of $E^i$ in the coordinate space $\boldsymbol{C}^i$ and $\boldsymbol{f} \stackrel{\text{def}}{=} (f_0\colon E^1 \longrightarrow E^2, f_1\colon \boldsymbol{C}^1 \longrightarrow \boldsymbol{C}^2)$ be a pair of continuous functions.

$\boldsymbol{f}$ is an $E^1$-$E^2$ m-atlas near morphism of $\boldsymbol{A}^1$ to $\boldsymbol{A}^2$ in the coordinate spaces $\boldsymbol{C}^1$, $\boldsymbol{C}^2$, abbreviated as isAtl$_{\text{Ar}}^{\text{near}}(\boldsymbol{A}^1, E^1, \boldsymbol{C}^1, \boldsymbol{A}^2, E^2, \boldsymbol{C}^2, f_0, f_1)$, and a $E^1$-$E^2$ m-atlas near morphism of $\boldsymbol{A}^1$ to $\boldsymbol{A}^2$ in the coordinate model categories $\mathscr{C}^1, \mathscr{C}^2$, abbreviated as

---

[10]This reverses the direction of the arrow $\hat{f}\colon U_m \rightarrowtail\!\!\!\xrightarrow{\tilde{=}}\!\!\!\twoheadrightarrow V_n$ from definition 1.24 (Nearly commutative diagrams in category $\mathscr{C}$) on page 14. This is permissible since $\hat{f}$ is an isomorphism

$$I = U^1 \cap f_0^{-1}[U^2] \xrightarrow{\quad f_0 \quad} U^2$$

$$\phi^1 \downarrow \qquad\qquad \tilde{=} \downarrow \phi^2$$

$$V^1 \xrightarrow{\quad f_1 \quad} C^2 \qquad\qquad V^2$$

Figure 4: Uncompleted m-atlas near morphism

$$u^1$$
$$\cap$$
$$\downarrow i$$

$$I = U^1 \cap f_0^{-1}[U^2] \xleftarrow{\;i\;} U'^1 \xrightarrow{\quad f_0 \quad} U'^2 \xhookrightarrow{\;i\;} U^2$$

$$\phi^1 \downarrow \qquad \tilde{=} \downarrow \phi^1 \qquad\qquad \tilde{=} \downarrow \phi^2 \qquad \tilde{=} \downarrow \phi^2$$

$$V^1 \xleftarrow{\;i\;} V'^1 \xrightarrow{f_1} \hat{V}'^2 \xleftarrow{\hat{f}}_{\tilde{=}} V'^2 \xhookrightarrow{\;i\;} V^2$$

Figure 5: Completed m-atlas near morphism

isAtl$_{\text{Ar}}^{\text{near}}(\boldsymbol{A}^1, E^1, \mathscr{C}^1, \boldsymbol{A}^2, E^2, \mathscr{C}^2, f_0, f_1)$, iff it is an $E^1$-$E^2$ m-atlas near morphism
$\text{triv}\ \text{triv}$
of $\boldsymbol{A}^1$ to $\boldsymbol{A}^2$ in the coordinate model spaces $\boldsymbol{C}^1, \boldsymbol{C}^2$.

The triple $\left(\boldsymbol{f}, (\boldsymbol{A}^1, E^1, \boldsymbol{C}^1), (\boldsymbol{A}^2, E^2, \boldsymbol{C}^2)\right)$ will refer to $\boldsymbol{f}$ considered as an $E^1$-$E^2$ m-atlas near morphism of $\boldsymbol{A}^1$ to $\boldsymbol{A}^2$ in the coordinate spaces $\boldsymbol{C}^1, \boldsymbol{C}^2$.

**Definition 11.4** (M-atlas morphisms for model spaces). Let $\mathscr{E}^i, \mathscr{C}^i, i = 1, 2$, be model categories, $\boldsymbol{E}^i \overset{\text{Ob}}{\in} \mathscr{E}^i, \boldsymbol{C}^i \overset{\text{Ob}}{\in} \mathscr{C}^i, \boldsymbol{A}^i$ be an m-atlas of $\boldsymbol{E}^i$ in the coordinate space $\boldsymbol{C}^i$ and $\boldsymbol{f} \overset{\text{def}}{=} (f_0 \colon \boldsymbol{E}^1 \longrightarrow \boldsymbol{E}^2, f_1 \colon \boldsymbol{C}^1 \longrightarrow \boldsymbol{C}^2)$ a pair of model functions.

$\boldsymbol{f}$ is an $\boldsymbol{E}^1$-$\boldsymbol{E}^2$ m-atlas morphism of $\boldsymbol{A}^1$ to $\boldsymbol{A}^2$ in the coordinate spaces $\boldsymbol{C}^1, \boldsymbol{C}^2$, abbreviated as isAtl$_{\text{Ar}}(\boldsymbol{A}^1, \boldsymbol{E}^1, \boldsymbol{C}^1, \boldsymbol{A}^2, \boldsymbol{E}^2, \boldsymbol{C}^2, f_0, f_1)$ and an $\mathscr{E}^1$-$\mathscr{E}^2$ m-atlas morphism of $\boldsymbol{A}^1$ to $\boldsymbol{A}^2$ in the coordinate model categories $\mathscr{C}^1, \mathscr{C}^2$, abbreviated as

isAtl$_{\text{Ar}}(\boldsymbol{A}^1, \mathscr{E}^1, \mathscr{C}^1, \boldsymbol{A}^2, \mathscr{E}^2, \mathscr{C}^2, f_0, f_1)$, iff for any charts $(U^i, V^i, \phi^i \colon U^i \xrightarrowtail{\tilde{=}} V^i) \in \boldsymbol{A}^i, i = 1, 2$, with $I \overset{\text{def}}{=} U^1 \cap f_0^{-1}[U^2] \neq \emptyset$ and any $u^1 \in I$, there exists a subchart $(U'^1, V'^1, \phi^1 \colon U'^1 \xrightarrowtail{\tilde{=}} V'^1) \in \boldsymbol{A}^1$ at $u^1$, a model neighborhood $U'^2 \subseteq U^2$ of $u^2 \stackrel{\text{def}}{=} f_0(u^1)$ and a chart $(U'^2, \hat{V}'^2, \phi'^2 \colon U'^2 \rightarrowtail\stackrel{\tilde{=}}{\twoheadrightarrow} \hat{V}'^2) \in \boldsymbol{A}^2$ at $u^2$ such that $f_0[U'^1] \subseteq U'^2$ and diagram (11.8) below is commutative, as shown in fig. 6 (Completed m-atlas morphism) on page 54.

$$D \stackrel{\text{def}}{=} \left\{ \begin{array}{l} i \colon \{u^1\} \rightarrowtail U'^1, \phi^1 \colon U'^1 \rightarrowtail\stackrel{\tilde{=}}{\twoheadrightarrow} V'^1, f_1 \colon V'^1 \longrightarrow \hat{V}'^2, \\ f_0 \colon U'^1 \longrightarrow U'^2, \phi'^2 \colon U'^2 \rightarrowtail\stackrel{\tilde{=}}{\twoheadrightarrow} \hat{V}'^2 \end{array} \right\} \quad (11.8)$$

Figure 6: Completed m-atlas morphism

The triple $\left(\boldsymbol{f}, (\boldsymbol{A}^1, \boldsymbol{E}^1, \boldsymbol{C}^1), (\boldsymbol{A}^2, \boldsymbol{E}^2, \boldsymbol{C}^2)\right)$ will refer to $\boldsymbol{f}$ considered as an $\boldsymbol{E}^1$-$\boldsymbol{E}^2$ m-atlas morphism of $\boldsymbol{A}^1$ to $\boldsymbol{A}^2$ in the coordinate spaces $\boldsymbol{C}^1$, $\boldsymbol{C}^2$.

$\boldsymbol{f}$ is also a constrained $\boldsymbol{E}^1$-$\boldsymbol{E}^2$ m-atlas morphism of $\boldsymbol{A}^1$ to $\boldsymbol{A}^2$ in the coordinate spaces $\boldsymbol{C}^1$, $\boldsymbol{C}^2$, abbreviated as $\text{isAtl}_{\text{Ar}}^{\text{const}}(\boldsymbol{A}^1, \boldsymbol{E}^1, \boldsymbol{C}^1, \boldsymbol{A}^2, \boldsymbol{E}^2, \boldsymbol{C}^2, f_0, f_1)$ and a constrained $\mathscr{C}^1$-$\mathscr{C}^2$ m-atlas morphism of $\boldsymbol{A}^1$ to $\boldsymbol{A}^2$ in the coordinate model categories $\mathscr{C}^1$, $\mathscr{C}^2$, abbreviated as $\text{isAtl}_{\text{Ar}}^{\text{const}}(\boldsymbol{A}^1, \mathscr{C}^1, \boldsymbol{C}^1, \boldsymbol{A}^2, \mathscr{C}^2, \boldsymbol{C}^2, f_0, f_1)$, iff $f_0$ and $f_1$ are constrained, i.e., $f_0[\boldsymbol{E}^1]$ is contained in a model neighborhood of $\boldsymbol{E}^2$ and $f_1[\boldsymbol{C}^1]$ is contained in a model neighborhood of $\boldsymbol{C}^2$.

*Remark* 11.5. If $\boldsymbol{f}$ is constrained then $\text{Top}(\boldsymbol{E}^1)$ is a model neighborhood of $\boldsymbol{E}^1$ and $\text{Top}(\boldsymbol{C}^1)$ is a model neighborhood of $\boldsymbol{C}^1$.

**Definition 11.6** (M-atlas morphisms for topological spaces). Let $\mathscr{C}^i, i = 1, 2$, be a model category, $\boldsymbol{C}^i \stackrel{\text{Ob}}{\in} \mathscr{C}^i$, $E^i$ be a topological space and $\boldsymbol{A}^i$ be an m-atlas of $E^i$ in the coordinate space $\boldsymbol{C}^i$.

A pair of functions $\boldsymbol{f} \stackrel{\text{def}}{=} (f_0 \colon E^1 \longrightarrow E^2, f_1 \colon \boldsymbol{C}^1 \longrightarrow \boldsymbol{C}^2)$ is an $E^1$-$E^2$ m-atlas morphism of $\boldsymbol{A}^1$ to $\boldsymbol{A}^2$ in the coordinate spaces $\boldsymbol{C}^i$, abbreviated as $\text{isAtl}_{\text{Ar}}(\boldsymbol{A}^1, E^1, \boldsymbol{C}^1, \boldsymbol{A}^2, E^2, \boldsymbol{C}^2, f_0, f_1)$, and an $E^1$-$E^2$ m-atlas morphism of $\boldsymbol{A}^1$ to $\boldsymbol{A}^2$ in the coordinate model categories $\mathscr{C}^i$, abbreviated as
isAtl$_{\text{Ar}}(A^1, E^1, \mathscr{C}^1, A^2, E^2, \mathscr{C}^2, f_0, f_1)$, iff it is an $\underset{\text{triv triv}}{E^1\text{-}E^2}$ m-atlas morphism of $A^1$ to $A^2$ in the coordinate model categories $\mathscr{C}^1, \mathscr{C}^2$.

The triple $(f, (A^1, E^1, C^1), (A^2, E^2, C^2))$ will refer to $f$ considered as an $E^1$-$E^2$ m-atlas (near) morphism of $A^1$ to $A^2$ in the coordinate spaces $C^1, C^2$.

$f$ is also a constrained $E^1$-$E^2$ m-atlas morphism of $A^1$ to $A^2$ in the coordinate spaces $C^1, C^2$, abbreviated as isAtl$_{\text{Ar}}^{\text{const}}(A^1, E^1, C^1, A^2, E^2, C^2, f_0, f_1)$, and a constrained $E^1$-$E^2$ m-atlas morphism of $A^1$ to $A^2$ in the coordinate model categories $\mathscr{C}^1, \mathscr{C}^2$, abbreviated as isAtl$_{\text{Ar}}^{\text{const}}(A^1, E^1, \mathscr{C}^1, A^2, E^2, \mathscr{C}^2, f_0, f_1)$, iff $f_1$ is constrained, i.e., $f_1[C^1]$ is contained in a model neighborhood of $C^2$.

**Definition 11.7** (Equivalence of m-atlas (near) morphisms). Let $E^i, C^i, i = 1, 2$, be model spaces, $A^i$ be an m-atlas of $E^i$ in the coordinate space $C^i$, and

$f \stackrel{\text{def}}{=} (f_0 \colon E^1 \longrightarrow E^2, f_1 \colon C^1 \longrightarrow C^2)$ and $g \stackrel{\text{def}}{=} (g_0 \colon E^1 \longrightarrow E^2, g_1 \colon C^1 \longrightarrow C^2)$ be m-atlas (near) morphisms of $A^1$ to $A^2$ in the coordinate spaces $C^1, C^2$. Then $f$ is equivalent to $g$ as an $E^1$-$E^2$ m-atlas (near) morphism of $A^1$ to $A^2$ in the coordinate spaces $C^1, C^2$, abbreviated $f \colon (A^1, E^1, C^1) \longrightarrow (A^2, E^2, C^2)$ is equivalent to $g \colon (A^1, E^1, C^1) \longrightarrow (A^2, E^2, C^2)$ and $(f, (A^1, E^1, C^1), (A^2, E^2, C^2))$ is equivalent to $(g, (A^1, E^1, C^1), (A^2, E^2, C^2))$ iff $f_0 = g_0$.

By abuse of language we write $f$ is equivalent to $g$ when the spaces and atlases are understood by context.

Let $E^i, i = 1, 2$, be a topological space, $C^i$ be a model space, $A^i$ be an m-atlas of $E^i$ in the coordinate space $C^i$ and $f \stackrel{\text{def}}{=} (f_0 \colon E^1 \longrightarrow E^2, f_1 \colon C^1 \longrightarrow C^2)$, and $g \stackrel{\text{def}}{=} (g_0 \colon E^1 \longrightarrow E^2, g_1 \colon C^1 \longrightarrow C^2)$ be m-atlas (near) morphisms. Then $f$ is equivalent to $g$ as an $E^1$-$E^2$ m-atlas (near) morphism, abbreviated $f \colon (A^1, E^1, C^1) \longrightarrow (A^2, E^2, C^2)$ is equivalent to $g \colon (A^1, E^1, C^1) \longrightarrow (A^2, E^2, C^2)$, iff $f_0 = g_0$.

By abuse of language we write $f$ is equivalent to $g$ when the spaces and atlases are understood by context.

### 11.1.2 Definitions of semi-strict and strict

**Definition 11.8** (Semi-strict m-atlas (near) morphisms). Let $\mathscr{C}^i, \mathscr{C}^i, i = 1, 2$, be model categories, $E^i \stackrel{\text{Ob}}{\in} \mathscr{C}^i, C^i \stackrel{\text{Ob}}{\in} \mathscr{C}^i, A^i$ be an m-atlas of $E^i$ in the coordinate space $C^i$ and $f \stackrel{\text{def}}{=} (f_0 \colon E^1 \longrightarrow E^2, f_1 \colon C^1 \longrightarrow C^2)$ an $E^1$-$E^2$ m-atlas (near) morphism of $A^1$ to $A^2$ in the coordinate spaces $C^1, C^2$.

$f$ is a semi-strict $E^1$-$E^2$ m-atlas (near) morphism of $A^1$ to $A^2$ in the coordinate spaces $C^1, C^2$, abbreviated as isAtl$_{\text{Ar}\atop\text{semi-strict}}^{(\text{near})}(A^1, E^1, C^1, A^2, E^2, C^2, f_0, f_1)$, iff $f_0$ is a strict local m-morphism of $E^1$ to $E^2$, $f_1$ is a strict local m-morphism of $C^1$ to $C^2$ and for every $(U^i, V^i, \phi^i \colon U^i \xrightarrow{\tilde{=}} V^i) \in A^i, i = 1, 2$, with $I \stackrel{\text{def}}{=} U^1 \cap f_0^{-1}[U^2] \neq \emptyset$, $\phi^2 \circ f_0 \circ \phi^{1-1} \colon \phi^1[I] \longrightarrow V^2$ is a local $C^1$-$C^2$ m-morphism of $\phi^1[I]$ to $V^2$.

$f$ is a semi-strict $E^1$-$E^2$-$\mathscr{C}^1$-$\mathscr{C}^2$ m-atlas (near) morphism of $A^1$ to $A^2$ in the coordinate spaces $C^1$, $C^2$, abbreviated as $\mathrm{isAtl}_{\mathrm{Ar}\ \mathrm{semi\text{-}strict}}^{(\mathrm{near})}(A^1, E^1, C^1, A^2, E^2, C^2, f_0, f_1, \mathscr{C}^1, \mathscr{C}^2)$, iff $f_0$ is a strict local $E^1$-$E^2$ m-morphism of $E^1$ to $E^2$, $f_1$ is a local $\mathscr{C}^1$-$\mathscr{C}^2$ m-morphism of $C^1$ to $C^2$ and for every $(U^i, V^i, \phi^i\colon U^i \rightarrowtail\!\!\xrightarrow{\tilde{=}}\!\!\twoheadrightarrow V^i) \in A^i$, $i = 1, 2$, with $I \stackrel{\mathrm{def}}{=} U^1 \cap f_0^{-1}[U^2] \neq \emptyset$, $\phi^2 \circ f_0 \circ \phi^{1-1}\colon \phi^1[I] \longrightarrow V^2$ is a local $\mathscr{C}^1$-$\mathscr{C}^2$ m-morphism of $\phi^1[I]$ to $V^2$.

It is a semi-strict $\mathscr{C}^1$-$\mathscr{C}^2$-$E^1$-$E^2$ m-atlas (near) morphism of $A^1$ to $A^2$ in the coordinate spaces $C^1$, $C^2$, abbreviated as $\mathrm{isAtl}_{\mathrm{Ar}\ \mathrm{semi\text{-}strict}}^{(\mathrm{near})}(A^1, E^1, C^1, A^2, E^2, C^2, f_0, f_1, \mathscr{C}^1, \mathscr{C}^2)$, iff $f_0$ is a local $\mathscr{C}^1$-$\mathscr{C}^2$ m-morphism of $E^1$ to $E^2$, $f_1$ is a strict local m-morphism of $C^1$ to $C^2$ and for every $(U^i, V^i, \phi^i\colon U^i \rightarrowtail\!\!\xrightarrow{\tilde{=}}\!\!\twoheadrightarrow V^i) \in A^i$, $i = 1, 2$, with $I \stackrel{\mathrm{def}}{=} U^1 \cap f_0^{-1}[U^2] \neq \emptyset$, $\phi^2 \circ f_0 \circ \phi^{1-1}\colon \phi^1[I] \longrightarrow U^2$ is a local $C^1$-$C^2$ m-morphism of $\phi^1[I]$ to $V^2$.

It is a semi-strict $\mathscr{C}^1$-$\mathscr{C}^2$-$E^1$-$E^2$-$\mathscr{C}^1$-$\mathscr{C}^2$ m-atlas (near) morphism of $A^1$ to $A^2$ in the coordinate spaces $C^1$, $C^2$, abbreviated as
$\mathrm{isAtl}_{\mathrm{Ar}\ \mathrm{semi\text{-}strict}}^{(\mathrm{near})}(A^1, E^1, C^1, A^2, E^2, C^2, f_0, f_1, \mathscr{C}^1, \mathscr{C}^2, \mathscr{C}^1, \mathscr{C}^2)$, iff $f_0$ is a local $\mathscr{C}^1$-$\mathscr{C}^2$ m-morphism of $E^1$ to $E^2$ and $f_1$ is a local $\mathscr{C}^1$-$\mathscr{C}^2$ m-morphism of $C^1$ to $C^2$ and for every $(U^i, V^i, \phi^i\colon U^i \rightarrowtail\!\!\xrightarrow{\tilde{=}}\!\!\twoheadrightarrow V^i) \in A^i$, $i = 1, 2$, with $I \stackrel{\mathrm{def}}{=} U^1 \cap f_0^{-1}[U^2] \neq \emptyset$, $\phi^2 \circ f_0 \circ \phi^{1-1}\colon \phi^1[I] \longrightarrow V^2$ is a local $\mathscr{C}^1$-$\mathscr{C}^2$ m-morphism of $\phi^1[I]$ to $V^2$.

**Definition 11.9** (Strict m-atlas (near) morphisms)**.** Let

1. $\mathscr{C}^i, \mathscr{C}^i, i = 1, 2$, be model categories

2. $E^i \stackrel{\mathrm{Ob}}{\in} \mathscr{C}^i$

3. $C^i \stackrel{\mathrm{Ob}}{\in} \mathscr{C}^i$

4. $A^i$ be an m-atlas of $E^i$ in the coordinate space $C^i$

5. $f \stackrel{\mathrm{def}}{=} (f_0\colon E^1 \longrightarrow E^2, f_1\colon C^1 \longrightarrow C^2)$ be an m-atlas (near) morphism

$f$ is a strict $E^1$-$E^2$ m-atlas (near) morphism of $A^1$ to $A^2$ in the coordinate spaces $C^1$, $C^2$, abbreviated as $\mathrm{isAtl}_{\mathrm{Ar}\ \mathrm{strict}}^{(\mathrm{near})}(A^1, E^1, C^1, A^2, E^2, C^2, f_0, f_1)$, iff $f_0$ is a strict m-morphism of $E^1$ to $E^2$, $f_1$ is a strict m-morphism of $C^1$ to $C^2$ and for every $(U^i, V^i, \phi^i\colon U^i \rightarrowtail\!\!\xrightarrow{\tilde{=}}\!\!\twoheadrightarrow V^i) \in A^i$, $i = 1, 2$, with $I \stackrel{\mathrm{def}}{=} U^1 \cap f_0^{-1}[U^2] \neq \emptyset$, $\phi^2 \circ f_0 \circ \phi^{1-1}\colon \phi^1[I] \longrightarrow U^2$ is a morphism of $\mathrm{Mod}(V^2, C^2)$.

$f$ is a strict $E^1$-$E^2$-$\mathscr{C}^1$-$\mathscr{C}^2$ m-atlas (near) morphism of $A^1$ to $A^2$ in the coordinate spaces $C^1$, $C^2$, abbreviated as $\mathrm{isAtl}_{\mathrm{Ar}\ \mathrm{strict}}^{(\mathrm{near})}(A^1, E^1, C^1, A^2, E^2, C^2, f_0, f_1, \mathscr{C}^1, \mathscr{C}^2)$, iff $f_0$ is a strict m-morphism of $E^1$ to $E^2$, $f_1$ is an m-morphism of $\mathscr{C}^1$ to $\mathscr{C}^2$ and for every

$(U^i, V^i, \phi^i\colon U^i \rightarrowtail\!\!\!\xrightarrow{\sim} V^i) \in \boldsymbol{A}^i$, $i = 1, 2$, with $I \stackrel{\text{def}}{=} U^1 \cap f_0^{-1}[U^2] \neq \emptyset$, $\phi^2 \circ f_0 \circ \phi^{1-1}\colon \phi^1[I] \longrightarrow V^2$ is a morphism of $\mathscr{C}^2$.

$\boldsymbol{f}$ is a strict $\mathscr{E}^1$-$\mathscr{E}^2$-$\boldsymbol{E}^1$-$\boldsymbol{E}^2$ m-atlas (near) morphism of $\boldsymbol{A}^1$ to $\boldsymbol{A}^2$ in the coordinate spaces $\boldsymbol{C}^1, \boldsymbol{C}^2$, abbreviated as $\text{isAtl}_{\text{Ar strict}}^{(\text{near})}(\boldsymbol{A}^1, \boldsymbol{E}^1, \boldsymbol{C}^1, \boldsymbol{A}^2, \boldsymbol{E}^2, \boldsymbol{C}^2, f_0, f_1, \mathscr{E}^1, \mathscr{E}^2)$, iff $f_0$ is an m-morphism of $\mathscr{E}^1$ to $\mathscr{E}^2$, $f_1$ is a strict m-morphism of $\boldsymbol{C}^1$ to $\boldsymbol{C}^2$ and for every

$(U^i, V^i, \phi^i\colon U^i \rightarrowtail\!\!\!\xrightarrow{\sim} V^i) \in \boldsymbol{A}^i$, $i = 1, 2$, with $I \stackrel{\text{def}}{=} U^1 \cap f_0^{-1}[U^2] \neq \emptyset$, $\phi^2 \circ f_0 \circ \phi^{1-1}\colon \phi^1[I] \longrightarrow U^2$ is a morphism of $\text{Mod}(\phi^1[I], \boldsymbol{C}^2)$.

$\boldsymbol{f}$ is a strict $\mathscr{E}^1$-$\mathscr{E}^2$-$\boldsymbol{E}^1$-$\boldsymbol{E}^2$-$\mathscr{C}^1$-$\mathscr{C}^2$ m-atlas (near) morphism of $\boldsymbol{A}^1$ to $\boldsymbol{A}^2$ in the coordinate spaces $\boldsymbol{C}^1, \boldsymbol{C}^2$, abbreviated as
$\text{isAtl}_{\text{Ar strict}}^{(\text{near})}(\boldsymbol{A}^1, \boldsymbol{E}^1, \boldsymbol{C}^1, \boldsymbol{A}^2, \boldsymbol{E}^2, \boldsymbol{C}^2, f_0, f_1, \mathscr{E}^1, \mathscr{E}^2, \mathscr{C}^1, \mathscr{C}^2)$, iff $f_0$ is an m-morphism of

$\mathscr{E}^1$ to $\mathscr{E}^2$, $f_1$ is an m-morphism of $\mathscr{C}^1$ to $\mathscr{C}^2$ and for every $(U^i, V^i, \phi^i\colon U^i \rightarrowtail\!\!\!\xrightarrow{\sim} V^i) \in \boldsymbol{A}^i$, $i = 1, 2$, with $I \stackrel{\text{def}}{=} U^1 \cap f_0^{-1}[U^2] \neq \emptyset$, $\phi^2 \circ f_0 \circ \phi^{1-1}\colon \phi^1[I] \longrightarrow V^2$ is a morphism of $\mathscr{C}^2$.

### 11.1.3 Abbreviated nomenclature for m-atlas morphisms

**Definition 11.10** (Abbreviated nomenclature for m-atlas (near) morphisms). Let $\mathscr{E}^i, \mathscr{C}^i$, $i = 1, 2$, be model categories, $\boldsymbol{E}^i \stackrel{\text{Ob}}{\in} \mathscr{E}^i$, $\boldsymbol{C}^i \stackrel{\text{Ob}}{\in} \mathscr{C}^i$, $\boldsymbol{A}^i$ be an m-atlas of $\boldsymbol{E}^i$ in the coordinate space $\boldsymbol{C}^i$ and $\boldsymbol{f} \stackrel{\text{def}}{=} (f_0\colon \boldsymbol{E}^1 \longrightarrow \boldsymbol{E}^2, f_1\colon \boldsymbol{C}^1 \longrightarrow \boldsymbol{C}^2)$ a pair of model functions.

$\boldsymbol{f}$ is a constrained semi-strict (strict) $\boldsymbol{E}^1$-$\boldsymbol{E}^2$ m-atlas (near) morphism of $\boldsymbol{A}^1$ to $\boldsymbol{A}^2$ in the coordinate spaces $\boldsymbol{C}^1, \boldsymbol{C}^2$, abbreviated as
$\text{isAtl}_{\text{Ar semi-strict (strict)}}^{(\text{const,const}-\text{near})}(\boldsymbol{A}^1, \boldsymbol{E}^1, \boldsymbol{C}^1, \boldsymbol{A}^2, \boldsymbol{E}^2, \boldsymbol{C}^2, f_0, f_1)$, iff it is a constrained $\boldsymbol{E}^1$-$\boldsymbol{E}^2$ m-atlas (near) morphism of $\boldsymbol{A}^1$ to $\boldsymbol{A}^2$ in the coordinate spaces $\boldsymbol{C}^1, \boldsymbol{C}^1$ and it is a semi-strict (trict) $\boldsymbol{E}^1$-$\boldsymbol{E}^2$ m-atlas (near) morphism of $\boldsymbol{A}^1$ to $\boldsymbol{A}^2$ in the coordinate spaces $\boldsymbol{C}^1, \boldsymbol{C}^1$.

$\boldsymbol{f}$ is a (constrained) (semi-strict, strict) $\boldsymbol{E}^1$-$\boldsymbol{E}^2$ m-atlas (near) morphism of $\boldsymbol{A}^1$ to $\boldsymbol{A}^2$ in the coordinate space $\boldsymbol{C}^1$, abbreviated as
$\text{isAtl}_{\text{Ar (semi-strict,strict)}}^{(\text{const,const}-\text{near,near})}(\boldsymbol{A}^1, \boldsymbol{E}^1, \boldsymbol{C}^1, \boldsymbol{A}^2, \boldsymbol{E}^2, f_0, f_1)$, iff it is a (constrained) (semi-strict, strict) $\boldsymbol{E}^1$-$\boldsymbol{E}^2$ m-atlas (near) morphism of $\boldsymbol{A}^1$ to $\boldsymbol{A}^2$ in the coordinate spaces $\boldsymbol{C}^1, \boldsymbol{C}^1$.

$\boldsymbol{f}$ is a (constrained) (semi-strict, strict) $\mathscr{E}^1$-$\mathscr{E}^2$-$\boldsymbol{E}^1$-$\boldsymbol{E}^2$ m-atlas (near) morphism of $\boldsymbol{A}^1$ to $\boldsymbol{A}^2$ in the coordinate space $\boldsymbol{C}^1$, abbreviated as
$\text{isAtl}_{\text{Ar (semi-strict,strict)}}^{(\text{const,const}-\text{near,near})}(\boldsymbol{A}^1, \boldsymbol{E}^1, \boldsymbol{C}^1, \boldsymbol{A}^2, \boldsymbol{E}^2, f_0, f_1, \mathscr{E}^1, \mathscr{E}^2)$, iff it is a (constrained) (semi-strict, strict) $\mathscr{E}^1$-$\mathscr{E}^2$-$\boldsymbol{E}^1$-$\boldsymbol{E}^2$ m-atlas (near) morphism of $\boldsymbol{A}^1$ to $\boldsymbol{A}^2$ in the coordinate spaces $\boldsymbol{C}^1, \boldsymbol{C}^1$.

$\boldsymbol{f}$ is a (constrained) (semi-strict, strict) $\boldsymbol{E}^1$ m-atlas (near) morphism of $\boldsymbol{A}^1$ to $\boldsymbol{A}^2$ in the coordinate spaces $C^1$, $C^2$, abbreviated as

$\text{isAtl}_{\text{Ar}\ (\text{semi-strict,strict})}^{(\text{const,const-near,near})}(A^1, E^1, C^1, A^2, C^2, f_0, f_1)$, iff it is a (constrained) (semi-strict, strict) $E^1$-$E^1$ m-atlas (near) morphism of $A^1$ to $A^2$ in the coordinate spaces $C^1$, $C^2$.

$f$ is a (constrained) (semi-strict, strict) $E^1$ m-atlas (near) morphism of $A^1$ to $A^2$ in the coordinate space $C^1$, abbreviated as

$\text{isAtl}_{\text{Ar}\ (\text{semi-strict,strict})}^{(\text{const,const-near,near})}(A^1, E^1, C^1, A^2, f_0, f_1)$, iff it is a (constrained) (semi-strict, strict) $E^1$-$E^1$ m-atlas (near) morphism of $A^1$ to $A^2$ in the coordinate spaces $C^1$, $C^1$.

$f$ is a (constrained) (semi-strict, strict) $E^1$-$\mathscr{C}^1$-$\mathscr{C}^2$ m-atlas (near) morphism of $A^1$ to $A^2$ in the coordinate spaces $C^1$, $C^2$, abbreviated as

$\text{isAtl}_{\text{Ar}\ (\text{semi-strict,strict})}^{(\text{const,const-near,near})}(A^1, E^1, C^1, A^2, C^2, f_0, f_1, \mathscr{C}^1, \mathscr{C}^2)$, iff it is a (constrained) (semi-strict, strict) $E^1$-$E^1$-$\mathscr{C}^1$-$\mathscr{C}^2$ m-atlas (near) morphism of $A^1$ to $A^2$ in the coordinate spaces $C^1$, $C^2$.

$f$ is a (constrained) (semi-strict, strict) $E^1$-$\mathscr{C}^1$ m-atlas (near) morphism of $A^1$ to $A^2$ in the coordinate space $C^1$, abbreviated as

$\text{isAtl}_{\text{Ar}\ (\text{semi-strict,strict})}^{(\text{const,const-near,near})}(A^1, E^1, C^1, A^2, f_0, f_1, \mathscr{C}^1)$, iff it is a (constrained) (semi-strict, strict) $E^1$-$E^1$-$\mathscr{C}^1$-$\mathscr{C}^1$ m-atlas (near) morphism of $A^1$ to $A^2$ in the coordinate spaces $C^1$, $C^1$.

Let $\mathscr{C}^i$, $i = 1, 2$, be a model category, $C^i \overset{\text{Ob}}{\in} \mathscr{C}^i$, $E^i$ be a topological space, $A^i$ be an m-atlas of $E^i$ in the coordinate space $C^i$, $f_0 \colon E^1 \longrightarrow E^2$ be continuous, $f_1 \colon C^1 \longrightarrow C^2$ be a model function and $f \overset{\text{def}}{=} (f_0 \colon E^1 \longrightarrow E^2, f_1 \colon C^1 \longrightarrow C^2)$.

$f$ is a (constrained) (semi-strict, strict) $E^1$-$E^2$ m-atlas (near) morphism of $A^1$ to $A^2$ in the coordinate spaces $C^1$, $C^2$, abbreviated as

$\text{isAtl}_{\text{Ar}\ (\text{semi-strict,strict})}^{(\text{const,const-near,near})}(A^1, E^1, C^1, A^2, E^2, C^2, f_0, f_1)$, and a (constrained) (semi-strict, strict) $E^1$-$E^2$ m-atlas (near) morphism of $A^1$ to $A^2$ in the coordinate model categories $\mathscr{C}^1$, $\mathscr{C}^2$, abbreviated as

$\text{isAtl}_{\text{Ar}\ (\text{semi-strict,strict})}^{(\text{const,const-near,near})}(A^1, E^1, \mathscr{C}^1, A^2, E^2, \mathscr{C}^2, f_0, f_1)$, iff it is a (constrained) (semi-strict, strict) $E^1_{\mathbf{triv}}$-$E^2_{\mathbf{triv}}$ m-atlas (near) morphism of $A^1$ to $A^2$ in the coordinate model spaces $C^1$, $C^2$.

$f$ is a (constrained) (semi-strict, strict) $E^1$-$E^2$ m-atlas (constrained) (near) morphism of $A^1$ to $A^2$ in the coordinate space $C^1$, abbreviated as

$\text{isAtl}_{\text{Ar}\ (\text{semi-strict,strict})}^{(\text{const,const-near,near})}(A^1, E^1, C^1, A^2, E^2, f_0, f_1)$, iff it is a (constrained) (semi-strict, strict) $E^1$-$E^2$ m-atlas (near) morphism of $A^1$ to $A^2$ in the coordinate spaces $C^1$, $C^1$.

$f$ is a (constrained) (semi-strict, strict) $E^1$-$E^2$-$\mathscr{C}^1$-$\mathscr{C}^2$ m-atlas (near) morphism of $A^1$ to $A^2$ in the coordinate spaces $C^1$, $C^2$, abbreviated as

$\text{isAtl}_{\text{Ar}\ (\text{semi-strict,strict})}^{(\text{const,const-near,near})}(A^1, E^1, C^1, A^2, E^2, C^2, f_0, f_1, \mathscr{C}^1, \mathscr{C}^2)$, iff it is a (constrained)

(semi-strict, strict) $E^1$-$E^2$-$\mathscr{C}^1$-$\mathscr{C}^2$ m-atlas (near) morphism of $A^1$ to $A^2$ in the co-
$\underset{\text{triv triv}}{}$
ordinate model spaces $C^1, C^2$.

$f$ is a (constrained) (semi-strict, strict) $E^1$-$E^2$-$\mathscr{C}^1$ m-atlas (near) morphism of $A^1$ to $A^2$ in the coordinate model category $\mathscr{C}^1$, abbreviated as
$\text{isAtl}_{\text{Ar}}^{(\text{const,const}-\text{near,near})}(A^1, E^1, C^1, A^2, E^2, f_0, f_1, \mathscr{C}^1)$, iff it is a (constrained) (semi-
$\underset{(\text{semi-strict,strict})}{}$
strict, strict) $E^1$-$E^2$-$\mathscr{C}^1$-$\mathscr{C}^1$ m-atlas (near) morphism of $A^1$ to $A^2$ in the coordinate model categories $\mathscr{C}^1, \mathscr{C}^1$.

Similar definitions apply with restrictions on the admissible atlases, i.e., $f$ is a (constrained) (semi-strict, strict) full (semi-maximal, maximal, full semi-maximal, full maximal) $\mathscr{C}^1$-$\mathscr{C}^2$-$E^1$-$E^2$-$\mathscr{C}^1$-$\mathscr{C}^2$ m-atlas (near) morphism of $A^1$ to $A^2$ in the coordinate spaces $C^1$, abbreviated as

$$\text{isAtl}_{\text{Ar}}^{(\text{const,const}-\text{near,near})}(A^1, E^1, C^1, A^2, E^2, C^2, f_0, f_1, \mathscr{C}^1, \mathscr{C}^2, \mathscr{C}^1, \mathscr{C}^2) \iff$$
full (S-max,max,S-max-full,max-full)
(semi-strict,strict)

$$\text{isAtl}_{\text{Ar}}^{(\text{const,const}-\text{near,near})}(A^1, E^1, C^1, A^2, E^2, C^2, f_0, f_1, \mathscr{C}^1, \mathscr{C}^2, \mathscr{C}^1, \mathscr{C}^2) \wedge$$
semi-strict (strict)

$$\text{isAtl}_{\text{Ob}} \quad (A^1, E^1, C^1) \wedge$$
full (S-max,max,S-max-full,max-full)

$$\text{isAtl}_{\text{Ob}} \quad (A^2, E^2, C^2)$$
full (S-max,max,S-max-full,max-full)

**Definition 11.11** (M-atlas identity morphisms). Let $\mathscr{C}^i, \mathscr{C}^i$, $i = 1, 2$, be model categories, $E^i \overset{\text{Ob}}{\in} \mathscr{C}^i$, $E^1 \overset{\text{mod}}{\subseteq} E^2$, $C^i \overset{\text{Ob}}{\in} \mathscr{C}^i$, $C^1 \overset{\text{mod}}{\subseteq} C^2$, $A^i$ be an m-atlas of $E^i$ in the coordinate space $C^i$ and $A^1 \subseteq A^2$, then the identity morphism of $(A^1, E^1, C^1)$ to $(A^2, E^2, C^2)$ is

$$\text{Id}_{(A^1,E^1,C^1),(A^2,E^2,C^2)} \overset{\text{def}}{=} \left(\text{ID}_{(E^1,C^1),(E^2,C^2)}, (A^1, E^1, C^1), (A^2, E^2, C^2)\right) \quad (11.9)$$

The identity morphism of $(A^i, E^i, C^i)$ is

$$\text{Id}_{(A^i,E^i,C^i)} \overset{\text{def}}{=} \left((\text{Id}_{E^i}, \text{Id}_{C^i}), (A^i, E^i, C^i), (A^i, E^i, C^i)\right) \quad (11.10)$$

Let $E^i$ be a topological space, $E^1 \subseteq E^2$, $\mathscr{C}^i$, $i = 1, 2$, be a model category, $C^i \overset{\text{Ob}}{\in} \mathscr{C}^i$, $C^1 \overset{\text{mod}}{\subseteq} C^2$, $A^i$ be an m-atlas of $E^i$ in the coordinate space $C^i$ and $A^1 \subseteq A^2$, then the identity morphism of $(A^1, E^1, C^1)$ to $(A^2, E^2, C^2)$ is

$$\text{Id}_{(A^1,E^1,C^1),(A^2,E^2,C^2)} \overset{\text{def}}{=} \left(\text{ID}_{(E^1,C^1),(E^2,C^2)}, (A^1, E^1, C^1), (A^2, E^2, C^2)\right) \quad (11.11)$$

The identity morphism of $(A^i, E^i, C^i)$ is

$$\text{Id}_{(A^i,E^i,C^i)} \overset{\text{def}}{=} \left((\text{Id}_{E^i}, \text{Id}_{C^i}), (A^i, E^i, C^i), (A^i, E^i, C^i)\right) \quad (11.12)$$

### 11.1.4 Proclamations on m-atlas (near) morphisms

**Lemma 11.12** (M-atlas (near) morphisms). *Let*

1. $\mathscr{E}^i, \mathscr{C}^i, \mathscr{E}'^i, \mathscr{C}'^i, i = 1, 2$, *be model categories*

2. $\mathscr{E}^i \overset{full-cat}{\subseteq} \mathscr{E}'^i, \mathscr{C}^i \overset{full-cat}{\subseteq} \mathscr{C}'^i$

3. $\boldsymbol{E}^i \overset{Ob}{\in} \mathscr{E}^i, \boldsymbol{C}^i \overset{Ob}{\in} \mathscr{C}^i, \boldsymbol{C}'^i \overset{Ob}{\in} \mathscr{C}'^i, \boldsymbol{C}^i \overset{mod}{\subseteq} \boldsymbol{C}'^i$

4. $f_1': \boldsymbol{C}'^1 \longrightarrow \boldsymbol{C}'^2$ *be a model function such that* $f_1'[\boldsymbol{C}^1] \subseteq \boldsymbol{C}^2$ *and* $f_1 \overset{def}{=} f_1'\!\restriction_{\boldsymbol{C}^1, \boldsymbol{C}^2}$ *is a model function*

5. $\boldsymbol{A}^i$ *be an M-atlas of* $\boldsymbol{E}^i$ *in the coordinate space* $\boldsymbol{C}^i$

6. $\boldsymbol{f} \overset{def}{=} (f_0: \boldsymbol{E}^1 \longrightarrow \boldsymbol{E}^2, f_1: \boldsymbol{C}^1 \longrightarrow \boldsymbol{C}^2)$ *be a (semi-strict, strict)* $(\mathscr{E}^1\text{-}\mathscr{E}^2\text{-})\boldsymbol{E}^1\text{-}\boldsymbol{E}^2(\text{-}\mathscr{C}^1\text{-}\mathscr{C}^2)$ *m-atlas (near) morphism of* $\boldsymbol{A}^1$ *to* $\boldsymbol{A}^2$ *in the coordinate spaces* $\boldsymbol{C}^1, \boldsymbol{C}^2$

7. $\boldsymbol{f}' \overset{def}{=} (f_0: \boldsymbol{E}^1 \longrightarrow \boldsymbol{E}^2, f_1': \boldsymbol{C}^1 \longrightarrow \boldsymbol{C}'^2)$,
   $\boldsymbol{f}'' \overset{def}{=} (f_0: \boldsymbol{E}^1 \longrightarrow \boldsymbol{E}^2, f_1': \boldsymbol{C}'^1 \longrightarrow \boldsymbol{C}'^2)$ *be pairs of model functions*

*Then*

1. *If* $\boldsymbol{f}$ *is an m-atlas morphism then* $\boldsymbol{f}$ *is an m-atlas near morphism.*

   *Proof.* Let $(U^i, V^i, \phi^i: U^i \overset{\tilde{=}}{\rightarrowtail\!\!\!\twoheadrightarrow} V^i) \in \boldsymbol{A}^i, i = 1, 2$, be charts with $I \overset{def}{=} U^1 \cap f_0^{-1}[U^2] \neq \emptyset, u^1 \in I$ and $u^2 \overset{def}{=} f_0(u^1)$. Let $(U'^1, V'^1, \phi^1: U'^1 \overset{\tilde{=}}{\rightarrowtail\!\!\!\twoheadrightarrow} V'^1) \in \boldsymbol{A}^1$ and $(U'^2, \hat{V}'^2, \phi'^2: U'^2 \overset{\tilde{=}}{\rightarrowtail\!\!\!\twoheadrightarrow} \hat{V}'^2) \in \boldsymbol{A}^2$ be charts as in definition 11.4 (M-atlas morphisms for model spaces) on page 53. $U'^2 \subseteq I$ and eqs. (11.2) to (11.7) of definition 11.1 (M-atlas near morphisms for model spaces) on page 51 hold. □

2. $\boldsymbol{f}'$ *is an* $\boldsymbol{E}^1\text{-}\boldsymbol{E}^2$ *m-atlas (near) morphism of* $\boldsymbol{A}^1$ *to* $\boldsymbol{A}^2$ *in the coordinate spaces* $\boldsymbol{C}^1$, $\boldsymbol{C}'^2$ *and* $\boldsymbol{f}''$ *is an* $\boldsymbol{E}^1\text{-}\boldsymbol{E}^2$ *m-atlas (near) morphism of* $\boldsymbol{A}^1$ *to* $\boldsymbol{A}^2$ *in the coordinate spaces* $\boldsymbol{C}'^1, \boldsymbol{C}'^2$.

   *Proof.* $f_0: \boldsymbol{E}^1 \longrightarrow \boldsymbol{E}^2$ is a model function by definition 11.1 (M-atlas near morphisms for model spaces) on page 51 ( definition 11.4 (M-atlas morphisms for model spaces) on page 53).

   $f_1': \boldsymbol{C}'^1 \longrightarrow \boldsymbol{C}'^2$ and $f_1': \boldsymbol{C}^1 \longrightarrow \boldsymbol{C}'^2$ are model functions by hypothesis.

   $\boldsymbol{C}^i \overset{mod}{\subseteq} \boldsymbol{C}'^i$ by hypothesis,

   $C^i$ is an open subset of $C'^i$ by lemma 2.8 (Model subspaces) on page 20.

   The expanded model space does not change the atlases or the commutation. □

3. *If $f$ is a strict ($\mathscr{E}^1$-$\mathscr{E}^2$-)$E^1$-$E^2$(-$\mathscr{C}^1$-$\mathscr{C}^2$) m-atlas morphism of $A^1$ to $A^2$ in the coordinate spaces $C^1$, $C^2$ then $f$ is a semi-strict ($\mathscr{E}^1$-$\mathscr{E}^2$-)$E^1$-$E^2$(-$\mathscr{C}^1$-$\mathscr{C}^2$) m-atlas morphism of $A^1$ to $A^2$ in the coordinate spaces $C^1$, $C^2$*

   *Proof.* There are four overlapping cases:

   **$E^i$ without category:** $f_0$ is an m-morphism of $E^1$ to $E^2$ by definition 11.9 (Strict m-atlas (near) morphisms) on page 56. $f_0$ is a local m-morphism of $E^1$ to $E^2$ by lemma 7.15 ((Local) m-morphisms) on page 32.

   **$C^i$ without category:** $f_1$ is an m-morphism of $C^1$ to $C^2$ by definition 11.9. $f_1$ is a local m-morphism of $C^1$ to $C^2$ by lemma 7.15.

   **$E^i$ with category:** $f_0$ is an $\mathscr{E}^1$-$\mathscr{E}^2$ m-morphism of $E^1$ to $E^2$ by definition 11.9. $f_0$ is a local $\mathscr{E}^1$-$\mathscr{E}^2$ m-morphism of $E^1$ to $E^2$ by lemma 7.15.

   **$C^i$ with category:** $f_1$ is a $\mathscr{C}^1$-$\mathscr{C}^2$ m-morphism of $C^1$ to $C^2$ by definition 11.9. $f_1$ is a local $\mathscr{C}^1$-$\mathscr{C}^2$ m-morphism of $C^1$ to $C^2$ by lemma 7.15.

   □

4. *If $f$ is a semi-strict (strict) ($\mathscr{E}^1$-$\mathscr{E}^2$-)$E^1$-$E^2$(-$\mathscr{C}^1$-$\mathscr{C}^2$) m-atlas (near) morphism of $A^1$ to $A^2$ in the coordinate spaces $C^1$, $C^2$ then $f'$ is a semi-strict (strict) ($\mathscr{E}'^1$-$\mathscr{E}'^2$-)$E^1$-$E^2$(-$\mathscr{C}'^1$-$\mathscr{C}'^2$) m-atlas (near) morphism of $A^1$ to $A^2$ in the coordinate spaces $C^1$, $C'^2$.*

   *Proof.* There are two cases.

   **Semi-strict:** There are four overlapping cases.

   > **$E^i$ without category:** $f_0$ is a (strict) local m-morphism of $E^1$ to $E^2$ by definition 11.8 (Semi-strict m-atlas (near) morphisms) on page 55.
   >
   > **$E^i$ with category:** $f_0$ is a (strict) local $\mathscr{E}^1$-$\mathscr{E}^2$ m-morphism of $E^1$ to $E^2$ by definition 11.8 (Semi-strict m-atlas (near) morphisms) on page 55. $\mathscr{E}^i \overset{\text{full-cat}}{\subseteq} \mathscr{E}'^i$ by hypothesis. $f_0$ is a (strict) local $\mathscr{E}'^1$-$\mathscr{E}'^2$ m-morphism of $E^1$ to $E^2$ by lemma 7.15 ((Local) m-morphisms) on page 32.
   >
   > **$C^i$ without category:** $f_1$ is a (strict) local m-morphism of $C^1$ to $C^2$ by definition 11.8 (Semi-strict m-atlas (near) morphisms) on page 55. $C^2 \overset{\text{mod}}{\subseteq} C'^2$ by hypothesis. $f_1$ is a (strict) local m-morphism of $C^1$ to $C'^2$ by lemma 7.15 ((Local) m-morphisms) on page 32.
   >
   > **$C^i$ with category:** $f_0$ is a local $\mathscr{C}^1$-$\mathscr{C}^2$ m-morphism of $C^1$ to $C^2$ by definition 11.8 (Semi-strict m-atlas (near) morphisms) on page 55. $\mathscr{C}^i \overset{\text{full-cat}}{\subseteq} \mathscr{C}'^i$ by hypothesis. $f_0$ is a (strict) local $\mathscr{C}'^1$-$\mathscr{C}'^2$ m-morphism of $C^1$ to $C^2$ by lemma 7.15 ((Local) m-morphisms) on page 32.

   **Strict:** There are four overlapping cases.

$E^i$ **without category:** $f_0$ is an m-morphism of $E^1$ to $E^2$ by definition 11.8 (Semi-strict m-atlas (near) morphisms) on page 55.

$E^i$ **with category:** $f_0$ is an $\mathscr{E}^1$-$\mathscr{E}^2$ m-morphism of $E^1$ to $E^2$ by definition 11.8 (Semi-strict m-atlas (near) morphisms) on page 55. $\mathscr{E}^i \stackrel{\text{full-cat}}{\subseteq} \mathscr{E}'^i$ by hypothesis. $f_0$ is an $\mathscr{E}'^1$-$\mathscr{E}'^2$ m-morphism of $E^1$ to $E^2$ by lemma 7.15 ((Local) m-morphisms) on page 32.

$C^i$ **without category:** $f_1$ is an m-morphism of $C^1$ to $C^2$ by definition 11.8 (Semi-strict m-atlas (near) morphisms) on page 55. $C^2 \stackrel{\text{mod}}{\subseteq} C'^2$ by hypothesis. $f_1$ is an m-morphism of $C^1$ to $C'^2$ by lemma 7.15 ((Local) m-morphisms) on page 32.

$C^i$ **with category:** $f_0$ is an $\mathscr{C}^1$-$\mathscr{C}^2$ m-morphism of $C^1$ to $C^2$ by definition 11.8 (Semi-strict m-atlas (near) morphisms) on page 55. $\mathscr{C}^i \stackrel{\text{full-cat}}{\subseteq} \mathscr{C}'^i$ by hypothesis. $f_0$ is a $\mathscr{C}'^1$-$\mathscr{C}'^2$ m-morphism of $C^1$ to $C^2$ by lemma 7.15 ((Local) m-morphisms) on page 32.

□

5. If $f$ is a semi-strict ($\mathscr{E}^1$-$\mathscr{E}^2$-)$E^1$-$E^2$(-$\mathscr{C}^1$-$\mathscr{C}^2$) m-atlas (near) morphism of $A^1$ to $A^2$ in the coordinate spaces $C^1$, $C^2$ and $f'_1$ is a (strict) local m-morphism (a $\mathscr{C}'^1$-$\mathscr{C}'^2$ m-morphism) of $C'^1$ to $C'^2$, then $f''$ is a semi-strict ($\mathscr{E}'^1$-$\mathscr{E}'^2$-)$E^1$-$E^2$(-$\mathscr{C}'^1$-$\mathscr{C}'^2$) m-atlas (near) morphism of $A^1$ to $A^2$ in the coordinate spaces $C'^1$, $C'^2$.

*Proof.* There are four overlapping cases:

$E^i$ **without categories:** $f_0$ is a (strict) local m-morphism of $E^1$ to $E^2$ by definition 11.8 (Semi-strict m-atlas (near) morphisms) on page 55.

$E^i$ **with categories:** $f_0$ is a (strict) local $\mathscr{E}^1$-$\mathscr{E}^2$ m-morphism of $E^1$ to $E^2$ by definition 11.8. $\mathscr{E}^i \stackrel{\text{full-cat}}{\subseteq} \mathscr{E}'^i$ by hypothesis. $f_0$ is a (strict) local $\mathscr{E}'^1$-$\mathscr{E}'^2$ m-morphism of $E^1$ to $E^2$ by lemma 7.15 ((Local) m-morphisms) on page 32.

$C'^i$ **without categories:** $f'_1$ is a (strict) local m-morphism of $C'^1$ to $C'^2$ by hypothesis.

$C'^i$ **with categories:** $f'_1$ is a local $\mathscr{C}'^1$-$\mathscr{C}'^2$ m-morphism of $C^1$ to $C^2$ by hypothesis.

□

6. If $f$ is a strict ($\mathscr{E}^1$-$\mathscr{E}^2$-)$E^1$-$E^2$(-$\mathscr{C}^1$-$\mathscr{C}^2$) m-atlas (near) morphism of $A^1$ to $A^2$ in the coordinate spaces $C^1$, $C^2$ and $f'_1$ and $f'_1$ is an m-morphism (a $\mathscr{C}'^1$-$\mathscr{C}'^2$ m-morphism) of $C'^1$ to $C'^2$, then $f''$ is a strict ($\mathscr{E}'^1$-$\mathscr{E}'^2$-)$E^1$-$E^2$(-$\mathscr{C}'^1$-$\mathscr{C}'^2$) m-atlas (near) morphism of $A^1$ to $A^2$ in the coordinate spaces $C'^1$, $C'^2$.

*Proof.* There are four overlapping cases:

$E^i$ **without categories:** $f_0$ is an m-morphism of $E^1$ to $E^2$ by definition 11.4 (M-atlas morphisms for model spaces) on page 53.

$E^i$ **with categories:** $f_0$ is an $\mathscr{C}^1$-$\mathscr{C}^2$ m-morphism of $E^1$ to $E^2$ by definition 11.4. $\mathscr{C}^i \overset{\text{full-cat}}{\subseteq} \mathscr{C}'^i$ by hypothesis. $f_0$ is an $\mathscr{C}'^1$-$\mathscr{C}'^2$ m-morphism of $E^1$ to $E^2$ by lemma 7.15 ((Local) m-morphisms) on page 32.

$C'^i$ **without categories:** $f_1'$ is an m-morphism of $C'^1$ to $C'^2$ by hypothesis.

$C'^i$ **with categories:** $f_1'$ is a $\mathscr{C}'^1$-$\mathscr{C}'^2$ m-morphism of $C^1$ to $C^2$ by hypothesis.

□

7. If $f$ is a strict $(\mathscr{C}^1$-$\mathscr{C}^2$-$)E^1$-$E^2(\text{-}\mathscr{C}^1$-$\mathscr{C}^2)$ m-atlas (near) morphism of $A^1$ to $A^2$ in the coordinate spaces $C^1$, $C^2$ then $f$ is a semi-strict $(\mathscr{C}^1$-$\mathscr{C}^2$-$)E^1$-$E^2(\text{-}\mathscr{C}^1$-$\mathscr{C}^2)$ m-atlas (near) morphism of $A^1$ to $A^2$ in the coordinate spaces $C^1$, $C^2$

*Proof.* There are four overlapping cases:

$E^i$ **without category:** $f_0$ is an m-morphism of $E^1$ to $E^2$ by definition 11.9 (Strict m-atlas (near) morphisms) on page 56. $f_0$ is a local m-morphism of $E^1$ to $E^2$ by lemma 7.15 ((Local) m-morphisms) on page 32.

$C^i$ **without category:** $f_1$ is an m-morphism of $C^1$ to $C^2$ by definition 11.9. $f_1$ is a local m-morphism of $C^1$ to $C^2$ by lemma 7.15.

$E^i$ **with category:** $f_0$ is an $\mathscr{C}^1$-$\mathscr{C}^2$ m-morphism of $E^1$ to $E^2$ by definition 11.9. $f_0$ is a local $\mathscr{C}^1$-$\mathscr{C}^2$ m-morphism of $E^1$ to $E^2$ by lemma 7.15.

$C^i$ **with category:** $f_1$ is a $\mathscr{C}^1$-$\mathscr{C}^2$ m-morphism of $C^1$ to $C^2$ by definition 11.9. $f_1$ is a local $\mathscr{C}^1$-$\mathscr{C}^2$ m-morphism of $C^1$ to $C^2$ by lemma 7.15.

□

8. If $A^i$, $i = 1, 2$, is semi-maximal and $f$ is a (semi-strict, strict) $(\mathscr{C}^1$-$\mathscr{C}^2$-$)E^1$-$E^2(\text{-}\mathscr{C}^1$-$\mathscr{C}^2)$ m-atlas near morphism of $A^1$ to $A^2$ in the coordinate spaces $C^1$, $C^2$ then it is a (semi-strict, strict) $(\mathscr{C}^1$-$\mathscr{C}^2$-$)E^1$-$E^2(\text{-}\mathscr{C}^1$-$\mathscr{C}^2)$ m-atlas morphism of $A^1$ to $A^2$ in the coordinate spaces $C^1$, $C^2$.

*Proof.*

$f$ **is an M-atlas morphism.** Let $(U^i, V^i, \phi^i\colon U^i \rightarrowtail\!\!\!\!\!\overset{\sim}{\twoheadrightarrow} V^i) \in A^i$, $i = 1, 2$. Diagram (11.1) in definition 11.1 (M-atlas near morphisms for model spaces) on page 51 is M-locally nearly commutative in $C^2$, i.e., for any $e^1 \in I \overset{\text{def}}{=} U^1 \cap f_0^{-1}[U^2]$ there are model neighborhoods $U'^1 \subseteq I, V'^1 \subseteq V^1, U'^2 \subseteq U^2, V'^2 \subseteq V^2, \hat{V}'^2 \subseteq C^2$ and an isomorphism $\hat{f}\colon V'^2 \rightarrowtail\!\!\!\!\!\overset{\sim}{\twoheadrightarrow} \hat{V}'^2$ such that eqs. (11.2) to (11.7) in definition 11.1 hold, as shown in figs. 4 to 5.

$f_0[U'^1] \subseteq U'^2$ by eq. (11.3) definition 11.1 . $f_1[V'^1] \subseteq \hat{V}'^2$ by eq. (11.5) definition 11.1 .

Since $A^i$ is semi-maximal, $(U'^i, V'^i, \phi^i \colon U'^i \overset{\tilde{=}}{\rightarrowtail\!\!\!\twoheadrightarrow} V'^i) \in A^i$ and $(U'^2, \hat{f}[V'^2], \hat{f} \circ \phi^2) \in A^2$.

Define $\phi'^2 \colon U'^2 \overset{\tilde{=}}{\rightarrowtail\!\!\!\twoheadrightarrow} \hat{f}[V'^2] \overset{\text{def}}{=} \hat{f} \circ \phi^2$. Then Diagram (11.8) is commutative.

**$f$ is (semi-strict, strict).** The definitions of strict and semi-strict are the same for m-atlas near morphisms and for m-atlas morphisms.

□

9. Let $E^i \overset{\text{def}}{=} \text{Top}(\boldsymbol{E}^i)$, $i = 1, 2$. If each $\boldsymbol{C}^i$ is fine grained then $\hat{f} \overset{\text{def}}{=} (f_0 \colon E^1 \longrightarrow E^2), f_1 \colon \boldsymbol{C}^1 \longrightarrow \boldsymbol{C}^2)$ is a (semi-strict, strict) $\boldsymbol{E}^1\text{-}\boldsymbol{E}^2$ m-atlas (near) morphism of $\boldsymbol{A}^1$ to $\boldsymbol{A}^2$ in the coordinate spaces $\boldsymbol{C}^1$, $\boldsymbol{C}^2$.

*Proof.*

(a) $\boldsymbol{A}^i$ is an m-atlas of $\boldsymbol{E}^i$ in the coordinate space $\boldsymbol{C}^i$ by lemma 10.2 (M-atlases) on page 47.

(b) $f_0$ is a model function, hence continuous. $f_1$ is a model function by hypothesis.

(c) Any model neighborhood of $\boldsymbol{E}^i$ is open, hence a model neighborhood of $\boldsymbol{E}^i_{\text{triv}}$.

(d) Let $(U^i, V^i, \phi^i) \in \boldsymbol{A}^i$, $i = 1, 2$. All of the model neighborhoods of $\boldsymbol{E}^i$ mentioned in definition 11.1 (M-atlas near morphisms for model spaces) on page 51 and definition 11.4 (M-atlas morphisms for model spaces) on page 53 are also model neighborhoods of $\boldsymbol{E}^i_{\text{triv}}$.

□

10. Let $\boldsymbol{f}$ be a (semi-strict, strict) $\mathcal{E}^1\text{-}\mathcal{E}^2\text{-}\boldsymbol{E}^1\text{-}\boldsymbol{E}^2\text{-}\mathcal{C}^1\text{-}\mathcal{C}^2$ ($\boldsymbol{E}^1\text{-}\boldsymbol{E}^2\text{-}\mathcal{C}^1\text{-}\mathcal{C}^2$) m-atlas (near) morphism of $\boldsymbol{A}^1$ to $\boldsymbol{A}^2$ in the coordinate spaces $\boldsymbol{C}^1$, $\boldsymbol{C}^2$.

$\boldsymbol{f}'$ is a (semi-strict, strict) $\mathcal{E}^1\text{-}\mathcal{E}'^2\text{-}\boldsymbol{E}^1\text{-}\boldsymbol{E}^2\text{-}\mathcal{C}^1\text{-}\mathcal{C}'^2$ ($\boldsymbol{E}^1\text{-}\boldsymbol{E}^2\text{-}\mathcal{C}^1\text{-}\mathcal{C}'^2$) m-atlas (near) morphism of $\boldsymbol{A}^1$ to $\boldsymbol{A}^2$ in the coordinate spaces $\boldsymbol{C}'^1$, $\boldsymbol{C}'^2$.

If $\boldsymbol{f}$ is semi-strict and $f'_1$ is a (strict) local m-morphism of $\boldsymbol{C}'^1$ to $\boldsymbol{C}'^2$, then $\boldsymbol{f}''$ is a semi-strict $(\mathcal{E}^1\text{-}\mathcal{E}^2\text{-})\boldsymbol{E}^1\text{-}\boldsymbol{E}^2(\text{-}\mathcal{C}'^1\text{-}\mathcal{C}'^2)$ m-atlas (near) morphism of $\boldsymbol{A}^1$ to $\boldsymbol{A}^2$ in the coordinate spaces $\boldsymbol{C}'^1$, $\boldsymbol{C}'^2$.

*Proof.*

- **$f'$ and $f''$ are M-atlas (near) morphisms.** The expanded categories and expanded model space do not change the atlases or the commutation.

- **$f'$ is (semi-strict, strict).** $\mathscr{C}^i \overset{\text{full-cat}}{\subseteq} \mathscr{C}'^i$ by hypothesis, $C^i \overset{\text{mod}}{\subseteq} C'^i$ by hypothesis, If $f_0$ is a local m-morphism of $\mathscr{C}^2$ then $f_0$ is a local m-morphism of $\mathscr{C}'^2$. If $f_1$ is a local m-morphism of $\mathscr{C}^2$ then $f_1$ is a local m-morphism of $\mathscr{C}'^2$. If $f_0 \overset{\text{Ar}}{\in} \mathscr{C}^2$ then $f_0 \overset{\text{Ar}}{\in} \mathscr{C}'^2$. If $f_1 \overset{\text{Ar}}{\in} \mathscr{C}^2$ then $f_1 \overset{\text{Ar}}{\in} \mathscr{C}'^2$.

$\square$

**Corollary 11.13** (M-atlas morphisms)**.** *Let*

1. *$\mathscr{C}^i, \mathscr{C}'^i, i = 1, 2,$ be model categories*

2. *$\mathscr{C}^i \overset{\text{full-cat}}{\subseteq} \mathscr{C}'^i$*

3. *$C^i \overset{\text{Ob}}{\in} \mathscr{C}^i$*

4. *$C'^i \overset{\text{Ob}}{\in} \mathscr{C}'^i$*

5. *$C^i \overset{\text{mod}}{\subseteq} C'^i$*

6. *$f'_1 \colon C'^1 \longrightarrow C'^2$ be a model function such that $f'_1[C^1] \subseteq C^2$ and $f_1 \overset{\text{def}}{=} f'_1 \restriction_{C^1, C^2}$ is a model function*

7. *$E^i$ be a topological space*

8. *$A^i$ be an m-atlas of $E^i$ in the coordinate space $C^i$*

9. *$f \overset{\text{def}}{=} (f_0 \colon E^1 \longrightarrow E^2, f_1 \colon C^1 \longrightarrow C^2)$ be a (semi-strict, strict) ($\mathscr{C}^1$-$\mathscr{C}^2$-)$E^1$-$E^2$(-$\mathscr{C}^1$-$\mathscr{C}^2$) m-atlas (near) morphism of $A^1$ to $A^2$ in the coordinate spaces $C^1, C^2$*

If $f \overset{\text{def}}{=} (f_0 \colon E^1 \longrightarrow E^2, f_1 \colon C^1 \longrightarrow C^2)$ is a (semi-strict, strict) $E^1$-$E^2$ m-atlas morphism of $A^1$ to $A^2$ in the coordinate spaces $C^1, C^2$ then it is a (semi-strict, strict) $E^1$-$E^2$ m-atlas near morphism of $A^1$ to $A^2$ in the coordinate spaces $C^1, C^2$

If $A^i$, $i = 1, 2$, is semi-maximal and $f$ is a (semi-strict, strict) $E^1$-$E^2$ m-atlas near morphism of $A^1$ to $A^2$ in the coordinate spaces $C^1, C^2$ then it is a (semi-strict, strict) $E^1$-$E^2$ m-atlas morphism of $A^1$ to $A^2$ in the coordinate spaces $C^1, C^2$.

If $f$ is a semi-strict (strict) $E^1$-$E^2$-$\mathscr{C}^1$-$\mathscr{C}^2$ m-atlas (near) morphism of $A^1$ to $A^2$ in the coordinate spaces $C^1, C^2$ then $f$ is a semi-strict (strict) $E^1$-$E^2$-$\mathscr{C}'^1$-$\mathscr{C}'^2$ m-atlas (near) morphism of $A^1$ to $A^2$ in the coordinate spaces $C^1, C^2$.

If a (near) m-atlas morphism is equivalent to a semi-strict (strict) (near) m-atlas morphism, it need not be semi-strict (strict). However:

**Lemma 11.14** (Equivalence of m-atlas (near) morphisms)**.** *Let*

1. $\mathcal{E}^i, \mathcal{C}^i, i = 1, 2$, be model categories

2. $E^i \overset{\text{Ob}}{\in} \mathcal{E}^i, C^i \overset{\text{Ob}}{\in} \mathcal{C}^i$,

3. $A^i$ be a full m-atlas of $E^i$ in the coordinate space $C^i$,

4. $f \overset{\text{def}}{=} (f_0 \colon E^1 \longrightarrow E^2, f_1 \colon C^1 \longrightarrow C^2)$ and $g \overset{\text{def}}{=} (g_0 \colon E^1 \longrightarrow E^2, g_1 \colon C^1 \longrightarrow C^2)$ be equivalent m-atlas near morphisms of $A^1$ to $A^2$ in the coordinate spaces $C^1$, $C^2$.

If $f$ is a semi-strict $(\mathcal{E}^1\text{-}\mathcal{E}^2\text{-})E^1\text{-}E^2(\text{-}\mathcal{C}^1\text{-}\mathcal{C}^2)$ m-atlas near morphism of $A^1$ to $A^2$ in the coordinate spaces $C^1$, $C^2$, then $g$ is a semi-strict $(\mathcal{E}^1\text{-}\mathcal{E}^2\text{-})E^1\text{-}E^2(\text{-}\mathcal{C}^1\text{-}\mathcal{C}^2)$ m-atlas near morphism of $A^1$ to $A^2$ in the coordinate spaces $C^1$, $C^2$.

*Proof.* Let $v^1 \in C^1$ be an arbitrary point. Let $(U^1, V^1, \phi^1) \in A^1$ be an arbitrary chart such that $v^1 \in V^1$, $u^1 \overset{\text{def}}{=} \phi^{1-1}(v^1)$, $(U^2, V^2, \phi^2) \in A^2$ be an arbitrary chart at $u^2 \overset{\text{def}}{=} f_0(u^1)$, $v^2 \overset{\text{def}}{=} \phi^2(u^2)$ and $I \overset{\text{def}}{=} U^1 \cap f_0^{-1}[U^2]$, as shown in fig. 7 (Uncompleted equivalent m-atlas near morphisms).

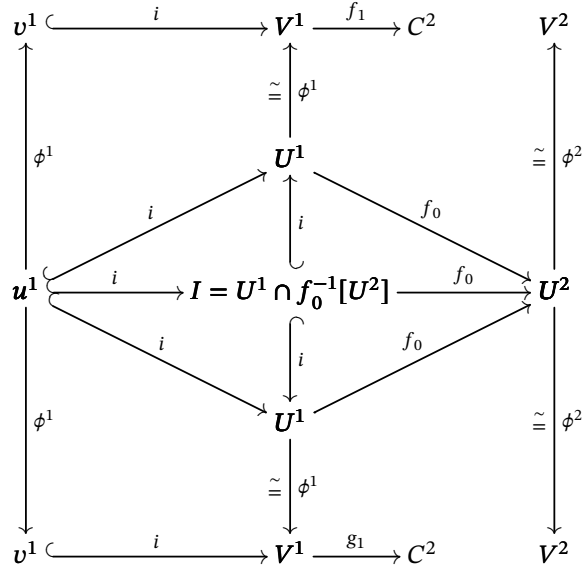

Figure 7: Uncompleted equivalent m-atlas near morphisms

Let ${}^f U'^1 \subseteq I$, ${}^f V'^1 \subseteq V^1$ ${}^f U'^2 \subseteq U^2$, ${}^f V'^2 \subseteq V^2$ and ${}^f \hat{V}'^2 \subseteq C^2$ be model neighborhoods and $\hat{f} \colon {}^f V'^2 \overset{\tilde{=}}{\rightarrowtail\!\!\!\twoheadrightarrow} {}^f \hat{V}'^2$ an isomorphism of $C^2$ ($\mathcal{C}^2$) as in definition 11.1 (M-atlas near morphisms for model spaces) on page 51.

Let ${}^gU'^1 \subseteq I$, ${}^gV'^1 \subseteq V^1$ ${}^gU'^2 \subseteq U^2$. ${}^gV'^2 \subseteq V^2$ and ${}^g\hat{V}'^2 \subseteq \boldsymbol{C}^2$ be model neighborhoods and $\hat{g}\colon {}^gV'^2 \rightarrowtail\overset{\tilde{=}}{\twoheadrightarrow} {}^g\hat{V}'^2$ an isomorphism of $\boldsymbol{C}^2$ ($\mathscr{C}^2$) as in definition 11.1.

${}^f\hat{V}'^2 \cap {}^g\hat{V}'^2$ is a model neighborhood of $\boldsymbol{C}^2$ by item 2 of definition 2.1 (Model spaces) on page 18.

There are four overlapping cases.

**$E^i$ without category:** $f_0 = g_0$ by definition 11.7 (Equivalence of m-atlas (near) morphisms) on page 55. $f_0$ is a (strict) local m-morphism of $\boldsymbol{E}^1$ to $\boldsymbol{E}^2$ by definition 11.8 (Semi-strict m-atlas (near) morphisms) on page 55.

**$C^i$ without category:** $f_1$ is a (strict) local m-morphism of $\boldsymbol{C}^1$ to $\boldsymbol{C}^2$ and $\phi^2 \circ f_0 \circ \phi^{1-1}\colon \phi^1[I] \longrightarrow V^2$ is a (strict) local m-morphism of $\phi^1[I]$ to $\boldsymbol{C}^2$ by definition 11.8 (Semi-strict m-atlas (near) morphisms) on page 55.

Let $\hat{V}'^2 \subseteq {}^f\hat{V}'^2 \cap {}^g\hat{V}'^2$ be a model neighborhod of $\boldsymbol{C}^2$ at $v^2$ such that $f_1\colon V'^1 \overset{\text{def}}{=} f_1^{-1}[\hat{V}'^2] \longrightarrow \hat{V}'^2$ is a morphism of $\boldsymbol{C}^2$.

The following, shown in figs. 8 to 9, are model neighborhoods of $\boldsymbol{C}^2$ by definition 7.1 (Model functions) on page 28:

$$V'^2 \overset{\text{def}}{=} \hat{f}^{-1}[\hat{V}'2]$$
$$U'^2 \overset{\text{def}}{=} \phi^{2-1}[V'^2]$$
$$U'^1 \overset{\text{def}}{=} f_0^{-1}[U'^2]$$

Then
$$\begin{aligned}
g_1\colon V'^1 \longrightarrow \hat{V}'^2 &= \hat{g}\colon V'^2 \longrightarrow \hat{V}'^2 \quad \circ \\
&\quad \phi^2\colon U'^2 \longrightarrow V'^2 \quad \circ \\
&\quad f_0\colon U'^1 \longrightarrow U'^2 \quad \circ \\
&\quad \phi^{1-1}\colon V'^1 \longrightarrow U'^1 \\
&= \hat{g}\colon V'^2 \longrightarrow \hat{V}'^2 \quad \circ \\
&\quad \hat{f}^{-1}\colon \hat{V}'^2 \longrightarrow V'^2 \quad \circ \\
&\quad \hat{f}\colon V'^2 \longrightarrow \hat{V}'^2 \quad \circ \\
&\quad \phi^2\colon U'^2 \longrightarrow '^2 \quad \circ \\
&\quad f_0\colon U'^1 \longrightarrow U'^2 \quad \circ \\
&\quad \phi^{1-1}\colon V'^1 \longrightarrow U'^1 \\
&= \hat{g}\colon V'^2 \longrightarrow \hat{V}'^2 \quad \circ \\
&\quad \hat{f}^{-1}\colon \hat{V}'^2 \longrightarrow V'^2 \quad \circ \\
&\quad f_1\colon V'^1 \longrightarrow \hat{V}'^2
\end{aligned}$$

is a composition of morphisms of $\boldsymbol{C}^2$, hence a morphism.

**$E^i$ with category:** $f_0 = g_0$ by definition 11.7. $f_0$ is a (strict) local $\mathscr{C}^1$-$\mathscr{C}^2$ m-morphism of $\boldsymbol{E}^1$ to $\boldsymbol{E}^2$ by definition 11.8.

**$C^i$ with category:** $f_1$ is a local $\mathscr{C}^1$-$\mathscr{C}^2$ m-morphism of $\boldsymbol{C}^1$ to $\boldsymbol{C}^2$ and $\phi^2 \circ f_0 \circ \phi^{1-1}\colon \phi^1[I] \longrightarrow V^2$ is a local $\mathscr{C}^1$-$\mathscr{C}^2$ m-morphism of $\operatorname{Mod}(\phi^1[I], \boldsymbol{C}^2)$ to $\mathrm{Mod}(V^2, \boldsymbol{C}^2)$ by definition 11.8 (Semi-strict m-atlas (near) morphisms) on page 55.

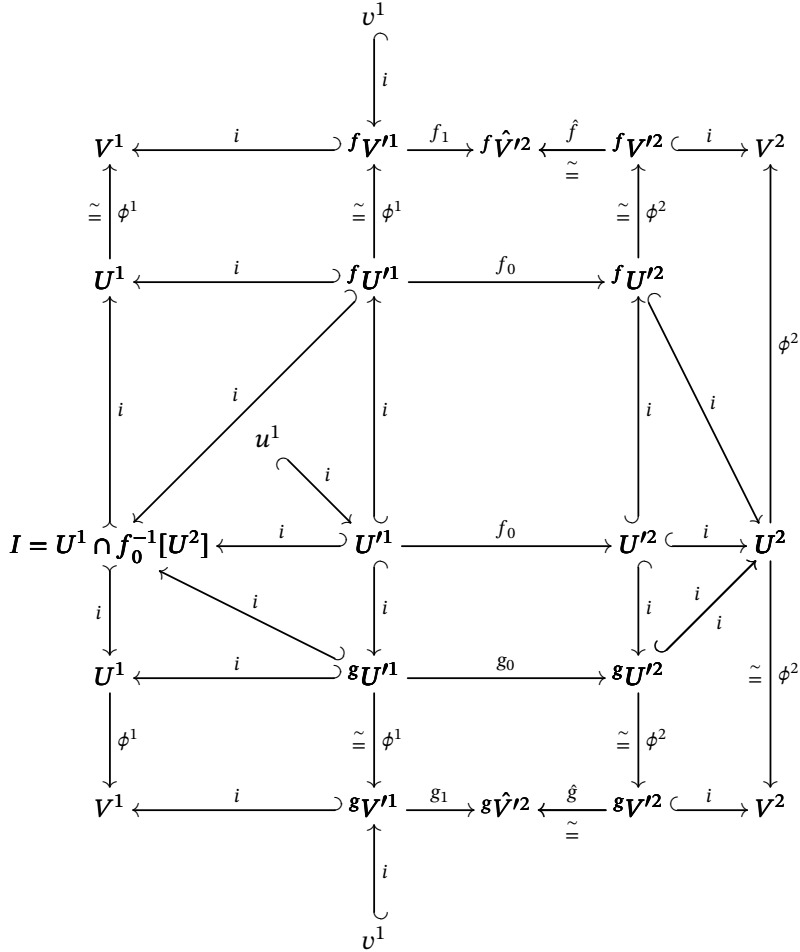

Figure 8: Completed equivalent m-atlas near morphisms

Let $U'^2 \subseteq {}^fU'^2 \cap {}^g U'^2$ be a model neighborhod of $\boldsymbol{C}^2$ at $u^2$ such that $f_0 \colon U'^1 \stackrel{\text{def}}{=} f_0^{-1}[U'^2] \longrightarrow U'^2$ is a morphism of $\mathscr{C}^2$.

Let $V'^2 \stackrel{\text{def}}{=} \phi^2[U'^2]$, ${}^f\hat{V}'^2 \stackrel{\text{def}}{=} \hat{f}[{}^fV'^2]$, ${}^g\hat{V}'^2 \stackrel{\text{def}}{=} \hat{g}[{}^gV'^2]$, ${}^fV'^1 \stackrel{\text{def}}{=} f_1^{-1}[{}^f\hat{V}'^2]$ and ${}^gV'^1 \stackrel{\text{def}}{=} g_1^{-1}[{}^g\hat{V}'^2]$. All of these are model neighborhoods by definition 7.1 (Model functions) on page 28 and the corresponding relative model subspaces of $\boldsymbol{V}^i$ are objects of $\mathscr{V}^i$ by item 5 of definition 7.7 (Model categories) on page 30. The restrictions of morphisms to them are morphisms by item 5 of definition 7.7.

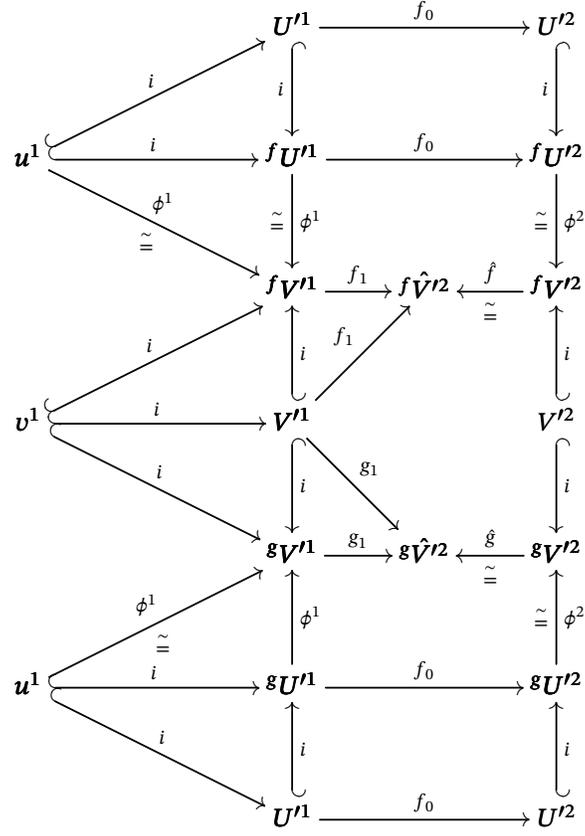

Figure 9: Completed equivalent m-atlas near morphisms

Then

$$
\begin{aligned}
g_1 \colon V'^1 \longrightarrow \hat{V}'^2 &= \hat{g} \colon V'^2 \longrightarrow \hat{V}'^2 \quad \circ \\
&\phantom{=} \phi^2 \colon U'^2 \longrightarrow V'^2 \quad \circ \\
&\phantom{=} f_0 \colon U'^1 \longrightarrow U'^2 \quad \circ \\
&\phantom{=} \phi^{1-1} \colon V'^1 \longrightarrow U'^1 \\
&= \hat{g} \colon V'^2 \longrightarrow \hat{V}'^2 \quad \circ \\
&\phantom{=} \hat{f}^{-1} \colon \hat{V}'^2 \longrightarrow V'^2 \quad \circ \\
&\phantom{=} \hat{f} \colon V'^2 \longrightarrow \hat{V}'^2 \quad \circ \\
&\phantom{=} \phi^2 \colon U'^2 \longrightarrow V'^2 \quad \circ \\
&\phantom{=} f_0 \colon U'^1 \longrightarrow U'^2 \quad \circ \\
&\phantom{=} \phi^{1-1} \colon V'^1 \longrightarrow U'^1 \\
&= \hat{g} \colon V'^2 \longrightarrow \hat{V}'^2 \quad \circ \\
&\phantom{=} \hat{f}^{-1} \colon \hat{V}'^2 \longrightarrow V'^2 \quad \circ \\
&\phantom{=} f_1 \colon V'^1 \longrightarrow \hat{V}'^2
\end{aligned}
$$

is a composition of morphisms of $\mathscr{C}^2$, hence a morphism.

$\square$

**Corollary 11.15** (Equivalence of m-atlas (near) morphisms)**.** *Let $\mathscr{E}^i, \mathscr{C}^i, i = 1, 2$, be model categories, $\mathscr{C}^2$ satisfy the restricted sheaf condition, $E^i \stackrel{Ob}{\in} \mathscr{E}^i$, $C^i \stackrel{Ob}{\in} \mathscr{C}^i$, $A^i$ be an m-atlas of $E^i$ in the coordinate space $C^i$, $A^1$ be full and $f \stackrel{def}{=} (f_0 \colon E^1 {\longrightarrow}\!\!\!\!\!\!\!\rightarrow E^2, f_1 \colon C^1 {\longrightarrow}\!\!\!\!\!\!\!\rightarrow C^2)$ and $g \stackrel{def}{=} (g_0 \colon E^1 \longrightarrow E^2, g_1 \colon C^1 \longrightarrow C^2)$ be equivalent m-atlas (near) morphisms of $A^1$ to $A^2$ in the coordinate spaces $C^1$, $C^2$.*

*If $f$ is a semi-strict $(\mathscr{E}^1\text{-}\mathscr{E}^2\text{-})E^1\text{-}E^2(\text{-}\mathscr{C}^1\text{-}\mathscr{C}^2)$ m-atlas morphism of $A^1$ to $A^2$ in the coordinate spaces $C^1$, $C^2$, then $g$ is a semi-strict $(\mathscr{E}^1\text{-}\mathscr{E}^2\text{-})E^1\text{-}E^2(\text{-}\mathscr{C}^1\text{-}\mathscr{C}^2)$ m-atlas morphism of $A^1$ to $A^2$ in the coordinate spaces $C^1$, $C^2$.*

*Proof.* $f$ and $g$ are $(\mathscr{E}^1\text{-}\mathscr{E}^2\text{-})E^1\text{-}E^2(\text{-}\mathscr{C}^1\text{-}\mathscr{C}^2)$ m-atlas near morphisms of $A^1$ to $A^2$ in the coordinate spaces $C^1$, $C^2$ by item 1 of lemma 11.12 (M-atlas (near) morphisms) on page 60. Then $g$ is a semi-strict $(\mathscr{E}^1\text{-}\mathscr{E}^2\text{-})E^1\text{-}E^2(\text{-}\mathscr{C}^1\text{-}\mathscr{C}^2)$ m-atlas near morphism of $A^1$ to $A^2$ in the coordinate spaces $C^1$, $C^2$ by lemma 11.14 (Equivalence of m-atlas (near) morphisms) on page 65. $g$ is a morphism by hypothesis. $\square$

*If $f$ is a strict $(\mathscr{E}^1\text{-}\mathscr{E}^2\text{-})E^1\text{-}E^2(\text{-}\mathscr{C}^1\text{-}\mathscr{C}^2)$ m-atlas (near) morphism of $A^1$ to $A^2$ in the coordinate spaces $C^1$, $C^2$, then $g$ is a strict $(\mathscr{E}^1\text{-}\mathscr{E}^2\text{-})E^1\text{-}E^2(\text{-}\mathscr{C}^1\text{-}\mathscr{C}^2)$ m-atlas near morphism of $A^1$ to $A^2$ in the coordinate spaces $C^1$, $C^2$.*

*Proof.* $f$ is a semi-strict $(\mathscr{E}^1\text{-}\mathscr{E}^2\text{-})E^1\text{-}E^2(\text{-}\mathscr{C}^1\text{-}\mathscr{C}^2)$ m-atlas (near) morphism of $A^1$ to $A^2$ in the coordinate spaces $C^1$, $C^2$ by item 7 of lemma 11.12 (M-atlas (near) morphisms) on page 63. $g$ is a semi-strict $(\mathscr{E}^1\text{-}\mathscr{E}^2\text{-})E^1\text{-}E^2(\text{-}\mathscr{C}^1\text{-}\mathscr{C}^2)$ m-atlas near morphism of $A^1$ to $A^2$ in the coordinate spaces $C^1$, $C^2$. The result follow from the appropriate rectricted sheaf condition of definition 2.1 (Model spaces) on page 18 and definition 7.7 (Model categories) on page 30. $\square$

**Lemma 11.16** (Composition of m-atlas (near) morphisms)**.** *Let*

1. *$\mathscr{E}^i, \mathscr{C}^i, i = 1, 2, 3$, be a model category, $E^i \stackrel{Ob}{\in} \mathscr{E}^i$, $C^i \stackrel{Ob}{\in} \mathscr{C}^i$,*

2. *$A^i$ an m-atlas of $E^i$ in the coordinate space $C^i$*

3. *$f^i \stackrel{def}{=} (f_0^i \colon E^i {\longrightarrow}\!\!\!\!\!\!\!\rightarrow E^{i+1}, f_1^i \colon C^i {\longrightarrow}\!\!\!\!\!\!\!\rightarrow C^{i+1})$ and $g^i \stackrel{def}{=} (g_0^i \colon E^i {\longrightarrow}\!\!\!\!\!\!\!\rightarrow E^{i+1}, g_1^i \colon C^i {\longrightarrow}\!\!\!\!\!\!\!\rightarrow C^{i+1})$ $i = 1, 2$, equivalent (semi-strict, strict) $E^i\text{-}E^{i+1}$ m-atlas (near) morphisms of $A^i$ to $A^{i+1}$ in the coordinate spaces $C^i$, $C^{i+1}$.*

*Then:*

1. *If $f^1$ is an $E^1\text{-}E^2$ m-atlas morphism of $A^1$ to $A^2$ in the coordinate spaces $C^1$, $C^2$ and If $f^2$ is an $E^2\text{-}E^3$ m-atlas near morphism of $A^2$ to $A^3$ in the coordinate spaces $C^2$, $C^3$ then $f^2 \stackrel{()}{\circ} f^1$ is an $E^1\text{-}E^3$ m-atlas near morphism of $A^1$ to $A^3$ in the coordinate spaces $C^1$, $C^3$.*

*Proof.* $f_0^2 \circ f_0^1$ and $f_1^2 \circ f_1^1$ are model functions by lemma 7.5 (Composition of model functions) on page 30.

Let $(U^i, V^i, \phi^i) \in \mathbf{A}^i$, $i = 1, 3$, be a chart such that $I \stackrel{\text{def}}{=} U^1 \cap (f_0^2 \circ f_0^1)^{-1}[U^3] \neq \emptyset$. It suffices to show that diagram (11.13) below is locally nearly commutative.

$$D \stackrel{\text{def}}{=} \left\{ \begin{array}{l} \phi^1 \colon I \rightarrowtail V^1, f_1^2 \circ f_1^1 \colon V^1 \longrightarrow \mathbf{C}^3, \\ f_0^2 \circ f_0^1 \colon I \longrightarrow U^3, \phi^3 \colon U^3 \xrightarrow{\tilde{=}} V^3 \end{array} \right\} \quad (11.13)$$

Let $u^1 \in I$, $u^2 \stackrel{\text{def}}{=} f_0^1(u^1)$, $u^3 \stackrel{\text{def}}{=} f_0^2(u^2)$. $(U^2, V^2, \phi^2) \in \mathbf{A}^2$ be a chart at $u^2$.

Since $\mathbf{f}^1$ is a morphism, there exists a subchart $(U'^1, V'^1, \phi^1 \colon U'^1 \xrightarrow{\tilde{=}} V'^1) \in \mathbf{A}^1$ at $u^1$, a model neighborhood $U'^2 \subseteq U^2$ of $u^2$ and a chart

$(U'^2, \hat{V}'^2, \phi'^2 \colon U'^2 \xrightarrow{\tilde{=}} \hat{V}'^2) \in \mathbf{A}^2$ at $u^2$ such that $f_0^1[U'^1] \subseteq U'^2$ and diagram (11.8) in definition 11.4 (M-atlas morphisms for model spaces) on page 54 is commutative, as shown in fig. 6 (Completed m-atlas morphism) on page 54.

Let $I^{2,3} \stackrel{\text{def}}{=} U'^2 \cap f_0^{2-1}[U^3]$. Since $\mathbf{f}^2$ is a near morphism, diagram (11.14) below is locally nearly commutative at $u^2$, and thus there exist model neighborhoods $U''^2 \subseteq I^{2,3}$, $U'^3 \subseteq U^3$, $V''^2 \subseteq V'^3$, $V'^3 \subseteq V^3$, $\hat{V}'^3 \subseteq \mathbf{C}^3$ and an isomorphism $\hat{f} \colon V'^3 \xrightarrow{\tilde{=}} \hat{V}'^3$ such that diagram (11.14) below is commutative, i.e., eqs. (11.15) to (11.20) below hold, as shown in figs. 4 to 5.

$$D \stackrel{\text{def}}{=} \left\{ \begin{array}{l} i \colon \{u^2\} \rightarrowtail I^{2,3}, \phi'^2 \colon I^{2,3} \rightarrowtail V'^2, f_1^2 \colon V'^2 \longrightarrow \mathbf{C}^3, \\ f_0^2 \colon I^{2,3} \longrightarrow U'^3, \phi^3 \colon U'^3 \xrightarrow{\tilde{=}} V'^3 \end{array} \right\} \quad (11.14)$$

$$u^2 \in U''^2 \quad (11.15)$$
$$f_0[U''^2] \subseteq U'^3 \quad (11.16)$$
$$\phi^1[U''^2] = V''^2 \quad (11.17)$$
$$f_1[V''^2] \subseteq \hat{V}'^3 \quad (11.18)$$
$$\phi^3[U'^3] = V'^3 \quad (11.19)$$
$$\hat{f} \circ \phi^3 \circ f_0^2 = f_1^2 \circ \phi'^2 \quad (11.20)$$

Define $U'''^1 \stackrel{\text{def}}{=} f_0^1[U''^2]$ and $V'''^1 \stackrel{\text{def}}{=} f_1^1[V''^2]$. Diagram (11.21) below is commutative, i.e., eqs. (11.22) to (11.27) below hold, as shown in figs. 10 and 11.

$$D \stackrel{\text{def}}{=} \{ \\ \quad i\colon \{u^1\} \rightarrowtail U'''^1, \phi^1\colon U'''^1 \rightarrowtail V'''^1, f_1^2 \circ f_1^1\colon V'''^1 \longrightarrow V'^3, \\ \quad f_0^2 \circ f_0^1\colon U'''^1 \longrightarrow U'^3 \phi^3\colon U'^3 \rightarrowtail\!\!\!\xrightarrow{\cong} V'^3 \\ \}$$

(11.21)

$$u^1 \in U'''^1 \subseteq U'^1 \tag{11.22}$$
$$f_0^2 \circ f_0^1[U'''^1] \subseteq U'^3 \tag{11.23}$$
$$\phi^1[U'^3] = V'^3 \tag{11.24}$$
$$f_1^2 \circ f_1^1[V'''^1] \subseteq \hat{V}'^3 \tag{11.25}$$
$$\phi^3[U'^3] = V'^3 \tag{11.26}$$
$$\hat{f} \circ \phi^3 \circ f_0^2 \circ f_0^1 = f_1^2 \circ \phi'^2 \circ f_0^1 \\ = f_1^2 \circ f_1^1 \circ \phi^1 \tag{11.27}$$

$u^1$

$U'^1 \xrightarrow{\phi^1}_{\cong} V'^1$

$f_0^1 \qquad f_1^1$

$U'^2 \xleftarrow{i} I^{2,3} = U'^2 \cap f_0^{-1}[U^3] \xrightarrow{\phi'^2} \hat{V}'^2$

$f_0^2 \qquad f_0^2 \qquad f_1^2$

$E^3 \xleftarrow{i} U^3 \xrightarrow{\phi^3}_{\cong} V^3 \qquad C^3$

Figure 10: Partially completed composition of M-atlas near morphisms

□

2. If each $\boldsymbol{f}^i$ is an $\boldsymbol{E}^1$-$\boldsymbol{E}^{i+1}$ m-atlas morphism of $\boldsymbol{A}^i$ to $\boldsymbol{A}^{i+1}$ in the coordinate spaces $\boldsymbol{C}^i$, $\boldsymbol{C}^{i+1}$ then $\boldsymbol{f}^2 \stackrel{()}{\circ} \boldsymbol{f}^1$ is an $\boldsymbol{E}^1$-$\boldsymbol{E}^3$ m-atlas morphism of $\boldsymbol{A}^1$ to $\boldsymbol{A}^3$ in the coordinate spaces $\boldsymbol{C}^1$, $\boldsymbol{C}^3$.

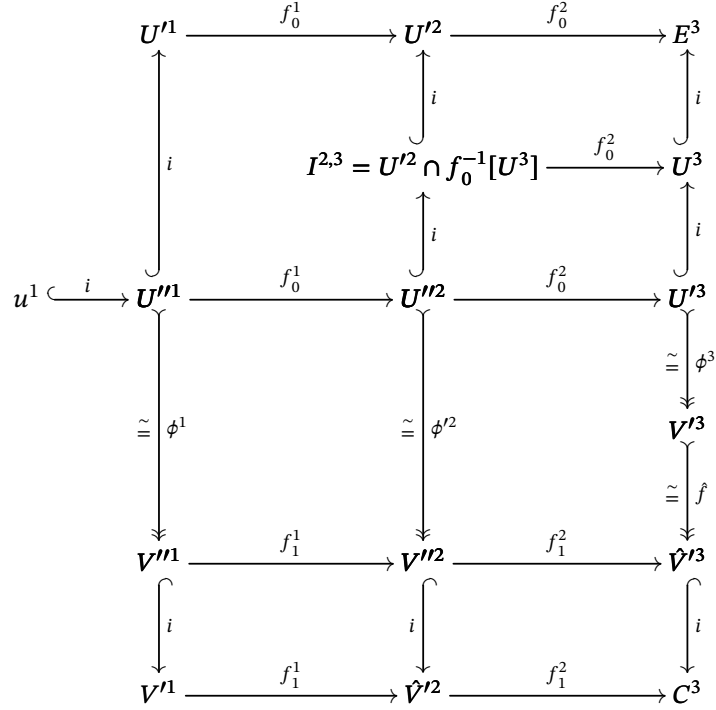

Figure 11: Completed composition of M-atlas near morphisms

*Proof.* $f_0^2 \circ f_0^1$ and $f_1^2 \circ f_1^1$ are model functions by lemma 7.5 (Composition of model functions) on page 30.

Let $(U^i, V^i, \phi^i) \in A^i$, $i = 1, 3$, be a chart, $I \stackrel{\text{def}}{=} U^1 \cap (f_0^2 \circ f_0^1)^{-1}[U^3] \neq \emptyset$.

Let $u^1 \in I$, $u^2 \stackrel{\text{def}}{=} f_0^1(u^1)$, $u^3 \stackrel{\text{def}}{=} f_0^2(u^2)$ and $(U^2, V^2, \phi^2) \in A^2$ be a chart at $u^2$. Since $\boldsymbol{f}^1$ is a morphism by hypothesis, there exists a subchart

$(U'^1, V'^1, \phi^1\colon U'^1 \stackrel{\tilde{=}}{\rightarrowtail\!\!\!\twoheadrightarrow} V'^1) \in A^1$ at $u^1$, a model neighborhood $U'^2 \subseteq U^2$ of $u^2$ and a chart $(U'^2, \hat{V}'^2, \phi'^2\colon U'^2 \stackrel{\tilde{=}}{\rightarrowtail\!\!\!\twoheadrightarrow} \hat{V}'^2) \in A^2$ at $u^2$ such that $f_0^1[U'^1] \subseteq U'^2$ and diagram (11.8) in definition 11.4 (M-atlas morphisms for model spaces) on page 54 is commutative, as shown in fig. 6 (Completed m-atlas morphism) on page 54.

Since $\boldsymbol{f}^2$ is a morphism by hypothesis, there exists a subchart

$(U''^2, V''^2, \phi^2\colon U''^2 \stackrel{\tilde{=}}{\rightarrowtail\!\!\!\twoheadrightarrow} V''^2) \in A^2$ at $u^2$, a model neighborhood $U'^3 \subseteq U^3$ of $u^3$ and a chart $(U'^3, \hat{V}'^3, \phi'^3\colon U'^3 \stackrel{\tilde{=}}{\rightarrowtail\!\!\!\twoheadrightarrow} \hat{V}'^3) \in A^3$ at $u^2$ such that $f_0^2[U''^2] \subseteq U'^3$ and diagram (11.28) below is commutative, as shown in figs. 12 and 13.

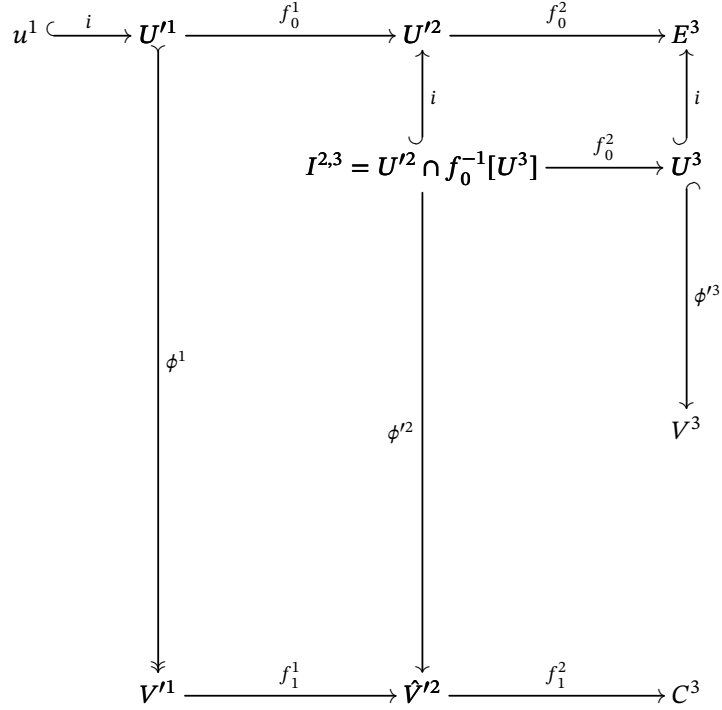

Figure 12: Uncompleted composition of m-atlas morphisms

$$D \stackrel{\text{def}}{=} \{$$
$$i\colon \{u^2\} \rightarrowtail U''^2, \phi^2\colon U''^2 \stackrel{\sim}{\underset{=}{\twoheadrightarrow}} V'^2, f_1\colon V''^2 \longrightarrow \hat{V}'^3, \quad (11.28)$$
$$f_0\colon U'^2 \longrightarrow U'^3, \phi'^3\colon U'^3 \stackrel{\sim}{\underset{=}{\twoheadrightarrow}} \hat{V}'^3$$
$$\}$$

Define $U'''^1 \stackrel{\text{def}}{=} f_0^{1-1}[U''^2]$, $V'''^1 \stackrel{\text{def}}{=} f_1^{1-1}[V''^2] = \phi^1[U'''^1]$, Figure 14 is commutative. □

3. *If $A^2$ is semi-maximal then $\boldsymbol{f}^2 \stackrel{()}{\circ} \boldsymbol{f}^1 \stackrel{\text{def}}{=} (f_0^2 \circ f_0^1, f_1^2 \circ f_1^1)$ is an $E^1$-$E^3$ m-atlas near morphism of $A^1$ to $A^3$ in the coordinate spaces $C^1$, $C^3$.*

*Proof.* Since $A^2$ is semi-maximal and $\boldsymbol{f}^1$ is a (near) morphism, $\boldsymbol{f}^1$ is a morphism either by hypothesis or by item 8 of lemma 11.12 (M-atlas (near) morphisms) on page 63.

The result follows from item 2 above, □

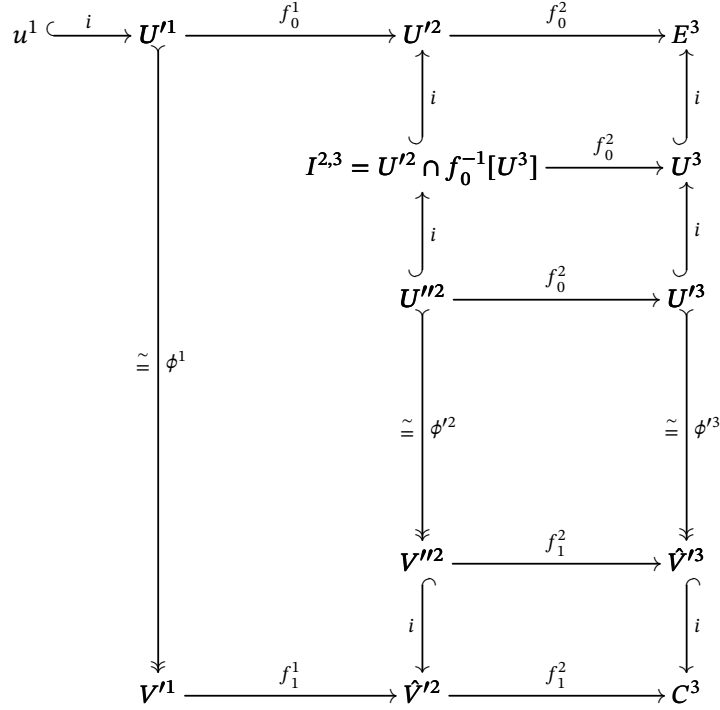

Figure 13: Partially completed composition of m-atlas morphisms

4. If $A^2$ and $A^3$ are semi-maximal then $f^2 \overset{()}{\circ} f^1 \overset{def}{=} (f_0^2 \circ f_0^1, f_1^2 \circ f_1^1)$ is an $E^1$-$E^3$ m-atlas morphism of $A^1$ to $A^3$ in the coordinate spaces $C^1, C^3$.

    *Proof.* Since $A^2$ is semi-maximal and $f^1$ is a (near) morphism, $f^1$ is a morphism either by hypothesis or by item 8 of lemma 11.12 (M-atlas (near) morphisms) on page 63. Similarly, $f^2$ is a morphism.
    
    The result follows from item 2 above, □

5. If each $f^i$, $i = 1, 2$, is a semi-strict ($\mathcal{E}^1$-$\mathcal{E}^2$-)$E^1$-$E^2$(-$\mathcal{C}^1$-$\mathcal{C}^2$) m-atlas (near) morphism then $f^2 \overset{()}{\circ} f^1$ is a semi-strict ($\mathcal{E}^1$-$\mathcal{E}^2$-)$E^1$-$E^2$(-$\mathcal{C}^1$-$\mathcal{C}^2$) m-atlas (near) morphism.

    *Proof.* There are four overlapping cases.

    **$E^i$ without category:** If each $E^i \overset{mod}{\subseteq} E^{i+1}$ and each $f_0^i$ is a local m-morphism of $E^i$ to $E^{i+1}$, then by lemma 7.19 (Composition of m-morphisms) on page 40 $f_0^2 \circ f_0^1$ is a local m-morphism of $E^1$ to $E^3$.

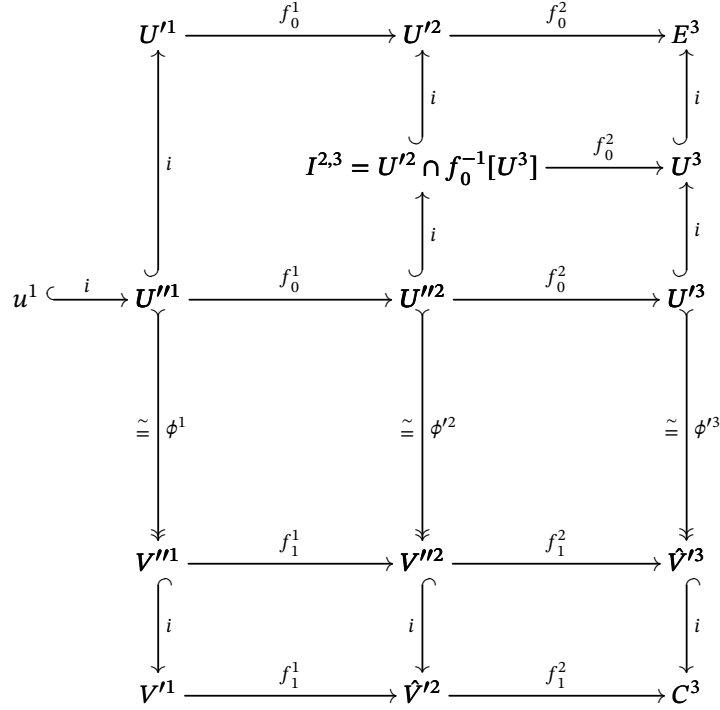

Figure 14: Completed composition of M-atlas morphisms

$C^i$ **without category:** If each $C^i \stackrel{\text{mod}}{\subseteq} C^{i+1}$ and each $f_1^i$ is a local m-morphism of $C^i$ to $C^{i+1}$, then by lemma 7.19 $f_1^2 \circ f_1^1$ is a local morphism of $C^1$ to $C^3$. Let $(U^i, V^i, \phi^i \colon U^i \stackrel{\tilde{=}}{\rightarrowtail\!\!\!\twoheadrightarrow} V^i) \in A^i$, $i = 1, 3$, with $I \stackrel{\text{def}}{=} U^1 \cap (f_0^2 \circ f_0^1)^{-1}[U^3] \neq \emptyset$. For every $u^1 \in I$ let $u^2 \stackrel{\text{def}}{=} f_0^1(u^1)$, $u^3 \stackrel{\text{def}}{=} f_0^2(u^2)$ and $(U^2, V^2, \phi^2 \colon U^2 \stackrel{\tilde{=}}{\rightarrowtail\!\!\!\twoheadrightarrow} V^2) \in A^2$ be a chart at $u^2$, $I^1 \stackrel{\text{def}}{=} U^1 \cap f_0^{1-1}[U^2] \neq \emptyset$ and $I^2 \stackrel{\text{def}}{=} U^2 \cap f_0^{2-1}[U^3] \neq \emptyset$. If $\phi^2 \circ f_0^1 \circ \phi^{1-1} \colon \phi^1[I^1] \longrightarrow V^2$ is a local m-morphism of $\text{Mod}(\phi^1[I^1], C^1)$ to $\text{Mod}(V^2, C^2)$ and $\phi^3 \circ f_0^2 \circ \phi^{2-1} \colon \phi^2[I^2] \longrightarrow V^3$ is a local m-morphism of $\text{Mod}(\phi^1[I^2], C^1)$ to $\text{Mod}(V^2, C^2)$ then by definition 7.10 there exist model neighborhoods $V'^i \subseteq V^i$, $i \in [1, 3]$, such that $f^i[V'^i] \subseteq V'^{i+1}$ and $\phi^{i+1} \circ f_0^i \circ \phi^{i-1} \colon V'^i \longrightarrow V'^{i+1}$ is an m-morphism of $C^{i+1}$, $i = 1, 2$. Then

$$\phi^3 \circ f_0^2 \circ f_0^1 \circ \phi^{1-1} \colon V'^1 \longrightarrow V'^3 \;=\; \phi^3 \circ f_0^2 \circ \phi^{2-1} \colon V'^2 \longrightarrow V'^3 \circ \phi^2 \circ f_0^1 \circ \phi^{1-1} \colon V'^1 \longrightarrow V'^2$$

is an m-morphism of $C^3$. Since $u^1 \in I$ is arbitray, $\phi^3 \circ f_0^2 \circ f_0^1 \circ \phi^{1-1} \colon \phi^1[I] \longrightarrow V^3$ is a local m-morphism of $\text{Mod}(\phi^1[I], \mathbf{C}^3)$ to $\text{Mod}(V^3, \mathbf{C}^3)$.

$E^i$ **with category:** If each $\mathscr{E}^i \overset{\text{full–cat}}{\subseteq} \mathscr{E}^{i+1}$ and each $f_0^i$ is a local $\mathscr{E}^i$-$\mathscr{E}^{i+1}$ m-morphism of $E^i$ to $E^{i+1}$, then by lemma 7.19 (Composition of m-morphisms) on page 40 $f_0^2 \circ f_0^1$ is a local $\mathscr{E}^1$-$\mathscr{E}^3$ m-morphism of $E^1$ to $E^3$.

$C^i$ **with category:** If each $\mathscr{E}^i \overset{\text{full–cat}}{\subseteq} \mathscr{E}^{i+1}$ and each $f_1^i$ is a local $\mathscr{E}^i$-$\mathscr{E}^{i+1}$ m-morphism of $C^i$ to $C^{i+1}$, then by lemma 7.19 (Composition of m-morphisms) on page 40 $f_1^2 \circ f_1^1$ is a local $\mathscr{E}^1$-$\mathscr{E}^3$ m-morphism of $C^1$ to $C^3$.

$\square$

6. *If each $f^i$, $i = 1, 2$, is strict then $f^2 \overset{()}{\circ} f^1$ is strict.*

*Proof.* There are four overlapping cases.

$E^i$ **without category:** If each $E^i \overset{\text{mod}}{\subseteq} E^{i+1}$ and each $f_0^i$ is an m-morphism of $E^i$ to $E^{i+1}$, then by lemma 7.19 (Composition of m-morphisms) on page 40 $f_0^2 \circ f_0^1$ is an m-morphism of $E^1$ to $E^3$.

$C^i$ **without category:** If each $C^i \overset{\text{mod}}{\subseteq} C^{i+1}$ and each $f_1^i$ is an m-morphism of $C^i$ to $C^{i+1}$, then by lemma 7.19 (Composition of m-morphisms) on page 40 $f_1^2 \circ f_1^1$ is an m-morphism of $C^1$ to $C^3$.

$E^i$ **with category:** If each $\mathscr{E}^i \overset{\text{full–cat}}{\subseteq} \mathscr{E}^{i+1}$ and each $f_0^i$ is a morphism of $\mathscr{E}^{i+1}$, then $f_0^2 \circ f_0^1$ is a morphism of $\mathscr{E}^3$.

$C^i$ **with category:** If each $\mathscr{E}^i \overset{\text{full–cat}}{\subseteq} \mathscr{E}^{i+1}$ and each $f_1^i$ is a morphism of $\mathscr{E}^i$, then $f_1^2 \circ f_1^1$ is a morphism of $\mathscr{E}^3$.

$\square$

7. *If each $f^i$ is constrained then $f^2 \overset{()}{\circ} f^1$ is constrained.*

*Proof.* This is a corollary of lemma 7.5 (Composition of model functions) on page 30. $\square$

8. $f^2 \overset{()}{\circ} f^1$ *is equivalent to* $g^2 \overset{()}{\circ} g^1$.

*Proof.* $f_0^1 = g_0^1$ and $f_0^2 = g_0^2$, hence $f_0^2 \circ f_0^1 = g_0^2 \circ g_0^1$. $\square$

**Corollary 11.17** (Composition of m-atlas (near) morphisms). *Let $\mathscr{E}^1$ be a model category, $E^1 \stackrel{Ob}{\in} \mathscr{E}^1$, $\mathscr{C}^i$, $i = 1, 2, 3$, be a model category, $C^i \stackrel{Ob}{\in} \mathscr{C}^i$, $A^i$ be an m-atlas of $E^1$ in the coordinate space $C^i$ and $f^i \stackrel{def}{=} (f_0^i \colon E^1 \longrightarrow E^1, f_1^i \colon C^i \longrightarrow C^{i+1})$, $i = 1, 2$, be a (semi-strict, strict) $(\mathscr{E}^1\text{-})E^1(\text{-}\mathscr{C}^i\text{-}\mathscr{C}^{i+1})$ m-atlas (near) morphism of $A^i$ to $A^{i+1}$ in the coordinate spaces $C^i$, $C^{i+1}$.*

*Then*

1. *If $A^2$ is semi-maximal then $f^2 \stackrel{()}{\circ} f^1 = (f_0^2 \circ f_0^1, f_1^2 \circ f_1^1)$ is a (semi-strict, strict) $(\mathscr{E}^1\text{-})E^1(\text{-}\mathscr{C}^1\text{-}\mathscr{C}^3)$ m-atlas near morphism of $A^1$ to $A^3$ in the coordinate spaces $C^1$, $C^3$.*

   *Proof.* Since each $f^i$ is a (semi-strict, strict) $E^1\text{-}E^1$ m-atlas (near) morphism of $A^i$ to $A^{i+1}$ in the coordinate spaces $C^i$, $C^{i+1}$ and $A^2$ is semi-maximal, $f^2 \stackrel{()}{\circ} f^1$ is a (semi-strict, strict) $E^1\text{-}E^1$ m-atlas near morphism of $A^1$ to $A^3$ in the coordinate spaces $C^1$, $C^3$. $\square$

2. *If each $f^i$ is a (semi-strict, strict) $(\mathscr{E}^1\text{-})E^1(\text{-}\mathscr{C}^i\text{-}\mathscr{C}^{i+1})$ m-atlas morphism of $A^i$ to $A^{i+1}$ in the coordinate spaces $C^i$, $C^{i+1}$ then $f^2 \stackrel{()}{\circ} f^1$ is a (semi-strict, strict) $(\mathscr{E}^1\text{-})E^1(\text{-}\mathscr{C}^1\text{-}\mathscr{C}^3)$ m-atlas morphism of $A^1$ to $A^3$ in the coordinate spaces $C^1$, $C^3$.*

   *Proof.* Since each $f^i$ is a (semi-strict, strict) $(\mathscr{E}^1\text{-}\mathscr{E}^1\text{-})E^1\text{-}E^1(\text{-}\mathscr{C}^i\text{-}\mathscr{C}^{i+1})$ m-atlas morphism of $A^i$ to $A^{i+1}$ in the coordinate spaces $C^i$, $C^{i+1}$, then $f^2 \stackrel{()}{\circ} f^1$ is a (semi-strict, strict) $(\mathscr{E}^1\text{-}\mathscr{E}^1\text{-})E^1\text{-}E^1(\text{-}\mathscr{C}^1\text{-}\mathscr{C}^3)$ m-atlas morphism of $A^1$ to $A^3$ in the coordinate spaces $C^1$, $C^3$. $\square$

*Let $\mathscr{E}^i$ $i = 1, 2, 3$, be a model category, $E^i \stackrel{Ob}{\in} \mathscr{E}^i$, $\mathscr{C}^1$ be a model category, $C^1 \stackrel{Ob}{\in} \mathscr{C}^1$, $A^i$ be an m-atlas of $E^i$ in the coordinate space $C^1$ and $f^i \stackrel{def}{=} (f_0^i \colon E^i \longrightarrow E^{i+1}, f_1^i \colon C^1 \longrightarrow C^1)$, $i = 1, 2$, be a (semi-strict, strict) $(\mathscr{E}^i\text{-}\mathscr{E}^{i+1}\text{-})E^i\text{-}E^{i+1}(\text{-}\mathscr{C}^1)$ m-atlas (near) morphism of $A^i$ to $A^{i+1}$ in the coordinate space $C^1$.*

3. *If $A^2$ is semi-maximal then $f^2 \stackrel{()}{\circ} f^1 = (f_0^2 \circ f_0^1, f_1^2 \circ f_1^1)$ is a (semi-strict, strict) $(\mathscr{E}^1\text{-}\mathscr{E}^3\text{-})E^1\text{-}E^3(\text{-}\mathscr{C}^1)$ m-atlas near morphism of $A^1$ to $A^3$ in the coordinate space $C^1$.*

   *Proof.* Since each $f^i$ is a (semi-strict, strict) $(\mathscr{E}^i\text{-}\mathscr{E}^{i+1}\text{-})E^i\text{-}E^{i+1}(\text{-}\mathscr{C}^1\text{-}\mathscr{C}^1)$ m-atlas (near) morphism of $A^i$ to $A^{i+1}$ in the coordinate spaces $C^1$, $C^1$ and $A^2$ is semi-maximal, $f^2 \stackrel{()}{\circ} f^1$ is a (semi-strict, strict) $(\mathscr{E}^1\text{-}\mathscr{E}^3\text{-})E^1\text{-}E^3(\text{-}\mathscr{C}^1\text{-}\mathscr{C}^1)$ m-atlas near morphism of $A^1$ to $A^3$ in the coordinate space $C^1$. $\square$

4. If each $\boldsymbol{f}^i$ is a (semi-strict, strict) $(\mathscr{E}^i\text{-}\mathscr{E}^{i+1}\text{-})\boldsymbol{E}^i\text{-}\boldsymbol{E}^{i+1}(\text{-}\mathscr{C}^1)$ m-atlas morphism of $\boldsymbol{A}^i$ to $\boldsymbol{A}^{i+1}$ in the coordinate spaces $\boldsymbol{C}^i$, $\boldsymbol{C}^{i+1}$ then $\boldsymbol{f}^2 \overset{()}{\circ} \boldsymbol{f}^1$ is a (semi-strict, strict) $(\mathscr{E}^1\text{-}\mathscr{E}^3\text{-})\boldsymbol{E}^1\text{-}\boldsymbol{E}^3(\text{-}\mathscr{C}^1\text{-}\mathscr{C}^1)$ m-atlas morphism of $\boldsymbol{A}^1$ to $\boldsymbol{A}^3$ in the coordinate spaces $\boldsymbol{C}^1$, $\boldsymbol{C}^1$.

*Proof.* Since each $\boldsymbol{f}^i$ is a (semi-strict, strict) $(\mathscr{E}^i\text{-}\mathscr{E}^{i+1}\text{-})\boldsymbol{E}^i\text{-}\boldsymbol{E}^{i+1}(\text{-}\mathscr{C}^1\text{-}\mathscr{C}^1)$ m-atlas morphism of $\boldsymbol{A}^i$ to $\boldsymbol{A}^{i+1}$ in the coordinate spaces $\boldsymbol{C}^1$, $\boldsymbol{C}^1$, then $\boldsymbol{f}^2 \overset{()}{\circ} \boldsymbol{f}^1$ is a (semi-strict, strict) $(\mathscr{E}^1\text{-}\mathscr{E}^3\text{-})\boldsymbol{E}^1\text{-}\boldsymbol{E}^3(\text{-}\mathscr{C}^1\text{-}\mathscr{C}^1)$ m-atlas morphism of $\boldsymbol{A}^1$ to $\boldsymbol{A}^3$ in the coordinate spaces $\boldsymbol{C}^1$, $\boldsymbol{C}^1$. $\square$

Let $\mathscr{E}^1$, $\mathscr{C}^1$ be model categories, $\boldsymbol{E}^1 \overset{\mathrm{Ob}}{\in} \mathscr{E}^1$, $\boldsymbol{C}^1 \overset{\mathrm{Ob}}{\in} \mathscr{C}^1$, $\boldsymbol{A}^i$, $i = 1, 2, 3$, be an m-atlas of $\boldsymbol{E}^1$ in the coordinate space $\boldsymbol{C}^1$ and $\boldsymbol{f}^i \overset{\mathrm{def}}{=} (f_0^i\colon \boldsymbol{E}^1 \longrightarrow \boldsymbol{E}^1, f_1^i\colon \boldsymbol{C}^1 \longrightarrow \boldsymbol{C}^1)$, $i = 1, 2$, be a (semi-strict, strict) $\boldsymbol{E}^1$ m-atlas (near) morphism of $\boldsymbol{A}^i$ to $\boldsymbol{A}^{i+1}$ in the coordinate space $\boldsymbol{C}^1$.

5. If $\boldsymbol{A}^2$ is semi-maximal then $\boldsymbol{f}^2 \overset{()}{\circ} \boldsymbol{f}^1 = (f_0^2 \circ f_0^1, f_1^2 \circ f_1^1)$ is a (semi-strict, strict) $\boldsymbol{E}^1$ m-atlas m-atlas near morphism of $\boldsymbol{A}^1$ to $\boldsymbol{A}^3$ in the coordinate space $\boldsymbol{C}^1$.

*Proof.* Since each $\boldsymbol{f}^i$ is a (semi-strict, strict) $\boldsymbol{E}^1$ m-atlas (near) morphism of $\boldsymbol{A}^i$ to $\boldsymbol{A}^{i+1}$ in the coordinate spaces $\boldsymbol{C}^1$, $\boldsymbol{C}^1$ and $\boldsymbol{A}^2$ is semi-maximal, $\boldsymbol{f}^2 \overset{()}{\circ} \boldsymbol{f}^1$ is a (semi-strict, strict) $\boldsymbol{E}^1$ m-atlas m-atlas near morphism of $\boldsymbol{A}^1$ to $\boldsymbol{A}^3$ in the coordinate spaces $\boldsymbol{C}^1$, $\boldsymbol{C}^1$. $\square$

6. If each $\boldsymbol{f}^i$ is a (semi-strict, strict) $\boldsymbol{E}^1$ m-atlas morphism of $\boldsymbol{A}^i$ to $\boldsymbol{A}^{i+1}$ in the coordinate space $\boldsymbol{C}^1$ then $\boldsymbol{f}^2 \overset{()}{\circ} \boldsymbol{f}^1$ is a (semi-strict, strict) $\boldsymbol{E}^1$ m-atlas morphism of $\boldsymbol{A}^1$ to $\boldsymbol{A}^3$ in the coordinate space $\boldsymbol{C}^1$.

*Proof.* Since each $\boldsymbol{f}^i$ is a (semi-strict, strict) $\boldsymbol{E}^1$ m-atlas morphism of $\boldsymbol{A}^i$ to $\boldsymbol{A}^{i+1}$ in the coordinate space $\boldsymbol{C}^1$, then $\boldsymbol{f}^2 \overset{()}{\circ} \boldsymbol{f}^1$ is a (semi-strict, strict) $\boldsymbol{E}^1$ m-atlas morphism of $\boldsymbol{A}^1$ to $\boldsymbol{A}^3$ in the coordinate space $\boldsymbol{C}^1$.

$\square$

Let $\mathscr{C}^i$, $i = 1, 2, 3$, be a model category, $\boldsymbol{C}^i \overset{\mathrm{Ob}}{\in} \mathscr{C}^i$, $E^i$ be a topological space, $\boldsymbol{A}^i$ be an m-atlas of $E^i$ in the coordinate space $\boldsymbol{C}^i$, $\boldsymbol{A}^2$ be semi-maximal and $\boldsymbol{f}^i \overset{\mathrm{def}}{=} (f_0^i\colon E^i \longrightarrow E^{i+1}, f_1^i\colon \boldsymbol{C}^i \longrightarrow \boldsymbol{C}^{i+1})$, $i = 1, 2$, be a (semi-strict, strict) $E^i\text{-}E^{i+1}$ m-atlas near morphism of $\boldsymbol{A}^i$ to $\boldsymbol{A}^{i+1}$ in the coordinate spaces $\boldsymbol{C}^i$, $\boldsymbol{C}^{i+1}$. Then

7. $\boldsymbol{f}^2 \overset{()}{\circ} \boldsymbol{f}^1 = (f_0^2 \circ f_0^1, f_1^2 \circ f_1^1)$ is a (semi-strict, strict) $E^1\text{-}E^3$ m-atlas near morphisms of $\boldsymbol{A}^1$ to $\boldsymbol{A}^3$ in the coordinate spaces $\boldsymbol{C}^1$, $\boldsymbol{C}^3$.

*Proof.* Since $f^i$ is a (semi-strict, strict) $E^i_{\text{triv}}$-$E^{i+1}_{\text{triv}}$ m-atlas near morphism of $A^i$ to $A^{i+1}$ in the coordinate spaces $C^i$, $C^{i+1}$, then $f^2 \overset{()}{\circ} f^1$ is a (semi-strict, strict) $E^1_{\text{triv}}$-$E^3_{\text{triv}}$ m-atlas near morphism of $A^1$ to $A^3$ in the coordinate spaces $C^1$, $C^3$. □

Let $\mathscr{E}^i, \mathscr{C}^i, i \in [1,3]$, be model categories, $E^i \overset{\text{Ob}}{\in} \mathscr{E}^i$, $C^i \overset{\text{Ob}}{\in} \mathscr{C}^i$ $A^i$ an m-atlas of $E^i$ in the coordinate space $C^i$, $A^2$ maximal and $f^i \overset{\text{def}}{=}$ $(f^i_0\colon E^i \longrightarrow E^{i+1}, f^i_1\colon C^i \longrightarrow C^{i+1})$ a (semi-strict, strict) $E^i$-$E^{i+1}$ m-atlas near morphism of $A^i$ to $A^{i+1}$ in the coordinate spaces $C^i$, $C^{i+1}$. Then

8. $f^2 \overset{()}{\circ} f^1 = (f^2_0 \circ f^1_0, f^2_1 \circ f^1_1)$ is a (semi-strict, strict) $E^1$-$E^3$ m-atlas near morphism of $A^1$ to $A^3$ in the coordinate spaces $C^1$, $C^3$.

   *Proof.* Since $A^2$ is maximal it is semi-maximal. □

Let $\mathscr{C}^i, i = 1, 2, 3$, be a model category, $C^i \overset{\text{Ob}}{\in} \mathscr{C}^i$, $E^i$ be a topological space, $A^i$ be an m-atlase of $E^i$ in the coordinate space $C^i$, $A^2$ maximal and $f^i \overset{\text{def}}{=}$ $(f^i_0\colon E^i \longrightarrow E^{i+1}, f^i_1\colon C^i \longrightarrow C^{i+1})$ be a (semi-strict, strict) $E^i$-$E^{i+1}$ m-atlas near morphism of $A^i$ to $A^{i+1}$ in the coordinate spaces $C^i$, $C^{i+1}$. Then

9. $f^2 \overset{()}{\circ} f^1 = (f^2_0 \circ f^1_0, f^2_1 \circ f^1_1)$ is a (semi-strict, strict) $E^1$-$E^3$ m-atlas near morphism of $A^1$ to $A^3$ in the coordinate spaces $C^1$, $C^3$.

   *Proof.* Since $A^2$ is maximal it is semi-maximal. □

Let $\mathscr{C}^i, i = 1, 2, 3$, be a model category, $C^i \overset{\text{Ob}}{\in} \mathscr{C}^i$, $E^i$ be a topological space, $A^i$ be an m-atlas of $E^i$ in the coordinate space $C^i$ and $f^i \overset{\text{def}}{=} (f^i_0\colon E^i \longrightarrow E^{i+1}, f^i_1\colon C^i \longrightarrow C^{i+1})$, $i = 1, 2$, be a (semi-strict, strict) $E^i$-$E^{i+1}$ m-atlas morphism of $A^i$ to $A^{i+1}$ in the coordinate spaces $C^i$, $C^{i+1}$. Then

10. $f^2 \overset{()}{\circ} f^1 = (f^2_0 \circ f^1_0, f^2_1 \circ f^1_1)$ is a (semi-strict, strict) $E^1$-$E^3$ m-atlas morphism of $A^1$ to $A^3$ in the coordinate spaces $C^1$, $C^3$.

    *Proof.* Since $(f^i_0\colon E^i \longrightarrow E^{i+1}, f^i_1\colon C^i \longrightarrow C^{i+1})$ is a (semi-strict, strict) $E^i_{\text{triv}}$-$E^{i+1}_{\text{triv}}$ m-atlas morphism of $A^i$ to $A^{i+1}$ in the coordinate spaces $C^i$, $C^{i+1}$, then $(f^2_0 \circ f^1_0, f^2_1 \circ f^1_1)$ is a (semi-strict, strict) $E^1_{\text{triv}}$-$E^3_{\text{triv}}$ m-atlas morphism of $A^1$ to $A^3$ in the coordinate spaces $C^1$, $C^3$. □

**Lemma 11.18** (M-atlas identity). *Let*

1. $\mathscr{E}^i, \mathscr{C}^i, i = 1, 2,$ *be model categories*

2. $\mathscr{E}^1 \overset{\text{full-cat}}{\subseteq} \mathscr{E}^2$

3. $\mathscr{C}^1 \overset{\text{full-cat}}{\subseteq} \mathscr{C}^2$

4. $\boldsymbol{E}^i \overset{\text{Ob}}{\in} \mathscr{E}^i$

5. $\boldsymbol{E}^1 \overset{\text{mod}}{\subseteq} \boldsymbol{E}^2$

6. $\boldsymbol{C}^i \overset{\text{Ob}}{\in} \mathscr{C}^i$

7. $\boldsymbol{C}^1 \overset{\text{mod}}{\subseteq} \boldsymbol{C}^2$

8. $\boldsymbol{A}^i$ *be an m-atlas of* $\boldsymbol{E}^i$ *in the coordinate space* $\boldsymbol{C}^i$

9. $\boldsymbol{A}^1 \subseteq \boldsymbol{A}^2$

*The identity morphism* $\mathrm{Id}_{(\boldsymbol{A}^1, \boldsymbol{E}^1, \boldsymbol{C}^1), (\boldsymbol{A}^2, \boldsymbol{E}^2, \boldsymbol{C}^2)}$ *is a semi-strict* $\boldsymbol{E}^1\text{-}\boldsymbol{E}^2\text{-}\mathscr{C}^1\text{-}\mathscr{C}^2$ *m-atlas near morphism of* $\boldsymbol{A}^1$ *to* $\boldsymbol{A}^2$ *in the coordinate spaces* $\boldsymbol{C}^1, \boldsymbol{C}^2$ *and a semi-strict* $\mathscr{E}^1\text{-}\mathscr{E}^2\text{-}\boldsymbol{E}^1\text{-}\boldsymbol{E}^2\text{-}\mathscr{C}^1\text{-}\mathscr{C}^2$ *m-atlas near morphism of* $\boldsymbol{A}^1$ *to* $\boldsymbol{A}^2$ *in the coordinate spaces* $\boldsymbol{C}^1, \boldsymbol{C}^2$. *It is a morphism if* $\boldsymbol{A}^1$ *and* $\boldsymbol{A}^2$ *are semi-maximal.*

*Proof.* $\mathrm{Id}_{(\boldsymbol{A}^1, \boldsymbol{E}^1, \boldsymbol{C}^1), (\boldsymbol{A}^2, \boldsymbol{E}^2, \boldsymbol{C}^2)}$ is an m-atlas near morphism:

Let $(U^i, V^i, \phi^i) \in \boldsymbol{A}^i, i = 1, 2,$ be charts such that $I \overset{\text{def}}{=} U^1 \cap U^2 \neq \emptyset$. To show that diagram $D$ below is M-locally nearly commutative in $\boldsymbol{C}^2$, let $u^1 \in I$. Since $\boldsymbol{E}^1 \overset{\text{mod}}{\subseteq} \boldsymbol{E}^2$ by hypothethis, $U^1$ is a model neighborhood of $\boldsymbol{E}^2$ and thus $I$ is a model neighborhood of $\boldsymbol{E}^2$ and the identity map $\mathrm{Id}_{I,U^2}$ is a morphism of $\boldsymbol{E}^2$. Then define $U'^1 \overset{\text{def}}{=} I$, $U'^2 \overset{\text{def}}{=} I$, $V'^1 \overset{\text{def}}{=} \phi^1[I]$, $V'^2 \overset{\text{def}}{=} \phi^2[I]$, $\hat{V}'^2 \overset{\text{def}}{=} V'^1$ and $\hat{f} \overset{\text{def}}{=} \phi^1 \circ \phi^{2-1} \colon V'^2 \overset{\tilde{=}}{\rightarrowtail\!\!\!\twoheadrightarrow} V'^1$. Equations (11.2) to (11.7) in definition 11.1 (M-atlas near morphisms for model spaces) hold.

$$D \overset{\text{def}}{=} \left\{ \begin{array}{l} \phi^1 \colon I \overset{\text{def}}{=} U^1 \cap U^2 \rightarrowtail V^1, \mathrm{Id}_{V^1, \boldsymbol{C}^2} \colon V^1 \longrightarrow \boldsymbol{C}^2, \\ \mathrm{Id}_{U^1, U^2} \colon U^1 \longrightarrow U^2, \phi^2 \colon U^2 \overset{\tilde{=}}{\rightarrowtail\!\!\!\twoheadrightarrow} V^2 \end{array} \right\}$$

$\mathrm{Id}_{(\boldsymbol{A}^1, \boldsymbol{E}^1, \boldsymbol{C}^1), (\boldsymbol{A}^2, \boldsymbol{E}^2, \boldsymbol{C}^2)}$ is semi-strict:

Let $u^1 \in \boldsymbol{A}$. Since the model neighborhoods of $\boldsymbol{E}^1$ form an open cover, there is a model neighborhood $U$ containing $u^1$. $\mathrm{Id}_{U, \boldsymbol{E}^1}$ is a morphism of $\boldsymbol{E}^1$ and of $\mathscr{E}^1$;

since $E^1 \overset{\text{mod}}{\subseteq} E^2$ and $\mathscr{C}^1 \overset{\text{full-cat}}{\subseteq} \mathscr{C}^2$, $\text{Id}_{U,E^1}$ is also a morphism of $E^2$ and of $\mathscr{C}^2$. Thus $\text{Id}_{(A^1,E^1,C^1),(A^2,E^2,C^2)}$ is a local morphism.

If $E^1 \overset{\text{Ob}}{\in} \mathscr{C}^2$ and $C^1 \overset{\text{Ob}}{\in} \mathscr{C}^2$ then $\text{Id}_{(A^1,E^1,C^1),(A^2,E^2,C^2)}$ is strict:
$\text{Id}_{E^1} \overset{\text{Ar}}{\in} \mathscr{C}^2 \overset{\text{full-cat}}{\subseteq} \mathscr{C}^2$ and $\text{Id}_{C^1} \overset{\text{Ar}}{\in} \mathscr{C}^1 \overset{\text{full-cat}}{\subseteq} \mathscr{C}^2$. $\square$

*If $A^i$ is semi-maximal, $i = 1, 2$, then $\boldsymbol{f} \overset{\text{def}}{=} \text{ID}_{(E^1,C^1),(E^2,C^2)}$ is an m-atlas morphism of $A^1$ to $A^2$ in the coordinate spaces $C^1$, $C^2$.*

*Proof.* Let $u^1 \in E^1$ and $(U^1, V^1, \phi^1\colon U^1 \rightarrowtail\!\!\!\!\!\!\overset{\sim}{\twoheadrightarrow} V^1) \in A^1 \subseteq A^2$ be a chart at $u^1$. Define $U'^1 \overset{\text{def}}{=} U^1$, $V'^1 \overset{\text{def}}{=} V^1$, $U'^2 \overset{\text{def}}{=} U^1$ and $V'^2 \overset{\text{def}}{=} V^1$. Diagram (11.8) in definition 11.4 is M-nearly commutative in $C^2$,
$\text{Id}_{E^i} \overset{\text{Ar}}{\in} \mathscr{C}^i \overset{\text{full-cat}}{\subseteq} \mathscr{C}^{i+1}$ and $\text{Id}_{C^i} \overset{\text{Ar}}{\in} \mathscr{C}^i \overset{\text{full-cat}}{\subseteq} \mathscr{C}^{i+1}$. $\square$

## 11.2 Categories of m-atlases and functors

This subsection defines categories of m-atlases and functors among them, and proves some basic results.

### 11.2.1 M-atlas categories

In the following, ... will refer to one of

$$A^1, E^1, C^1, A^2, E^2, C^2, f_0, f_1$$

$$A^1, E^1, C^1, A^2, E^2, C^2, f_0, f_1$$

depending on context.

**Definition 11.19** (Sets of m-atlas morphisms)**.** Let $E^i$ and $C^i$. $i = 1, 2$, be model spaces. Then the set of (constrained) (semi-strict, strict) (full, semi-maximal, maximal, full semi-maximal, full maximal) $E^1$-$E^2$ m-atlas (near) morphisms in the coordinate spaces $C^1$, $C^2$ is

$$\mathscr{A}t\ell_{\text{Ar}\ (\text{full,S-max,max,S-max-full,max-full})\ (\text{semi-strict,strict})}^{(\text{const,const-near,near})}(E^1, C^1, E^2, C^2) \overset{\text{def}}{=}$$

$$\left\{ ((f_0, f_1), (A^1, E^1, C^1), (A^2, E^2, C^2)) \,\middle|\, \text{isAtl}^{(\text{const,const-near,near})}_{\substack{(\text{full,S-max,max,S-max-full,max-full})\\(\text{semi-strict,strict})}}(A^1, E^1, C^1, A^2, E^2, C^2, f_0, f_1) \right\} \quad (11.29)$$

Let $E^i$, $i = 1, 2$, be a topological space and $\boldsymbol{C}^1 = (C^i, \mathscr{C}^i)$, $i = 1, 2$, be a model space. Then the set of (constrained) (semi-strict, strict) (full, semi-maximal, maximal, full semi-maximal, full maximal) $E^1$-$E^2$ m-atlas (near) morphisms in the coordinate spaces $\boldsymbol{C}^1$, $\boldsymbol{C}^2$ is

$$\mathscr{A}t\ell_{\mathrm{Ar}\ (\text{full,S-max,max,S-max-full,max-full})\ (\text{semi-strict,strict})}^{(\text{const,const–near,near})}(E^1, \boldsymbol{C}^1, E^2, \boldsymbol{C}^2) \stackrel{\text{def}}{=}$$

$$\left\{ ((f_0, f_1), (\boldsymbol{A}^1, E^1, \boldsymbol{C}^1), (\boldsymbol{A}^2, E^2, \boldsymbol{C}^2)) \;\middle|\; \mathrm{isAtl}_{(\text{full,S-max,max,S-max-full,max-full})\ (\text{semi-strict,strict})}^{(\text{const,const–near,near})}(\boldsymbol{A}^1, E^1, \boldsymbol{C}^1, \boldsymbol{A}^2, E^2, \boldsymbol{C}^2, f_0, f_1) \right\} \quad (11.30)$$

Let

1. $\boldsymbol{E}$ and $\boldsymbol{C}$ be sets of model spaces

2. $\boldsymbol{P} \stackrel{\text{def}}{=} \boldsymbol{E} \times \boldsymbol{C}$

$$\mathscr{A}t\ell_{\mathrm{Ar}\ (\text{full,S-max,max,S-max-full,max-full})\ (\text{semi-strict,strict})}^{(\text{const,const–near,near})}(\boldsymbol{E}, \boldsymbol{C}) \stackrel{\text{def}}{=}$$

$$\left\{ ((f_0, f_1), (\boldsymbol{A}^1, E^1, \boldsymbol{C}^1), (\boldsymbol{A}^2, E^2, \boldsymbol{C}^2)) \;\middle|\; (E^i, \boldsymbol{C}^i) \in \boldsymbol{P} \wedge \mathrm{isAtl}_{(\text{full,S-max,max,S-max-full,max-full})\ (\text{semi-strict,strict})}^{(\text{const,const–near,near})}(\ldots) \right\} \quad (11.31)$$

Let

1. $\mathscr{E}$ be a model category

2. $\boldsymbol{C}$ be a set of model spaces

3. $\boldsymbol{P} \stackrel{\text{def}}{=} \mathrm{Ob}(\mathscr{E}) \times \boldsymbol{C}$

$$\mathcal{A}t\ell_{\text{Ar}\,(\text{full},\text{S-max},\text{max},\text{S-max-full},\text{max-full})}^{(\text{const},\text{const}-\text{near},\text{near})}(\boldsymbol{E},\mathscr{C}) \stackrel{\text{def}}{=}$$
$$\left\{ ((f_0, f_1), (\boldsymbol{A}^1, \boldsymbol{E}^1, \boldsymbol{C}^1), (\boldsymbol{A}^2, \boldsymbol{E}^2, \boldsymbol{C}^2)) \,\middle|\, (\boldsymbol{E}^i, \boldsymbol{C}^i) \in \boldsymbol{P} \land \text{isAtl}_{\substack{(\text{full},\text{S-max},\text{max},\text{S-max-full},\text{max-full}) \\ (\text{semi-strict},\text{strict})}}^{(\text{const},\text{const}-\text{near},\text{near})}(..., \mathscr{C}, \mathscr{C}) \right\} \quad (11.32)$$

Let

1. $\boldsymbol{E}$ be a set of model spaces

2. $\mathscr{C}$ be a model category

3. $\boldsymbol{P} \stackrel{\text{def}}{=} \boldsymbol{E} \times Ob(\mathscr{C})$

$$\mathcal{A}t\ell_{\text{Ar}\,(\text{full},\text{S-max},\text{max},\text{S-max-full},\text{max-full})}^{(\text{const},\text{const}-\text{near},\text{near})}(\mathscr{E},\mathscr{C}) \stackrel{\text{def}}{=}$$
$$\left\{ ((f_0, f_1), (\boldsymbol{A}^1, \boldsymbol{E}^1, \boldsymbol{C}^1), (\boldsymbol{A}^2, \boldsymbol{E}^2, \boldsymbol{C}^2)) \,\middle|\, (\boldsymbol{E}^i, \boldsymbol{C}^i) \in \boldsymbol{P} \land \text{isAtl}_{\substack{(\text{full},\text{S-max},\text{max},\text{S-max-full},\text{max-full}) \\ (\text{semi-strict},\text{strict})}}^{(\text{const},\text{const}-\text{near},\text{near})}(..., \mathscr{C}, \mathscr{C}) \right\} \quad (11.33)$$

Let

1. $\mathscr{E}$ and $\mathscr{C}$ be model categories

2. $\boldsymbol{P} \stackrel{\text{def}}{=} Ob(\mathscr{E}) \times Ob(\mathscr{C})$

$$\mathcal{A}t\ell_{\text{Ar}\,(\text{full},\text{S-max},\text{max},\text{S-max-full},\text{max-full})}^{(\text{const},\text{const}-\text{near},\text{near})}(\mathscr{E},\mathscr{C}) \stackrel{\text{def}}{=}$$
$$\left\{ ((f_0, f_1), (\boldsymbol{A}^1, \boldsymbol{E}^1, \boldsymbol{C}^1), (\boldsymbol{A}^2, \boldsymbol{E}^2, \boldsymbol{C}^2)) \,\middle|\, (\boldsymbol{E}^i, \boldsymbol{C}^i) \in \boldsymbol{P} \land \text{isAtl}_{\substack{(\text{full},\text{S-max},\text{max},\text{S-max-full},\text{max-full}) \\ (\text{semi-strict},\text{strict})}}^{(\text{const},\text{const}-\text{near},\text{near})}(..., \mathscr{E}, \mathscr{E}, \mathscr{C}, \mathscr{C}) \right\} \quad (11.34)$$

**Lemma 11.20** (Sets of m-atlas morphisms). *Let $E$ and $C$ be sets of model spaces.*

$$\mathscr{A}t\ell_{\text{Ar}\atop(\text{full,S-max,max,S-max-full,max-full})\atop(\text{semi-strict,strict})}^{(\text{const,const}-\text{near,near})}(E, C) \subseteq \mathscr{A}t\ell_{\text{Ar}}^{(\text{near})}(E, C).$$

$$\mathscr{A}t\ell_{\text{Ar}\atop\text{max-full}\atop(\text{semi-strict,strict})}^{(\text{const,const}-\text{near,near})}(E, C) \subseteq \mathscr{A}t\ell_{\text{Ar}\atop\text{S-max-full}\atop(\text{semi-strict,strict})}^{(\text{const,const}-\text{near,near})}(E, C)$$

$$\subseteq \mathscr{A}t\ell_{\text{Ar}\atop\text{full}\atop(\text{semi-strict,strict})}^{(\text{const,const}-\text{near,near})}(E, C)$$

*Proof.* The result follows from definition 11.10 (Abbreviated nomenclature for m-atlas (near) morphisms) on page 57 and definition 11.19 above. □

$$\mathscr{A}t\ell_{\text{Ar}\atop(\text{full,S-max,max,S-max-full,max-full})\atop(\text{semi-strict,strict})}(E, C) \subseteq \mathscr{A}t\ell_{\text{Ar}\atop(\text{full,S-max,max,S-max-full,max-full})\atop(\text{semi-strict,strict})}^{\text{near}}(E, C)$$

$$\mathscr{A}t\ell_{\text{Ar}\atop(\text{S-max,max})\atop(\text{semi-strict,strict})}(E, C) = \mathscr{A}t\ell_{\text{Ar}\atop(\text{S-max,max})\atop(\text{semi-strict,strict})}^{\text{near}}(E, C)$$

*Proof.* The result follows from item 1 of lemma 11.12 (M-atlas (near) morphisms) on page 60, definition 11.10 (Abbreviated nomenclature for m-atlas (near) morphisms) on page 57 and definition 11.19 above. □

*Proof.* The result follows from lemma 11.12 (M-atlas (near) morphisms) on page 60, definition 11.10 (Abbreviated nomenclature for m-atlas (near) morphisms) on page 57 and definition 11.19 above. □

**Corollary 11.21** (Sets of m-atlas morphisms). *Let $E$ and $C$ be sets of topological spaces.*

$$\mathscr{A}t\ell_{\text{Ar}\atop(\text{full,S-max,max,S-max-full,max-full})\atop(\text{semi-strict,strict})}^{(\text{const,const}-\text{near,near})}(E, C) \subseteq \mathscr{A}t\ell_{\text{Ar}}^{(\text{near})}(E, C)$$

$$\mathscr{A}t\ell_{\text{Ar}\atop\text{max-full}\atop(\text{semi-strict,strict})}^{(\text{const,const}-\text{near,near})}(E, C) \subseteq \mathscr{A}t\ell_{\text{Ar}\atop\text{S-max-full}\atop(\text{semi-strict,strict})}^{(\text{const,const}-\text{near,near})}(E, C)$$

$$\subseteq \mathscr{A}t\ell_{\text{Ar}\atop\text{full}\atop(\text{semi-strict,strict})}^{(\text{const,const}-\text{near,near})}(E, C)$$

$$\mathscr{A}t\ell_{\text{Ar}\atop(\text{full,S-max,max,S-max-full,max-full})\atop(\text{semi-strict,strict})}(E, C) \subseteq \mathscr{A}t\ell_{\text{Ar}\atop(\text{full,S-max,max,S-max-full,max-full})\atop(\text{semi-strict,strict})}^{\text{near}}(E, C)$$

*Proof.* The result follows from definition 11.10 (Abbreviated nomenclature for m-atlas (near) morphisms) on page 57 and definition 11.19 above. □

**Definition 11.22** (Categories $\mathcal{A}t\ell(E,C)$). Let $E$ and $C$ be sets of model spaces. Let $P \stackrel{\text{def}}{=} E \times C$.

$$\mathcal{A}t\ell^{(\text{const})}_{\substack{(\text{full,S-max,max,S-max-full,max-full}) \\ (\text{semi-strict,strict})}}(E,C) \stackrel{\text{def}}{=}$$

$$\left( \mathcal{A}t\ell_{\text{Ob} \substack{(\text{full,S-max,max,S-max-full,max-full})}}(E,C), \mathcal{A}t\ell^{(\text{const})}_{\text{Ar} \substack{(\text{full,S-max,max,S-max-full,max-full}) \\ (\text{semi-strict,strict})}}(E,C), \stackrel{A}{\circ} \right) \quad (11.35)$$

Let $E^i$ and $C^i$, $i = 1,2$, be model spaces, $A^i$ be an atlas of $E^i$ in the coordinate space $C^i$ and $f \stackrel{\text{def}}{=} (f_0 \colon E^1 \longrightarrow E^2, f_0 \colon C^1 \longrightarrow C^2)$ be a (strict, semi-strict) $E^1$-$E^2$ M-atlas (near) morphism of $A^1$ to $A^2$ in the coordinate spaces $C^1, C^2$.

The identity functor $\mathcal{F}_{M,\text{Id}}$ is

$$\mathcal{F}_{M,\text{Id}}(A^i, E^i, C^i) \stackrel{\text{def}}{=} (A^i, E^i, C^i) \quad (11.36)$$

This nomenclature will be justified below.

$$\mathcal{F}_{M,\text{Id}}\left((f_0 \colon E^1 \longrightarrow E^2, f_1 \colon C^1 \longrightarrow C^2,), (A^1, E^1, C^1), (A^2, E^2, C^2)\right) \stackrel{\text{def}}{=}$$

$$\left((f_0 \colon E^1 \longrightarrow E^2, f_1 \colon C^1 \longrightarrow C^2,), (A^1, E^1, C^1), (A^2, E^2, C^2)\right) \quad (11.37)$$

These are fine grained iff $C$ is fine grained.

**Lemma 11.23** ($\mathcal{A}t\ell(E,C)$ is a category). *Let $E$ and $C$ be sets of model spaces. Then*

1. $\mathcal{A}t\ell^{(\text{const-near,near})}_{\substack{(\text{S-max,max}) \\ (\text{semi-strict,strict})}}(E,C)$ *is a category and the identity morphism for an object* $(A^i, E^i, C^i)$ *of* $\mathcal{A}t\ell^{(\text{const-near,near})}_{\text{Ob} \substack{(\text{S-max,max}) \\ (\text{semi-strict,strict})}}(E,C)$ *is* $\text{Id}_{(A^i, E^i, C^i)}$.

2. $\mathcal{A}t\ell^{(\text{const})}_{\substack{(\text{full,S-max,max,S-max-full,max-full}) \\ (\text{semi-strict,strict})}}(E,C)$ *is a category and the identity morphism for an object* $(A^i, E^i, C^i)$ *of* $\mathcal{A}t\ell^{(\text{const})}_{\text{Ob} \substack{(\text{full,S-max,max,S-max-full,max-full}) \\ (\text{semi-strict,strict})}}(E,C)$ *is* $\text{Id}_{(A^i, E^i, C^i)}$.

*Proof.* Let $(A^i, E^i, C^i)$, $i \in [1,3]$, be an object of $\mathcal{A}t\ell(E,C)$ and let $m^i \stackrel{\text{def}}{=} \left((f_0^i, f_1^i), (A^i, E^i, C^i), (A^{i+1}, E^{i+1}, C^{i+1})\right)$ be a morphism of $\mathcal{A}t\ell(E,C)$. Then

**Composition:**
$m^2 \overset{A}{\circ} m^1 = \left((f_0^2 \circ f_0^1, f_1^2 \circ f_1^1), (A^1, E^1, C^1), (A^3, E^3, C^3)\right)$ is a morphism of $\mathcal{A}t\ell(E, C)$ by item 2 of lemma 11.16 (Composition of m-atlas (near) morphisms) on page 72.

**Associativity:**
Composition is associative by lemma 1.19 (Tuple composition for labeled morphisms) on page 12.

**Identity:**
$\mathrm{Id}_{(A^i, E^i, C^i)}$ is an identity morphism by lemma 1.19.

$\square$

**Corollary 11.24** (Subcategories of $\mathcal{A}t\ell(E, C)$)**.** *Let $E$ and $C$ be sets of model spaces. Then* $\mathcal{A}t\ell^{(\mathrm{const})}_{(\mathrm{full,S\text{-}max,max,S\text{-}max\text{-}full,max\text{-}full})\,(\mathrm{semi\text{-}strict,strict})}(E, C)$ *is a subcategory of* $\mathcal{A}t\ell^{(\mathrm{const})}(E, C)$.

*Remark* 11.25. They are not, in general, full subcategories.

*The identity functor $\mathcal{F}_{\mathrm{M,Id}}$ is a functor from each of the subcategories to each of the containing categories.*

**Definition 11.26** (Categories $\mathcal{A}t\ell(\mathscr{E}, \mathscr{C})$)**.** Let $\mathscr{E}$ and $\mathscr{C}$ be model categories. Let $P \overset{\mathrm{def}}{=} Ob(\mathscr{E}) \times Ob(\mathscr{C})$. Then

$$\mathcal{A}t\ell_{\mathrm{Ob}\,(\mathrm{full,S\text{-}max,max,S\text{-}max\text{-}full,max\text{-}full})}(\mathscr{E}, \mathscr{C}) \overset{\mathrm{def}}{=} \bigcup_{\substack{E \in \mathscr{E}^{\mathrm{Ob}} \\ C \in \mathscr{C}^{\mathrm{Ob}}}} \mathcal{A}t\ell_{\mathrm{Ob}\,(\mathrm{full,S\text{-}max,max,S\text{-}max\text{-}full,max\text{-}full})}(E, C) \quad (11.38)$$

$$\mathcal{A}t\ell^{(\mathrm{const,const-near,near})}_{\mathrm{Ar}\,(\mathrm{full,S\text{-}max,max,S\text{-}max\text{-}full,max\text{-}full})\,(\mathrm{semi\text{-}strict,strict})}(\mathscr{E}, \mathscr{C}) \overset{\mathrm{def}}{=}$$

$$\left\{ ((f_0, f_1), (A^1, E^1, C^1), (A^2, E^2, C^2)) \,\middle|\, (E^i, C^i) \in P \wedge \mathrm{isAtl}^{(\mathrm{const,const-near,near})}_{\mathrm{Ar}\,(\mathrm{full,S\text{-}max,max,S\text{-}max\text{-}full,max\text{-}full})\,(\mathrm{semi\text{-}strict,strict})}(A^1, E^1, C^1, A^2, E^2, C^2, f_0, f_1, \mathscr{E}, \mathscr{C}) \right\} \quad (11.39)$$

$$\mathcal{A}t\ell^{(\text{const,const-near,near})}_{\substack{(\text{full,S-max,max,S-max-full,max-full})\\(\text{semi-strict,strict})}}(\mathcal{E},\mathcal{C}) \stackrel{\text{def}}{=}$$

$$\left(\mathcal{A}t\ell_{\text{Ob}}(\mathcal{E},\mathcal{C}), \mathcal{A}t\ell^{(\text{const,const-near,near})}_{\text{Ar}}(\mathcal{E},\mathcal{C}), \stackrel{A}{\circ}\right) \quad (11.40)$$

**Lemma 11.27** ($\mathcal{A}t\ell(\mathcal{E},\mathcal{C})$ is a category). *Let $\mathcal{E}$ and $\mathcal{C}$ be model categories. Then*

1. $\mathcal{A}t\ell^{(\text{const})}_{\substack{(\text{full,S-max,max,S-max-full,max-full})\\(\text{semi-strict,strict})}}(\mathcal{E},\mathcal{C})$ *is a category and the identity morphism for an object $(A^i, E^i, C^i)$ of $\mathcal{A}t\ell^{(\text{const})}_{\substack{(\text{full,S-max,max,S-max-full,max-full})\\(\text{semi-strict,strict})}}(\mathcal{E},\mathcal{C})$ is $\text{Id}_{(A^i,E^i,C^i)}$.*

   *Proof.* Let $(A^i, E^i, C^i)$, $i \in [1,3]$, be an object of $\mathcal{A}t\ell(\mathcal{E},\mathcal{C})$ and let $m^i \stackrel{\text{def}}{=} \big((f^i_0, f^i_1), (A^i, E^i, C^i), (A^{i+1}, E^{i+1}, C^{i+1})\big)$ be a morphism of $\mathcal{A}t\ell(\mathcal{E},\mathcal{C})$. Then

   **Composition:**
   $f^2_0 \circ f^1_0$ is a morphism of $\mathcal{E}$, $f^2_1 \circ f^1_1$ is a morphism of $\mathcal{C}$ and $\boldsymbol{f}^2 \stackrel{()}{\circ} \boldsymbol{f}^1 = (f^2_0 \circ f^1_0, f^2_1 \circ f^1_1)$ is an $E^1$-$E^3$ m-atlas morphism of $A^1$ to $A^3$ in the coordinate spaces $C^1, C^2$ by item 2 of lemma 11.16 (Composition of m-atlas (near) morphisms) on page 72.

   **Associativity:**
   Composition is associative by lemma 1.19 (Tuple composition for labeled morphisms) on page 12.

   **Identity:** $\text{Id}_{(A^i,E^i,C^i)}$ is an identity morphism by lemma 1.19.

   □

2. $\mathcal{A}t\ell^{(\text{const-near,near})}_{\substack{(\text{S-max,max})\\(\text{semi-strict,strict})}}(\mathcal{E},\mathcal{C})$ *is a category and the identity morphism for an object $(A^i, E^i, C^i)$ of $\mathcal{A}t\ell^{(\text{const-near,near})}_{\substack{(\text{S-max,max})\\(\text{semi-strict,strict})}}(\mathcal{E},\mathcal{C})$ is $\text{Id}_{(A^i,E^i,C^i)}$.*

**Definition 11.28** (Categories $\mathcal{A}t\ell(E, C)$ of topological spaces). Let $E$ be a set of topological spaces and $C$ a set of model spaces. Let $P \stackrel{\text{def}}{=} E \times C$. Then

$$\mathcal{A}t\ell^{(\text{const})}_{\substack{(\text{full,S-max,max,S-max-full,max-full})\\(\text{semi-strict,strict})}}(E,C) \stackrel{\text{def}}{=}$$

$$\left(\mathcal{A}t\ell_{\text{Ob}\,(\text{full,S-max,max,S-max-full,max-full})}(E,C), \mathcal{A}t\ell^{(\text{const})}_{\substack{\text{Ar}\,(\text{full,S-max,max,S-max-full,max-full})\\(\text{semi-strict,strict})}}(E,C), \stackrel{A}{\circ}\right) \quad (11.41)$$

Let $(A^1, E^1, C^1) \in \mathcal{A}t\ell_{\text{Ob}}(E, C)$.

$$\text{Id}_{(A^1,E^1,C^1)} \stackrel{\text{def}}{=} \big((\text{Id}_{E^1}, \text{Id}_{C^1}), (A^1, E^1, C^1), (A^1, E^1, C^1)\big) \quad (11.42)$$

**Lemma 11.29** ($\mathcal{A}t\ell(E, C)$ of topological spaces is a category)**.** *Let $E$ be a set of topological spaces and $C$ a set of model spaces. Then $\mathcal{A}t\ell(E, C)$ is a category and the identity morphism for an object $(A^1, E^1, C^1)$ is $\text{Id}_{(A^1,E^1,C^1)}$.*

*Proof.* Let $(A^i, E^i, C^i)$, $i \in [1, 3]$, be an object of $\mathcal{A}t\ell(E, C)$. and let $m^i \stackrel{\text{def}}{=} \big((f_0^i, f_1^i), (A^i, E^i, C^i), (A^{i+1}, E^{i+1}, C^{i+1})\big)$ be a morphism of $\mathcal{A}t\ell(E, C)$. Then

**Composition:**
$\big((f_0^2 \circ f_0^1, f_1^2 \circ f_1^1), (A^1, E^1, C^1), (A^3, E^3, C^3)\big)$ is a morphism of $\mathcal{A}t\ell(E, C)$ by corollary 11.17 (Composition of m-atlas (near) morphisms) on page 78.

**Associativity:**
Composition is associative by lemma 1.19 (Tuple composition for labeled morphisms) on page 12.

**Identity:** $\text{Id}_{(A^i, E^i, C^i)}$ is an identity morphism by lemma 1.19.

□

*Similar results follow with restrictions on the admissible atlases, the admissible morphisms, or both:*

1. $\mathcal{A}t\ell(\mathcal{E}, \mathcal{C})_{\text{full}}$

2. $\mathcal{A}t\ell(\mathcal{E}, \mathcal{C})_{\text{S-max}}$

3. $\mathcal{A}t\ell(\mathcal{E}, \mathcal{C})_{\text{max}}$

4. $\mathcal{A}t\ell_{\text{S-max-full}}(\mathcal{E}, \mathcal{C})$

5. $\mathcal{A}t\ell_{\text{max-full}}(\mathcal{E}, \mathcal{C})$

6. $\mathcal{A}t\ell_{\text{full}}^{\text{semi-strict}}(\mathcal{E}, \mathcal{C})$

7. $\mathcal{A}t\ell_{\text{S-max}}^{\text{semi-strict}}(\mathcal{E}, \mathcal{C})$

8. $\mathcal{A}t\ell_{\text{max}}^{\text{semi-strict}}(\mathcal{E}, \mathcal{C})$

9. $\mathcal{A}t\ell_{\text{S-max-full}}^{\text{semi-strict}}(\mathcal{E}, \mathcal{C})$

10. $\mathcal{A}t\ell_{\text{semi-strict}}^{\text{max-full}}(\mathcal{E}, \mathcal{C})$

11. $\mathcal{A}t\ell_{\text{strict}}^{\text{full}}(\mathcal{E}, \mathcal{C})$

12. $\mathcal{A}t\ell_{\text{strict}}^{\text{S-max}}(\mathcal{E}, \mathcal{C})$

13. $\mathcal{A}t\ell_{\text{strict}}^{\text{max}}(\mathcal{E}, \mathcal{C})$

14. $\mathcal{A}t\ell_{\text{strict}}^{\text{S-max-full}}(\mathcal{E}, \mathcal{C})$

15. $\mathcal{A}t\ell_{\text{strict}}^{\text{max-full}}(\mathcal{E}, \mathcal{C})$

**Definition 11.30** ($\boldsymbol{E}_{\text{triv}}$ and $\mathcal{A}t\ell_{\text{triv}}(\boldsymbol{E}, \boldsymbol{C})$). Let $\boldsymbol{E}$ be a set of topological spaces and $\boldsymbol{C}$ a set of model spaces. Then

$$\boldsymbol{E}_{\text{triv}} \overset{\text{def}}{=} \left\{ E^\mu_{\text{triv}} \,\middle|\, E^\mu \in \boldsymbol{E} \right\} \tag{11.43}$$

$$\mathcal{A}t\ell_{\text{triv}}(\boldsymbol{E}, \boldsymbol{C}) \overset{\text{def}}{=} \mathcal{A}t\ell\left(\boldsymbol{E}_{\text{triv}}, \boldsymbol{C}\right) \tag{11.44}$$

### 11.2.2 M-atlas functors

**Definition 11.31** ($\mathcal{F}_{\text{M,Top}}$). Let $\boldsymbol{E}^i \overset{\text{def}}{=} (E^i, \mathcal{E}^i)$, $i = 1, 2$, and $\boldsymbol{C}^i \overset{\text{def}}{=} (C^i, \mathcal{C}^i)$ be model spaces and $\boldsymbol{A}^i$ be an m-atlas of $\boldsymbol{E}^i$ in the coordinate space $\boldsymbol{C}^i$. Then

$$\mathcal{F}_{\text{M,Top}}(\boldsymbol{A}^i, \boldsymbol{E}^i, \boldsymbol{C}^i) \overset{\text{def}}{=} (\boldsymbol{A}^i, E^i, \boldsymbol{C}^i) \tag{11.45}$$

Let $f_0 \colon \boldsymbol{E}^1 \longrightarrow \boldsymbol{E}^2$ and $f_1 \colon \boldsymbol{C}^1 \longrightarrow \boldsymbol{C}^2$ be model functions. Then

$$\mathcal{F}_{\text{M,Top}}\big((f_0 \colon \boldsymbol{E}^1 \longrightarrow \boldsymbol{E}^2, f_1 \colon \boldsymbol{C}^1 \longrightarrow \boldsymbol{C}^2,), (\boldsymbol{A}^1, \boldsymbol{E}^1, \boldsymbol{C}^1), (\boldsymbol{A}^2, \boldsymbol{E}^2, \boldsymbol{C}^2)\big) \overset{\text{def}}{=}$$
$$\big((f_0 \colon E^1 \longrightarrow E^2, f_1 \colon \boldsymbol{C}^1 \longrightarrow \boldsymbol{C}^2), (\boldsymbol{A}^1, E^1, \boldsymbol{C}^1), (\boldsymbol{A}^2, E^2, \boldsymbol{C}^2)\big) \tag{11.46}$$

**Lemma 11.32** ($\mathcal{F}_{\text{M,Top}}$ maps M-atlases to M-atlases if the coordinate model space is fine grained.)**.** *Let $\boldsymbol{E} \overset{\text{def}}{=} (E, \mathcal{E})$ and $\boldsymbol{C}$ be model spaces, $\boldsymbol{C}$ be fine grained and $\boldsymbol{A}$ be an M-atlas of $\boldsymbol{E}$ in the model space $\boldsymbol{C}$. Then $\boldsymbol{A}$ is an M-atlas of $E$ in the model space $\boldsymbol{C}$.*

*Proof.* Let $(U, V, \phi)$ be an M-chart in $\boldsymbol{A}$, Every open subset of $V$ is a model neighborhood of $\boldsymbol{C}$ since it is fine grained by hypothesis. Then $(U, V, \phi)$ is an M-chart of $E$ in the model space $\boldsymbol{C}$ by lemma 9.3 (M-charts) on page 43.

The definition of mutually compatible charts neighborhood of $C$. $(U, V, \phi)$ is also an M-chart of $E$ in the coordinate space $C$:

Since $U$ is a model neighborhood of $E$, it is open and hence a model neighborhood of $E_{\text{triv}}$.

$V$ is a model neighborhood of $C$ by hypothesis.

Since $U'$ is a model neighborhood of $E_{\text{triv}}$ and hence open, and $\phi$ is a homeomorphism, $\phi[U']$ is open. Since $V$ is a model neighborhood of $C$, $C$ is fine grained and $\phi[U'] \subseteq V$ is open, $\phi[U']$ is a model neighborhood of $C$,

Since $V'$ is a model neighborhood of $C$, it is open, $\phi^{-1}[V']$ is open and hence a model neighborhood of $E_{\text{triv}}$. □

**Theorem 11.33** ($\mathcal{F}_{\text{M,Top}}$ is a functor if $C$ is fine grained.). *Let $E$ and $C$ be sets of model spaces and $C$ be fine grained. Then $\mathcal{F}_{\text{M,Top}}$ is a functor from $\mathcal{A}t\ell(E, C)$ to $\mathcal{A}t\ell(\text{Top}[E], C)$*

*Proof.* Let $o^i \stackrel{\text{def}}{=} (A^i, E^i \stackrel{\text{def}}{=} (E^i, \mathcal{E}^i), C^i), i \in [1,3]$, be an object of $\mathcal{A}t\ell(E, C)$, $o'^i \stackrel{\text{def}}{=} \mathcal{F}_{\text{M,Top}} o^i = (A^i, E^i, C^i)$ the corresponding object of $\mathcal{A}t\ell(\text{Top}[E], C)$, $m^i \stackrel{\text{def}}{=} ((f_0^i \colon E^i \longrightarrow E^{i+1}, f_1^i \colon C^i \longrightarrow C^{i+1},), o^i, o^{i+1}), i = 1, 2$, be a morphism of $\mathcal{A}t\ell(E, C)$ and $m'^i \stackrel{\text{def}}{=} \mathcal{F}_{\text{M,Top}} m^i = ((f_0^i \colon E^i \longrightarrow E^{i+1}, f_1^i \colon C^i \longrightarrow C^{i+1}), o'^i, o'^{i+1})$ the corresponding morphism of $\mathcal{A}t\ell(\text{Top}[E], C)$.

$\mathcal{F}_{\text{M,Top}}$ satisfies these criteria:

**Preservation of objects:**

$A^i$ is an M-atlas of $E^i$ in the coordinate space $C^i$ by lemma 11.32 ($\mathcal{F}_{\text{M,Top}}$ maps M-atlases to M-atlases if the coordinate model space is fine grained) above.

**Preservation of morphisms:**

$f_0^i \colon E^i \longrightarrow E^{i+1}$ is continuous and $f_1^i \colon C^i \longrightarrow C^{i+1}$ is a model function.

**Preservation of endpoints:**

$$\mathcal{F}_{\text{M,Top}} o^i = o'^i \stackrel{\text{Ob}}{\in} \mathcal{A}t\ell(\text{Top}[E], C)$$

$$\mathcal{F}_{\text{M,Top}} m^i = ((f_0^i \colon E^i \longrightarrow E^{i+1}, f_1^i \colon C^i \longrightarrow C^{i+1}), o'^i, o'^{i+1}) \stackrel{\text{Ar}}{\in} \mathcal{A}t\ell(\text{Top}[E], C)$$

is a morphism from $o'^i$ to $o'^{i+1}$.

**Identity:**

Let $((\text{Id}_{E^i}, \text{Id}_{C^i}), o^i, o^i)$ be an identity morphism of $o^i$ in $\mathcal{A}t\ell(E, C)$. Then $((\text{Id}_{E^i} \colon \text{Top}(E^i) \longrightarrow \text{Top}(E^i) \, \text{Id}_{C^i}), o'^i, o'^i)$ is an identity morphism of $o'^i$ in $\mathcal{A}t\ell(\text{Top}[E], C)$.

**Composition:**

$$\begin{aligned}
\mathscr{F}_{M,\text{Top}}(m^2 \overset{A}{\circ} m^1) &= \mathscr{F}_{M,\text{Top}}((f_0^2 \circ f_0^1, f_1^2 \circ f_1^1), o^1, o^3) \\
&= ((f_0^2 \circ f_0^1, f_1^2 \circ f_1^1), o'^1, o'^3) \\
&= ((f_0^2, f_1^2), o'^1, o'^2) \overset{A}{\circ} ((f_0^1, f_1^1), o'^2, o'^3) \\
&= \mathscr{F}_{M,\text{Top}}((f_0^2, f_1^2), o^2, o^3) \overset{A}{\circ} \\
&\quad \mathscr{F}_{M,\text{Top}}((f_0^1, f_1^1), o^1, o^2) \\
&= \mathscr{F}_{M,\text{Top}} m^2 \overset{A}{\circ} \mathscr{F}_{M,\text{Top}} m^1
\end{aligned}$$

$\square$

**Definition 11.34** ($\mathscr{F}_{\text{Top},M}$). Let $E^i$, $i = 1, 2$, be a topological space, $\boldsymbol{C}^i$ a model space and $\boldsymbol{A}^i$ an m-atlas of $E^i$ and $\boldsymbol{A}$ an m-atlas of $E$ in the coordinate space $\boldsymbol{C}^i$. Then

$$\mathscr{F}_{\text{Top},M}(\boldsymbol{A}^i, E^i, \boldsymbol{C}^i) \overset{\text{def}}{=} \left(\boldsymbol{A}, \underset{\text{triv}}{E^i}, \boldsymbol{C}^i\right) \quad (11.47)$$

Let $f_0 \colon E^1 \longrightarrow E^2$ be continuous and $f_1 \colon \boldsymbol{C}^1 \longrightarrow \boldsymbol{C}^2$ be a model function. Then

$$\mathscr{F}_{\text{Top},M}((f_0 \colon E^1 \longrightarrow E^2, f_1 \colon \boldsymbol{C}^1 \longrightarrow \boldsymbol{C}^2,), (\boldsymbol{A}^1, E^1, \boldsymbol{C}^1), (\boldsymbol{A}^2, E^2, \boldsymbol{C}^2)) \overset{\text{def}}{=}$$
$$\left(\left(f_0 \colon \underset{\text{triv}}{E^1} \longrightarrow \underset{\text{triv}}{E^2}, f_1 \colon \boldsymbol{C}^1 \longrightarrow \boldsymbol{C}^2\right), \left(\boldsymbol{A}^1, \underset{\text{triv}}{E^1}, \boldsymbol{C}^1\right), \left(\boldsymbol{A}^2, \underset{\text{triv}}{E^2}, \boldsymbol{C}^2\right)\right) \quad (11.48)$$

**Lemma 11.35** ($\mathscr{F}_{\text{Top},M}$ maps M-atlases to M-atlases). *Let $E$ be a topological space, $\boldsymbol{C}$ be a model space and $\boldsymbol{A}$ be an M-atlas of $E$ in the model space $\boldsymbol{C}$. Then $\boldsymbol{A}$ is an M-atlas of $\underset{\text{triv}}{E}$ in the model space $\boldsymbol{C}$.*

*Proof.* Let $(U, V, \phi)$ be a M-chart in $\boldsymbol{A}$, $U' \subseteq U$ a model neighborhood of $\underset{\text{triv}}{E}$ and $V' \subseteq V$ a model neighborhood of $\boldsymbol{C}$. $(U, V, \phi)$ is also an M-chart of $E$ in the coordinate space $\boldsymbol{C}$: $\square$

**Theorem 11.36** ($\mathscr{F}_{\text{Top},M}$ is a functor). *Let $\boldsymbol{E}$ be a set of topological spaces and $\boldsymbol{C}$ a set of model spaces, Then $\mathscr{F}_{\text{Top},M}$ is a functor from $\mathscr{A}t\ell(\boldsymbol{E}, \boldsymbol{C})$ to $\mathscr{A}t\ell(\underset{\text{triv}}{\boldsymbol{E}}, \boldsymbol{C})$.*

*Proof.* Let $o^i \overset{\text{def}}{=} (\boldsymbol{A}^i, E^i, \boldsymbol{C}^i)$, $i \in [1, 3]$, be an object of $\mathscr{A}t\ell(\boldsymbol{E}, \boldsymbol{C})$, $o'^i \overset{\text{def}}{=} (\boldsymbol{A}^i, \underset{\text{triv}}{E^i}, \boldsymbol{C}^i)$ the corresponding object of $\mathscr{A}t\ell(\underset{\text{triv}}{\boldsymbol{E}}, \boldsymbol{C})$, $m^i \overset{\text{def}}{=} ((f_0^i \colon E^i \longrightarrow E^{i+1}, f_1^i \colon \boldsymbol{C}^i \longrightarrow \boldsymbol{C}^{i+1}), o^i, o^{i+1})$, $i = 1, 2$, be a morphism of $\mathscr{A}t\ell(\boldsymbol{E}, \boldsymbol{C})$ and $m'^i \overset{\text{def}}{=} ((f_0^i, f_1^i), o'^i, o'^{i+1})$ the corresponding morphism of $\mathscr{A}t\ell(\underset{\text{triv}}{\boldsymbol{E}}, \boldsymbol{C})$.

$\mathscr{F}_{\text{Top},M}$ satisfies these criteria:

**Preservation of endpoints:** $A^i$ is an M-atlas of $E^i_{\text{triv}}$ in the coordinate space $\boldsymbol{C}^i$ by lemma 11.32 (<span style="color:red">$\mathscr{F}_{\text{M,Top}}$ maps M-atlases to M-atlases if the coordinate model space is fine grained</span>) above. $f_0^i\colon E^i \longrightarrow E^{i+1}$ is continuous and $f_1^i\colon \boldsymbol{C}^i \longrightarrow \boldsymbol{C}^{i+1}$ is a model function.

$$\mathscr{F}_{\text{Top,M}}\, o^i = o'^i \overset{\text{Ob}}{\in} \mathscr{A}t\ell\big(\underset{\text{triv}}{\boldsymbol{E}}, \boldsymbol{C}\big)$$

$$\mathscr{F}_{\text{Top,M}}\, m^i = m'^i \overset{\text{Ar}}{\in} \mathscr{A}t\ell\big(\underset{\text{triv}}{\boldsymbol{E}}, \boldsymbol{C}\big)$$

is a morphism from $o'^i$ to $o'^{i+1}$.

**Identity:** $\text{Id}_{o^i} \overset{\text{def}}{=} \big((\text{Id}_{E^i}, \text{Id}_{\boldsymbol{C}^i}), o^i, o^i\big)$. Then $\mathscr{F}_{\text{Top,M}}\, \text{Id}_{o^i} = \big((\text{Id}_{E^i_{\text{triv}}}, \text{Id}_{\boldsymbol{C}^i}), o'^i, o'^i\big)$ is an identity morphism of $o'^i$.

**Composition:**

$$\begin{aligned}
\mathscr{F}_{\text{Top,M}}(m^2 \overset{A}{\circ} m^1) &= \mathscr{F}_{\text{Top,M}}\big((f_0^2 \circ f_0^1, f_1^2 \circ f_1^1), o^1, o^3\big) \\
&= \big((f_0^2 \circ f_0^1, f_1^2 \circ f_1^1), o'^1, o'^3\big) \\
&= \big((f_0^2, f_1^2), o'^1, o'^2\big) \overset{A}{\circ} \big((f_0^1, f_1^1), o'^2, o'^3\big) \\
&= \mathscr{F}_{\text{Top,M}}\big((f_0^2, f_1^2), o^2, o^3\big) \overset{A}{\circ} \\
&\quad \mathscr{F}_{\text{Top,M}}\big((f_0^1, f_1^1), o^1, o^2\big) \\
&= \mathscr{F}_{\text{Top,M}}\, m^2 \overset{A}{\circ} \mathscr{F}_{\text{Top,M}}\, m^1
\end{aligned}$$

□

$\mathscr{F}_{\text{M,Top}} \circ \mathscr{F}_{\text{Top,M}}$ *is the identity functor on* $\mathscr{A}t\ell(\boldsymbol{E}, \boldsymbol{C})$.

*Proof.*

$$\begin{aligned}
\mathscr{F}_{\text{M,Top}} \circ \mathscr{F}_{\text{Top,M}}\, o^i &= \mathscr{F}_{\text{M,Top}}\, o'^i \\
&= o^i
\end{aligned}$$

$$\begin{aligned}
\mathscr{F}_{\text{M,Top}} \circ \mathscr{F}_{\text{Top,M}}\, m^i &= \mathscr{F}_{\text{M,Top}}\, m'^i \\
&= m^i
\end{aligned}$$

□

# 12 Associated model spaces and functors

**Definition 12.1** (Coordinate model spaces associated with m-atlases). Let $A^i$, $i = 1, 2$, be an m-atlas of $E^i$ in the coordinate space $C^i$, $f \stackrel{\text{def}}{=} (f_0, f_1)$ be an $E^1$-$E^2$ m-atlas (near) morphism $(f_0, f_1)$ of $A^1$ to $A^2$ in the coordinate spaces $C^1$, Then:

The minimal coordinate model space with neighborhoods in the m-atlas $A^i$ of $E^i$ in the coordinate space $C^i$ is

$$\mathcal{F}_2^{\min \text{M}}(A^i, E^i, C^i) \stackrel{\text{def}}{=}$$
$$\operatorname*{Mod}_{\min}\left(C^i, \pi_2[A^i], \left\{\phi' \circ \phi^{-1} \;\middle|\; \left(\exists_{\substack{(U,V,\phi) \in A^i \\ (U',V',\phi') \in A^i}}\right) U \cap U' \neq \emptyset\right\}\right) \quad (12.1)$$

The coordinate mapping associated with the $E^1$-$E^2$ m-atlas (near) morphism $(f_0, f_1)$ of $A^1$ to $A^2$ in the coordinate spaces $C^1$, $C^2$ is

$$\mathcal{F}_2^{\min \text{M}}((f_0, f_1), (A^1, E^1, C^1), (A^2, E^2, C^2)) \stackrel{\text{def}}{=}$$
$$f_1 \colon \mathcal{F}_2^{\min \text{M}}(A^1, E^1, C^1) \longrightarrow \mathcal{F}_2^{\min \text{M}}(A^2, E^2, C^2) \quad (12.2)$$

If it is a model function then it is also the coordinate m-morphism associated with the $E^1$-$E^2$ m-atlas (near) morphism $(f_0, f_1)$ of $A^1$ to $A^2$ in the coordinate spaces $C^1$, $C^2$.

**Lemma 12.2** (Coordinate model spaces associated with m-atlases). *Let $A^i$ be an m-atlas of $E^i$ in the coordinate space $C^i$, $i = 1, 2$.*
$\mathcal{F}_2^{\min \text{M}}(A^i, E^i, C^i)$ *is a model space.*

*Proof.* $\mathcal{F}_2^{\min \text{M}}(A^i, E^i, C^i)$ satisfies the conditions for a model space.

1. Since $\pi_2[A^i]$ is an open cover of $\bigcup \pi_2[A^i]$, the set of finite intersections is also an open cover.

2. Finite intersections of finite intersections are finite intersections

3. Restrictions of continuous functions are continuous

4. If $f \colon A \longrightarrow B$ is a morphism of $\mathcal{F}_2^{\min \text{M}}(A^i, E^i, C^i)$, $A, A', B, B'$ model meighborhoods of $\mathcal{F}_2^{\min \text{M}}(A^i, E^i, C^i)$, $A' \subseteq A$, $B' \subseteq B$ and $f[A'] \subseteq B'$ then since $f \colon A \longrightarrow B$ is a morphism it is a restriction of a transition function between sets in $\pi_2[A^i]$ and its restrictions are also, hence morphisms, and thus $f\restriction_{A'} A' \longrightarrow B'$ is a morphism.

5. If $(U, V, \phi) \in A^i$ then $\text{Id}_V = \phi \circ \phi^{-1}$ is a transition function and hence a morphism of $\mathcal{F}_2^{\overset{M}{\min}}(A^i, E^i, C^i)$. If $A, A'$ are objects of $\pi_2(\mathcal{F}_2^{\overset{M}{\min}}(A^i, E^i, C^i))$ and $A' \subseteq A$ then the inclusion map $i\colon A' \hookrightarrow A$ is a restriction of an identity morphism of $\mathcal{F}_2^{\overset{M}{\min}}(A^i, E^i, C^i)$ and hence a morphism.

6. Restricted sheaf condition: let

    (a) $U_\alpha, V_\alpha, \alpha \prec A$, be model neighborhoods of $\mathcal{F}_2^{\overset{M}{\min}}(A^i, E^i, C^i)$

    (b) $f_\alpha \colon U_\alpha \longrightarrow V_\alpha$ be a morphism of $\mathcal{F}_2^{\overset{M}{\min}}(A^i, E^i, C^i)$

    (c) $U \overset{\text{def}}{=} \bigcup_{\alpha \prec A} U_\alpha$

    (d) $V \overset{\text{def}}{=} \bigcup_{\alpha \prec A} V_\alpha$

    (e) $f\colon U \longrightarrow V$ be continuous and $\left(\underset{\substack{\alpha \prec A \\ x \in U_\alpha}}{\forall}\right) f(x) = f_\alpha(x)$

    Then $f$ is a morphism of $C^i$ and hence a morphism of $\mathcal{F}_2^{\overset{M}{\min}}(A^i, E^i, C^i)$.

    $\square$

Let $\boldsymbol{f} \overset{\text{def}}{=} (f_0, f_1)$ be an $E^1$-$E^2$ m-atlas (near) morphism from $A^1$ to $A^2$ in the coordinate spaces $C^1, C^2$. If $\boldsymbol{f}$ is a morphism or each $A^i$ is semi-maximal then $f_1 \colon \mathcal{F}_2^{\overset{M}{\min}}(A^1, E^1, C^1) \longrightarrow \mathcal{F}_2^{\overset{M}{\min}}(A^2, E^2, C^2)$ is well defined, i.e., $f_1\left[\mathcal{F}_2^{\overset{M}{\min}}(A^1, E^1, C^1)\right] \subseteq \mathcal{F}_2^{\overset{M}{\min}}(A^2, E^2, C^2)$.

*Proof.* $\boldsymbol{f}$ is a morphism either by hypothesis or by lemma 11.12 (M-atlas (near) morphisms) on page 60. Let $v^1 \in \mathcal{F}_2^{\overset{M}{\min}}(A^1, E^1, C^1)$, $(U^1, V^1, \phi^1) \in A^1$ be a chart with $v^1 \in V^1$, $u^1 \overset{\text{def}}{=} \phi^{1-1}(v^1)$, $u^2 \overset{\text{def}}{=} f_0(u^1)$ and $(U^2, V^2, \phi^2) \in A^2$ be a chart at $u^2$. Let $(U'^1, V'^1, \phi^1 \colon U'^1 \overset{\sim}{\twoheadrightarrow} V'^1) \in A^1$, $(U'^2, \hat{V}'^2, \phi'^2 \colon U'^2 \overset{\sim}{\twoheadrightarrow} \hat{V}'^2) \in A^2$ be as in definition 11.4 (M-atlas morphisms for model spaces) on page 53. Then $f_1(v^1) \in \hat{V}'^2 \subseteq \mathcal{F}_2^{\overset{M}{\min}}(A^1, E^1, C^1)$.  $\square$

**Definition 12.3** (Model spaces associated with m-atlases). Let $A^i$, $i = 1, 2$, be an m-atlas of $E^i$ in the coordinate space $C^i$ and $\boldsymbol{f} \overset{\text{def}}{=} (f_0, f_1)$ be an $E^1$-$E^2$ m-atlas (near) morphism $(f_0, f_1)$ of $A^1$ to $A^2$ in the coordinate spaces $C^1$, Then:

The minimal model space with neighborhoods in the m-atlas $A^i$ of $E^i$ in the coordinate space $C^i$ is

$$\overset{\text{M}}{\mathcal{F}}{}_1^{\min}(A^i, E^i, C^i) \overset{\text{def}}{=}$$

$$\underset{\min}{\text{Mod}}\left(E^i, \pi_1[A^i], \left\{\phi'^{-1} \circ \phi \;\middle|\; \left(\underset{(U',V',\phi') \in A^i}{\exists_{(U,V,\phi) \in A^i}}\right) U \cap U' \neq \emptyset\right\}\right) \quad (12.3)$$

The mapping associated with the $E^1$-$E^2$ m-atlas morphism $(f_0, f_1)$ of $A^1$ to $A^2$ in the coordinate spaces $C^1, C^2$ is

$$\overset{\text{M}}{\mathcal{F}}{}_1^{\min}\big((f_0, f_1), (A^1, E^1, C^1), (A^2, E^2, C^2)\big) \overset{\text{def}}{=}$$

$$f_0 \colon \overset{\text{M}}{\mathcal{F}}{}_1^{\min}(A^1, E^1, C^1) \longrightarrow \overset{\text{M}}{\mathcal{F}}{}_1^{\min}(A^2, E^2, C^2) \quad (12.4)$$

If it is a model function then the it is also the m-atlas morphism associated with the $E^1$-$E^2$ m-atlas morphism $(f_0, f_1)$ of $A^1$ to $A^2$ in the coordinate spaces $C^1, C^2$.

**Lemma 12.4** (Model spaces associated with m-atlases). *Let $A$ be an m-atlas of $E$ in the coordinate space $C$. Then $\overset{\text{M}}{\mathcal{F}}{}_1^{\min}(A, E, C)$ is a model space.*

*Proof.* The result follows from Lemma 5.4 (Minimal model spaces are model spaces) on page 27. □

**Theorem 12.5** (Functors from m-atlases to model spaces). *Let $E$ and and $C$ be sets of model spaces.*

*$\overset{\text{M}}{\mathcal{F}}{}_1^{\min}$ is a functor from $\mathcal{A}t\ell(E, C)$ to $\underset{\text{full-submod-}\mathcal{T}riv}{E} \overset{\text{fullcat}}{\subseteq} \underset{\text{submod-}\mathcal{T}riv}{E}$*

*$\overset{\text{M}}{\mathcal{F}}{}_1^{\min}$ is a functor from $\underset{\text{full}}{\mathcal{A}t\ell(E, C)}$ to $\underset{\text{full-submod-}\mathcal{T}riv}{E} \overset{\text{fullcat}}{\subseteq} \underset{\text{submod-}\mathcal{T}riv}{E}$*

*$\overset{\text{M}}{\mathcal{F}}{}_2^{\min}$ is a functor from $\underset{\text{S-max}}{\mathcal{A}t\ell(E, C)}$ to $\underset{\text{submod-}\mathcal{T}riv}{C}$.*

*$\overset{\text{M}}{\mathcal{F}}{}_2^{\min}$ is a functor from $\underset{\text{full}}{\mathcal{A}t\ell(E, C)}$ to $\underset{\text{full-submod-}\mathcal{T}riv}{C} \overset{\text{fullcat}}{\subseteq} \underset{\text{submod-}\mathcal{T}riv}{C}$*

*Proof.* $\underset{\text{full-submod-}\mathcal{T}riv}{E} \overset{\text{fullcat}}{\subseteq} \underset{\text{submod-}\mathcal{T}riv}{E}$ and $\underset{\text{full-submod-}\mathcal{T}riv}{C} \overset{\text{fullcat}}{\subseteq} \underset{\text{submod-}\mathcal{T}riv}{C}$ by lemma 4.2 (Trivial model spaces) on page 23.

Let $o^i \overset{\text{def}}{=} (A^i, E^i, C^i)$, $i \in [1, 3]$, be objects in $\mathcal{A}t\ell(E, C)$ and let $m^i \overset{\text{def}}{=} \big((f_0^i, f_1^i), o^i, o^{i+1}\big)$ be a morphism in $\mathcal{A}t\ell(E, C)$.

$\overset{\text{M}}{\mathcal{F}}{}_1^{\min} \colon \mathcal{A}t\ell(E, C) \longrightarrow \underset{\text{submod-}\mathcal{T}riv}{E} :$

**Preservation of endpoints:** $\overset{\text{M}}{\mathcal{F}}{}_1^{\min}(m^i) = f_0^i \colon \overset{\text{M}}{\mathcal{F}}{}_1^{\min} o^i \longrightarrow \overset{\text{M}}{\mathcal{F}}{}_1^{\min} o^{i+1}$

**Composition:**

$$\mathcal{F}_1^{\mathrm{min}\,M}(m^2 \overset{A}{\circ} m^1) = \mathcal{F}_1^{\mathrm{min}\,M}((f_0^2 \circ f_0^1, f_1^2 \circ f_1^1)(\boldsymbol{A}^1, \boldsymbol{E}^1, \boldsymbol{C}^1), (\boldsymbol{A}^3, \boldsymbol{E}^3, \boldsymbol{C}^3))$$

$$= f_0^2 \circ f_0^1 \colon \mathcal{F}_1^{\mathrm{min}\,M} o^1 \longrightarrow \mathcal{F}_1^{\mathrm{min}\,M} o^3$$

$$= (f_0^2 \colon \mathcal{F}_1^{\mathrm{min}\,M} o^2 \longrightarrow \mathcal{F}_1^{\mathrm{min}\,M} o^3) \circ (f_0^1 \colon \mathcal{F}_1^{\mathrm{min}\,M} o^1 \longrightarrow \mathcal{F}_1^{\mathrm{min}\,M} o^2)$$

$$= \mathcal{F}_1^{\mathrm{min}\,M}((f_0^2, f_1^2)(\boldsymbol{A}^2, \boldsymbol{E}^2, \boldsymbol{C}^2), (\boldsymbol{A}^3, \boldsymbol{E}^3, \boldsymbol{C}^3)) \circ$$

$$\mathcal{F}_1^{\mathrm{min}\,M}((f_0^1, f_1^1)(\boldsymbol{A}^1, \boldsymbol{E}^1, \boldsymbol{C}^1), (\boldsymbol{A}^2, \boldsymbol{E}^2, \boldsymbol{C}^2))$$

$$= \mathcal{F}_1^{\mathrm{min}\,M}(m^2) \circ \mathcal{F}_1^{\mathrm{min}\,M}(m^1)$$

**Identity:**

1. $\mathcal{F}_1^{\mathrm{min}\,M}(\mathrm{Id}_{o^i}) = \mathcal{F}_1^{\mathrm{min}\,M}((\mathrm{Id}_{E^i}, \mathrm{Id}_{C^i}), (\boldsymbol{A}^i, \boldsymbol{E}^i, \boldsymbol{C}^i), (\boldsymbol{A}^i, \boldsymbol{E}^i, \boldsymbol{C}^i)) =$
   $\mathrm{Id}_{E^i} \colon \mathcal{F}_1^{\mathrm{min}\,M} o^i \longrightarrow \mathcal{F}_1^{\mathrm{min}\,M} o^i$

2. $\mathrm{Id}_{\mathcal{F}_1^{\mathrm{min}\,M}\, o^i} = \mathrm{Id}_{E^i} \colon \mathcal{F}_1^{\mathrm{min}\,M} o^i \longrightarrow \mathcal{F}_1^{\mathrm{min}\,M} o^i$

The proof for $\mathcal{F}_1^{\mathrm{min}\,M} \colon \mathcal{A}t\ell(\boldsymbol{E}, \boldsymbol{C}) \underset{\mathrm{full}}{\longrightarrow} \underset{\mathrm{full-submod-}\mathcal{T}riv}{\boldsymbol{E}} \overset{\mathrm{fullcat}}{\subseteq}$ is identical.

$\mathcal{F}_2^{\mathrm{min}\,M} \colon \mathcal{A}t\ell(\boldsymbol{E}, \boldsymbol{C}) \longrightarrow \underset{\mathrm{op-triv}}{\boldsymbol{C}}$:

**Preservation of endpoints:**

$$\mathcal{F}_2^{\mathrm{min}\,M}(m^i) = f_1^i \colon \mathcal{F}_2^{\mathrm{min}\,M} o^i \longrightarrow \mathcal{F}_2^{\mathrm{min}\,M} o^{i+1}$$

$\mathcal{F}(g \circ f) = \mathcal{F}(g) \circ \mathcal{F}(f) \colon \mathcal{F}_2^{\mathrm{min}\,M}(m^2 \overset{A}{\circ} m^1) =$
$\mathcal{F}_2^{\mathrm{min}\,M}((f_0^2 \circ f_0^1, f_1^2 \circ f_1^1)(\boldsymbol{A}^1, \boldsymbol{E}^1, \boldsymbol{C}^1), (\boldsymbol{A}^3, \boldsymbol{E}^3, \boldsymbol{C}^3)) =$
$f_0^2 \circ f_0^1 \colon \mathcal{F}_2^{\mathrm{min}\,M} o^1 \longrightarrow \mathcal{F}_2^{\mathrm{min}\,M} o^3 =$
$(f_1^2 \colon \mathcal{F}_2^{\mathrm{min}\,M} o^2 \longrightarrow \mathcal{F}_2^{\mathrm{min}\,M} o^3) \circ (f_1^1 \colon \mathcal{F}_2^{\mathrm{min}\,M} o^1 \longrightarrow \mathcal{F}_2^{\mathrm{min}\,M} o^2) =$
$\mathcal{F}_2^{\mathrm{min}\,M}((f_0^2, f_1^2)(\boldsymbol{A}^2, \boldsymbol{E}^2, \boldsymbol{C}^2), (\boldsymbol{A}^3, \boldsymbol{E}^3, \boldsymbol{C}^3)) \circ$
$\mathcal{F}_2^{\mathrm{min}\,M}((f_0^1, f_1^1)(\boldsymbol{A}^1, \boldsymbol{E}^1, \boldsymbol{C}^1), (\boldsymbol{A}^2, \boldsymbol{E}^2, \boldsymbol{C}^2)) =$
$\mathcal{F}_2^{\mathrm{min}\,M}(m^2) \circ \mathcal{F}_2^{\mathrm{min}\,M}(m^1)$

$\mathcal{F}(\mathrm{Id}_A) = \mathrm{Id}_{\mathcal{F}(A)}$:

1. $\mathcal{F}_2^{\mathrm{min}\,M}(\mathrm{Id}_{o^i}) = \mathcal{F}_2^{\mathrm{min}\,M}((\mathrm{Id}_{E^i}, \mathrm{Id}_{C^i}), (\boldsymbol{A}^i, \boldsymbol{E}^i, \boldsymbol{C}^i), (\boldsymbol{A}^i, \boldsymbol{E}^i, \boldsymbol{C}^i)) =$
   $\mathrm{Id}_{C^i} \colon \mathcal{F}_2^{\mathrm{min}\,M} o^i \longrightarrow \mathcal{F}_2^{\mathrm{min}\,M} o^i$

2. $\text{Id}_{\mathscr{F}_2^{\text{M min}}} o^i = \text{Id}_{C^i} \colon \mathscr{F}_2^{\text{M min}} o^i \longrightarrow \mathscr{F}_2^{\text{M min}} o^i$

The proof for $\mathscr{F}_2^{\text{M min}} \colon \underset{\text{full}}{\mathscr{Atl}(\boldsymbol{E}, \boldsymbol{C})} \longrightarrow \underset{\text{full-submod-}\mathscr{T}riv}{\boldsymbol{C}} \overset{\text{fullcat}}{\subseteq}$ is identical. $\square$

# 13 Classic m-atlas morphisms and functors

This subsection introduces an alternate definition of and taxonomy for morphisms between m-atlases, defines categories of m-atlases, defines functors amomg them and proves some basic reasults. It introduces the notions of classic m-atlas morphisms, which are very similar to the conventional definitions for a manifold, although it does not require the space to be locally Euclidean.

While it is easy to define satisfactory notions of classic m-morphisms, there are technical difficulties in defining a satisfactory notion of classic near m-morphisms.

## 13.1 Classic m-atlas morphisms

This subsection introduces the notions of classic m-atlas morphisms, and proves some basic results.

### 13.1.1 Definitions of classic m-atlas morphisms

**Definition 13.1** (Classic m-atlas morphisms for model spaces). Let $\mathscr{E}^i, \mathscr{C}^i, i = 1, 2$, be a model category, $\boldsymbol{E}^i \overset{\text{Ob}}{\in} \mathscr{E}^i$, $\boldsymbol{C}^i \overset{\text{Ob}}{\in} \mathscr{C}^i$ and $\boldsymbol{A}^i$ be an m-atlas of $\boldsymbol{E}^i$ in the coordinate space $\boldsymbol{C}^i$.

$f \colon \boldsymbol{E}^1 \longrightarrow \boldsymbol{E}^2$ is a classic $\boldsymbol{E}^1$-$\boldsymbol{E}^2$ m-atlas morphism of $\boldsymbol{A}^1$ to $\boldsymbol{A}^2$ in the coordinate spaces $\boldsymbol{C}^1, \boldsymbol{C}^2$, abbreviated as $\text{isAtl}_{\text{Ar}}^{\text{classic}}(\boldsymbol{A}^1, \boldsymbol{E}^1, \boldsymbol{C}^1, \boldsymbol{A}^2, \boldsymbol{E}^2, \boldsymbol{C}^2, f)$. and a classic $\mathscr{E}^1$-$\mathscr{E}^2$ m-atlas morphism of $\boldsymbol{A}^1$ to $\boldsymbol{A}^2$ in the cordinate model categories $\mathscr{E}^1, \mathscr{E}^2$, abbreviated as $\text{isAtl}_{\text{Ar}}^{\text{classic}}(\boldsymbol{A}^1, \mathscr{E}^1, \mathscr{C}^1, \boldsymbol{A}^2, \mathscr{E}^2, \mathscr{C}^2, f)$, iff for any $(U^i, V^i, \phi^i \colon U^i \overset{\tilde{=}}{\rightarrowtail\!\!\!\twoheadrightarrow} V^i) \in \boldsymbol{A}^i$, $i = 1, 2$, with $I \overset{\text{def}}{=} U^1 \cap f^{-1}[U^2] \neq \emptyset$, $\phi^2 \circ f \circ \phi^{1-1} \colon \phi^1[I] \longrightarrow V^2$ is a $\boldsymbol{C}^1$-$\boldsymbol{C}^2$ m-morphism of $\phi^1[I]$ to $V^2$.

$\boldsymbol{f}$ is also a constrained classic $E^1$-$E^2$ m-atlas morphism of $\boldsymbol{A}^1$ to $\boldsymbol{A}^2$ in the coordinate spaces $\boldsymbol{C}^1, \boldsymbol{C}^2$, abbreviated as $\text{isAtl}_{\text{Ar}}^{\text{constrained,classic}}(\boldsymbol{A}^1, \boldsymbol{E}^1, \boldsymbol{C}^1, \boldsymbol{A}^2, \boldsymbol{E}^2, \boldsymbol{C}^2, f_0, f_1)$, and a constrained classic $E^1$-$E^2$ m-atlas morphism of $\boldsymbol{A}^1$ to $\boldsymbol{A}^2$ in the coordinate model categories $\mathscr{E}^1, \mathscr{E}^2$, abbreviated as $\text{isAtl}_{\text{Ar}}^{\text{constrained,classic}}(\boldsymbol{A}^1, E^1, \mathscr{E}^1, \boldsymbol{A}^2, E^2, \mathscr{E}^2, f_0, f_1)$, iff $f$ is constrained, i.e., $f[\boldsymbol{E}]$ is contained in a model neighborhood of $\boldsymbol{E}^2$.

$f \colon \boldsymbol{E}^1 \longrightarrow \boldsymbol{E}^2$ is a classic $\boldsymbol{E}^1$-$\boldsymbol{E}^2$-$\mathscr{E}^1$-$\mathscr{E}^2$ m-atlas morphism of $\boldsymbol{A}^1$ to $\boldsymbol{A}^2$ in the coordinate model spaces $\mathscr{E}^1, \mathscr{E}^2$. abbreviated as $\text{isAtl}_{\text{Ar}}^{\text{classic}}(\boldsymbol{A}^1, \boldsymbol{E}^1, \boldsymbol{C}^1, \boldsymbol{A}^2, \boldsymbol{E}^2, \boldsymbol{C}^2, f, \mathscr{E}^1, \mathscr{E}^2)$, and a classic $\mathscr{E}^1$-$\mathscr{E}^2$ m-atlas morphism of $\boldsymbol{A}^1$ to $\boldsymbol{A}^2$ in the cordinate model categories $\mathscr{E}^1$, $\mathscr{E}^2$, abbreviated as $\text{isAtl}_{\text{Ar}}^{\text{classic}}(\boldsymbol{A}^1, \mathscr{E}^1, \mathscr{C}^1, \boldsymbol{A}^2, \mathscr{E}^2, \mathscr{C}^2, f)$, iff $\mathscr{E}^1 \overset{\text{cat}}{\subseteq} \mathscr{E}^2$ and for any $(U^i, V^i, \phi^i\colon U^i \rightarrowtail\!\!\!\!\!\xrightarrow{\sim} V^i) \in \boldsymbol{A}^i$, $i = 1, 2$, with $I \stackrel{\text{def}}{=} U^1 \cap f^{-1}[U^2] \neq \emptyset$,
$\phi^2 \circ f \circ \phi^{1-1}\colon \phi^1[I] \longrightarrow V^2$ is a morphism of $\mathscr{C}^2$.

Let $\boldsymbol{f}$ be a classic $\boldsymbol{E}^1$-$\boldsymbol{E}^2$ m-atlas morphism of $\boldsymbol{A}^1$ to $\boldsymbol{A}^2$ in the coordinate spaces $\boldsymbol{C}^1$, $\boldsymbol{C}^2$.

The triple $\big((f), (\boldsymbol{A}^1, \boldsymbol{E}^1, \boldsymbol{C}^1), (\boldsymbol{A}^2, \boldsymbol{E}^2, \boldsymbol{C}^2)\big)$ will refer to $f$ considered as a classic $\boldsymbol{E}^1$-$\boldsymbol{E}^2$ ($\boldsymbol{E}^1$-$\boldsymbol{E}^2$-$\mathscr{C}^1$-$\mathscr{C}^2$) m-atlas morphism of $\boldsymbol{A}^1$ to $\boldsymbol{A}^2$ in the coordinate spaces $\boldsymbol{C}^1$, $\boldsymbol{C}^2$.

The classic m-atlas identity morphism of $(\boldsymbol{A}^i, \boldsymbol{E}^i, \boldsymbol{C}^i)$ is

$$\text{Id}_{(\boldsymbol{A}^i, \boldsymbol{E}^i, \boldsymbol{C}^i)} \stackrel{\text{def}}{=} \big(\text{Id}_{\boldsymbol{E}^i}(\boldsymbol{A}^i, \boldsymbol{E}^i, \boldsymbol{C}^i), (\boldsymbol{A}^i, \boldsymbol{E}^i, \boldsymbol{C}^i)\big) \tag{13.1}$$

**Definition 13.2** (Classic m-atlas morphisms for topological spaces). Let $\mathscr{C}^i$, $i = 1, 2$, be a model category, $\mathscr{E}^i$ be a topological category, $E^i \stackrel{\text{Ob}}{\in} \mathscr{E}^i$, $\boldsymbol{C}^i \stackrel{\text{Ob}}{\in} \mathscr{C}^i$ and $\boldsymbol{A}^i$ be an m-atlas of $E^i$ in the coordinate space $\boldsymbol{C}^i$.

$f\colon E^1 \longrightarrow E^2$ is a classic $E^1$-$E^2$ m-atlas morphism of $\boldsymbol{A}^1$ to $\boldsymbol{A}^2$ in the coordinate spaces $\boldsymbol{C}^1$, $\boldsymbol{C}^2$, abbreviated as $\text{isAtl}_{\text{Ar}}^{\text{classic}}(\boldsymbol{A}^1, E^1, \boldsymbol{C}^1, \boldsymbol{A}^2, E^2, \boldsymbol{C}^2, f)$. and a classic $\mathscr{E}^1$-$\mathscr{E}^2$ m-atlas morphism of $\boldsymbol{A}^1$ to $\boldsymbol{A}^2$ in the cordinate model categories $\mathscr{C}^1$, $\mathscr{C}^2$, abbreviated as $\text{isAtl}_{\text{Ar}}^{\text{classic}}(\boldsymbol{A}^1, \mathscr{E}^1, \mathscr{C}^1, \boldsymbol{A}^2, \mathscr{E}^2, \mathscr{C}^2, f)$, iff $f$ is a classic $\underset{\textbf{triv}}{E^1}$-$\underset{\textbf{triv}}{E^2}$ m-atlas morphism of $\boldsymbol{A}^1$ to $\boldsymbol{A}^2$ in the coordinate spaces $\boldsymbol{C}^1$, $\boldsymbol{C}^2$.

$f\colon E^1 \longrightarrow E^2$ is a classic $E^1$-$E^2$-$\mathscr{C}^1$-$\mathscr{C}^2$ m-atlas morphism of $\boldsymbol{A}^1$ to $\boldsymbol{A}^2$ in the coordinate spaces $\boldsymbol{C}^1$, $\boldsymbol{C}^2$, abbreviated as $\text{isAtl}_{\text{Ar}}^{\text{classic}}(\boldsymbol{A}^1, E^1, \boldsymbol{C}^1, \boldsymbol{A}^2, E^2, \boldsymbol{C}^2, f, \mathscr{C}^1, \mathscr{C}^2)$. and a classic $\mathscr{E}^1$-$\mathscr{E}^2$-$\mathscr{C}^1$-$\mathscr{C}^2$ m-atlas morphism of $\boldsymbol{A}^1$ to $\boldsymbol{A}^2$ in the cordinate model categories $\mathscr{C}^1$, $\mathscr{C}^2$, abbreviated as $\text{isAtl}_{\text{Ar}}^{\text{classic}}(\boldsymbol{A}^1, \mathscr{E}^1, \mathscr{C}^1, \boldsymbol{A}^2, \mathscr{E}^2, \mathscr{C}^2, f, \mathscr{C}^1, \mathscr{C}^2)$, iff $f$ is a classic $\underset{\textbf{triv}}{E^1}$-$\underset{\textbf{triv}}{E^2}$-$\mathscr{C}^1$-$\mathscr{C}^2$ m-atlas morphism of $\boldsymbol{A}^1$ to $\boldsymbol{A}^2$ in the coordinate spaces $\boldsymbol{C}^1$, $\boldsymbol{C}^2$.

The triple $\big((f), (\boldsymbol{A}^1, E^1, \boldsymbol{C}^1), (\boldsymbol{A}^2, E^2, \boldsymbol{C}^2)\big)$ will refer to $f$ considered as a classic $E^1$-$E^2$ m-atlas morphism of $\boldsymbol{A}^1$ to $\boldsymbol{A}^2$ in the coordinate spaces $\boldsymbol{C}^1$, $\boldsymbol{C}^2$.

### 13.1.2 Definitions of semi-strict and strict

**Definition 13.3** (Semi-strict classic m-atlas (near) morphisms). Let $\mathscr{E}^i, \mathscr{C}^i, i = 1, 2$, be model categories, $\boldsymbol{E}^i \stackrel{\text{Ob}}{\in} \mathscr{E}^i$, $\boldsymbol{C}^i \stackrel{\text{Ob}}{\in} \mathscr{C}^i$ and $\boldsymbol{A}^i$ be an m-atlas of $\boldsymbol{E}^i$ in the coordinate space $\boldsymbol{C}^i$.

Let $f\colon \boldsymbol{E}^1 \longrightarrow \boldsymbol{E}^2$ be a classic $\boldsymbol{E}^1$-$\boldsymbol{E}^2$ m-atlas (near) morphism of $\boldsymbol{A}^1$ to $\boldsymbol{A}^2$ in the coordinate spaces $\boldsymbol{C}^1$, $\boldsymbol{C}^2$.

$f$ is a semi-strict classic $\boldsymbol{E}^1$-$\boldsymbol{E}^2$ m-atlas morphism of $\boldsymbol{A}^1$ to $\boldsymbol{A}^2$ in the coordinate spaces $\boldsymbol{C}^1$, $\boldsymbol{C}^2$, abbreviated as $\text{isAtl}_{\underset{\text{semi-strict}}{\text{Ar}}}^{\text{classic}}(\boldsymbol{A}^1, \boldsymbol{E}^1, \boldsymbol{C}^1, \boldsymbol{A}^2, \boldsymbol{E}^2, \boldsymbol{C}^2, f)$, iff $f$ is a local m-morphism of $\boldsymbol{E}^1$ to $\boldsymbol{E}^2$.

*Remark* 13.4. Definitions of semi-strict classic m-atlas morphisms for topological spaces would be pointless, as any semi-strict classic m-atlas morphism would be strict.

It is a semi-strict classic $\mathscr{E}^1$-$\mathscr{E}^2$-$E^1$-$E^2$ m-atlas morphism of $A^1$ to $A^2$ in the coordinate spaces $C^1, C^2$, abbreviated as
$\text{isAtl}_{\text{Ar}\atop\text{semi-strict}}^{\text{classic}}(A^1, E^1, C^1, A^2, E^2, C^2, f, \mathscr{E}^1, \mathscr{E}^2)$, iff $f$ is a local $\mathscr{E}^1$-$\mathscr{E}^2$ m-morphism of $E^1$ to $E^2$.

Let $f\colon E^1 \longrightarrow E^2$ be a classic $E^1$-$E^2$-$\mathscr{C}^1$-$\mathscr{C}^2$ m-atlas morphism of $A^1$ to $A^2$ in the coordinate spaces $C^1, C^2$.

$f$ is a semi-strict classic $E^1$-$E^2$-$C^1$-$C^2$ m-atlas morphism of $A^1$ to $A^2$ in the coordinate spaces $C^1, C^2$, abbreviated as
$\text{isAtl}_{\text{Ar}\atop\text{semi-strict}}^{\text{classic}}(A^1, E^1, C^1, A^2, E^2, C^2, f, \mathscr{C}^1, \mathscr{C}^2)$, iff $f$ is a local m-morphism of $E^1$ to $E^2$.

It is a semi-strict classic $\mathscr{E}^1$-$\mathscr{E}^2$-$E^1$-$E^2$-$\mathscr{C}^1$-$\mathscr{C}^2$ m-atlas morphism of $A^1$ to $A^2$ in the coordinate spaces $C^1, C^2$, abbreviated as
$\text{isAtl}_{\text{Ar}\atop\text{semi-strict}}^{\text{classic}}(A^1, E^1, C^1, A^2, E^2, C^2, f, \mathscr{E}^1, \mathscr{E}^2, \mathscr{C}^1, \mathscr{C}^2)$, iff $f$ is a local $\mathscr{E}^1$-$\mathscr{E}^2$ m-morphism of $E^1$ to $E^2$.

**Definition 13.5** (Strict classic m-atlas morphisms). Let $\mathscr{E}^i, \mathscr{C}^i, i = 1, 2$, be model categories, $E^i \stackrel{\text{Ob}}{\in} \mathscr{E}^i$, $C^i \stackrel{\text{Ob}}{\in} \mathscr{C}^i$, $A^i$ be an m-atlas of $E^i$ in the coordinate space $C^i$.

Let $f\colon E^1 \longrightarrow E^2$ be a classic $E^1$-$E^2$ m-atlas morphism of $A^1$ to $A^2$ in the coordinate spaces $C^1, C^2$.

$f$ is a strict classic $E^1$-$E^2$ m-atlas morphism of $A^1$ to $A^2$ in the coordinate spaces $C^1, C^2$, abbreviated as $\text{isAtl}_{\text{Ar}\atop\text{strict}}^{\text{classic}}(A^1, E^1, C^1, A^2, E^2, C^2, f)$, iff $f$ is an m-morphism of $E^1$ to $E^2$.

It is a strict classic $\mathscr{E}^1$-$\mathscr{E}^2$-$E^1$-$E^2$ m-atlas morphism of $A^1$ to $A^2$ in the coordinate spaces $C^1, C^2$, abbreviated as
$\text{isAtl}_{\text{Ar}\atop\text{strict}}^{\text{classic}}(A^1, E^1, C^1, A^2, E^2, C^2, f, \mathscr{E}^1, \mathscr{E}^2)$, iff $f$ is an $\mathscr{E}^1$-$\mathscr{E}^2$ m-morphism of $E^1$ to $E^2$.

Let $f\colon E^1 \longrightarrow E^2$ be a classic $E^1$-$E^2$-$\mathscr{C}^1$-$\mathscr{C}^2$ m-atlas morphism of $A^1$ to $A^2$ in the coordinate spaces $C^1, C^2$.

$f$ is a strict classic $E^1$-$E^2$-$\mathscr{C}^1$-$\mathscr{C}^2$ m-atlas morphism of $A^1$ to $A^2$ in the coordinate spaces $C^1, C^2$, abbreviated as $\text{isAtl}_{\text{Ar}\atop\text{strict}}^{\text{classic}}(A^1, E^1, C^1, A^2, E^2, C^2, f, \mathscr{C}^1, \mathscr{C}^2)$, iff $f$ is an m-morphism of $E^1$ to $E^2$.

It is a strict classic $\mathscr{E}^1$-$\mathscr{E}^2$-$E^1$-$E^2$-$\mathscr{C}^1$-$\mathscr{C}^2$ m-atlas morphism of $A^1$ to $A^2$ in the coordinate spaces $C^1, C^2$, abbreviated as
$\text{isAtl}_{\text{Ar}\atop\text{strict}}^{\text{classic}}(A^1, E^1, C^1, A^2, E^2, C^2, f, \mathscr{E}^1, \mathscr{E}^2, \mathscr{C}^1, \mathscr{C}^2)$, iff $f$ is an $\mathscr{E}^1$-$\mathscr{E}^2$ m-morphism of $E^1$ to $E^2$.

### 13.1.3 Abbreviated nomenclature for classic m-atlas morphisms

**Definition 13.6** (Abbreviated nomenclature for classic m-atlas morphisms). Let $\mathscr{E}^i, \mathscr{C}^i, i = 1, 2$, be model categories, $E^i \stackrel{\text{Ob}}{\in} \mathscr{E}^i$, $C^i \stackrel{\text{Ob}}{\in} \mathscr{C}^i$, $A^i$ be an m-atlas of $E^i$ in the coordinate space $C^i$ and $f\colon E^1 \longrightarrow E^2$ be a classic $E^1$-$E^2$ m-atlas morphism of $A^1$ to $A^2$ in the coordinate spaces $C^1$, $C^2$.

$f$ is a (full) (maximal) (constrained) (semi-strict, strict) classic $E^1$-$E^2$ m-atlas morphism of $A^1$ to $A^2$ in the coordinate space $C^1$, abbreviated as
$$\text{isAtl}_{\text{Ar}\,(\text{full,S-max,max,S-max-full,max-full})\,(\text{semi-strict,strict})}^{(\text{constrained,})\text{classic}} (A^1, E^1, C^1, A^2, E^2, f),$$
iff it is a (full) (maximal) (constrained) (semi-strict, strict) classic $E^1$-$E^2$ m-atlas morphism of $A^1$ to $A^2$ in the coordinate spaces $C^1$, $C^1$.

$f$ is a (full) (maximal) (constrained) (semi-strict, strict) classic $\mathscr{E}^1$-$\mathscr{E}^2$-$E^1$-$E^2$ m-atlas morphism of $A^1$ to $A^2$ in the coordinate space $C^1$, abbreviated as
$$\text{isAtl}_{\text{Ar}\,(\text{full,S-max,max,S-max-full,max-full})\,(\text{semi-strict,strict})}^{(\text{constrained,})\text{classic}} (A^1, E^1, C^1, A^2, E^2, f, \mathscr{E}^1, \mathscr{E}^2),$$
iff it is a (full) (maximal) (constrained) (semi-strict, strict) classic $\mathscr{E}^1$-$\mathscr{E}^2$-$E^1$-$E^2$ m-atlas morphism of $A^1$ to $A^2$ in the coordinate spaces $C^1$, $C^1$.

Let $\mathscr{C}^i$, $i = 1,2$, be a model category, $C^i \stackrel{\text{Ob}}{\in} \mathscr{C}^i$, $E^i$ be a topological space, $A^i$ be an m-atlas of $E^i$ in the coordinate space $C^i$ and $f\colon E^1 \longrightarrow E^2$ be a classic $E^1$-$E^2$ m-atlas morphism of $A^1$ to $A^2$ in the coordinate spaces $C^1$, $C^2$.

$f$ is a (full) (maximal) (constrained) (semi-strict, strict) classic $E^1$-$E^2$ m-atlas morphism of $A^1$ to $A^2$ in the coordinate spaces $C^1$, $C^2$, abbreviated as
$$\text{isAtl}_{\text{Ar}\,(\text{full,S-max,max,S-max-full,max-full})\,(\text{semi-strict,strict})}^{(\text{constrained,})\text{classic}} (A^1, E^1, C^1, A^2, E^2, C^2, f),$$
iff it is a (full) (maximal) (constrained) (semi-strict, strict) classic $E^1_{\text{triv}}$-$E^2_{\text{triv}}$ m-atlas morphism of $A^1$ to $A^2$ in the coordinate spaces $C^1$, $C^2$.

$f$ is a (full) (maximal) (constrained) (semi-strict, strict) classic $E^1$-$E^2$ m-atlas morphism of $A^1$ to $A^2$ in the coordinate space $C^1$, abbreviated as
$$\text{isAtl}_{\text{Ar}\,(\text{full,S-max,max,S-max-full,max-full})\,(\text{semi-strict,strict})}^{(\text{constrained,})\text{classic}} (A^1, E^1, C^1, A^2, E^2, f),$$
iff it is a (full) (maximal) (constrained) (semi-strict, strict) classic $E^1$-$E^2$ m-atlas morphism of $A^1$ to $A^2$ in the coordinate spaces $C^1$, $C^1$.

### 13.1.4 Proclamations on classic m-atlas (near) morphisms

**Lemma 13.7** (Classic m-atlas (near) morphisms). *Let $E$ be a set of topological spaces, $E^i \in E$, $i = 1, 2$, $C^i$ be a model space, $A^i$ be an m-atlas of $E^i_{\text{triv}}$ in the coordinate space $C^i$.*

  1. *Let $f\colon E^1_{\text{triv}} \longrightarrow E^2_{\text{triv}}$ be a classic $E^1_{\text{triv}}$-$E^2_{\text{triv}}$ m-atlas morphism of $A^1$ to $A^2$ in the coordinate spaces $C^1$, $C^2$ and $C^1 \stackrel{\text{mod}}{\subseteq} C^2$. Then $f$ is a strict constrained classic $E_{\mathcal{T}riv}$-$E_{\mathcal{T}riv}$-$E^1_{\text{triv}}$-$E^2_{\text{triv}}$ m-atlas morphism of $A^1$ to $A^2$ in the coordinate spaces $C^1$, $C^2$.*

*Proof.*

**$f$ is constrained:** $f[E^1] \subseteq E^2$. $E^2$ is a model neighborhood of $E^2_{\mathbf{triv}}$ by definition 4.1 (Trivial model spaces and categories) on page 22.

**$f$ is strict:** $\boldsymbol{E}_{\mathcal{T}riv} \overset{\text{fullcat}}{\subseteq} \boldsymbol{E}_{\mathcal{T}riv}$. $f$ is continuous by definition 7.1 (Model functions) on page 28. $E^i \in \boldsymbol{E}$ by hypothesis, $E^i \overset{\text{Ob}}{\in} \boldsymbol{E}_{\mathcal{T}riv}$ and $f \overset{\text{Ar}}{\in} \boldsymbol{E}_{\mathcal{T}riv}$ by definition 4.1. $\boldsymbol{C}^1 \overset{\text{mod}}{\subseteq} \boldsymbol{C}^2$ by definition 13.1 (Classic m-atlas morphisms for model spaces) on page 98.

□

2. *If $E^1$ is an open subspace of $E^2$ then $f \colon E^1_{\mathbf{triv}} \longrightarrow E^2_{\mathbf{triv}}$ is also a strict constrained classic $E^1_{\mathbf{triv}}$-$E^2_{\mathbf{triv}}$ m-atlas morphism of $\boldsymbol{A}^1$ to $\boldsymbol{A}^2$ in the coordinate spaces $\boldsymbol{C}^1$, $\boldsymbol{C}^2$.*

*Proof.* $E^1_{\mathbf{triv}} \overset{\text{strictmod}}{\subseteq} E^2_{\mathbf{triv}}$ by lemma 4.2 (Trivial model spaces) on page 23. $E^i \overset{\text{Ob}}{\in} \boldsymbol{E}_{\mathcal{T}riv}$ and $f \overset{\text{Ar}}{\in} \boldsymbol{E}_{\mathcal{T}riv}$ by definition 4.1. □

Let

1. $\mathscr{E}^i, \mathscr{C}^i, \mathscr{E}'^i, \mathscr{C}'^i$, $i = 1, 2$, *be model categories*
2. $\mathscr{E}^i \overset{\text{full-cat}}{\subseteq} \mathscr{E}'^i$
3. $\mathscr{C}^i \overset{\text{full-cat}}{\subseteq} \mathscr{C}'^i$
4. $\mathscr{C}'^1 \overset{\text{full-cat}}{\subseteq} \mathscr{C}'^2$
5. $\boldsymbol{E}^i \overset{\text{Ob}}{\in} \mathscr{E}^i$
6. $\boldsymbol{C}^i \overset{\text{Ob}}{\in} \mathscr{C}^i$
7. $\boldsymbol{C}'^i \overset{\text{Ob}}{\in} \mathscr{C}'^i$
8. $\boldsymbol{C}^i \overset{\text{mod}}{\subseteq} \boldsymbol{C}'^i$
9. $\boldsymbol{A}^i$ *be an m-atlas of $\boldsymbol{E}^i$ in the coordinate space $\boldsymbol{C}^i$*
10. $f \colon \boldsymbol{E}^1 \longrightarrow \boldsymbol{E}^2$ *be a (semi-strict, strict) ($\mathscr{E}^1$-$\mathscr{E}^2$-)$\boldsymbol{E}^1$-$\boldsymbol{E}^2$(-$\mathscr{C}^1$-$\mathscr{C}^2$) m-atlas (neae) morphism of $\boldsymbol{A}^1$ to $\boldsymbol{A}^2$ in the coordinate spaces $\boldsymbol{C}^1$, $\boldsymbol{C}^2$ (model function?)*

*Then*

1. If $f$ is a classic $\boldsymbol{E}^1$-$\boldsymbol{E}^2$ m-atlas morphism of $\boldsymbol{A}^1$ to $\boldsymbol{A}^2$ in the coordinate spaces $\boldsymbol{C}^1$, $\boldsymbol{C}^2$ then $f$ is a classic $\boldsymbol{E}^1$-$\boldsymbol{E}^2$ m-atlas near morphism of $\boldsymbol{A}^1$ to $\boldsymbol{A}^2$ in the coordinate spaces $\boldsymbol{C}^1$, $\boldsymbol{C}^2$.

   *Proof.* Let $(U^i, V^i, \phi^i \colon U^i \rightarrowtail\!\!\!\xrightarrow{\tilde{=}} V^i) \in \boldsymbol{A}^i$, $i = 1, 2$, be charts with $I \stackrel{\text{def}}{=} U^1 \cap f_0^{-1}[U^2] \neq \emptyset$. Then $\phi^2 \circ f \circ \phi^{1-1} \colon \phi^1[I] \longrightarrow V^2$ is a $\boldsymbol{C}^1$-$\boldsymbol{C}^2$ m-morphism of $\phi^1[I]$ to $V^2$ by definition 13.1 (Classic m-atlas morphisms for model spaces) on page 98 and thus a local $\boldsymbol{C}^1$-$\boldsymbol{C}^2$ m-morphism of $\phi^1[I]$ to $V^2$ by item 1 of lemma 7.15 ((Local) m-morphisms) on page 32. . □

2. If $f$ is a classic $\boldsymbol{E}^1$-$\boldsymbol{E}^2$-$\mathscr{C}^1$-$\mathscr{C}^2$ m-atlas morphism of $\boldsymbol{A}^1$ to $\boldsymbol{A}^2$ in the coordinate spaces $\boldsymbol{C}^1$, $\boldsymbol{C}^2$ then $f$ is a classic $\boldsymbol{E}^1$-$\boldsymbol{E}^2$-$\mathscr{C}^1$-$\mathscr{C}^2$ m-atlas near morphism of $\boldsymbol{A}^1$ to $\boldsymbol{A}^2$ in the coordinate spaces $\boldsymbol{C}^1$, $\boldsymbol{C}^2$.

   *Proof.* Let $(U^i, V^i, \phi^i \colon U^i \rightarrowtail\!\!\!\xrightarrow{\tilde{=}} V^i) \in \boldsymbol{A}^i$, $i = 1, 2$, be charts with $I \stackrel{\text{def}}{=} U^1 \cap f_0^{-1}[U^2] \neq \emptyset$. Then $\phi^2 \circ f \circ \phi^{1-1} \colon \phi^1[I] \longrightarrow V^2$ is a $\boldsymbol{C}^1$-$\boldsymbol{C}^2$-$\mathscr{C}^1$-$\mathscr{C}^2$ m-morphism of $\phi^1[I]$ to $V^2$ by definition 13.1 (Classic m-atlas morphisms for model spaces) on page 98 and thus a local $\boldsymbol{C}^1$-$\boldsymbol{C}^2$-$\mathscr{C}^1$-$\mathscr{C}^2$ m-morphism of $\phi^1[I]$ to $V^2$ by item 3 of lemma 7.15. . □

3. If $f$ is a (semi-strict, strict) $\mathscr{C}^1$-$\mathscr{C}^2$-$\boldsymbol{E}^1$-$\boldsymbol{E}^2$ ($\boldsymbol{E}^1$-$\boldsymbol{E}^2$) classic m-atlas morphism of $\boldsymbol{A}^1$ to $\boldsymbol{A}^2$ in the coordinate spaces $\boldsymbol{C}^1$, $\boldsymbol{C}^2$, then $f$ is a (semi-strict, strict) $\mathscr{C}'^1$-$\mathscr{C}'^2$-$\boldsymbol{E}^1$-$\boldsymbol{E}^2$ ($\boldsymbol{E}^1$-$\boldsymbol{E}^2$) classic m-atlas morphism of $\boldsymbol{A}^1$ to $\boldsymbol{A}^2$ in the coordinate spaces $\boldsymbol{C}'^1$, $\boldsymbol{C}'^2$.

   *Proof.*

   **$f$ remains a classic m-atlas morphism with the expanded categories.**
   $\mathscr{C}'^1 \stackrel{\text{full-cat}}{\subseteq} \mathscr{C}'^2$ by hypothesis.

   Let $(U^i, V^i, \phi^i \colon U^i \rightarrowtail\!\!\!\xrightarrow{\tilde{=}} V^i) \in \boldsymbol{A}^i$, $i = 1, 2$, with $I \stackrel{\text{def}}{=} U^1 \cap f^{-1}[U^2] \neq \emptyset$.

   $\phi^2 \circ f \circ \phi^{-1}$ is a morphism of $\mathscr{C}^2 \stackrel{\text{full-cat}}{\subseteq} \mathscr{C}'^2$ by definition 13.1.

   **$f$ remains semi-strict (strict).** If $f$ is a local morphism of $\mathscr{C}^2$ then $f$ is a local morphism of $\mathscr{C}'^2$. If $f$ is a local morphism of $\boldsymbol{E}^2$ then $f$ is a local morphism of $\boldsymbol{E}'^2$. If $f \stackrel{\text{Ar}}{\in} \mathscr{C}^2$ then $f \stackrel{\text{Ar}}{\in} \mathscr{C}'^2$. If $f$ is a morphism of $\boldsymbol{E}^2$ then $f$ is a morphism of $\boldsymbol{E}'^2$.

   □

4. If $f$ is a classic $\boldsymbol{E}^1$-$\boldsymbol{E}^2$ m-atlas near morphism of $\boldsymbol{A}^1$ to $\boldsymbol{A}^2$ in the coordinate spaces $\boldsymbol{C}^1$, $\boldsymbol{C}^2$ and $\boldsymbol{C}^2$ is semi-maximal, then $f$ is a classic $\boldsymbol{E}^1$-$\boldsymbol{E}^2$ m-atlas morphism of $\boldsymbol{A}^1$ to $\boldsymbol{A}^2$ in the coordinate spaces $\boldsymbol{C}^1$, $\boldsymbol{C}^2$.

*Proof.* Let $(U^i, V^i, \phi^i\colon U^i \rightarrowtail\overset{\sim}{\twoheadrightarrow} V^i) \in \boldsymbol{A}^i$, $i = 1, 2$, be charts with $I \overset{\text{def}}{=} U^1 \cap f_0^{-1}[U^2] \neq \emptyset$. Then $\phi^2 \circ f \circ \phi^{1-1}\colon \phi^1[I] \longrightarrow V^2$ is a local $\boldsymbol{C}^1$-$\boldsymbol{C}^2$ m-morphism of $\phi^1[I]$ to $V^2$ by definition 13.1 (Classic m-atlas morphisms for model spaces) on page 98. □

5. *If $\mathscr{C}^2$ satisfies the restricted sheaf condition, $f$ is a classic $\boldsymbol{E}^1$-$\boldsymbol{E}^2$-$\mathscr{C}^1$-$\mathscr{C}^2$ m-atlas near morphism of $\boldsymbol{A}^1$ to $\boldsymbol{A}^2$ in the coordinate spaces $\boldsymbol{C}^1$, $\boldsymbol{C}^2$ and $\boldsymbol{A}^2$ is semi-maximal, then $f$ is a classic $\boldsymbol{E}^1$-$\boldsymbol{E}^2$-$\mathscr{C}^1$-$\mathscr{C}^2$ m-atlas morphism of $\boldsymbol{A}^1$ to $\boldsymbol{A}^2$ in the coordinate spaces $\boldsymbol{C}^1$, $\boldsymbol{C}^2$.*

   *Proof.* Let $(U^i, V^i, \phi^i\colon U^i \rightarrowtail\overset{\sim}{\twoheadrightarrow} V^i) \in \boldsymbol{A}^i$, $i = 1, 2$, be charts with $I \overset{\text{def}}{=} U^1 \cap f_0^{-1}[U^2] \neq \emptyset$. Then $\phi^2 \circ f \circ \phi^{1-1}\colon \phi^1[I] \longrightarrow V^2$ is a local $\boldsymbol{C}^1$-$\boldsymbol{C}^2$-$\mathscr{C}^1$-$\mathscr{C}^2$ m-morphism of $\phi^1[I]$ to $V^2$ by definition 13.1 (Classic m-atlas morphisms for model spaces) on page 98 □

**Lemma 13.8** (Composition of classic m-atlas morphisms). *Let $\mathscr{E}^i$, $\mathscr{C}^i$, $i \in [1, 3]$, be a model category, $\boldsymbol{E}^i \overset{\text{Ob}}{\in} \mathscr{E}^i$, $\boldsymbol{C}^i \overset{\text{Ob}}{\in} \mathscr{C}^i$, $\boldsymbol{A}^i$ be an m-atlas of $\boldsymbol{E}^i$ in the coordinate space $\boldsymbol{C}^i$ and $f^i\colon \boldsymbol{E}^i \longrightarrow \boldsymbol{E}^{i+1}$, $i = 1, 2$, be a (semi-strict, strict) classic $\boldsymbol{E}^i$-$\boldsymbol{E}^{i+1}$ m-atlas morphism of $\boldsymbol{A}^i$ to $\boldsymbol{A}^{i+1}$ in the coordinate spaces $\boldsymbol{C}^i$, $\boldsymbol{C}^{i+1}$.*

*$f^2 \circ f^1\colon \boldsymbol{E}^1 \longrightarrow \boldsymbol{E}^3$ is a (semi-strict, strict) classic $\boldsymbol{E}^1$-$\boldsymbol{E}^3$ m-atlas morphism of $\boldsymbol{A}^1$ to $\boldsymbol{A}^3$ in the coordinate spaces $\boldsymbol{C}^1$, $\boldsymbol{C}^3$.*

*Proof.*

**$f^2 \circ f^1$ is a classic $\boldsymbol{E}^1$-$\boldsymbol{E}^3$ m-atlas morphism.** If $\mathscr{C}^1 \overset{\text{cat}}{\subseteq} \mathscr{C}^2$ and $\mathscr{C}^2 \overset{\text{cat}}{\subseteq} \mathscr{C}^3$ then $\mathscr{C}^1 \overset{\text{cat}}{\subseteq} \mathscr{C}^3$. Let $(U^1, V^1, \phi^1\colon U^1 \rightarrowtail\overset{\sim}{\twoheadrightarrow} V^1) \in \boldsymbol{A}^1$, $(U^3, V^3, \phi^3\colon U^3 \rightarrowtail\overset{\sim}{\twoheadrightarrow} V^3) \in \boldsymbol{A}^3$ with $I \overset{\text{def}}{=} U^1 \cap (f^2 \circ f^1)^{-1}[U^3] \neq \emptyset$, $(f^2 \circ f^1)^{-1}[U^3]$ is a model neighborhood by definition 7.1 (Model functions) on page 28 and $I$ is a model neighborhood by item 2 of definition 2.1 (Model spaces) on page 18.

For any $u^1 \in I$ and any chart $(U^2, V^2, \phi^2\colon U^2 \rightarrowtail\overset{\sim}{\twoheadrightarrow} V^2) \in \boldsymbol{A}^2$ at $f^1(u^1)$, define $I^{1,2} \overset{\text{def}}{=} U^1 \cap f^{1-1}[U^2]$, $I^{2,3} \overset{\text{def}}{=} U^2 \cap f^{2-1}[U^3]$ and $I^{1,3}_{u^1} \overset{\text{def}}{=} U^1 \cap f^{1-1}[I^{2,3}] \subseteq I$. $I^{1,3}_{u^1}$ is a model neighborhood by definition 7.1 and item 2 of definition 2.1.

$\phi^2 \circ f^1 \circ \phi^{1-1}\colon \phi^1[I^{1,2}] \longrightarrow V^2$ is a morphism of $\mathscr{C}^2 \overset{\text{cat}}{\subseteq} \mathscr{C}^3$ by hypothesis and $\phi^2 \circ f^1 \circ \phi^{1-1}\colon \phi^1[I^{1,3}_{u^1}] \longrightarrow V^2$ is a morphism of $\mathscr{C}^2 \overset{\text{cat}}{\subseteq} \mathscr{C}^3$ by item 4. of definition 2.1. $\phi^3 \circ f^2 \circ \phi^{2-1}\colon \phi^2[I^{2,3}] \longrightarrow V^3$ is a morphism of $\mathscr{C}^3$ by hypothesis and and thus $(\phi^3 \circ f^2 \circ f^1 \circ \phi^{1-1}\colon \phi^1[I^{1,3}_{u^1}] \longrightarrow V^3) = (\phi^3 \circ f^2 \circ \phi^{2-1}\colon \phi^2[I^{2,3}] \longrightarrow V^3) \circ (\phi^2 \circ f^1 \circ \phi^{1-1}\colon \phi^1[I^{1,3}_{u^1}] \longrightarrow V^2)$ is a morphism of $\mathscr{C}^3$.

$(\phi^3 \circ f^2 \circ f^1 \circ \phi^{1-1} \colon \phi^1[I^{1,3}_{u^1}] \longrightarrow V^3)$ agrees with $(\phi^3 \circ f^2 \circ f^1 \circ \phi^{1-1} \colon \phi^1[I] \longrightarrow V^3)$ on $\phi^1[I^{1,3}_{u^1}]$ and $\phi^1[I] = \bigcup_{u^1 \in I} \phi^1[I^{1,3}_{u^1}]$, thus $(\phi^3 \circ f^2 \circ f^1 \circ \phi^{1-1} \colon \phi^1[I] \longrightarrow V^3)$ is a morphism of $\mathscr{C}^3$ by item 6 of definition 2.1.

$f^2 \circ f^1$ **is (semi-strict, strict).** If $E^1 \overset{\text{mod}}{\subseteq} E^2$ and $E^2 \overset{\text{mod}}{\subseteq} E^3$ then $E^1 \overset{\text{mod}}{\subseteq} E^3$. If $C^1 \overset{\text{mod}}{\subseteq} C^2$ and $C^2 \overset{\text{mod}}{\subseteq} C^3$ then $C^1 \overset{\text{mod}}{\subseteq} C^3$. If $\mathscr{C}^1 \overset{\text{cat}}{\subseteq} \mathscr{C}^2$ and $\mathscr{C}^2 \overset{\text{cat}}{\subseteq} \mathscr{C}^3$ then $\mathscr{C}^1 \overset{\text{cat}}{\subseteq} \mathscr{C}^3$.

If $f^1$ is a local morphism of $E^2$ and $f^2$ is a local morphism of $E^3$ then $f^2 \circ f^1$ is a local morphism of $E^3$ by lemma 7.17 (Composition of local m-morphisms) on page 39.

If $f^1$ is a local $\mathscr{C}^1$-$\mathscr{C}^2$ morphism of $E^2$ and $f^2$ is a local $\mathscr{C}^2$-$\mathscr{C}^3$ morphism of $E^3$ then $f^2 \circ f^1$ is a local $\mathscr{C}^1$-$\mathscr{C}^3$ morphism of $E^3$ by lemma 7.17.

If $f^1$ is a morphism of $E^2 E^3$ and $f^2$ is a morphism of $E^3$ then $f^2 \circ f^1$ is a morphism of $E^3$.

□

*Similar results apply with restrictions on the admissible atlases.*

**Lemma 13.9** (Classic m-atlas identity)**.** *Let $\mathscr{C}^i$, $\mathscr{C}^i$, $i = 1, 2$, be a model category, $\mathscr{C}^1 \overset{\text{full-cat}}{\subseteq} \mathscr{C}^2$, $\mathscr{C}^1 \overset{\text{full-cat}}{\subseteq} \mathscr{C}^2$. $E^i \overset{\text{Ob}}{\in} \mathscr{C}^i$, $E^1 \subseteq E^2$, $C^i \overset{\text{Ob}}{\in} \mathscr{C}^i$, $C^1 \subseteq C^2$, $A^i$ a classic m-atlas of $E^i$ in the coordinate space $C^i$ and $A^1 \subseteq A^2$.*

*The identity morphism $f \overset{\text{def}}{=} \mathrm{ID}_{(E^1, C^1), (E^2, C^2)}$ is a strict classic $\mathscr{C}^1$-$\mathscr{C}^2$-$E^1$-$E^2$-$\mathscr{C}^1$-$\mathscr{C}^2$ m-atlas morphism of $A^1$ to $A^2$ in the coordinate spaces $C^1$, $C^2$.*

*Proof.* Let $u^1 \in E^1$ and $(U^1, V^1, \phi^1 \colon U^1 \rightarrowtail\mathrel{\mkern-14mu}\xrightarrow{\sim} V^1) \in A^1 \subseteq A^2$ be a chart at $u^1$. Define $U'^1 \overset{\text{def}}{=} U^1$, $V'^1 \overset{\text{def}}{=} V^1$, $U'^2 \overset{\text{def}}{=} U^1$ and $V'^2 \overset{\text{def}}{=} V^1$. Diagram (11.8) in definition 11.4 is M-nearly commutative in $C^2$,

$\mathrm{Id}_{E^i} \overset{\text{Ar}}{\in} \mathscr{C}^i \overset{\text{full-cat}}{\subseteq} \mathscr{C}^{i+1}$ and $\mathrm{Id}_{C^i} \overset{\text{Ar}}{\in} \mathscr{C}^i \overset{\text{full-cat}}{\subseteq} \mathscr{C}^{i+1}$. □

**Corollary 13.10** (Classic m-atlas identity)**.** *Let $\mathscr{C}$, $\mathscr{C}$, be model categories, $E \overset{\text{Ob}}{\in} \mathscr{C}$, $C \overset{\text{Ob}}{\in} \mathscr{C}$, $A$ an m-atlas of $E$ in the coordinate spaces $C$.*

*The identity morphism $f \overset{\text{def}}{=} \mathrm{Id}_{(A, E, C)}$ is a strict m-atlas morphism of $A$ to $A$ in the coordinate spaces $C$, $C$.*

## 13.2 Categories of classic m-atlases and functors

This subsection defines categories of m-atlases with classic m-atlas morphisms and functors among them, and proves some basic results.

### 13.2.1 Classic m-atlas categories

**Definition 13.11** (Sets of classic m-atlas morphisms)**.** Let $E^i$ and $C^i$, $i = 1, 2$, be model spaces. Then

$$\mathcal{A}t\ell_{\text{Ar}\,(\text{full,S-max,max,S-max-full,max-full})}^{\text{classic}}(E^1, C^1, E^2, C^2) \stackrel{\text{def}}{=}$$
$$_{(\text{semi-strict,strict})}$$

$$\left\{ ((f), (A^1, E^1, C^1), (A^2, E^2, C^2)) \;\middle|\; \text{isAtl}_{\text{Ar}\,(\text{full,S-max,max,S-max-full,max-full})}^{\text{classic}}(A^1, E^1, C^1, A^2, E^2, C^2, f) \right\} \quad (13.2)$$
$$_{(\text{semi-strict,strict})}$$

**Definition 13.12** (Categories $\mathcal{A}t\ell^{\text{Classic}}(E, C)$). Let $E$ and $C$ be sets of model spaces. Let $P \stackrel{\text{def}}{=} E \times C$. Then

$$\mathcal{A}t\ell_{\text{Ar}\,(\text{full,S-max,max,S-max-full,max-full})}^{\text{classic}}(E, C) \stackrel{\text{def}}{=} \bigcup_{\substack{(E^\mu, C^\mu) \in P \\ (E^\nu, C^\nu) \in P}} \mathcal{A}t\ell_{\text{Ar}\,(\text{full,S-max,max,S-max-full,max-full})}^{\text{classic}}(E^\mu, C^\mu, E^\nu, C^\nu)$$
$$_{(\text{semi-strict,strict})} \qquad\qquad _{(\text{semi-strict,strict})}$$
$$(13.3)$$

$$\mathcal{A}t\ell^{\text{classic}}_{(\text{full,S-max,max,S-max-full,max-full})}(E, C) \stackrel{\text{def}}{=} \left( \mathcal{A}t\ell_{\text{Ob}\,(\text{full,S-max,max,S-max-full,max-full})}(E, C),\; \mathcal{A}t\ell_{\text{Ar}\,(\text{full,S-max,max,S-max-full,max-full})}(E, C),\; \stackrel{A}{\circ} \right),$$
$$_{(\text{semi-strict,strict})}$$
$$(13.4)$$

Let $E^i$ and $C^i$, $i = 1, 2$, be model spaces, $A^i$ be an atlas of $E^i$ in the coordinate space $C^i$ and $f\colon E^1 \longrightarrow E^2$ be a (strict, semi-strict) $E^1$-$E^2$ classic M-atlas morphism of $A^1$ to $A^2$ in the coordinate spaces $C^1$, $C^2$.

The identity functor $\mathcal{F}_{\text{Classic,Id}}$ is

$$\mathcal{F}_{\text{Classic,Id}}(A^i, E^i, C^i) \stackrel{\text{def}}{=} (A^i, E^i, C^i) \quad (13.5)$$

$$\mathcal{F}_{\text{Classic,Id}}\big((f\colon E^1 \longrightarrow E^2), (A^1, E^1, C^1), (A^2, E^2, C^2)\big) \stackrel{\text{def}}{=}$$
$$\big((f\colon E^1 \longrightarrow E^2), (A^1, E^1, C^1), (A^2, E^2, C^2)\big) \quad (13.6)$$

This nomenclature will be justified below.

**Lemma 13.13** ($\mathcal{A}t\ell^{\text{Classic}}(E, C)$ is a category)**.** *Let $E$ and $C$ be sets of model spaces. Then $\mathcal{A}t\ell^{\text{Classic}}(E, C)$ is a category and the identity morphism for an object $(A^i, E^i, C^i)$ of $\mathcal{A}t\ell_{\text{Ob}}^{\text{Classic}}(E, C)$ is $\text{Id}_{(A^i, E^i, C^i)}$.*

*Proof.* Let $(A^i, E^i, C^i)$, $i \in [1, 3]$, be an object of $\mathcal{A}t\ell^{\text{Classic}}(E, C)$ and let $m^i \stackrel{\text{def}}{=} ((f^i), (A^i, E^i, C^i), (A^{i+1}, E^{i+1}, C^{i+1}))$ be a morphism of $\mathcal{A}t\ell^{\text{Classic}}(E, C)$. Then

**Composition:** $f^2 \circ f^1$ is a classic $\boldsymbol{E}^1$-$\boldsymbol{E}^3$ m-atlas morphism of $\boldsymbol{A}^1$ to $\boldsymbol{A}^3$ in the coordinate spaces $\boldsymbol{C}^1$, $\boldsymbol{C}^3$ by lemma 13.8 (Composition of classic m-atlas morphisms) on page 104 and $(f^2) \overset{()}{\circ} (f^1) = (f^2 \circ f^1)$.

**Associativity:** Composition is associative by lemma 1.19 (Tuple composition for labeled morphisms) on page 12.

**Identity:** $\mathrm{Id}_{(\boldsymbol{A}^i, \boldsymbol{E}^i, \boldsymbol{C}^i)}$ is an identity morphism by lemma 1.19.

$\square$

**Corollary 13.14** (Subcategories of $\mathcal{A}t\ell^{\mathrm{Classic}}(\boldsymbol{E}, \boldsymbol{C})$)**.** *Let $\boldsymbol{E}$ and $\boldsymbol{C}$ be sets of model spaces. Then all of the following are subcategories of $\mathcal{A}t\ell^{\mathrm{Classic}}(\boldsymbol{E}, \boldsymbol{C})$.*

1. $\mathcal{A}t\ell^{\mathrm{Classic}}_{\mathrm{full}}(\boldsymbol{E}, \boldsymbol{C})$
2. $\mathcal{A}t\ell^{\mathrm{Classic}}_{\mathrm{S\text{-}max}}(\boldsymbol{E}, \boldsymbol{C})$
3. $\mathcal{A}t\ell^{\mathrm{Classic}}_{\mathrm{max}}(\boldsymbol{E}, \boldsymbol{C})$
4. $\mathcal{A}t\ell^{\mathrm{Classic}}_{\mathrm{S\text{-}max\text{-}full}}(\boldsymbol{E}, \boldsymbol{C})$
5. $\mathcal{A}t\ell^{\mathrm{Classic}}_{\mathrm{max\text{-}full}}(\boldsymbol{E}, \boldsymbol{C})$
6. $\mathcal{A}t\ell^{\mathrm{Classic}}_{\mathrm{semi\text{-}strict}}(\boldsymbol{E}, \boldsymbol{C})$
7. $\mathcal{A}t\ell^{\mathrm{Classic}}_{\substack{\mathrm{full} \\ \mathrm{semi\text{-}strict}}}(\boldsymbol{E}, \boldsymbol{C})$
8. $\mathcal{A}t\ell^{\mathrm{Classic}}_{\substack{\mathrm{S\text{-}max} \\ \mathrm{semi\text{-}strict}}}(\boldsymbol{E}, \boldsymbol{C})$
9. $\mathcal{A}t\ell^{\mathrm{Classic}}_{\substack{\mathrm{max} \\ \mathrm{semi\text{-}strict}}}(\boldsymbol{E}, \boldsymbol{C})$
10. $\mathcal{A}t\ell^{\mathrm{Classic}}_{\substack{\mathrm{S\text{-}max\text{-}full} \\ \mathrm{semi\text{-}strict}}}](\boldsymbol{E}, \boldsymbol{C})$
11. $\mathcal{A}t\ell^{\mathrm{Classic}}_{\substack{\mathrm{max\text{-}full} \\ \mathrm{semi\text{-}strict}}}(\boldsymbol{E}, \boldsymbol{C})$
12. $\mathcal{A}t\ell^{\mathrm{Classic}}_{\mathrm{strict}}(\boldsymbol{E}, \boldsymbol{C})$
13. $\mathcal{A}t\ell^{\mathrm{Classic}}_{\substack{\mathrm{full} \\ \mathrm{strict}}}(\boldsymbol{E}, \boldsymbol{C})$

14. $\mathcal{A}t\ell^{\text{Classic}}_{\substack{\text{S-max}\\\text{strict}}}(E, C)$

15. $\mathcal{A}t\ell^{\text{Classic}}_{\substack{\text{max}\\\text{strict}}}(E, C)$

16. $\mathcal{A}t\ell^{\text{Classic}}_{\substack{\text{S-max-full}\\\text{strict}}}(E, C)$

17. $\mathcal{A}t\ell^{\text{Classic}}_{\substack{\text{max-full}\\\text{strict}}}(E, C)$

*Remark* 13.15. They are not, in general, full subcategories.

The identity functor $\mathcal{F}_{\text{Classic,Id}}$ is a functor from each of the subcategories to each of the containing categories.

**Definition 13.16** (Categories $\mathcal{A}t\ell^{\text{Classic}}(\mathcal{E}, \mathcal{C})$). Let $\mathcal{E}$ and $\mathcal{C}$ be model categories. Let $P \stackrel{\text{def}}{=} \text{Ob}(\mathcal{E}) \times \text{Ob}(\mathcal{C})$. Then

$$\mathcal{A}t\ell^{\text{Classic}}_{\text{Ob}}(\mathcal{E}, \mathcal{C}) \stackrel{\text{def}}{=} \bigcup_{\substack{E \in \mathcal{E}^{\text{Ob}} \\ C \in \mathcal{C}^{\text{Ob}}}} \mathcal{A}t\ell^{\text{Classic}}_{\text{Ob}}(E, C) \tag{13.7}$$

$$\mathcal{A}t\ell^{\text{Classic}}_{\text{Ar}}(\mathcal{E}, \mathcal{C}) \stackrel{\text{def}}{=} \left\{ ((f), (A^1, E^1, C^1), (A^2, E^2, C^2)) \,\middle|\, (E^i, C^i) \in P \land \underset{\text{strict}}{\text{isAtl}_{\text{Ar}}}(A^1, E^1, C^1, A^2, E^2, C^2, f, \mathcal{E}, \mathcal{C}) \right\} \tag{13.8}$$

$$\mathcal{A}t\ell^{\text{Classic}}(\mathcal{E}, \mathcal{C}) \stackrel{\text{def}}{=} \left( \mathcal{A}t\ell^{\text{Classic}}_{\text{Ob}}(\mathcal{E}, \mathcal{C}), \mathcal{A}t\ell^{\text{Classic}}_{\text{Ar}}(\mathcal{E}, \mathcal{C}), \stackrel{A}{\circ} \right) \tag{13.9}$$

Similar definitions apply with restrictions on the admissible atlases, the admissible morphisms, or both:

1. $\mathcal{A}t\ell^{\text{Classic}}_{\text{full}}$

2. $\mathcal{A}t\ell^{\text{Classic}}_{\text{S-max}}$

3. $\mathcal{A}t\ell^{\text{Classic}}_{\text{max}}$

4. $\mathcal{A}t\ell^{\text{Classic}}_{\text{S-max-full}}$

5. $\mathcal{A}t\ell^{\text{Classic}}_{\text{max-full}}$

6. $\mathcal{A}t\ell^{\text{Classic}}_{\substack{\text{full}\\\text{semi-strict}}}$

7. $\mathcal{A}t\ell^{\text{Classic}}_{\substack{\text{S-max}\\\text{semi-strict}}}$

8. $\mathcal{A}t\ell^{\text{Classic}}_{\substack{\text{max}\\\text{semi-strict}}}$

9. $\mathcal{A}t\ell^{\text{Classic}}_{\substack{\text{S-max-full}\\\text{semi-strict}}}$

10. $\mathcal{A}t\ell^{\text{Classic}}_{\substack{\text{max-full}\\\text{semi-strict}}}$

11. $\mathcal{A}t\ell^{\text{Classic}}_{\substack{\text{full}\\\text{strict}}}$

12. $\mathcal{A}t\ell^{\text{Classic}}_{\substack{\text{S-max}\\\text{strict}}}$

13. $\mathcal{A}t\ell^{\text{Classic}}_{\substack{\text{max}\\\text{strict}}}$

14. $\mathcal{A}t\ell^{\text{Classic}}_{\substack{\text{S-max-full}\\\text{strict}}}$

15. $\mathcal{A}t\ell^{\text{Classic}}_{\substack{\text{max-full}\\\text{strict}}}$

e.g.,

$$\mathcal{A}t\ell^{\text{Classic}}_{\substack{\text{Ob}\\\text{S-max-full}}}l(\mathcal{E},\mathcal{C}) \stackrel{\text{def}}{=} \bigcup_{(\boldsymbol{E}^\mu,\boldsymbol{C}^\mu)\in\boldsymbol{P}} \mathcal{A}t\ell^{\text{Classic}}_{\substack{\text{Ob}\\\text{S-max-full}}}(\boldsymbol{E}^\mu,\boldsymbol{C}^\mu) \tag{13.10}$$

$$\mathcal{A}t\ell^{\text{Classic}}_{\substack{\text{Ar}\\\text{S-max-full}\\\text{semi-strict}}}(\mathcal{E},\mathcal{C}) \stackrel{\text{def}}{=} \bigcup_{\substack{(\boldsymbol{E}^\mu,\boldsymbol{C}^\mu)\in\boldsymbol{P}\\(\boldsymbol{E}^\nu,\boldsymbol{C}^\nu)\in\boldsymbol{P}}} \mathcal{A}t\ell^{\text{Classic}}_{\substack{\text{Ar}\\\text{S-max-full}\\\text{semi-strict}}}(\boldsymbol{E}^\mu,\boldsymbol{C}^\mu,\boldsymbol{E}^\nu,\boldsymbol{C}^\nu) \tag{13.11}$$

$$\mathcal{A}t\ell^{\text{Classic}}_{\substack{\text{S-max-full}\\\text{semi-strict}}}(\mathcal{E},\mathcal{C}) \stackrel{\text{def}}{=} \left(\mathcal{A}t\ell^{\text{Classic}}_{\substack{\text{Ob}\\\text{S-max-full}}}l(\mathcal{E},\mathcal{C}), \mathcal{A}t\ell^{\text{Classic}}_{\substack{\text{Ar}\\\text{S-max-full}\\\text{semi-strict}}}(\boldsymbol{E},\boldsymbol{C}), \stackrel{A}{\circ}\right) \tag{13.12}$$

Let $(\boldsymbol{A}^1,\boldsymbol{E}^1,\boldsymbol{C}^1) \in \mathcal{A}t\ell^{\text{Classic}}_{\text{Ob}}(\mathcal{E},\mathcal{C})$.

$$\text{Id}_{(\boldsymbol{A}^1,\boldsymbol{E}^1,\boldsymbol{C}^1)} \stackrel{\text{def}}{=} \left((\text{Id}_{\boldsymbol{E}^1},\text{Id}_{\boldsymbol{C}^1}),(\boldsymbol{A}^1,\boldsymbol{E}^1,\boldsymbol{C}^1),(\boldsymbol{A}^1,\boldsymbol{E}^1,\boldsymbol{C}^1)\right) \tag{13.13}$$

**Lemma 13.17** ($\mathcal{A}t\ell^{\text{Classic}}(\mathcal{E}, \mathcal{C})$ is a category). *Let $\mathcal{E}$ and $\mathcal{C}$ be model categories. Then $\mathcal{A}t\ell^{\text{Classic}}(\mathcal{E}, \mathcal{C})$ is a category and the identity morphism for an object $(A^i, E^i, C^i)$ of $\mathcal{A}t\ell^{\text{Classic}}_{\text{Ob}}(\mathcal{E}, \mathcal{C})$ is $\text{Id}_{(A^i, E^i, C^i)}$.*

*Proof.* Let $(A^i, E^i, C^i), i \in [1, 3]$, be an object of $\mathcal{A}t\ell^{\text{Classic}}(\mathcal{E}, \mathcal{C})$ and let $m^i \stackrel{\text{def}}{=} ((f^i), (A^i, E^i, C^i), (A^{i+1}, E^{i+1}, C^{i+1})), i = 1, 2$, be a morphism of $\mathcal{A}t\ell^{\text{Classic}}(\mathcal{E}, \mathcal{C})$. Then

**Composition:** $f^2 \circ f^1$ is a morphism of $\mathcal{E}$ and is an $E^1$-$E^3$ m-atlas morphism of $A^1$ to $A^3$ in the coordinate spaces $C^1, C^3$ by lemma 13.8 (Composition of classic m-atlas morphisms) on page 104.

**Associativity:** Composition is associative by lemma 1.19 (Tuple composition for labeled morphisms) on page 12.

**Identity:** $\text{Id}_{(A^i, E^i, C^i)}$ is an identity morphism by lemma 1.19.

$\square$

*Similar results follow with restrictions on the admissible atlases, the admissible morphisms, or both.*

### 13.2.2 Classic m-atlas functors

**Definition 13.18** ($\mathcal{F}_{\text{M,Classic}}$). Let $E^i, i = 1, 2, C^i$ be model spaces and $A^i$ an m-atlas of $E^i$ in the coordinate model space $C^i$. Then

$$\mathcal{F}_{\text{M,Classic}}(A^i, E^i, C^i) \stackrel{\text{def}}{=} (A^i, E^i, C^i) \tag{13.14}$$

Let $f \stackrel{\text{def}}{=} (f_0 \colon E^1 \longrightarrow E^2, f_1 \colon C^1 \longrightarrow C^2,)$ be an M-atlas morphism of $A^1$ to $A^2$ in the coordinate spaces $C^1, C^2$. Then

$$\mathcal{F}_{\text{M,Classic}}((f_0 \colon E^1 \longrightarrow E^2, f_1 \colon C^1 \longrightarrow C^2,), (A^1, E^1, C^1), (A^2, E^2, C^2)) \stackrel{\text{def}}{=}$$
$$((f_0 \colon E^1 \longrightarrow E^2), (A^1, E^1, C^1), (A^2, E^2, C^2)) \tag{13.15}$$

**Theorem 13.19** ($\mathcal{F}_{\text{M,Classic}}$ is a functor). *Let $E$ and $C$ be sets of model spaces, Then $\mathcal{F}_{\text{M,Classic}}$ is a functor from $\mathcal{A}t\ell(E, C)$ to $\mathcal{A}t\ell^{\text{classic}}(E, C)$.*

*Proof.* Let $o^i \stackrel{\text{def}}{=} (A^i, E^i, C^i), i \in [1, 3]$, be an object of $\mathcal{A}t\ell(E, C)$ and $m^i \stackrel{\text{def}}{=} ((f_0^i, f_1^i), o^i, o^{i+1})$, $i = 1, 2$, be a morphism of $\mathcal{A}t\ell(E, C)$.

$\mathcal{F}_{\text{M,Classic}}$ satisfies these criteria:

**Preservation of endpoints:**

$$\mathcal{F}_{\text{M,Classic}} \, o^i = o^i \stackrel{\text{Ob}}{\in} \mathcal{A}t\ell^{\text{classic}}(E, C)$$
$$\mathcal{F}_{\text{Top,M}} \, m^i = ((f_0 \colon E^i \longrightarrow E^{i+1}), (A^i, E^i, C^i), (A^{i+1}, E^{i+1}, C^{i+1}))$$
$$\text{is a morphism from } o^i \text{ to } o^{i+1}.$$

**Identity:** Let $\big((\mathrm{Id}_{E^i}, \mathrm{Id}_{C^i}), o^i, o^i\big)$ be an identity morphism of $o^i$ in $\mathcal{A}t\ell(E, C)$. Then $\big((\mathrm{Id}_{E^i}\colon E^i \longrightarrow E^i), o^i, o^i\big)$ is an identity morphism of $o^i$ in $\mathcal{A}t\ell^{\mathrm{classic}}(E, C)$.

**Composition:**

$$\begin{aligned}
\mathscr{F}_{\mathrm{M,Classic}}(m^2 \overset{A}{\circ} m^i) &= \mathscr{F}_{\mathrm{M,Classic}}\big((f_0^2 \circ f_0^1, f_1^2 \circ f_1^1), o^1, o^3\big) \\
&= \big((f_0^2 \circ f_0^1), o^1, o^3\big) \\
&= \big((f_0^2), o^1, o^2\big) \overset{A}{\circ} \big((f_0^1), o^2, o^3\big) \\
&= \mathscr{F}_{\mathrm{M,Classic}}\big((f_0^2, f_1^2), o^2, o^3\big) \overset{A}{\circ} \\
&\quad \mathscr{F}_{\mathrm{M,Classic}}\big((f_0^1, f_1^1), o^1, o^2\big) \\
&= \mathscr{F}_{\mathrm{M,Classic}}\, m^2 \overset{A}{\circ} \mathscr{F}_{\mathrm{M,Classic}}\, m^1
\end{aligned}$$

□

# Part VII
# Equivalence of manifolds

For a manifold[11] the coordinate category is open subsets of a Banach space or more generally a Fréchet space, with an appropriate choice of morphisms. Choosing a separating hyperplane and half space, with open sets in the chosen half space, allows manifolds with boundary. Similarly, choosing a ball with tails allows a manifold with tails.

For differentiable manifolds the coordinate category is similar, but the morphisms are limited to those sufficiently differentiable, in order to impose a diferentiability constraint on the transition functions $t_\beta^\alpha = \phi_\beta \circ \phi_\alpha^{-1}$.

This part of the paper defines $C^k$-atlases, $C^k$-manifolds, local coordinate spaces equivalent to $C^k$-manifolds, categories of them and functors, and gives some basic results.

## 14 Linear spaces and linear model spaces

This subsection defines spaces and related model spaces suitable for use as the coordinate spaces of generalized $C^k$ manifolds.

**Definition 14.1** (Linear spaces). Let $C$ be a locally arcwise connected topological subspace of a real (complex) Banach or Fréchet space. Then $C$ is a real (complex) linear space. If $C$ has a non-void interior, i.e., contains a ball, then it is non-degenerate.

---

[11]The literature has several different definitions of a manifold. This paper uses one chosen for ease of exposition.

*Remark* 14.2. This paper uses the term *ball* to refer to balls of the underlying space but uses the terms *open set* and *neighborhoods* to refer to the relative topology.

Let $\mathscr{C}$ be a small category whose objects are real (complex) linear spaces and whose morphisms are $C^k$ functions. Then $\mathscr{C}$ is a $C^k$ linear category. If each object of $\mathscr{C}$ is non-degenerate then $\mathscr{C}$ is non-degenerate.

Let $C$ be a real (complex) linear space. Then $C_{C^k-\mathcal{T}riv}$, the category of all $C^k$ functions between open subspaces of $C$, is the trivial $C^k$ linear category of $C$.

**Lemma 14.3** (Open subsets of linear spaces are locally arcwise connected). *Let $U$ be an open subset of the real (complex) linear space $S$. Then $U$ is locally arcwise connected.*

*Proof.* An open subset of a locally arcwise connected space is locally arcwise connected. □

**Definition 14.4** (Linear model spaces). Let $S$ be a real (complex) linear space and $\mathbf{S} \stackrel{\text{def}}{=} (S, \mathcal{S})$ a model space for $S$. Then $\mathbf{S}$ is a real (complex) linear model space.

Let $\mathbf{S} \stackrel{\text{def}}{=} (S, \mathcal{S})$ be a real (complex) linear model space such that every morphism of $\mathcal{S}$ is a $C^k$ function. Then $\mathbf{S}$ is a real (complex) $C^k$ linear model space.

Let $\mathcal{S}$ be a small category whose objects are real (complex) $C^k$ linear model spaces and whose morphisms are $C^k$ model functions. Then $\mathcal{S}$ is a $C^k$ linear model category.

**Definition 14.5** (Trivial $C^k$ linear model spaces). Let $C$ be a real (complex) linear space and $\mathscr{C}$ the category of all $C^k$ functions between open sets of $C$. Then $C_{C^k-\mathbf{triv}} \stackrel{\text{def}}{=} (C, \mathscr{C})$ is the trivial $C^k$ linear model space of $C$ and $C_{C^k-\mathbf{triv}}$ is a real (complex) trivial $C^k$ linear model space.

Let $\mathbf{C}$ be a set of real (complex) linear spaces.

$\mathbf{C}_{C^k-\mathbf{triv}} \stackrel{\text{def}}{=} \left\{ C'_{C^k-\mathbf{triv}} \middle| C' \in \mathbf{C} \right\}$ is the set of trivial $C^k$ linear model spaces of $\mathbf{C}$.

The category of trivial $C^k$ linear model spaces of $\mathbf{C}$, abbreviated $\mathbf{C}_{C^k-\mathcal{T}riv}$, is the category whose objects are $\mathbf{C}_{C^k-\mathbf{triv}}$ and whose morphisms are all the $C^k$ functions among them.

$\mathbf{C}_{C^k-\text{op-triv}}$, the category of open trivial $C^k$ model spaces in $\mathbf{C}$, is the category whose objects are $\left\{ U_{C^k-\mathbf{triv}} \middle| U \in \mathbf{C}_{\text{op}} \right\}$, the trivial $C^k$ linear model spaces of non-null open sets of spaces in $\mathbf{C}$, and whose morphisms are all the $C^k$ functions among them.

**Lemma 14.6** (Trivial $C^k$ linear model spaces). *Let $C^i$, $i = 1, 2$, be a real (complex) linear space and $\mathscr{C}^i \stackrel{\text{def}}{=} \text{Cat}\left( C^i_{C^k-\mathbf{triv}} \right)$.*

1. $C^i_{C^k-\mathbf{triv}}$ is a $C^k$ linear model space.

   *Proof.*

   (a) $\mathrm{Ob}(\mathscr{C}^i)$ is an open cover for $C^i$.
   (b) $\mathrm{Ob}(\mathscr{C}^i)$ is closed under finite intersections.
   (c) The morphisms of $\mathscr{C}^i$ are $C^k$, hence continuous.
   (d) If $f\colon A \longrightarrow B$ is a morphism, $A' \overset{\mathrm{Ob}}{\in} \mathscr{C}^i \subseteq A \overset{\mathrm{Ob}}{\in} \mathscr{C}^i$, $B' \overset{\mathrm{Ob}}{\in} \mathscr{C}^i \subseteq B \overset{\mathrm{Ob}}{\in} \mathscr{C}^i$ and $f[A'] \subseteq B'$ then $f\restriction_{A'}\colon A' \longrightarrow B'$ is $C^k$ and thus a morphism.
   (e) If $A' \overset{\mathrm{Ob}}{\in} \mathscr{C}^i \subseteq A \overset{\mathrm{Ob}}{\in} \mathscr{C}^i$ then the inclusion map $i\colon A' \hookrightarrow A$ is $C^k$ and thus a morphism.
   (f) Restricted sheaf condition: Whenever
      i. $U_\alpha$ and $V_\alpha$, $\alpha \prec A$, are objects of $\mathscr{C}^i$.
      ii. $f_\alpha\colon U_\alpha \longrightarrow V_\alpha$ are morphisms of $\mathscr{C}^i$.
      iii. $U \overset{\mathrm{def}}{=} \bigcup_{\alpha \prec A} U_\alpha \overset{\mathrm{Ob}}{\in} \mathscr{C}^i$,
      iv. $V \overset{\mathrm{def}}{=} \bigcup_{\alpha \prec A} V_\alpha \overset{\mathrm{Ob}}{\in} \mathscr{C}^i$
      v. $f\colon U \longrightarrow V$ is a continuous function and for every $\alpha \prec A$, $f$ agrees with $f_\alpha$ on $U_\alpha$

      then $f$ is $C^k$ and thus a morphism of $\mathscr{C}^i$.

      □

2. Let $f\colon C^1 \longrightarrow C^2$ be $C^k$. Then $f$ is a model function from $C^1_{C^k-\mathbf{triv}}$ to $C^2_{C^k-\mathbf{triv}}$.

   *Proof.* Let $U^i$ be a model neighborhood of $C^i_{C^k-\mathbf{triv}}$, $i = 1, 2$.

   (a) Since $U^2$ is a model neighborhood of $C^2_{C^k-\mathbf{triv}}$, $U^2$ is open, $f^{-1}[U^2]$ is open and thus a model neighborhood of $C^1_{C^k-\mathbf{triv}}$.
   (b) $f[U^1] \subseteq C^2$. Since $C^2$ is open, it is a model neighborhood of $C^2_{C^k-\mathbf{triv}}$.

   □

3. Let $C^1$ be an open subspace of $C^2$. Then $C^1_{C^k-\mathbf{triv}} \overset{\text{strict}-\text{mod}}{\subseteq} C^2_{C^k-\mathbf{triv}}$.

   *Proof.* By definition 14.5 above, $C^i_{C^k-\mathbf{triv}} = (C^i, \mathscr{C}^i)$, where $\mathscr{C}^i)$ is the category of all $C^k$ functions between open sets of $C^i$.

$$C^1_{C^k-\mathbf{triv}} \overset{\text{mod}}{\subseteq} C^2_{C^k-\mathbf{triv}}:$$

(a) $C^1$ is a subspace of $C^2$ by hypothesis.

(b) Every object of $\mathscr{C}^1$ is open in $C^1$, thus open in $C^2$ and thus an object of $\mathscr{C}^2$. If $f\colon U \longrightarrow V \overset{\text{Ar}}{\in} \mathscr{C}^1$ then $U$ and $V$ are open in $C^1$ and $f$ is $C^k$ in $C^1$, thus $U$ and $V$ are open in $C^2$ and $f$ is $C^k$ in $C^2$ and $f\colon U \longrightarrow V \overset{\text{Ar}}{\in} \mathscr{C}^2$. If $f\colon U \longrightarrow V$ is a morphism of $\mathscr{C}^2$, $U \overset{\text{Ob}}{\in} \mathscr{C}^1$ and $V \overset{\text{Ob}}{\in} \mathscr{C}^1$, then $U$ and $V$ are open in $C^1$ and $f$ is $C^k$ in $C^2$, thus $C^k$ in $C^1$, and $f\colon U \longrightarrow V \overset{\text{Ar}}{\in} \mathscr{C}^1$.

(c) If $U \overset{\text{Ob}}{\in} \mathscr{C}^2$ then $U$ is open in $C^2$, hence $U \cap C^2$ is open in $C^2$, hence $U \cap C^2$ is open in $C^2$ and $U \overset{\text{Ob}}{\in} \mathscr{C}^2$.

(d) If $U \overset{\text{Ob}}{\in} \mathscr{C}^2$ and $U \subseteq C^1$ then U is open is $C^2$, hence open is $C^1$, hence $U \overset{\text{Ob}}{\in} \mathscr{C}^1$.

$$C^1_{C^k-\mathbf{triv}} \overset{\text{strict-mod}}{\subseteq} C^2_{C^k-\mathbf{triv}}:$$ $\text{Id}_{C^1,C^2}$ is open and continuous because $C^1$ is an open subspace of $C^2$.

□

**Definition 14.7** ($C^k$ singleton categories). Let $C$ be a $C^k$ linear model space. Then the $C^k$ singleton category of $C$, abbreviated $C_{C^k-\mathscr{S}ing}$, is the category whose sole object is $C$ and whose morphisms are all of the $C^k$ model functions from $C$ to itself.

# 15 $C^k$-nearly commutative diagrams

Let $C$ be a linear space, $\mathscr{C} \overset{\text{def}}{=} S_{C^k-\mathscr{T}riv}$, $C \overset{\text{def}}{=} C_{C^k-\mathbf{triv}}$ and $D$ a tree with two branches, whose nodes are topological spaces $U_i$ and $V^j$ and whose links are continuous functions $f_i\colon U_i \longrightarrow U_{i+1}$ and $f'_j\colon U_j \longrightarrow U_{j+1}$ between the spaces:

$$D = \{f_0\colon U_0 = V_0 \longrightarrow U_1, \ldots, f_{m-1}\colon U_{m-1} \longrightarrow U_m,$$
$$f'_0\colon U_0 = V_0 \longrightarrow V_1, \ldots, f'_{m-1}\colon V_{m-1} \longrightarrow V_n\}$$

with $U_0 = V_0$, $U_m \subseteq C$ and $V_n \subseteq C$ open, as shown in fig. 1 (Uncompleted nearly commutative diagram) on page 15.

**Definition 15.1** ($C^k$-nearly commutative diagrams). $D$ is (left,right,strongly) $C^k$-nearly commutative in linear space $C$ iff $D$ is (left,right,strongly) M-nearly commutative in model space $C$.

**Definition 15.2** ($C^k$-nearly commutative diagrams at a point). Let $C$, $\mathscr{C}$, $\boldsymbol{C}$ and $D$ be as above and $x$ be an element of the initial node. $D$ is (left,right,strongly) $C^k$-nearly commutative in $C$ at $x$ iff $D$ is (left,right,strongly) M-nearly commutative in model space $\boldsymbol{C}$ at $x$.

**Definition 15.3** ($C^k$-locally nearly commutative diagrams). Let $C$, $\mathscr{C}$, $\boldsymbol{C}$ and $D$ be as above. $D$ is (left,right,strongly) $C^k$-locally nearly commutative in $C$ iff $D$ is (left,right,strongly) M-locally nearly commutative in model space $\boldsymbol{C}$.

**Lemma 15.4.** *Let $C$, $\mathscr{C}$, $\boldsymbol{C}$ and $D$ be as above. $D$ is $C^k$-nearly commutative in linear space $C$ iff There is a $C^k$ diffeomorphism $\hat{f}\colon U_m \rightarrowtail\!\!\!\!\xrightarrow{\tilde{=}}\!\!\!\twoheadrightarrow V_n$ making the graph a commutative diagram, as shown in fig. 2.*

*Proof.* A set is a model neighborhood of $\boldsymbol{C}$ (an object of $\underset{C^k-\mathcal{T}riv}{S}$) iff it is an open subset of $A$. A function is a morphism of $\boldsymbol{C}$ (a morphism of $\underset{C^k-\mathcal{T}riv}{S}$) iff it is a $C^k$ function between open subsets of $S$. □

# 16   $C^k$ charts

**Definition 16.1** ($C^k$ charts). Let $C$ be a linear space and $E$ a topological space. A $C^{k}$[12] chart $(U, V, \phi)$ of $E$ in the coordinate space $C$ consists of

1. A nonvoid open subset $U \subseteq E$, known as a coordinate patch

2. An open subset $V \subseteq C$

3. A homeomorphism $\phi \colon U \rightarrowtail\!\!\!\!\xrightarrow{\tilde{=}}\!\!\!\twoheadrightarrow V$, known as a coordinate function

*Remark* 16.2. I consider it clearer to explicate the range, rather than the conventional usage of specifying only the domain and function or the minimalist usage of specifying only the function.

A chart $(U, V, \phi)$ is non-degenerate iff $V$ contains a ball of the underlying Banach or Fréchet space.

**Lemma 16.3** ($C^k$ charts). *Let $C$ be a linear space and $E$ a topological space. The triple $(U, V, \phi)$ is a $C^k$ chart of $E$ in the coordinate space $C$ iff it is an m-chart of $\underset{\mathbf{triv}}{E}$ in the coordinate space $\underset{C^k-\mathbf{triv}}{C}$.*

*Proof.* $U$ is a model neighborhood of $\underset{\mathbf{triv}}{E}$ iff it is a nonvoid open set of $E$. $V$ is a model neighborhood of $\underset{C^k-\mathbf{triv}}{C}$ iff it is a nonvoid open set of $C$. $\phi$ is required to be a homeomorphis from $U$ to $V$ in either case. □

---
[12]With $k \in \mathbb{N} \cup \{\infty, \omega\}$.

**Definition 16.4** (C$^k$ subcharts)**.** Let $(U, V, \phi)$ be a C$^k$ chart of $E$ in the coordinate space $C$ and $U'$ be a novoid open subset of $U$. Then $(U', V', \phi') \stackrel{\text{def}}{=} (U', \phi[U'], \phi\!\restriction_{U',V'})$ is a subchart of $(U, V, \phi)$.

By abuse of language we will write $(U', V', \phi)$ for $(U', V', \phi')$.

**Lemma 16.5** (C$^k$ subcharts)**.** *Let $(U, V, \phi)$ be a C$^k$ chart of $E$ in the coordinate space $C$ and $(U', V', \phi')$ be a triple. Then $(U', V', \phi')$ is a subchart of $(U, V, \phi)$ iff it is a subchart of $(U, V, \phi)$ considered as a chart of $E_{\text{triv}}$ in the coordinate space $C_{\text{C}^k-\text{triv}}$.*

*Proof.* $(U', V', \phi')$ satisfies the conditions of definition 16.1;

1. $U$ is open by definition 16.1

2. $U'$ is required to be a nonvoid open subset of $U$ in either case, and thus a model neighborhood of $E_{\text{triv}}$.

3. $\phi$ is a homeomorphism, so $V' = \phi[U']$ is open, and thus a model neighborhood of $C_{\text{C}^k-\text{triv}}$.

4. $U'$ is a model neighborhood of $E_{\text{triv}}$ iff it is an open set of $E$.

5. $\phi$ is a homeomorphism, so $\phi\!\restriction_{U'}\colon U' \underset{\sim}{\rightarrowtail\!\!\!\twoheadrightarrow} \phi[U']$ is also.

□

**Corollary 16.6** (C$^k$ subcharts)**.** *Let $(U, V, \phi)$ be a C$^k$ chart of $E$ in the coordinate space $C$ and $(U', V', \phi')$ be a subchart of $(U, V, \phi)$. Then $(U', V', \phi')$ is a C$^k$ chart of $E$ in the coordinate space $C$.*

*Proof.* The result follows from lemma 16.3 (C$^k$ charts) above, lemma 9.3 (M-charts) on page 43 and definition 9.5 (Subcharts) on page 44. □

**Definition 16.7** (C$^k$ compatibility)**.** Let $(U, V, \phi)$ and $(U', V', \phi')$ be C$^k$ charts of $E$ in the coordinate space $C$. Then $(U, V, \phi)$ is C$^k$ compatible with $(U', V', \phi')$ iff either

1. $U$ and $U'$ are disjoint

2. The transition function $t = \phi' \circ \phi^{-1}\!\restriction_{\phi[U \cap U']}$ is a C$^k$ diffeomorphism.

**Lemma 16.8** (Symmetry of C$^k$ compatibility)**.** *Let $(U, V, \phi)$ and $(U', V', \phi')$ be C$^k$ charts of $E$ in the coordinate space $C$. Then $(U, V, \phi)$ is C$^k$ compatible with $(U', V', \phi')$ iff $(U', V', \phi')$ is C$^k$ compatible with $(U, V, \phi)$.*

*Proof.* It suffices to prove the implication in only one direction.

1. $U \cap U' = U' \cap U$.

2. Since the transition function $t = \phi' \circ \phi^{-1}\restriction_{\phi[U \cap U']}$ is a $C^k$ diffeomorphism of $C$, so is $t^{-1} = \phi \circ \phi'^{-1}\restriction_{\phi'[U \cap U']}$.

□

**Lemma 16.9** ($C^k$ compatibility of subcharts). *Let $(U_i, V_i, \phi_i)$, $i = 1, 2$, be $C^k$ charts of $E$ in the coordinate space $C$, $(U'_i, V'_i, \phi'_i)$ be subcharts and $(U_1, V_1, \phi_1)$ be $C^k$ compatible with $(U_2, V_2, \phi_2)$. Then $(U'_1, V'_1, \phi'_1)$ is $C^k$ compatible with $(U'_2, V'_2, \phi'_2)$.*

*Proof.* If $U_1 \cap U_2 = \emptyset$ then $U'_1 \cap U'_2 = \emptyset$. If $U'_1 \cap U'_2 = \emptyset$ then $(U'_1, V'_1, \phi'_1)$ is $C^k$ compatible with $(U'_2, V'_2, \phi'_2)$. Otherwise, the transition function $t^1_2 \stackrel{\text{def}}{=} \phi_2 \circ \phi_1^{-1}\restriction_{\phi_1[U_1 \cap U_2]}$ is a $C^k$ diffeomorphism and hence $t^1_2\restriction_{\phi_1[U'_1 \cap U'_2]}: \phi_1[U'_1 \cap U'_2] \xrightarrow{\tilde{=}} \phi_2[U'_1 \cap U'_2]$ is a $C^k$ diffeomorphism.

□

**Corollary 16.10** ($C^k$ compatibility with subcharts). *Let $(U, V, \phi)$ be a $C^k$ charts of $E$ in the coordinate space $C$ and $(U', V', \phi')$ a subchart. Then $(U', V', \phi')$ is $C^k$ compatible with $(U, V, \phi)$.*

*Proof.* $(U, V, \phi)$ is $C^k$ compatible with itself and is a subchart of itself,

□

**Definition 16.11** (Covering by $C^k$ charts). *Let $A$ be a set of charts of $E$ in the coordinate space $C$. $A$ covers $E$ iff $\pi_1[A]$ covers $E$, i.e., $E = \bigcup \pi_1[A]$.*

# 17   $C^k$-atlases

A set of charts can be atlases for different coordinate spaces even if it is for the same total space. In order to aggregate them into categories, there must be a way to distinguish them. Including the two[13] spaces in the definitions of the categories serves the purpose.

**Definition 17.1** ($C^k$-atlases). *Let $A$ be a set of mutually $C^k$ compatible charts of $E$ in the coordinate space $C$. $A$ is a $C^k$-atlas of $E$ in the coordinate space $C$, abbreviated $\text{isAtl}^{C^k}_{\text{Ob}}(A, E, C)$, iff $A$ covers $E$.*

*$A$ is a full $C^k$-atlas of $E$ in the coordinate space $C$, abbreviated $\text{isAtl}^{C^k}_{\text{Ob}\atop\text{full}}(A, E, C)$, iff*

1. $\pi_1[A]$ covers $E$
2. $\pi_2[A]$ covers $C$.

---
[13]The total space is redundant, but convenient.

**A** is non-degenerate iff it contains a non-degenerate chart.
By abuse of language we write $U \in A$ for $U \in \pi_1[A]$.
Let $E$ be a topological space and $C$ a linear space. Then

$$\mathcal{A}t\ell^{C^k}_{Ob}(E, C) \stackrel{def}{=} \left\{ (A, E, C) \,\middle|\, \text{isAtl}^{C^k}_{Ob}(A, E, C) \right\} \tag{17.1}$$

$$\mathcal{A}t\ell^{C^k}_{Ob\,\text{full}}(E, C) \stackrel{def}{=} \left\{ (A, E, C) \,\middle|\, \text{isAtl}^{C^k}_{Ob\,\text{full}}(A, E, C) \right\} \tag{17.2}$$

Let $\boldsymbol{E}$ be a set of topological spaces and $\boldsymbol{C}$ a set of linear spaces. Then

$$\mathcal{A}t\ell^{C^k}_{Ob}(\boldsymbol{E}, \boldsymbol{C}) \stackrel{def}{=} \bigcup_{\substack{E_\mu \in \boldsymbol{E} \\ C_\mu \in \boldsymbol{C}}} \mathcal{A}t\ell_{Ob}(E_\mu, C_\mu) \tag{17.3}$$

$$\mathcal{A}t\ell^{C^k}_{Ob\,\text{full}}(\boldsymbol{E}, \boldsymbol{C}) \stackrel{def}{=} \left\{ (A, E, C) \,\middle|\, \text{isAtl}^{C^k}_{Ob\,\text{full}}(A, E, C) \right\} \tag{17.4}$$

**Lemma 17.2** ($C^k$-atlases). *Let $E$ be a topological space, $C$ be a linear space and $A$ be a $C^k$-atlas of $E$ in the coordinate space $C$.*

*If $A$ is non-degenerate then $C$ is non-degenerate.*

*Proof.* Let $(U, V, \phi)$ be a non-degenerate chart of $A$. $V$ contains a ball and is contained in $C$. $\square$

*If $C$ is non-degenerate and $A$ is full then $A$ is non-degenerate.*

*Proof.* Let $B$ be a ball in $C$ with center $v$. Since $A$ is full, it contains a chart $(U, V, \phi)$ with $v \in V$. Then $V$ contains a ball with center $v$. $\square$

**Definition 17.3** (Compatibility of charts with $C^k$-atlases). A chart $(U, V, \phi)$ of $E$ in the coordinate space $C$ is $C^k$ compatible with a $C^k$-atlas $A$ iff it is $C^k$ compatible with every chart in the atlas.

**Lemma 17.4** (Compatibility of subcharts with $C^k$-atlases). *Let $A$ be a $C^k$-atlas of $E$ in the coordinate space $C$ and $\boldsymbol{C}_1 = (U_1, V_1, \phi_1)$ a $C^k$ chart in $A$. Then any subchart of $\boldsymbol{C}_1$ is $C^k$ compatible with $A$.*

*Proof.* Let $\boldsymbol{C}' = (U', V', \phi')$ be a subchart of $\boldsymbol{C}_1$ and $\boldsymbol{C}_2 = (U_2, V_2, \phi_2)$ another chart in $A$.

1. If $U_1 \cap U_2 = \emptyset$, then $U' \cap U_2 = \emptyset$.

2. If $U' \cap U_2 = \emptyset$ then $\boldsymbol{C}'$ is $C^k$ compatible with $\boldsymbol{C}_2$.

3. Otherwise the transition function $t_2^1 \stackrel{def}{=} \phi_2 \circ \phi_1^{-1} \restriction_{\phi_1[U_1 \cap U_2]}$ is a $C^k$ diffeomorphism and thus $t_2^1 \restriction_{\phi_1[U' \cap U_2]}$ is a $C^k$ diffeomorphism.

**Lemma 17.5** (Extensions of $C^k$-atlases). *Let $A$ be a $C^k$ atlas of $E$ in the coordinate space $C$ and $(U_i, V_i, \phi_i)$, $i = 1, 2$, be a $C^k$ chart of $E$ in the coordinate space $C$ $C^k$ compatible with $A$ in the coordinate space $C$. Then $(U_1, V_1, \phi_1)$ is $C^k$ compatible with $(U_2, V_2, \phi_2)$ in the coordinate space $C$.*

*Proof.* If $U_1 \cap U_2 = \emptyset$ then $(U_1, V_1, \phi_1)$ is $C^k$ compatible with $(U_2, V_2, \phi_2)$. Otherwise, $\phi_2 \circ \phi_1^{-1} \restriction_{\phi_1[U_1 \cap U_2]} \colon \phi_1[U_1 \cap U_2] \xrightarrow{\tilde{=}} \phi_2[U_1 \cap U_2]$ is a homeomorphism. It remains to show that $\phi_2 \circ \phi_1^{-1} \restriction_{\phi_1[U_1 \cap U_2]}$ is a $C^k$ diffeomorphism. Let $(U'_\alpha, V'_\alpha, \phi'_\alpha)$, $\alpha \prec A$, be charts in $A$ such that $U_1 \cap U_2 \subseteq \bigcup_{\alpha \prec A} U'_\alpha$ and $U_1 \cap U_2 \cap U'_\alpha \neq \emptyset$, $\alpha \prec A$. Since the charts are $C^k$ compatible with $(U'_\alpha, V'_\alpha, \phi'_\alpha)$, $\phi_2 \circ \phi'^{-1}_\alpha \restriction_{U_1 \cap U_2 \cap U'_\alpha}$ and $\phi'_\alpha \circ \phi_1^{-1} \restriction_{U_1 \cap U_2 \cap U'_\alpha}$ are $C^k$ diffeomorphisms and thus $\phi_2 \circ \phi_1^{-1} = \phi_2 \circ \phi'^{-1}_\alpha \circ \phi'_\alpha \circ \phi_1^{-1}$ is a $C^k$ diffeomorphism. □

**Definition 17.6** (Maximal $C^k$-atlases). Let $E$ be a topological space, $C$ be a linear space and $A$ be a non-degenerate $C^k$-atlas of $E$ in the coordinate space $C$.

$A$ is a maximal $C^k$-atlas of $E$ in the coordinate space $C$, abbreviated $\mathrm{isAtl}^{C^k}_{\mathrm{Ob}\ \max}(A, E, C)$, iff $A$ is a $C^k$-atlas that cannot be extended by adding an additional $C^k$ compatible chart.

$$\underset{\text{max-full}}{\mathrm{isAtl}^{C^k}_{\mathrm{Ob}}(A, E, C)} \overset{\text{def}}{=} \underset{\text{full}}{\mathrm{isAtl}^{C^k}_{\mathrm{Ob}}(A, E, C)} \wedge \underset{\text{max}}{\mathrm{isAtl}^{C^k}_{\mathrm{Ob}}(A, E, C)} \tag{17.5}$$

$A$ is a semi-maximal $C^k$-atlas of $E$ in the coordinate space $C$, abbreviated $\mathrm{isAtl}^{C^k}_{\mathrm{Ob}\ \text{S-max}}(A, E, C)$, iff whenever $(U, V, \phi) \in A$, $U' \subseteq U$, $V' \subseteq V$ and $V'' \subseteq C$ are open, $\phi[U'] = V'$ and $\phi' \colon V' \xrightarrow{\tilde{=}} V''$ is a $C^k$ diffeomorphism then $(U', V'', \phi' \circ \phi) \in A$.

$$\underset{\text{S-max-full}}{\mathrm{isAtl}^{C^k}_{\mathrm{Ob}}(A, E, C)} \overset{\text{def}}{=} \underset{\text{full}}{\mathrm{isAtl}^{C^k}_{\mathrm{Ob}}(A, E, C)} \wedge \underset{\text{S-max}}{\mathrm{isAtl}^{C^k}_{\mathrm{Ob}}(A, E, C)} \tag{17.6}$$

**Lemma 17.7** (Maximal $C^k$-atlases are semi-maximal $C^k$-atlases). *Let $E$ be a topological space, $C$ a $C^k$-linear space and $A$ a maximal $C^k$-atlas of $E$ in the coordinate space $C$. Then $A$ is a semi-maximal $C^k$-atlas of $E$ in the coordinate space $C$.*

*Proof.* Let $(U, V, \phi) \in A$, $U' \subseteq U$, $V' \subseteq V$ and $V'' \subseteq C$ be open, $\phi[U'] = V'$ and $\phi' \colon V' \xrightarrow{\tilde{=}} V''$ be a $C^k$ diffeomorphism. $(U', V', \phi)$ is a subchart of $(U, V, \phi)$ and by lemma 17.4 (Compatibility of subcharts with $C^k$-atlases) on page 118 is $C^k$ compatible with the charts of $A$. Since $\phi'$ is a $C^k$ diffeomorphism, $(U', V'', \phi' \circ \phi)$ is $C^k$ compatible with the charts of $A$. Since $A$ is maximal, $(U', V'', \phi' \circ \phi)$ is a chart of $A$. □

**Theorem 17.8** (Existence and uniqueness of maximal $C^k$-atlases). *Let $A$ be a $C^k$-atlas of $E$ in the coordinate space $C$. Then there exists a unique maximal $C^k$-atlas $\text{Atlas}^{C^k}_{\max}(A, E, C)$ of $E$ in the coordinate space $C$ compatible with $A$.*

*Proof.* Let $P$ be the set of all $C^k$-atlases of $E$ in the coordinate space $C$ containing $A$ and $C^k$ compatible in the coordinate space $C$ with all of the $C^k$ charts in $A$. Let $P_{\max}$ be a maximal chain of $A$. Then $A' = \bigcup P_{\max}$ is a maximal $C^k$ atlas of $E$ in the coordinate space $C$ $C^k$ compatible with $A$. Uniqueness follows from lemma 17.5 (Extensions of $C^k$-atlases) on page 119. $\square$

**Definition 17.9** (Sets of $C^k$-atlases). Let $E$ be a topological space and $C$ be a linear space.

$$\mathcal{Atl}^{C^k}_{\text{Ob}\,(\text{full,S-max,max,S-max-full,max-full})}(E, C) \stackrel{\text{def}}{=} \left\{ (A, E, C) \,\middle|\, \text{isAtl}^{C^k}_{\text{Ob}\,(\text{full,S-max,max,S-max-full,max-full})}(A, E, C) \right\} \quad (17.7)$$

Let $\boldsymbol{E}$ be a set of topological spaces and $\boldsymbol{C}$ a set of linear spaces. Then

$$\mathcal{Atl}^{C^k}_{\text{Ob}\,(\text{full,S-max,max,S-max-full,max-full})}(\boldsymbol{E}, \boldsymbol{C}) \stackrel{\text{def}}{=} \bigcup_{\substack{E \in \boldsymbol{E} \\ C \in \boldsymbol{C}}} \mathcal{Atl}^{C^k}_{\text{Ob}\,(\text{full,S-max,max,S-max-full,max-full})}(E, C) \quad (17.8)$$

# 18 $C^k$-atlas (near) morphisms and functors

This section introduces a taxonomy for morphisms between $C^k$-atlases, defines categories of $C^k$-atlases, defines functors amomg them, defines functors between them and categories of m-atlases, defines inverse functors and proves some basic reasults.

## 18.1 $C^k$-atlas (near) morphisms

This subsection introduces the notions of classic $C^k$-atlas near morphisms and of classic $C^k$-atlas morphisms, and proves some basic results. Classic $C^k$-atlas near morphisms between maximal atlases will be proven to be classic $C^k$-atlas morphisms.

### 18.1.1 Definitions of $C^k$-atlas (near) morphisms

**Definition 18.1** ($C^k$-atlas near morphisms). Let $E^i$, $i = 1, 2$, be a topological space, $C^i$ a linear space, $A^i$ a $C^k$-atlases of $E^i$ in the coordinate space $C^i$ and $\boldsymbol{f} \stackrel{\text{def}}{=} (f_0, f_1)$ a pair[14] of functions.

---

[14]The conventional definition uses only the first of the two functions and a slightly different compatibility condition.

$f$ is an $E^1$-$E^2$ $C^k$ near morphism of $A^1$ to $A^2$ in the coordinate spaces $C^1$, $C^2$, abbreviated as $\text{isAtl}_{\text{Ar}}^{C^k,\text{near}}(A^1, E^1, C^1, A^2, E^2, C^2, f_0, f_1)$, iff

1. $f_0\colon E^1 \longrightarrow E^2$ is a continuous function.

2. $f_1\colon C^1 \longrightarrow C^2$ is a $C^k$ function.

3. for any $(U^i, V^i, \phi^i\colon U^i \underset{\sim}{\rightarrowtail\!\!\!\twoheadrightarrow} V^i) \in A^i$, $i = 1, 2$, diagram (11.1) in definition 11.1 (M-atlas near morphisms for model spaces) on page 51 is $C^k$-locally nearly commutative in $C^2$, i.e., for any $u^1 \in I$ there are open sets $U'^1 \subseteq I$, $V'^1 \subseteq V^1$, $U'^2 \subseteq U^2$, $V'^2 \subseteq V^2$, $\hat{V}'^2 \subseteq C^2$ and a $C^k$ diffeomorphism $\hat{f}\colon \hat{V}'^2 \underset{\sim}{\rightarrowtail\!\!\!\twoheadrightarrow} V'^2$ such that eqs. (11.2) to (11.7) on pages 52 to 52 in definition 11.4 (M-atlas morphisms for model spaces) on page 53 hold.

$f$ is also a full (semi-maximal, maximal, full semi-maximal, full maximal) $E^1$-$E^2$ $C^k$ near morphism of $A^1$ to $A^2$ in the coordinate spaces $C^1$, $C^2$, abbreviated as $\text{isAtl}_{\text{Ar full (S-max,max,S-max-full,max-full)}}^{C^k,\text{near}}(A^1, E^1, C^1, A^2, E^2, C^2, f_0, f_1)$, iff

$$\text{isAtl}_{\text{Ob full (S-max,max,S-max-full,max-full)}}^{C^k,\text{near}}(A^1, E^1, C^1) \quad \wedge \quad \text{isAtl}_{\text{Ob full (S-max,max,S-max-full,max-full)}}^{C^k,\text{near}}(A^2, E^2, C^2)$$

The identity morphism of $(A^i, E^i, C^i)$ is

$$\text{Id}_{(A^i, E^i, C^i)} \overset{\text{def}}{=} ((\text{Id}_{E^i}, \text{Id}_{C^i}), (A^i, E^i, C^i), (A^i, E^i, C^i)) \tag{18.1}$$

This nomenclature will be justified later.

Let $E^i$, $i = 1, 2$, be topological spaces and $C^i$ be linear spaces. Then

$$\mathcal{A}tl_{\text{Ar (full,S-max,max,S-max-full,max-full)}}^{C^k,\text{near}}(E^1, C^1, E^2, C^2) \overset{\text{def}}{=} \left\{ ((f_0, f_1), (A^1, E^1, C^1), (A^2, E^2, C^2)) \,\middle|\, \right.$$
$$\left. \text{isAtl}_{\text{Ar (full,S-max,max,S-max-full,max-full)}}^{C^k,\text{near}}(A^1, E^1, C^1, A^2, E^2, C^2, f_0, f_1) \right\} \tag{18.2}$$

The identity morphism of $(A^i, E^i, C^i)$ is

$$\text{Id}_{(A^i, E^i, C^i)} \overset{\text{def}}{=} ((\text{Id}_{E^i}, \text{Id}_{C^i}), (A^i, E^i, C^i), (A^i, E^i, C^i)) \tag{18.3}$$

This nomenclature will be justified later.

**Definition 18.2** ($C^k$-atlas morphisms). Let $E^i$, $i = 1, 2$, be a topological space, $C^i$ a linear space, $A^i$ a $C^k$-atlases of $E^i$ in the coordinate space $C^i$ and $f \overset{\text{def}}{=} (f_0, f_1)$ a pair[15] of functions.

$f$ is an $E^1$-$E^2$ $C^k$-morphism of $A^1$ to $A^2$ in the coordinate spaces $C^1$, $C^2$, abbreviated as $\text{isAtl}_{\text{Ar}}^{C^k}(A^1, E^1, C^1, A^2, E^2, C^2, f_0, f_1)$, iff

---

[15]The conventional definition uses only the first of the two functions and a slightly different compatibility condition.

1. $f_0 \colon E^1 \longrightarrow E^2$ is a continuous function.
2. $f_1 \colon C^1 \longrightarrow C^2$ is a $C^k$ function.
3. for any $u^1 \in E^1$, any chart $(U^1, V^1, \phi^1 \colon U^1 \overset{\tilde{=}}{\rightarrowtail\!\!\!\twoheadrightarrow} V^1) \in \boldsymbol{A}^1$ at $u^1$ and any chart $(U^2, V^2, \phi^2 \colon U^2 \overset{\tilde{=}}{\rightarrowtail\!\!\!\twoheadrightarrow} V^2) \in \boldsymbol{A}^2$ at $f_0(u^1)$ there exists a subchart $(U'^1, V'^1, \phi^1 \colon U'^1 \overset{\tilde{=}}{\rightarrowtail\!\!\!\twoheadrightarrow} V'^1) \in \boldsymbol{A}^1$ at $u^1$, an open set $U'^2 \subseteq U^2$ and a chart $(U'^2, \hat{V}'^2, \phi'^2 \colon U'^2 \overset{\tilde{=}}{\rightarrowtail\!\!\!\twoheadrightarrow} \hat{V}'^2) \in \boldsymbol{A}^2$ at $f_0(u^1)$ such that $f_0[U'^1] \subseteq U'^2$ and diagram (11.8) in definition 11.4 (M-atlas morphisms for model spaces) on page 53 is commutative, as shown in fig. 6 (Completed m-atlas morphism) on page 54.

$\boldsymbol{f}$ is also a full (semi-maximal, maximal, full semi-maximal, full maximal) $E^1$-$E^2$ $C^k$ morphism of $\boldsymbol{A}^1$ to $\boldsymbol{A}^2$ in the coordinate spaces $C^1, C^2$, abbreviated as

$\mathrm{isAtl}_{\mathrm{Ar}\ \text{full (S-max,max,S-max-full,max-full)}}^{C^k}(\boldsymbol{A}^1, E^1, C^1, \boldsymbol{A}^2, E^2, C^2, f_0, f_1)$, iff

$\mathrm{isAtl}_{\mathrm{Ob}\ \text{full (S-max,max,S-max-full,max-full)}}^{C^k}(\boldsymbol{A}^1, E^1, C^1) \quad \wedge \quad \mathrm{isAtl}_{\mathrm{Ob}\ \text{full (S-max,max,S-max-full,max-full)}}^{C^k}(\boldsymbol{A}^2, E^2, C^2)$

The identity morphism of $(\boldsymbol{A}^i, E^i, C^i)$ is

$$\mathrm{Id}_{(\boldsymbol{A}^i, E^i, C^i)} \overset{\mathrm{def}}{=} ((\mathrm{Id}_{E^i}, \mathrm{Id}_{C^i}), (\boldsymbol{A}^i, E^i, C^i), (\boldsymbol{A}^i, E^i, C^i)) \tag{18.4}$$

This nomenclature will be justified later.

Let $E^i$, $i = 1, 2$, be topological spaces and $C^i$ be linear spaces. Then

$$\mathcal{A}t\ell_{\mathrm{Ar}}^{C^k}(E^1, C^1, E^2, C^2) \overset{\mathrm{def}}{=} \Big\{((f_0, f_1), (\boldsymbol{A}^1, E^1, C^1), (\boldsymbol{A}^2, E^2, C^2)) \Big|$$
$$\mathrm{isAtl}_{\mathrm{Ar}}^{C^k}(\boldsymbol{A}^1, E^1, C^1, \boldsymbol{A}^2, E^2, C^2, f_0, f_1)\Big\} \tag{18.5}$$

$$\mathcal{A}t\ell_{\mathrm{Ar}\ \text{full}}^{C^k}(E^1, C^1, E^2, C^2) \overset{\mathrm{def}}{=} \Big\{((f_0, f_1), (\boldsymbol{A}^1, E^1, C^1), (\boldsymbol{A}^2, E^2, C^2)) \Big|$$
$$\mathrm{isAtl}_{\mathrm{Ar}\ \text{full}}^{C^k}(\boldsymbol{A}^1, E^1, C^1, \boldsymbol{A}^2, E^2, C^2, f_0, f_1)\Big\} \tag{18.6}$$

$$\mathcal{A}t\ell_{\mathrm{Ar}\ \text{max}}^{C^k}(E^1, C^1, E^2, C^2) \overset{\mathrm{def}}{=} \Big\{((f_0, f_1), (\boldsymbol{A}^1, E^1, C^1), (\boldsymbol{A}^2, E^2, C^2)) \Big|$$
$$\mathrm{isAtl}_{\mathrm{Ar}\ \text{max}}^{C^k}(\boldsymbol{A}^1, E^1, C^1, \boldsymbol{A}^2, E^2, C^2, f_0, f_1)\Big\} \tag{18.7}$$

$$\mathcal{A}t\ell_{\mathrm{Ar}\ \text{S-max}}^{C^k}(E^1, C^1, E^2, C^2) \overset{\mathrm{def}}{=} \Big\{((f_0, f_1), (\boldsymbol{A}^1, E^1, C^1), (\boldsymbol{A}^2, E^2, C^2)) \Big|$$
$$\mathrm{isAtl}_{\mathrm{Ar}\ \text{S-max}}^{C^k}(\boldsymbol{A}^1, E^1, C^1, \boldsymbol{A}^2, E^2, C^2, f_0, f_1)\Big\} \tag{18.8}$$

$$\mathcal{A}t\ell^{C^k}_{\substack{Ar \\ \text{max-full}}}(E^1,C^1,E^2,C^2) \stackrel{\text{def}}{=} \left\{ ((f_0,f_1),(A^1,E^1,C^1),(A^2,E^2,C^2)) \,\middle|\, \right.$$
$$\left. \text{isAtl}^{C^k}_{\substack{Ar \\ \text{max-full}}}(A^1,E^1,C^1,A^2,E^2,C^2,f_0,f_1) \right\} \quad (18.9)$$

$$\mathcal{A}t\ell^{C^k}_{\substack{Ar \\ \text{S-max-full}}}(E^1,C^1,E^2,C^2) \stackrel{\text{def}}{=} \left\{ ((f_0,f_1),(A^1,E^1,C^1),(A^2,E^2,C^2)) \,\middle|\, \right.$$
$$\left. \text{isAtl}^{C^k}_{\substack{Ar \\ \text{S-max-full}}}(A^1,E^1,C^1,A^2,E^2,C^2,f_0,f_1) \right\} \quad (18.10)$$

**Definition 18.3** (Equivalence of $C^k$-atlas (near) morphisms). Let $E^i$, $i = 1,2$, be a topological space, $C^i$ a linear space, $A^i$ a $C^k$-atlases of $E^i$ in the coordinate space $C^i$ and $f \stackrel{\text{def}}{=} (f_0,f_1)$ and $g \stackrel{\text{def}}{=} (g_0,g_1)$ $E^1$-$E^2$ $C^k$ (near) morphisms of $A^1$ to $A^2$ in the coordinate spaces $C^1, C^2$.

$f$ is $C^k$-equivalent to $g$ iff $f_0 = g_0$.

### 18.1.2 Proclamations on $C^k$-atlas (near) morphisms

**Lemma 18.4** ($C^k$-atlas (near) morphisms). *Let $E$ be a set of topological spaces, $C$ a set of $C^k$ linear spaces, $E^i \in E$, $i = 1,2$, $C^i \in C$, $A^i$ a $C^k$-atlas of $E^i$ in the coordinate space $C^i$ and $f \stackrel{\text{def}}{=} (f_0\colon E^1 \longrightarrow E^2, f_1\colon C^1 \longrightarrow C^2)$ a pair of functions.*

*$f$ is an $E^1$-$E^2$ $C^k$ (near) morphism of $A^1$ to $A^2$ in the coordinate spaces $C^1, C^2$ iff $f$ is a strict $\underset{\mathcal{T}riv}{E}$-$\underset{\mathcal{T}riv}{E}$-$\underset{\text{triv}}{E^1}$-$\underset{\text{triv}}{E^2}$-$\underset{C^k-\text{op--triv}}{C}$-$\underset{C^k-\text{op--triv}}{C}$ m-atlas (near) morphism of $A^1$ to $A^2$ in the coordinate spaces $\underset{C^k-\textbf{triv}}{C^1}, \underset{C^k-\textbf{triv}}{C^2}$.*

*Proof.* The model neighborhoods of $\underset{C^k-\textbf{triv}}{C^i}$ are the open sets of $C^i$, $\underset{\mathcal{T}riv}{E} \overset{\text{full--cat}}{\subseteq} \underset{\mathcal{T}riv}{E}$, the morphisms of $\underset{\mathcal{T}riv}{E}$ are the continuous functions between spaces in $E$, $\underset{C^k-\mathcal{T}riv}{C} \overset{\text{full--cat}}{\subseteq} \underset{C^k-\mathcal{T}riv}{C}$, and the morphisms of $\underset{C^k-\mathcal{T}riv}{C}$ are the $C^k$ functions between spaces in $C$. □

*$f$ is an $E^1$-$E^2$ $C^k$-morphism of $A^1$ to $A^2$ in the coordinate spaces $C^1, C^2$ iff $f$ is a strict $E^1$-$E^2$-$\underset{C^k-\mathcal{T}riv}{C}$-$\underset{C^k-\mathcal{T}riv}{C}$ m-atlas morphism of $A^1$ to $A^2$ in the coordinate spaces $\underset{C^k-\textbf{triv}}{C^1}, \underset{C^k-\textbf{triv}}{C^2}$.*

*Proof.* The model neighborhoods of $\underset{C^k-\textbf{triv}}{C^i}$ are the open sets of $C^i$, $\underset{\mathcal{T}riv}{E} \overset{\text{full--cat}}{\subseteq} \underset{\mathcal{T}riv}{E}$, $\underset{C^k-\mathcal{T}riv}{C} \overset{\text{full--cat}}{\subseteq} \underset{C^k-\mathcal{T}riv}{C}$ and the morphisms of $\underset{C^k-\mathcal{T}riv}{C}$ are the $C^k$ functions between spaces in $C$. □

**Corollary 18.5** ($C^k$-atlas (near) morphisms)**.** *Let $E^i$, $i = 1, 2$, be a topological space, $C^i$ a linear space, $A^i$ a semi-maximal $C^k$-atlas of $E^i$ in the coordinate space $C^i$ and $\boldsymbol{f} \stackrel{\text{def}}{=} (f_0, f_1)$ an $E^1$-$E^2$ $C^k$ near morphism of $\boldsymbol{A}^1$ to $\boldsymbol{A}^2$ in the coordinate spaces $C^1$, $C^2$. Then $\boldsymbol{f}$ is an $E^1$-$E^2$ $C^k$ morphism of $\boldsymbol{A}^1$ to $\boldsymbol{A}^2$ in the coordinate spaces $C^1$, $C^2$.*

*Proof.* The result follows from lemma 11.12 (M-atlas (near) morphisms) on page 60
. □

Let $E^i$, $i = 1, 2, 3$, be a topological space, $C^i$ a linear space, $A^i$ a $C^k$-atlas of $E^i$ in the coordinate space $C^i$ and $(f_0^i, f_1^i)$ an $E^i$-$E^{i+1}$ $C^k$-morphism of $\boldsymbol{A}^i$ to $\boldsymbol{A}^{i+1}$ in the coordinate spaces $C^i$, $C^{i+1}$. Then $(f_0^2 \circ f_0^1, f_1^2 \circ f_1^1)$ is a $E^1$-$E^3$ $C^k$-morphism of $\boldsymbol{A}^1$ to $\boldsymbol{A}^3$ in the coordinate spaces $C^1$, $C^3$.

*Proof.* The result follows from item 2 of lemma 11.16 (Composition of m-atlas (near) morphisms) on page 72. □

**Lemma 18.6** (Composition of equivalent $C^k$-atlas (near) morphisms)**.** *Let $E^i$, $i = 1, 2, 3$, be a topological space, $C^i$ be a $C^k$ linear space, $A^i$ be a $C^k$-atlas of $E^i$ in the coordinate space $C^i$ and $\boldsymbol{f}^i \stackrel{\text{def}}{=} (f_0^i \colon \boldsymbol{E}^i \longrightarrow \boldsymbol{E}^{i+1}, f_1^i \colon \boldsymbol{C}^i \longrightarrow \boldsymbol{C}^{i+1})$ and $\boldsymbol{g}^i \stackrel{\text{def}}{=} (g_0^i \colon \boldsymbol{E}^i \longrightarrow \boldsymbol{E}^{i+1}, g_1^i \colon \boldsymbol{C}^i \longrightarrow \boldsymbol{C}^{i+1})$ be $C^k$-equivalent $E^i$-$E^{i+1}$ m-atlas (near) morphisms of $\boldsymbol{A}^i$ to $\boldsymbol{A}^{i+1}$ in the coordinate spaces $C^i$, $C^{i+1}$.*

*Then $\boldsymbol{f}^2 \stackrel{()}{\circ} f^1$ is $C^k$-equivalent to $\boldsymbol{g}^2 \stackrel{()}{\circ} g^1$.*

*Proof.* $f_0^1 = g_0^1$ and $f_0^2 = g_0^2$, hence $f_0^2 \circ f_0^1 = g_0^2 \circ g_0^1$. □

## 18.2 Categories of $C^k$ atlases and functors

### 18.2.1 Categories of $C^k$ atlases

This subsubsection defines categories of $C^k$-atlases with $C^k$-atlas morphisms as morphisms. It does not define categories of $C^k$-atlases with $C^k$-atlas near morphisms as morphisms becaue the composition of two $C^k$-atlas near morphisms is not in general a $C^k$-atlas near morphism, and requiring the atlasses to be semi-maximal would cause all $C^k$-atlas near morphisms to be $C^k$-atlas morphisms.

**Definition 18.7** (Categories of $C^k$ atlases)**.** Let $\boldsymbol{E}$ be a set of topological spaces and $\boldsymbol{C}$ a set of linear spaces. Let $\boldsymbol{P} \stackrel{\text{def}}{=} \boldsymbol{E} \times \boldsymbol{C}$. Then

$$\mathcal{A}t\ell^{C^k}_{\text{Ar}}(\boldsymbol{E}, \boldsymbol{C})_{\text{(full,S-max,max,S-max-full,max-full)}} \stackrel{\text{def}}{=} \bigcup_{\substack{(E^\mu, C^\mu) \in \boldsymbol{P} \\ (E^\nu, C^\nu) \in \boldsymbol{P}}} \mathcal{A}t\ell^{C^k}_{\text{Ar}}(E^\mu, C^\mu, E^\nu, C^\nu)_{\text{(full,S-max,max,S-max-full,max-full)}} \tag{18.11}$$

$$\mathcal{A}t\ell^{C^k}(\boldsymbol{E}, \boldsymbol{C})_{\text{(full,S-max,max,S-max-full,max-full)}} \stackrel{\text{def}}{=} \left( \mathcal{A}t\ell^{C^k}_{\text{Ob}}(\boldsymbol{E}, \boldsymbol{C})_{\text{(full,S-max,max,S-max-full,max-full)}}, \mathcal{A}t\ell^{C^k}_{\text{Ar}}(\boldsymbol{E}, \boldsymbol{C})_{\text{(full,S-max,max,S-max-full,max-full)}}, \stackrel{A}{\circ} \right) \tag{18.12}$$

*Remark* 18.8. It is pointless to define categories of near morphisms, since
$$is\mathcal{A}t\ell_{Ar\ \text{S-max (max)}}^{C^k,near}(E^\mu, C^\mu, E^\nu, C^\nu) \iff is\mathcal{A}t\ell_{Ar\ \text{S-max (max)}}^{C^k}(E^\mu, C^\mu, E^\nu, C^\nu).$$

**Lemma 18.9** ($\mathcal{A}t\ell^{C^k}(E, C)$ is a category). *Let $E$ be a set of topological spaces and $C$ a set of linear spaces. Then each of $\mathcal{A}t\ell_{\text{(full,S-max,max,S-max-full,max-full)}}^{C^k}(E, C)$ is a category.*

*Let $(A^i, E^i, C^i) \in \mathcal{A}t\ell_{Ob}^{C^k}(E, C)$. Then $\text{Id}_{(A^i, E^i, C^i)}$ is the identity morphism for $(A^i, E^i, C^i)$.*

*Proof.* Let $(A^i, E^i, C^i)$, $i \in [1, 3]$, be an object of $\mathcal{A}t\ell^{C^k}(E, C)$ and let $m^i \stackrel{\text{def}}{=} ((f_0^i, f_1^i), (A^i, E^i, C^i), (A^{i+1}, E^{i+1}, C^{i+1}))$, $i = 1, 2$, be a morphism of $\mathcal{A}t\ell^{C^k}(E, C)$. Then

**Composition:** $((f_0^2 \circ f_0^1, f_1^2 \circ f_1^1), (A^1, E^1, C^1), (A^3, E^3, C^3))$ is a morphism of $\mathcal{A}t\ell^{C^k}(E, C)$ by corollary 18.5 ($C^k$-atlas (near) morphisms) on page 124.

**Associativity:** Composition is associative by lemma 1.19 (Tuple composition for labeled morphisms) on page 12.

**Identity:** $\text{Id}_{(A^i, E^i, C^i)}$ is an identity morphism by lemma 1.19.

□

### 18.2.2 $C^k$ atlas functors

**Definition 18.10** (Functors from $C^k$ atlases to m-atlases). Let $E^i$, $i = 1, 2$, be a topological space, $C^i$ a linear space, $A^i$ $C^k$-atlases of $E^i$ in the coordinate space $C^i$, $f_0: E^1 \longrightarrow E^2$ continuous and $f_1: C^1 \longrightarrow C^2$ $C^k$. Then

$$\mathcal{F}_{C^k,M}^{C^k}(A^i, E^i, C^i) \stackrel{\text{def}}{=} \left(A^i, E^i_{\text{triv}}, C^i_{C^k-\text{triv}}\right) \tag{18.13}$$

$$\mathcal{F}_{C^k,M}^{C^k}((f_0, f_1), (A^1, E^1, C^1), (A^2, E^2, C^2)) \stackrel{\text{def}}{=}$$
$$\left((f_0, f_1), \left(A^1, E^1_{\text{triv}}, C^1_{C^k-\text{triv}}\right), \left(A^2, E^2_{\text{triv}}, C^2_{C^k-\text{triv}}\right)\right) \tag{18.14}$$

**Theorem 18.11** (Functors from $C^k$ atlases to m-atlases). *Let $E$ be a set of topological spaces and $C$ a set of linear spaces. Then $\mathcal{F}_{C^k,M}^{C^k}$ is a functor from $\mathcal{A}t\ell^{C^k}(E, C)$ to $\mathcal{A}t\ell(E_{\text{triv}}, C_{C^k-\text{triv}})$*

*Proof.* Let $o^i \stackrel{\text{def}}{=} (A^i, E^i, C^i)$, $i \in [1,3]$, be an object of $\mathcal{A}t\ell^{C^k}(E, C)$,
$o'^i \stackrel{\text{def}}{=} \mathcal{F}_{C^k,M}^{C^k} o^i = (A^i, \underset{\text{triv}}{E^i}, \underset{C^k-\text{triv}}{C^i})$, $m^i \stackrel{\text{def}}{=} ((f_0^i, f_1^i), o^i, o^{i+1})$, $i = 1, 2$, be a morphism from $o^i$ to $o^{i+1}$ and $m'^i \stackrel{\text{def}}{=} \mathcal{F}_{C^k,M}^{C^k} m^i = ((f_0^i, f_1^i).o'^i, o'^{i+1})$, be the corresponding morphism in $\mathcal{A}t\ell(\underset{\text{triv}}{E}, \underset{C^k-\text{triv}}{C})$.

**Preservation of endpoints:**

$$\mathcal{F}_{C^k,M}^{C^k} o^i = (A^i, \underset{\text{triv}}{E^i}, \underset{C^k-\text{triv}}{C^i})$$
$$= o'^i$$
$$\mathcal{F}_{C^k,M}^{C^k} m^i = m'^i$$
$$= ((f_0^i, f_1^i), o'^i, o'^{i+1})$$

$\mathcal{F}_{C^k,M}^{C^k} m^i$ is a morphism from $\mathcal{F}_{C^k,M}^{C^k} o^i$ to $\mathcal{F}_{C^k,M}^{C^k} o^{i+1}$:

**Identity:**

$$\mathcal{F}_{C^k,M}^{C^k} \text{Id}_{o^i} = \mathcal{F}_{C^k,M}^{C^k}((\text{Id}_{E^1}, \text{Id}_{C^1}), o^i, o^i)$$
$$= ((\text{Id}_{\underset{\text{triv}}{E^i}}, \text{Id}_{C^i}), o'^i, o'^i)$$
$$= ((\text{Id}_{\underset{\text{triv}}{E^i}}, \text{Id}_{C^i}), \mathcal{F}_{C^k,M}^{C^k} o^i, \mathcal{F}_{C^k,M}^{C^k} o^i)$$
$$= \text{Id}_{\mathcal{F}_{C^k,M}^{C^k} o^i}$$

**Composition:**

$$m^2 \overset{A}{\circ} m^1 = ((f_0^2 \circ f_0^1, f_1^2 \circ f_1^1), o^1, o^3)$$
$$\mathcal{F}_{C^k,M}^{C^k}(m^2 \overset{A}{\circ} m^1) = ((f_0^2 \circ f_0^1, f_1^2 \circ f_1^1), o'^1, o'^3)$$
$$= ((f_0^2, f_1^2), o'^2, o'^3) \overset{A}{\circ} ((f_0^1, f_1^1), o'^1, o'^2)$$
$$= \mathcal{F}_{C^k,M}^{C^k}((f_0^2, f_1^2), o^2, o^3) \overset{A}{\circ} \mathcal{F}_{C^k,M}^{C^k}((f_0^1, f_1^1), o^1, o^2)$$
$$= \mathcal{F}_{C^k,M}^{C^k} m^2 \overset{A}{\circ} \mathcal{F}_{C^k,M}^{C^k} m^1$$

□

**Definition 18.12** (Functors from m-atlases to $C^k$ atlases). Let $E^i$, $i = 1, 2$, be model spaces, $C^i$ be linear model spaces, $A^i$ a maximal m-atlas of $E^i$ in the coordinate space $C^i$ and $(f_0, f_1)$ an $E^1$-$E^2$ m-atlas morphism of $A^1$ to $A^2$ in the coordinate spaces $C^1, C^2$, Then

$$\mathcal{F}_{M,C^k}^{C^k}(A^i, E^i, C^i) \stackrel{\text{def}}{=} (A^i, \pi_1(E^i), \pi_1(C^i)) \tag{18.15}$$

$$\mathcal{F}_{M,C^k}^{C^k}((f_0, f_1), (A^1, E^1, C^1), (A^2, E^2, C^2)) \stackrel{\text{def}}{=}$$
$$\left((f_0, f_1), (A^1, \pi_1(E^1), \pi_1(C^1)), (A^2, \pi_1(E^2), \pi_1(C^2))\right) \quad (18.16)$$

**Theorem 18.13** (Functors from m-atlases to $C^k$ atlases). *Let $E$ be a set of model spaces and $C$ a set of $C^k$ linear model spaces. Then $\mathcal{F}_{M,C^k}^{C^k}$ is a functor from $\mathcal{A}t\ell(E, C)$ to $\mathcal{A}t\ell^{C^k}(\pi_1[E], \pi_1[C])$.*

*Proof.* Let $o^i \stackrel{\text{def}}{=} (A^i, (E^i, \mathcal{E}^i), (C^i, \mathcal{C}^i))$, $i \in [1, 3]$, be an object of $\mathcal{A}t\ell^{C^k}(E, C)$ and $m^i \stackrel{\text{def}}{=} ((f_0^i, f_1^i).o^i, o^{i+1})$, $i = 1, 2$, be a morphism from $o^i$ to $o^{i+1}$.

$\mathcal{F}_{M,C^k}^{C^k}(m^1)$ is a morphism from $\mathcal{F}_{M,C^k}^{C^k} o^1$ to $\mathcal{F}_{M,C^k}^{C^k} o^2$:

$$\mathcal{F}_{M,C^k}^{C^k}(m^1) =$$
$$\mathcal{F}_{M,C^k}^{C^k}((f_0^1, f_1^1), (A^1, E^1, C^1), (A^2, E^2, C^2)) =$$
$$\left((f_0^1, f_1^1), (A^1, \pi_1(E^1), \pi_1(C^1)), (A^2, \pi_1(E^2), \pi_1(C^2))\right)$$

$\mathcal{F}_{M,C^k}^{C^k}$ maps identity functions to identity functions:

$$\mathcal{F}_{M,C^k}^{C^k} \text{Id}_{(A^i, E^i, C^i)} =$$
$$\mathcal{F}_{M,C^k}^{C^k}((\text{Id}_{E^i}, \text{Id}_{C^i}), (A^i, E^i, C^i), (A^i, E^i, C^i)) =$$
$$((\text{Id}_{\pi_1(E^i)}, \text{Id}_{\pi_1(C^i)}), (A^i, \pi_1(E^i), \pi_1(C^i)), (A^i, \pi_1(E^i), \pi_1(C^i))) =$$
$$((\text{Id}_{E^i}, \text{Id}_{C^i}), \mathcal{F}_{M,C^k}^{C^k}(A^i, E^i, C^i), \mathcal{F}_{M,C^k}^{C^k}(A^i, E^i, C^i)) =$$
$$\text{Id}_{\mathcal{F}_{M,C^k}^{C^k}(A^i, E^i, C^i)}$$

$\mathcal{F}_{M,C^k}^{C^k} m^2 \stackrel{A}{\circ} \mathcal{F}_{M,C^k}^{C^k} m^1 = \mathcal{F}_{M,C^k}^{C^k}(m^2 \stackrel{A}{\circ} m^1)$:

$$m^2 \stackrel{A}{\circ} m^1 = ((f_0^2 \circ f_0^1, f_1^2 \circ f_1^1), (A^1, E^1, C^1), (A^3, E^3, C^3))$$
$$\mathcal{F}_{M,C^k}^{C^k}((A^i, E^i, C^i)) = (A^i, \pi_1(E^i), \pi_1(C^i))$$
$$\mathcal{F}_{M,C^k}^{C^k}(m^i) = \left((f_0^i, f_1^i), (A^i, \pi_1(E^i), \pi_1(C^i)), (A^{i+1}, \pi_1(E^{i+1}), \pi_1(C^{i+1}))\right)$$
$$\mathcal{F}_{M,C^k}^{C^k} m^2 \stackrel{A}{\circ} \mathcal{F}_{M,C^k}^{C^k} m^1 = \left((f_0^2 \circ f_0^1, f_1^2 \circ f_1^2), (A^1, \pi_1(E^1), \pi_1(C^1)), (A^3, \pi_1(E^3), \pi_1(C^3))\right)$$
$$\mathcal{F}_{M,C^k}^{C^k}(m^2 \circ m^1) = \left((f_0^2 \circ f_0^1, f_1^2 \circ f_1^2), (A^1, \pi_1(E^1), \pi_1(C^1)), (A^3, \pi_1(E^3), \pi_1(C^3))\right)$$

□

# 19 Associated model spaces and functors

$\mathcal{F}_{C^k,M}^{C^k}$ is an obvious functor of $C^k$-atlases to m-atlases, but $\underset{\mathbf{triv}}{E}$ and $\underset{C^k-\mathbf{triv}}{C}$ may have more model neighborhoods or more morphisms than needed for consistency with the atlas. There exist, however, model spaces with the minimum model neighborhoods and morphisms needed.

**Definition 19.1** (Coordinate model spaces associated with $C^k$-atlases). Let $A^i$, $i = 1, 2$, be a $C^k$-atlas of $E^i$ in the coordinate space $C^i$, $f_0 \colon E^1 \longrightarrow E^2$ a continuous function and $f_1 \colon C^1 \longrightarrow C^2$ a $C^k$ function. Then

$$\mathcal{F}_2^{\min\,C^k}(A^i, E^i, C^i) \stackrel{\text{def}}{=}$$

$$\underset{\min}{\mathrm{Mod}}\left(C^i, \pi_2[A^i], \left\{\phi' \circ \phi^{-1} \,\middle|\, \left(\exists_{\substack{(U,V,\phi) \in A^i \\ (U',V',\phi') \in A^i}}\right) U \cap U' \neq \emptyset\right\}\right) \quad (19.1)$$

$$\mathcal{F}_2^{\min\,C^k}\bigl((f_0, f_1), (A^1, E^1, C^1), (A^2, E^2, C^2)\bigr) \stackrel{\text{def}}{=}$$

$$f_1 \colon \mathcal{F}_2^{\min\,C^k}(A^1, E^1, C^1) \longrightarrow \mathcal{F}_2^{\min\,C^k}(A^2, E^2, C^2) \quad (19.2)$$

The minimal coordinate $C^k$ model space with neighborhoods in the $C^k$-atlas $A^i$ of $E^i$ in the coordinate space $C^i$ is $\mathcal{F}_2^{\min\,C^k}(A^1, E^i, C^i)$.

The coordinate mapping associated with the $E^1$-$E^2$ $C^k$-atlas morphism $(f_0, f_1)$ of $A^1$ to $A^2$ in the coordinate spaces $C^1$, $C^2$ is $f_1 \colon \mathcal{F}_2^{\min\,C^k}(A^1, E^1, C^1) \longrightarrow \mathcal{F}_2^{\min\,C^k}(A^2, E^2, C^2)$. If it is a model function then it is also the coordinate m-atlas morphism associated with the $E^1$-$E^2$ $C^k$-atlas morphism $(f_0, f_1)$ of $A^1$ to $A^2$ in the coordinate spaces $C^1$, $C^2$.

**Lemma 19.2** (Coordinate model spaces associated with $C^k$-atlases). *Let $A^i$ be a $C^k$-atlas of $E^i$ in the coordinate space $C^i$, $i = 1, 2$. Then $\mathcal{F}_2^{\min\,C^k}(A^i, E^i, C^i)$ is a $C^k$ linear model space.*

*Proof.* $\mathcal{F}_2^{\min\,C^k}(A^i, E^i, C^i)$ satisfies the conditions for a model space.

1. Since $\pi_2[A^i]$ is an open cover of $\bigcup \pi_2[A^i]$, the set of finite intersections is also an open cover.

2. Finite intersections of finite intersections are finite intersections

3. Restrictions of continuous functions are continuous

4. If $f\colon A \longrightarrow B$ is a morphism of $\mathcal{F}_2^{\min\;C^k}(\boldsymbol{A}^i, E^i, C^i)$, $A, A', B, B'$ model meighborhoods of $\mathcal{F}_2^{\min\;C^k}(\boldsymbol{A}^i, E^i, C^i)$, $A' \subseteq A$, $B' \subseteq B$ and $f[A'] \subseteq B'$ then since $f\colon A \longrightarrow B$ is a morphism it is a restriction of a transition function between its restrictions to sets in $\pi_2[\boldsymbol{A}^i]$ and its restrictions are also, hence morphisms, and thus $f\!\restriction_{A'}\colon A' \longrightarrow B'$ is a morphism.

5. If $(U, V, \phi) \in \boldsymbol{A}^i$ then $\mathrm{Id}_V = \phi \circ \phi^{-1} A$ is a transition function and hence a morphism of $\mathcal{F}_2^{\min\;C^k}(\boldsymbol{A}^i, E^i, C^i)$. If $A, A'$ are objects of $\pi_2(\mathcal{F}_2^{\min\;C^k}(\boldsymbol{A}^i, E^i, C^i))$ and $A' \subseteq A$ then the inclusion map $i\colon A' \hookrightarrow A$ is a restriction of an identity morphism of $\mathcal{F}_2^{\min\;C^k}(\boldsymbol{A}^i, E^i, C^i)$ and hence a morphism.

6. Restricted sheaf condition: let

   (a) $U_\alpha, V_\alpha, \alpha \prec A$, be an object of $\pi_2\!\left(\mathcal{F}_2^{\min\;C^k}(\boldsymbol{A}^i, E^i, C^i)\right)$

   (b) $f_\alpha\colon U_\alpha \longrightarrow V_\alpha$ be a morphism of $\pi_2\!\left(\mathcal{F}_2^{\min\;C^k}(\boldsymbol{A}^i, E^i, C^i)\right)$

   (c) $U \stackrel{\text{def}}{=} \bigcup_{\alpha \prec A} U_\alpha$

   (d) $V \stackrel{\text{def}}{=} \bigcup_{\alpha \prec A} V_\alpha$

   (e) $f\colon U \longrightarrow V$ be continuous and $\left(\forall_{\substack{\alpha \prec A \\ x \in U_\alpha}}\right) f(x) = f_\alpha(x)$

   Then $f$ is $C^k$ and hence a morphism of $\pi_2\!\left(\mathcal{F}_2^{\min\;C^k}(\boldsymbol{A}^i, E^i, C^i)\right)$.

$\square$

$\mathcal{F}_2^{\min\;C^k}(\boldsymbol{A}^i, E^i, C^i)$ *is non-degenerate iff $\boldsymbol{A}^i$ is non-degenerate.*

*Proof.* If $\boldsymbol{A}^i$ is non-degenerate then by definition 17.1 ($C^k$-atlases) on page 117 there is at least one chart $(U, V, \phi) \in \boldsymbol{A}$ such that $V$ contains a ball of the underlying Banach or Fréchet space. $U \subseteq \bigcup \pi_2[\boldsymbol{A}^i]$, so by definition 19.1 (Coordinate model spaces associated with $C^k$-atlases) on page 128 that ball is contained in $\mathcal{F}_2^{\min\;C^k}(\boldsymbol{A}^i, E^i, C^i)$.

Conversly, if $\mathcal{F}_2^{\min\;C^k}(\boldsymbol{A}^i, E^i, C^i)$ is non-degenerate then by definition 14.1 (Linear spaces) on page 111 it contains a ball $B$ at $v \in \mathcal{F}_2^{\min\;C^k}(\boldsymbol{A}^i, E^i, C^i)$ with radius $r$. By definition 19.1 there is a chart $(U, V, \phi) \in \boldsymbol{A}^i$ at $v$. Since $V$ is open in $C^i$, there is a ball $B'$ at $v$ with rafius $r' > 0$ such that $B' \cap C^i \subseteq V$. Then the ball at $v$ with radius $\min(r, r')$ is contained in $V$. $\square$

Let $A^i$, $i = 1, 2$, be a $C^k$-atlas of $E^i$ in the coordinate space $C^i$, $f_0\colon E^1 \longrightarrow E^2$ a continuous function, $f_1\colon C^1 \longrightarrow C^2$ a $C^k$ function and $f \stackrel{\text{def}}{=} (f_0, f_1)$ either a semi-maximal $C^k$-atlas near morphism or a $C^k$-atlas morphism from $A^1$ to $A^2$. Then $f_1\colon \mathcal{F}_2^{\min\,C^k}(A^1, E^1, C^1) \longrightarrow \mathcal{F}_2^{\min\,C^k}(A^2, E^2, C^2)$ is well defined, i.e., $f_1\left[\mathcal{F}_2^{\min\,C^k}(A^1, E^1, C^1)\right] \subseteq \mathcal{F}_2^{\min\,C^k}(A^2, E^2, C^2)$.

*Proof.* $f$ is a morphism either by corollary 18.5 ($C^k$-atlas (near) morphisms) on page 124 or by hypothesis. By definition 18.2 ($C^k$-atlas morphisms) on page 121 there exists a subchart $(U'^1, V'^1, \phi^1\colon U'^1 \xrightarrow{\tilde{=}} V'^1) \in A^1$ at $u^1$, an open set $U'^2 \subseteq U^2$ and a chart $(U'^2, \hat{V}'^2, \phi'^2\colon U'^2 \xrightarrow{\tilde{=}} \hat{V}'^2) \in A^2$ at $u^2$ such that $f_0[U'^1] \subseteq U'^2$ and diagram (11.8) in definition 11.4 (M-atlas morphisms for model spaces) on page 53 is commutative, as shown in fig. 6 (Completed m-atlas morphism) on page 54. Then $f_1(v^1) \in \hat{V}'^2$. □

**Definition 19.3** (Model spaces associated with $C^k$-atlases). Let $A^i$, $i = 1, 2$, be a $C^k$-atlas of $E^i$ in the coordinate space $C^i$ and $f \stackrel{\text{def}}{=} (f_0, f_1)$ an $E^1$-$E^2$ $C^k$-atlas morphism of $A^1$ to $A^2$ in the coordinate spaces $C^1$, $C^2$. Then:

The minimal $C^k$ model space with neighborhoods in the $C^k$-atlas $A^i$ of $E^i$ in the coordinate space $C^i$ is

$$\mathcal{F}_1^{\min\,C^k}(A^i, E^i, C^i) \stackrel{\text{def}}{=}$$
$$\operatorname*{Mod}_{\min}\left(E^i, \pi_1[A^i], \left\{ \phi'^{-1} \circ \phi \;\middle|\; \left(\exists_{\substack{(U,V,\phi) \in A^i \\ (U',V',\phi') \in A^i}}\right) U \cap U' \neq \emptyset \right\}\right) \quad (19.3)$$

The mapping associated with the $E^1$-$E^2$ $C^k$-atlas morphism $(f_0, f_1)$ of $A^1$ to $A^2$ in the coordinate spaces $C^1$, $C^2$ is

$$\mathcal{F}_1^{\min\,C^k}\bigl((f_0, f_1), (A^1, E^1, C^1), (A^2, E^2, C^2)\bigr) \stackrel{\text{def}}{=}$$
$$f_0\colon \mathcal{F}_1^{\min\,C^k}(A^1, E^1, C^1) \longrightarrow \mathcal{F}_1^{\min\,C^k}(A^2, E^2, C^2) \quad (19.4)$$

If it is a model function then the it is also the m-atlas morphism associated with the $E^1$-$E^2$ $C^k$-atlas morphism $(f_0, f_1)$ of $A^1$ to $A^2$ in the coordinate spaces $C^1$, $C^2$.

**Lemma 19.4** (Model spaces associated with $C^k$-atlases). *Let $A$ be a $C^k$-atlas of $E$ in the coordinate space $C$. Then $\mathcal{F}_1^{\min\,C^k}(A, E, C)$ is a model space.*

*Proof.* The result follows from Lemma 5.4 (Minimal model spaces are model spaces) on page 27. □

**Theorem 19.5** (Functors from $C^k$ atlases to model spaces). *Let $E$ be a set of topological spaces and $C$ a set of linear spaces. Then $\mathcal{F}_1^{\min^{C^k}}$ is a functor from $\mathcal{A}t\ell^{C^k}(E, C)$ to $\underset{\mathcal{T}riv}{E}$, $\mathcal{F}_1^{\min^{C^k}}$ is a functor from $\mathcal{A}t\ell^{C^k}_{\text{full}}(E, C)$ to $\underset{\mathcal{T}riv}{E}$, $\mathcal{F}_2^{\min^{C^k}}$ is a functor from $\mathcal{A}t\ell^{C^k}_{\text{S-max}}(E, C)$ to $\underset{C^k-\text{op-triv}}{C}$ and $\mathcal{F}_2^{\min^{C^k}}$ is a functor from $\mathcal{A}t\ell^{C^k}_{\text{full}}(E, C)$ to $\underset{C^k-\textbf{triv}}{C}$.*

*Proof.* Let $o^i \stackrel{\text{def}}{=} (A^i, E^i, C^i)$, $i \in [1, 3]$, be objects in $\mathcal{A}t\ell^{C^k}(E, C)$ and let $m^i \stackrel{\text{def}}{=} ((f_0^i, f_1^i), o^i, o^{i+1})$ be a morphism in $\mathcal{A}t\ell^{C^k}(E, C)$.

$$\mathcal{F}_1^{\min^{C^k}}: \mathcal{A}t\ell^{C^k}(E, C) \longrightarrow \underset{\textbf{triv}}{E}:$$

**Preservation of endpoints:** $\mathcal{F}_1^{\min^{C^k}}(m^i) = f_0^i: \mathcal{F}_1^{\min^{C^k}} o^i \longrightarrow \mathcal{F}_1^{\min^{C^k}} o^{i+1}$

**Composition:**

$$\mathcal{F}_1^{\min^{C^k}}(m^2 \stackrel{A}{\circ} m^1) = \mathcal{F}_1^{\min^{C^k}}((f_0^2 \circ f_0^1, f_1^2 \circ f_1^1)(A^1, E^1, C^1), (A^3, E^3, C^3))$$
$$= f_0^2 \circ f_0^1: \mathcal{F}_1^{\min^{C^k}} o^1 \longrightarrow \mathcal{F}_1^{\min^{C^k}} o^3$$
$$= (f_0^2: \mathcal{F}_1^{\min^{C^k}} o^2 \longrightarrow \mathcal{F}_1^{\min^{C^k}} o^3) \circ (f_0^1: \mathcal{F}_1^{\min^{C^k}} o^1 \longrightarrow \mathcal{F}_1^{\min^{C^k}} o^2)$$
$$= \mathcal{F}_1^{\min^{C^k}}((f_0^2, f_1^2)(A^2, E^2, C^2), (A^3, E^3, C^3)) \circ$$
$$\mathcal{F}_1^{\min^{C^k}}((f_0^1, f_1^1)(A^1, E^1, C^1), (A^2, E^2, C^2))$$
$$= \mathcal{F}_1^{\min^{C^k}}(m^2) \circ \mathcal{F}_1^{\min^{C^k}}(m^1)$$

**Identity:**

1. $\mathcal{F}_1^{\min^{C^k}}(\text{Id}_{o^i}) = \mathcal{F}_1^{\min^{C^k}}((\text{Id}_{E^i}, \text{Id}_{C^i}), (A^i, E^i, C^i), (A^i, E^i, C^i)) = \text{Id}_{E^i}: \mathcal{F}_1^{\min^{C^k}} o^i \longrightarrow \mathcal{F}_1^{\min^{C^k}} o^i$

2. $\text{Id}_{\mathcal{F}_1^{\min^{C^k}}} o^i = \text{Id}_{E^i}: \mathcal{F}_1^{\min^{C^k}} o^i \longrightarrow \mathcal{F}_1^{\min^{C^k}} o^i$

   The proof for $\mathcal{F}_1^{\min^{C^k}}: \mathcal{A}t\ell^{C^k}_{\text{full}}(E, C) \longrightarrow \underset{\textbf{triv}}{E}$ is identical.

$$\mathcal{F}_2^{\min^{C^k}}: \mathcal{A}t\ell^{C^k}(E, C) \longrightarrow \underset{C^k-\text{op-triv}}{C}:$$

**Preservation of endpoints:**

$$\mathcal{F}_2^{\mathrm{min}\,C^k}(m^i) = f_1^i \colon \mathcal{F}_2^{\mathrm{min}\,C^k} o^i \longrightarrow \mathcal{F}_2^{\mathrm{min}\,C^k} o^{i+1}$$

$$\mathcal{F}(g \circ f) = \mathcal{F}(g) \circ \mathcal{F}(f) \colon \mathcal{F}_2^{\mathrm{min}\,C^k}(m^2 \overset{A}{\circ} m^1) =$$
$$\mathcal{F}_2^{\mathrm{min}\,C^k}((f_0^2 \circ f_0^1, f_1^2 \circ f_1^1)(A^1, E^1, C^1), (A^3, E^3, C^3)) =$$
$$f_0^2 \circ f_0^1 \colon \mathcal{F}_2^{\mathrm{min}\,C^k} o^1 \longrightarrow \mathcal{F}_2^{\mathrm{min}\,C^k} o^3 =$$
$$(f_1^2 \colon \mathcal{F}_2^{\mathrm{min}\,C^k} o^2 \longrightarrow \mathcal{F}_2^{\mathrm{min}\,C^k} o^3) \circ (f_1^1 \colon \mathcal{F}_2^{\mathrm{min}\,C^k} o^1 \longrightarrow \mathcal{F}_2^{\mathrm{min}\,C^k} o^2) =$$
$$\mathcal{F}_2^{\mathrm{min}\,C^k}((f_0^2, f_1^2)(A^2, E^2, C^2), (A^3, E^3, C^3)) \circ$$
$$\mathcal{F}_2^{\mathrm{min}\,C^k}((f_0^1, f_1^1)(A^1, E^1, C^1), (A^2, E^2, C^2)) =$$
$$\mathcal{F}_2^{\mathrm{min}\,C^k}(m^2) \circ \mathcal{F}_2^{\mathrm{min}\,C^k}(m^1)$$

$\mathcal{F}(\mathrm{Id}_A) = \mathrm{Id}_{\mathcal{F}(A)}$:

1. $\mathcal{F}_2^{\mathrm{min}\,C^k}(\mathrm{Id}_{o^i}) = \mathcal{F}_2^{\mathrm{min}\,C^k}((\mathrm{Id}_{E^i}, \mathrm{Id}_{C^i}), (A^i, E^i, C^i), (A^i, E^i, C^i)) =$
   $\mathrm{Id}_{C^i} \colon \mathcal{F}_2^{\mathrm{min}\,C^k} o^i \longrightarrow \mathcal{F}_2^{\mathrm{min}\,C^k} o^i$

2. $\mathrm{Id}_{\mathcal{F}_2^{\mathrm{min}\,C^k} o^i} = \mathrm{Id}_{C^i} \colon \mathcal{F}_2^{\mathrm{min}\,C^k} o^i \longrightarrow \mathcal{F}_2^{\mathrm{min}\,C^k} o^i$

The proof for $\mathcal{F}_2^{\mathrm{min}\,C^k} \colon \mathcal{A}t\ell_{\mathrm{full}}^{C^k}(E, C) \longrightarrow \underset{C^k-\mathbf{triv}}{C}$ is identical. $\square$

# 20 Classic $C^k$-atlas morphisms and functors

This section introduces an alternate definition of and taxonomy for morphisms between $C^k$-atlases, defines categories of $C^k$-atlases, defines functors amomg them and proves some basic reasults. It introduces the notion of classic $C^k$-atlas morphisms, which are equivalent to the conventional definitions for a manifold.

## 20.1 Classic $C^k$-atlas morphisms

This subsection introduces the notion of classic $C^k$-atlas morphisms, and proves some basic results.

### 20.1.1 Definitions of classic $C^k$-atlas morphisms

**Definition 20.1** (Classic $C^k$-atlas morphisms). Let $E^i$, $i = 1, 2$, be a topological space, $C^i$ be a linear space, $A^i$ be a $C^k$-atlases of $E^i$ in the coordinate space $C^i$ and $f$ be a continuous function.

$f$ is a $E^1$-$E^2$ classic $C^k$ morphism of $A^1$ to $A^2$ in the coordinate spaces $C^1, C^2$, abbreviated as isAtl$_{\text{Ar}}^{\text{classic},C^k}(A^1, E^1, C^1, A^2, E^2, C^2, f)$, iff for any $(U^i, V^i, \phi^i\colon U^i \overset{\sim}{=}\!\!\twoheadrightarrow V^i) \in A^i$, $i = 1, 2$, with $I \overset{\text{def}}{=} U^1 \cap f^{-1}[U^2] \neq \emptyset$, $\phi^2 \circ f \circ \phi^{1-1}\colon \phi^1[I] \longrightarrow V^2$ is a $C^k$ function.

*Remark* 20.2. Definitions of constrained, semistrict and strict classic $C^k$-atlas morphisms would be pointless, as any classic $C^k$-atlas morphism would be constrained and strict.

### 20.1.2 Proclamations on classic $C^k$-atlas morphisms

**Lemma 20.3** (Classic $C^k$-atlas morphisms). *Let $E$ be a set of topological spaces, $C$ a set of $C^k$ linear spaces, $E^i \in E$, $i = 1, 2$, $C^i \in C$, $A^i$ a $C^k$-atlas of $E^i$ in the coordinate space $C^i$ and $f\colon E^1 \longrightarrow E^2$ a function.*

*$f$ is an $E^1$-$E^2$ classic $C^k$-atlas near morphism of $A^1$ to $A^2$ in the coordinate spaces $C^1$, $C^2$ iff $f$ is an $E^1$-$E^2$-$C_{C^k-\text{op}-\text{triv}}$-$C_{C^k-\text{op}-\text{triv}}$ classic m-atlas near morphism of $A^1$ to $A^2$ in the coordinate spaces $C^1_{C^k-\text{triv}}, C^2_{C^k-\text{triv}}$.*

*Let $E^1$ be an open subspace of $E^2$. Then $f$ is an $E^1$-$E^2$ classic $C^k$-atlas morphism of $A^1$ to $A^2$ in the coordinate spaces $C^1, C^2$ iff $f$ is a strict $E^1$-$E^2$ classic m-atlas morphism of $A^1$ to $A^2$ in the coordinate spaces $C^1_{C^k-\text{triv}}, C^2_{C^k-\text{triv}}$.*

*$f$ is an $E^1$-$E^2$ classic $C^k$-atlas morphism of $A^1$ to $A^2$ in the coordinate spaces $C^1$, $C^2$ iff $f$ is a strict $E^1$-$E^2$-$C_{C^k-\text{op}-\text{triv}}$-$C_{C^k-\text{op}-\text{triv}}$ classic m-atlas morphism of $A^1$ to $A^2$ in the coordinate spaces $C^1_{C^k-\text{triv}}, C^2_{C^k-\text{triv}}$.*

# 21 $C^k$ manifolds

Conventionally a manifold is different from its atlases, but $\mathcal{A}t\ell^{C^k}_{\max}(E, C)$ in definition 18.7 (Categories of $C^k$ atlases) on page 124 encourages treating them on an equal footing. All of the results for maximal atlases carry directly over to results for manifolds.

**Definition 21.1** ($C^k$ manifolds). Let $E$ be a topological space, $C$ a linear space and $A$ a maximal[16] $C^k$-atlas of $E$ in the coordinate space $C$. Then $(E, C, A)$ is a $C^k$ manifold.

Let $E$ be a set of topological spaces and $C$ be a set of linear spaces. Then

$$\text{Man}^{C^k}_{\text{Ob}}(E, C) \overset{\text{def}}{=} \mathcal{A}t\ell^{C^k}_{\max}(E, C) \tag{21.1}$$

*Remark* 21.2. The manifold $(E, C, A)$ corresponds to the object $(A, E, C)$.

---

[16]Requiring that the atlas be full would eliminate some pathologies.

**Definition 21.3** ($C^k$ manifold morphisms)**.** Let $S^i \stackrel{\text{def}}{=} (E^i, C^i, A^i)$, $i = 1, 2$, be $C^k$ manifolds and $(f_0\colon E^1 \longrightarrow E^2, f_1\colon C^1 \longrightarrow C^2)$ be an $E^1$-$E^2$ $C^k$-morphism of $A^1$ to $A^2$ in the coordinate spaces $C^1$, $C^2$. Then $(f_0, f_1)$ is a $C^k$ morphism of $S^1$ to $S^2$.

Let $E^i$, $i = 1, 2$, be a topological space and $C^i$ be a linear space. Then

$$\mathrm{Man}^{C^k}_{\mathrm{Ar}}(E^1, C^1, E^2, C^2) \stackrel{\text{def}}{=} \mathcal{Atl}^{C^k}_{\mathrm{Ar}}(E^1, C^1, E^2, C^2) \quad (21.2)$$

Let $\boldsymbol{E}$ be a set of topological spaces and $\boldsymbol{C}$ a set of linear spaces. Then

$$\mathrm{Man}^{C^k}_{\mathrm{Ar}}(\boldsymbol{E}, \boldsymbol{C}) \stackrel{\text{def}}{=} \mathcal{Atl}^{C^k}_{\mathrm{Ar}\ \mathrm{max}}(\boldsymbol{E}, \boldsymbol{C}) \quad (21.3)$$

$$\mathrm{Man}^{C^k}(\boldsymbol{E}, \boldsymbol{C}) \stackrel{\text{def}}{=} \mathcal{Atl}^{C^k}_{\mathrm{max}}(\boldsymbol{E}, \boldsymbol{C}) \quad (21.4)$$

**Theorem 21.4** (Categories of $C^k$ manifolds)**.** *Let $\boldsymbol{E}$ be a set of topological spaces and $\boldsymbol{C}$ a set of linear spaces. Then* $\mathrm{Man}^{C^k}(\boldsymbol{E}, \boldsymbol{C})$ *is a category and the identity morphism of* $(A, E, C)$ *is an identity morphism.*

*Proof.* The result follows directly from definitions 21.1 and 21.3 above and lemma 18.9 ($\mathcal{Atl}^{C^k}(\boldsymbol{E}, \boldsymbol{C}$ is a category) on page 125. □

# Part VIII
# Equivalence of fiber bundles

For fiber bundles[17], the adjunct spaces are the base space $\boldsymbol{X} = (X, \mathcal{X})$, the fiber $\boldsymbol{Y} = (Y, \mathcal{Y})$ and the group $G$; the category of the coordinate space is the category of Cartesian products $\{U \times Y \mid U \stackrel{\mathrm{Ob}}{\in} \mathcal{X}\}$ of model neighborhoods in the base space with the entire fiber, with morphisms $t\colon U \times Y \longrightarrow U \times Y$ that preserve the fibers, i.e., $\pi_1 \circ t = \pi_1$, and are generated by the group action on the fiber (Equation (24.1)).

The sole adjunct functions are the projection $\pi\colon E \twoheadrightarrow X$, the group operation and the group action on the fiber.

This part of the paper defines bundle atlases, fiber bundles, local coordinate spaces equivalent to fiber bundles, categories of them and functors, and gives basic results.

**Definition 21.5** (Trivial group category of groups)**.** Let $\boldsymbol{G}$ be a set of topological groups. The trivial group category of $\boldsymbol{G}$, abbreviated $\boldsymbol{G}_{\text{group}-\mathcal{T}riv}$, is the category of all continuous homomorphisms between groups of $\boldsymbol{G}$. By abuse of language it will be shortened to $\boldsymbol{G}_{\mathcal{T}riv}$ when the meaning is clear from context.

---
[17]The literature has several definitions of fiber bundle. This paper uses one chosen for clarity of exposition. It differs from [Steenrod, 1999, p. 8] in that, e.g., it uses the machinery of maximal atlases rather than equivalence classes of coordinate bundles, the nomenclature differs in several minor regards.

**Definition 21.6** (Group actions)**.** Let $Y$ be a topological space, $G$ a topological group, $\rho$ an effective group action of $G$ on $Y$, $y \in Y$ and $g \in G$. Then $y \star g \stackrel{\text{def}}{=} \rho(y, g)$.

Let $X$ be a topological space and $x \in X$. Then $(x, y) \star g \stackrel{\text{def}}{=} (x, y \star g)$.

A $\star$ with a subscript, superset, underset or overset refers to the group action $\rho$ with the same subscript, superset, underset or overset.

*Remark* 21.7. This notation is only used when it is clear from context what the group action is.

**Definition 21.8** (Protobundles)**.** Let $\boldsymbol{B} \stackrel{\text{def}}{=} (E, X, Y, \pi, G, \rho)$, where $E$, $X$ and $Y$ are topological spaces, $G$ a topological group, $\pi\colon E \twoheadrightarrow X$ a continuous surjection and $\rho$ an efective group action of $G$ on $Y$. Then $\boldsymbol{B}$ is a protobundle.

*Remark* 21.9. While this definition does not itself require $E$ to have a local product structure nor require $\pi$ to have the Covering Homotopy Property, only those protobundles having an atlas are of interest, and for them definition 25.1 (Bundle atlases) on page 144 imposes additional constraints.

## 22 $G$-*rho* **model spaces**

**Definition 22.1** ($G$-$\rho$-model spaces)**.** Let $\mathbf{XY} \stackrel{\text{def}}{=} (X \times Y, \mathcal{XY})$ be a model space, $G$ a topological group and $\rho$ an effective group action of $G$ on $Y$ such that the objects of $\mathcal{XY}$ are products of open sets with $Y$ and the morphisms are fiber-preserving automorphisms generated by the group action, i.e.,

$$\left(\forall_{U \in \mathcal{XY}}^{\text{Ob}}\right)\left(\exists_{V \in X}^{\text{op}}\right) U = V \times Y \tag{22.1}$$

$$\left(\forall_{f \in \mathcal{XY}\colon V \times Y \stackrel{\tilde{=}}{\rightarrowtail\!\!\!\twoheadrightarrow} V \times Y}^{\text{Ar}}\right)\left(\exists_{g \in G^V}\right)\left(\forall_{(x,y) \in V \times Y}\right) f(x, y) = (x, y \star g(x)) \tag{22.2}$$

Then $\mathbf{XY}$ is a $G$-$\rho$ model space of $X \times Y$, abbreviated isG$\rho(\mathbf{XY}, Y, G, \rho)$, $G_{\mathbf{XY},f} \stackrel{\text{def}}{=} \text{range}(g)$ is the set of group elements associated with $f$ and $G_{\mathbf{XY}} \stackrel{\text{def}}{=} \bigcup_{f \in \mathcal{XY}}^{\text{Ar}} G_{\mathbf{XY},f}$ is the set of group elements associated with $\mathbf{XY}$ morphisms.

**Lemma 22.2** ($G$-$\rho$-model spaces)**.** *Let $\mathbf{XY}$ be a $G$-$\rho$ model space of $X \times Y$ and $f \stackrel{\text{Ar}}{\in} \mathcal{XY}\colon V \times Y \stackrel{\tilde{=}}{\rightarrowtail\!\!\!\twoheadrightarrow} V \times Y$.*

*The function $g$ in eq. (22.2) is unique.*

*Proof.* The group action is effective. □

*$G_{\mathbf{XY}}$ is unique.*

$$
\begin{CD}
((x^1,y^1),g^1) @>{f_C\times f_G}>> ((x^2,y^2),g^2) \\
@V{\star^1}VV @VV{\star^2}V \\
(x^1,y^1\star g^1) @>>{f_C}> (x^2,y^2\star g^2)
\end{CD}
$$

Figure 15: Preserving group action

*Proof.* The function $g$ is unique. □

**Definition 22.3** (Morphisms of $G$-$\rho$-model spaces). Let $X^i.Y^i$, $i = 1, 2$, be topological spaces, $G^i$ a topological group, $\rho^i$ an effective group action on $Y^i$ and $\mathbf{XY}^i \stackrel{\text{def}}{=} (X^i \times Y^i, \mathcal{XY}^i)$ a $G^i$-$\rho^i$ model space of $X^i \times Y^i$. Then a model function $f_C \colon \mathbf{XY}^1 \longrightarrow \mathbf{XY}^2$ is a $G^1$-$G^2$-$\rho^1$-$\rho^2$ morphism of $X^1 \times Y^1$ to $X^2 \times Y^2$, abbreviated isG$\rho$morph($\mathbf{XY}^1, G^1, \rho^1, \mathbf{XY}^2, G^2, \rho^2, f_C$) iff it preserves the group action, i.e., there is a continuous homomorphism $f_G \colon G^1 \longrightarrow G^2$ such that fig. 15 (Preserving group action) is commutative, i.e., eq. (22.3) holds.

$$\left(\exists_{f_G\colon G_1\longrightarrow G_2}\right)\left(\forall_{\substack{(x,y)\in X^1\times Y^1 \\ g\in G^1}}\right) f_C((x,y) \star^1 g) = f_C((x,y)) \star^2 f_G(g) \quad (22.3)$$

**Lemma 22.4** (Morphisms of $G$-$\rho$-model spaces). *Let $X^i.Y^i$, $i = 1, 2, 3$, be topological spaces, $G^i$ a topological group, $\rho^i$ an effective group action on $Y^i$, $\mathbf{XY}^i \stackrel{\text{def}}{=} (X^i \times Y^i, \mathcal{XY}^i)$ a $G^i$-$\rho^i$ model space of $X^i \times Y^i$ and $f_C^i \colon X^i \times Y^i \longrightarrow X^{i+1} \times Y^{i+1}$ a $G^i$-$G^{i+1}$-$\rho^i$-$\rho^{i+1}$ morphism of $X^i \times Y^i$ to $X^{i+1} \times Y^{i+1}$.*

1. *The function $f_G$ in eq. (22.3) is unique.*

   *Proof.* The group action is effective. □

2. *$f_C^i$ preserves fibers, i.e., $\pi_1(f_C^i(x,y)) = \pi_1(f_C^i(x,y'))$ for $x$ in $X^i$ and $y, y'$ in $Y^i$,*

   *Proof.* Since $\rho^i$ is effective, there is a $g \in G^i$ such that $y' = y \star^i g$ and thus

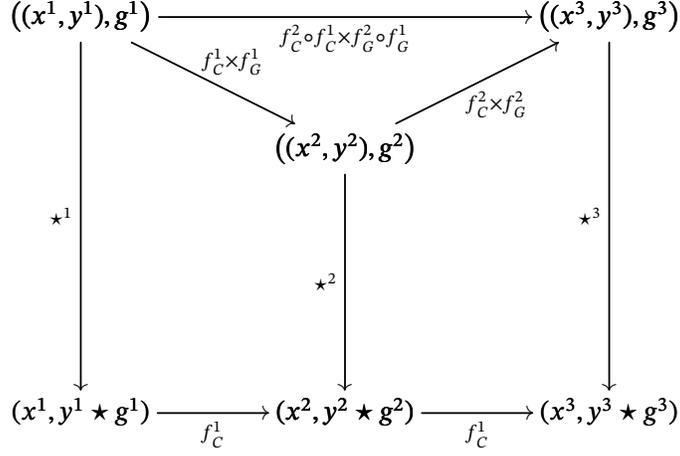

Figure 16: Preserving group actions

by definition 22.3

$$\begin{aligned}\pi_1(f_C^i(x,y')) &= \\ \pi_1(f_C^i(x,y\star^i g)) &= \\ \pi_1(f_C^i(x,y)\star^{i+1} f_G^i(g)) &= \\ \pi_1(f_C^i(x,y))\end{aligned} \qquad (22.4)$$

□

3. There exists a unique function $f_X^i\colon X^i \longrightarrow X^{i+1}$ such that $f_X^i \circ \pi_1 = \pi_1 \circ f_C^i$.

   *Proof.* Define $f_X^i(x) = \pi_1(f_C^i(x,y))$, where $y$ is an arbitrary point of $Y^i$. It does not depend on the choice of $y$ because $f_C^i$ preserves fibers. □

   *Remark* 22.5. $f_C^i$ may be a twisted product: there need not exist $f_Y^i\colon Y^i \longrightarrow Y^{i+1}$ such that $f_C^i = f_X^i \times f_Y^i$.

4. $f_C^2 \circ f_C^1$ is a $G^1$-$G^3$-$\rho^1$-$\rho^3$ morphism of $X^1 \times Y^1$ to $X^3 \times Y^3$.

   *Proof.* Let $f_G^i\colon G^i \longrightarrow G^{i+1}$ be a continuous homomorphism such that

$$\left(\forall_{\substack{(x^i,y^i)\in X^i\times Y^i \\ g^i\in G^i}}\right) f_C^i((x^i,y^i)\star g^i) = f_C^i((x^i,y^i))\star f_G^i(g^i) \qquad (22.5)$$

Let $(x^1, y^1) \in X^1 \times Y^1$ and $g^1 \in G^1$. Then fig. 16 (Preserving group actions) is commutative and

$$f_C^2 \circ f_C^1((x^1, y^1) \star g^1) = \tag{22.6}$$
$$= f_C^2\left(f_C^1((x^1, y^1)) \star f_G^1(g^1)\right) \tag{22.7}$$
$$= f_C^2\left(f_C^1((x^1, y^1))\right) \star f_G^2 \circ f_G^1(g^1) = \tag{22.8}$$
$$\tag{22.9}$$
$$= f_C^2 \circ f_C^1((x^1, y^1)) \star f_G^2 \circ f_G^1(g^1) \tag{22.10}$$

$\square$

**Definition 22.6** (Categories of $G$-$\rho$-model spaces). Let $X^\alpha. Y^\alpha$, $\alpha \prec A$, be topological spaces, $G^\alpha$ a topological group, $\rho^\alpha$ an effective group action on $Y^\alpha$, $\mathbf{XY}^\alpha \stackrel{\text{def}}{=} (X^\alpha \times Y^\alpha, \mathcal{XY}^\alpha)$ a $G^\alpha$-$\rho^\alpha$ model space of $X^\alpha \times Y^\alpha$, $\mathcal{XYG\rho}_{\text{Ob}} \stackrel{\text{def}}{=} \{(X^\alpha, Y^\alpha, G^\alpha, \rho^\alpha)|\alpha \prec A\}$, $\mathcal{XYG\rho}_{\text{Ar}} \stackrel{\text{def}}{=} \{f_C\colon X^\alpha \times Y^\alpha \longrightarrow X^\beta \times Y^\beta \mid$ isG$\rho$morph$(\mathbf{XY}^\alpha, G^\alpha, \rho^\alpha, \mathbf{XY}^\beta, G^\beta, \rho^\beta, f_C) \wedge \alpha \prec A \wedge \beta \prec A\}$ and $\mathcal{XYG\rho} \stackrel{\text{def}}{=} (\mathcal{XYG\rho}_{\text{Ob}}, \mathcal{XYG\rho}_{\text{Ar}})$. Then any subcategory of $\mathcal{XYG\rho}$ is a $G$-$\rho$ model category.

**Definition 22.7** (Trivial $G$-$\rho$-model spaces). Let $X$ and $Y$ be topological spaces, $G$ a topological group, $\rho$ an effective group action of $G$ on $Y$ and $\mathcal{XY}$ the category of all products of open subsets of $X$ with $Y$ and all homeomorphisms induced by the group action, i.e,

$$\text{Ob}(\mathcal{XY}) \stackrel{\text{def}}{=} \left\{V \times Y \,\middle|\, V \in \underset{\text{op}}{X}\right\} \tag{22.11}$$

$$\text{Ar}(\mathcal{XY}) \stackrel{\text{def}}{=} \left\{ f\colon V \times Y \xrightarrow{\tilde{=}} V \times Y \,\middle|\, V \in \underset{\text{op}}{X} \wedge \pi_1 \circ f = \pi_1 \wedge \left(\exists_{g \in G^V}\right)\left(\forall_{\substack{x \in V \\ y \in Y}}\right) f(x, y) = (x, y \star g(x)) \right\} \tag{22.12}$$

Then the trivial $G$-$\rho$ model space of $X, Y$, abbreviated $\underset{G-\rho-\mathbf{triv}}{X, Y}$, is $(X \times Y, \mathcal{XY})$ and $\underset{G-\rho-\mathbf{triv}}{X, Y}$ is a trivial $G$-$\rho$ model space of $X, Y$.

*Remark* 22.8. Let $G'$ be a topological group and $\rho'$ an effective group action of $G'$ on $Y$ such that $\underset{G-\rho-\mathbf{triv}}{X, Y} = \underset{G'-\rho'-\mathbf{triv}}{X, Y}$. Although $G'$ must be isomorphic to $G$. it need not have the same topology.

The identity morphism of $\underset{G-\rho-\mathbf{triv}}{X, Y}$ is $\text{Id}_{\underset{G-\rho-\mathbf{triv}}{X,Y}} \stackrel{\text{def}}{=} \text{Id}_{X \times Y}$.

Let $\boldsymbol{B}^\alpha \stackrel{\text{def}}{=} (E^\alpha, X^\alpha, Y^\alpha, \pi^\alpha, G^\alpha, \rho^\alpha)$, $\alpha \prec A$, be a protobundle and $\boldsymbol{B} \stackrel{\text{def}}{=} \{\boldsymbol{B}^\alpha | \alpha \prec A\}$ be a set of protobundles.

The trivial coordinate model category of $\boldsymbol{B}$, $\boldsymbol{B}_{\text{Bun}-\mathcal{T}riv}$, is the category with objects all trivial $G_\alpha$-$\rho_\alpha$ model spaces of $X_\alpha, Y_\alpha$, $\alpha \prec A$ and morphisms all continuous functions compatible with the group action:

$$\boldsymbol{B}_{\text{Bun}-\mathcal{T}riv_{\text{Ob}}} \stackrel{\text{def}}{=} \left\{ X, Y_{\mathbf{G}-\rho-\mathbf{triv}} \,\middle|\, (E, X, Y, \pi, G, \rho) \in \boldsymbol{B} \right\} \tag{22.13}$$

$$\boldsymbol{B}_{\text{Bun}-\mathcal{T}riv_{\text{Ar}}} \stackrel{\text{def}}{=} \left\{ f_C \colon X^1 \times Y^1 \longrightarrow X^2 \times Y^2 \,\middle|\, \left( \exists_{\substack{(E^i, X^i, Y^i, \pi^i, G^i, \rho^i) \in \boldsymbol{B} \\ f_G \colon G^1 \longrightarrow G^2}} \right) \text{isG}\rho\text{morph}(\underset{\mathbf{G}-\rho-\mathbf{triv}}{X^1, Y^1}, G^1, \rho^1, \underset{\mathbf{G}-\rho-\mathbf{triv}}{X^2, Y^2}, G^2, \rho^2, f^C) \right\} \tag{22.14}$$

$$\boldsymbol{B}_{\text{Bun}-\mathcal{T}riv} \stackrel{\text{def}}{=} \left( \boldsymbol{B}_{\text{Bun}-\mathcal{T}riv_{\text{Ob}}}, \boldsymbol{B}_{\text{Bun}-\mathcal{T}riv_{\text{Ar}}} \right) \tag{22.15}$$

*Remark* 22.9. The morphisms $f_C$ may be twisted products: there need not exist $f_Y \colon Y^1 \longrightarrow Y^2$ such that $f_C = f_X \times f_Y$.

The trivial product coordinate model category of $\boldsymbol{B}$, $\boldsymbol{B}_{\text{Bun}-\text{prod}-\mathcal{T}riv}$, is the category with objects all trivial $G^\alpha$-$\rho^\alpha$ model spaces of $X^\alpha, Y^\alpha$, $\alpha \prec A$ and morphisms all products of continuous functions compatible with the group action:

$$\boldsymbol{B}_{\text{Bun}-\text{prod}-\mathcal{T}riv_{\text{Ob}}} \stackrel{\text{def}}{=} \boldsymbol{B}_{\text{Bun}-\mathcal{T}riv_{\text{Ob}}} \tag{22.16}$$

$$\boldsymbol{B}_{\text{Bun}-\text{prod}-\mathcal{T}riv_{\text{Ar}}} \stackrel{\text{def}}{=} \left\{ f_C \stackrel{\text{def}}{=} f_X \times f_Y \colon X^1 \times Y^1 \longrightarrow X^2 \times Y^2 \,\middle|\, \left( \exists_{\substack{(E^i, X^i, Y^i, \pi^i, G^i, \rho^i) \in \boldsymbol{B} \\ f_G \colon G^1 \longrightarrow G^2}} \right) \text{isG}\rho\text{morph}(\underset{\mathbf{G}-\rho-\mathbf{triv}}{X^1, Y^1}, G^1, \rho^1, \underset{\mathbf{G}-\rho-\mathbf{triv}}{X^2, Y^2}, G^2, \rho^2, f_C) \right\} \tag{22.17}$$

$$\boldsymbol{B}_{\text{Bun}-\text{prod}-\mathcal{T}riv} \stackrel{\text{def}}{=} \left( \boldsymbol{B}_{\text{Bun}-\text{prod}-\mathcal{T}riv_{\text{Ob}}}, \boldsymbol{B}_{\text{Bun}-\text{prod}-\mathcal{T}riv_{\text{Ar}}} \right) \tag{22.18}$$

Any subcategory of $\boldsymbol{B}_{\text{Bun}-\mathcal{T}riv}$ is a trivial coordinate model category and any subcategory of $\boldsymbol{B}_{\text{Bun}-\text{prod}-\mathcal{T}riv}$ is a trivial product coordinate model category.

**Lemma 22.10** (The trivial coordinate model category of $\boldsymbol{B}$ is a category). *Let $\boldsymbol{B}^\alpha \stackrel{\text{def}}{=} (E^\alpha, X^\alpha, Y^\alpha, \pi^\alpha, G^\alpha, \rho^\alpha)$, $\alpha \prec A$, be a protobundle and $\boldsymbol{B} \stackrel{\text{def}}{=} \{\boldsymbol{B}^\alpha | \alpha \prec A\}$ be a set of protobundles.*

*Then $\boldsymbol{B}_{\text{Bun}-\mathcal{T}riv}$ is a category and the identity morphism for object $\boldsymbol{B}^\alpha$ is $\text{Id}_{X^\alpha \times Y^\alpha}$.*

*Proof.*

**Composition:** The composition of $G$-$\rho$ morphisms is a $G$-$\rho$ morphism by lemma 22.4 (Morphisms of $G$-$\rho$-model spaces) on page 136.

**Associativity:** Morphisms are simply functions and composition of morphisms is simply composition of functions.

**Unit:** The identity morphisms are simply identity functions and composition of morphisms is simply composition of functions.

□

$\boldsymbol{B}_{\text{Bun}-\text{prod}-\mathcal{T}riv}$ *is a category and the identity morphism for object $\boldsymbol{B}^\alpha$ is $\text{Id}_{X^\alpha \times Y^\alpha}$.*

*Proof.*

**Composition:** The composition of $G$-$\rho$ morphisms is a $G$-$\rho$ morphism by lemma 22.4 (Morphisms of $G$-$\rho$-model spaces) on page 136. Let $X^i, Y^i \underset{\mathbf{G^i}-\boldsymbol{\rho^i}-\mathbf{triv}}{\overset{\text{Ob}}{\in}} \boldsymbol{B}_{\text{Bun}-\text{prod}-\mathcal{T}riv}$, $i = 1, 2, 3$, $f_C^i = f_X^i \times f_Y^i \colon X^i \times Y^i \longrightarrow X^{i+1} \times Y^{i+1} \overset{\text{Ar}}{\in} \boldsymbol{B}_{\text{Bun}-\text{prod}-\mathcal{T}riv}$. Then $f_C^{i+1} \circ f_C^i = (f_X^{i+1} \circ f_X^i) \times (f_Y^{i+1} \circ f_Y^i)$.

**Associativity:** Morphisms are simply functions and composition of morphisms is simply composition of functions.

**Unit:** The identity morphisms are simply identity functions and composition of morphisms is simply composition of functions.

□

**Lemma 22.11** (The trivial $G$-$\rho$ model space of $X, Y$ is a $G$-$\rho$ model space of $X, Y$). *Let $X$ and $Y$ be topological spaces, $G$ a topological group and $\rho$ an effective group action on $Y$. Then $X, Y_{\mathbf{G}-\boldsymbol{\rho}-\mathbf{triv}}$ is a $G$-$\rho$ model space of $X \times Y$.*

*Proof.* $X, Y_{\mathbf{G}-\boldsymbol{\rho}-\mathbf{triv}}$ satisfies the conditions of definition 22.1 ($G$-$\rho$-model spaces) on page 135:

1. Product with fiber:

$$\left(\forall_{U \in \pi_2^{\text{Ob}}\left(\underset{\mathbf{G}-\rho-\mathbf{triv}}{X,Y}\right)}\right)\left(\exists_{V \in X_{\text{op}}}\right) U = V \times Y$$

   by definition 22.7.

2. Generated by $\rho$: Let $f \in \pi_2^{\text{Ar}}\left(\underset{\mathbf{G}-\rho-\mathbf{triv}}{X,Y}\right): V \times Y \rightarrowtail\twoheadrightarrow V \times Y$. Then

$$\left(\exists_{g \in G^V}\right)\left(\forall_{(x,y) \in V \times Y}\right) f(x,y) = (x, y \star g(x))$$

   by definition 22.7.

$\square$

## 23 $G$-$rho$-nearly commutative diagrams

Let

1. $X$ and $Y$ be topological spaces

2. $G$ be a topological group,

3. $\rho$ an effective group action on $Y$

4. $\boldsymbol{C} = (C, \mathscr{C}) \stackrel{\text{def}}{=} \underset{\mathbf{G}-\rho-\mathbf{triv}}{X,Y}$

5. $D$ be, as in diagram (1.20) of 14, a tree with two branches, whose nodes are topological spaces $U_i$ and $V^j$ and whose links are continuous functions $f_i: U_i \longrightarrow U_{i+1}$ and $f'_j: U_j \longrightarrow U_{j+1}$ between the spaces, with $U_0 = V_0$, $U_m \subseteq C$ open and $V_n \subseteq C$ open, as shown in fig. 1 (Uncompleted nearly commutative diagram) on page 15.

**Definition 23.1** ($G$-$\rho$-nearly commutative diagrams). $D$ is nearly commutative in $X, Y, \pi, \rho$ iff $D$ is nearly commutative in category $\mathscr{C}$.

**Definition 23.2** ($G$-$\rho$-nearly commutative diagrams at a point). Let $\mathscr{C}$ and $D$ be as above and $x$ be an element of the initial node. $D$ is nearly commutative in $X, Y, \pi, \rho$ at $x$ iff $D$ is nearly commutative in $\mathscr{C}$ at $x$.

**Definition 23.3** ($G$-$\rho$-locally nearly commutative diagrams). Let $\mathscr{C}$ and $D$ be as above. $D$ is locally nearly commutative in $X, Y, \pi, \rho$ iff $D$ is locally nearly commutative in $\mathscr{C}$.

# 24 Bundle charts

**Definition 24.1** ($Y$-$\pi$-bundle charts). Let $E, X, Y$ be topological spaces and $\pi\colon E \longrightarrow X$. A $Y$-$\pi$-bundle chart $C \stackrel{\text{def}}{=} (U, V \times Y, \phi)$ of $E$ in the coordinate space $X \times Y$ consists of

1. a nonvoid open set $U \subseteq E$, known as a coordinate patch
2. An open set $V \times Y \subseteq X \times Y$
3. a homeomorphism $\phi\colon U \rightarrowtail\!\!\!\stackrel{\tilde{=}}{\twoheadrightarrow} V \times Y$, known as a coordinate function, that preserves fibers. i.e., $\pi_1 \circ \phi = \pi$.

**Lemma 24.2** (Properties of projection). *Let $(U, V \times Y, \phi)$ be a $Y$-$\pi$-bundle chart of $E$ in the coordinate space $X \times Y$ and $v \in V$. Then*
   $\pi\!\restriction_U$ *is a surjection.*

*Proof.* Let $v \in V$, $y \in Y$ and $u \stackrel{\text{def}}{=} \phi^{-1}\big((v, y)\big) \in U$. Then $\pi(u) = v$. $\square$

   $\pi^{-1}[\{v\}]$ *is homeomorphic to $Y$.*

*Proof.* $\phi$ and $\phi^{-1}$ are homeomorphisms, so their restrictions are homeomorphisms and thus $\phi^{-1}[\{v\}]$ is homeomorphic to $\{v\} \times Y$, which is homeomorphic to $Y$. $\square$

**Definition 24.3** ($Y$-$\pi$ subcharts). Let $(U, V \times Y, \phi)$ be a $Y$-$\pi$-bundle chart of $E$ in the coordinate space $X \times Y$ and $U'$ be a nonvoid open subset of $U$. Then $(U', V' \times Y, \phi') \stackrel{\text{def}}{=} (U', \phi[U'], \phi\!\restriction_{U'})$ is a subchart of $(U, V \times Y, \phi)$.

**Lemma 24.4** ($Y$-$\pi$ subcharts). *Let $(U, V \times Y, \phi)$ be a $Y$-$\pi$-bundle chart of $E$ in the coordinate space $X \times Y$ and $(U', V' \times Y, \phi')$ a subchart of $(U, V \times Y, \phi)$. Then $(U', V' \times Y, \phi')$ is a $Y$-$\pi$-bundle chart of $E$ in the coordinate space $X \times Y$.*

*Proof.* $(U', V' \times Y, \phi')$ satisfies the conditions of definition 24.1:

1. $U' \subseteq E$ is open.
2. $\phi[U']$ os open since $\phi$ is a homeomorphism.
3. $\phi\!\restriction_{U'}\colon U' \longrightarrow V' \times Y$ is the restriction of a homeomorphism and thus a homeomorphism. $\phi\!\restriction_{U'}$ preserves fibers because $\phi$ does.

$\square$

**Definition 24.5** ($G$-$\rho$-compatibility). Let $E, X, Y$ be topological spaces, $G$ a topological group, $\rho\colon Y \times G \longrightarrow Y$ an effective right action of $G$ on $Y$, $\pi\colon E \twoheadrightarrow X$ surjective and $(U^\mu, V^\mu \times Y, \phi^\mu)$, $(U^\nu, V^\nu \times Y, \phi^\nu)$ $Y$-$\pi$-bundle charts. $(U^\mu, V^\mu \times Y, \phi^\mu)$ and $(U^\nu, V^\nu \times Y, \phi^\nu)$ are $G$-$\rho$-compatible if either

1. $U^\mu$ and $U^\nu$ are disjoint

2. The transition function $t^\mu_\nu \stackrel{\text{def}}{=} \phi^\mu \circ \phi^{\nu-1}\restriction_{\phi^\nu[U^\mu \cap U^\nu]}$ is generated by the group action, i.e., there is a continuous function $g^\mu_\nu \colon \pi_1[\phi^\nu[U^\mu \cap U^\nu]] \longrightarrow G$ such that

$$\left(\forall_{(x,y) \in \phi^\nu[U^\mu \cap U^\nu]}\right) t^\mu_\nu(x,y) = (x, y \star g^\mu_\nu(x)) \qquad (24.1)$$

**Lemma 24.6** (Symmetry of $G$-$\rho$ compatibility). *Let $(U^\mu, V^\mu \times Y, \phi^\mu)$ and $(U^\nu, V^\nu \times Y, \phi^\nu)$ be $Y$-$\pi$-bundle charts of $E$ in the coordinate space $X \times Y$. Then $(U^\mu, V^\mu \times Y, \phi^\mu)$ is $G$-$\rho$-compatible with $(U^\nu, V^\nu \times Y, \phi^\nu)$ iff $(U^\nu, V^\nu \times Y, \phi^\nu)$ is $G$-$\rho$-compatible with $(U^\mu, V^\mu \times Y, \phi^\mu)$.*

*Proof.* It suffices to prove the result in one direction. If $(U^\mu, V^\mu \times Y, \phi^\mu) \cap (U^\nu, V^\nu \times Y, \phi^\nu)$ then $(U^\nu, V^\nu \times Y, \phi^\nu) \cap (U^\mu, V^\mu \times Y, \phi^\mu)$. Otherwise, let $g^\mu_\nu \colon \pi_1[\phi^\nu[U^\mu \cap U^\nu]] \longrightarrow G$ be a continuous function such that

$$\left(\forall_{(x,y) \in \phi^\nu[U^\mu \cap U^\nu]}\right) t^\mu_\nu(x,y) = (x, y \star g^\mu_\nu(x))$$

and the inverse transition function is also generated by the group action:

$$\left(\forall_{(x,y) \in \phi^\nu[U^\mu \cap U^\nu]}\right) t^\nu_\mu(x,y) = \left(x, y \star g^\mu_\nu(x)^{-1}\right)$$

$\square$

**Lemma 24.7** ($G$-$\rho$-compatibility of subcharts). *Let $(U^i, V^i \times Y, \phi^i)$ be a $Y$-$\pi$-bundle chart of $E$ in the coordinate space $X \times Y$, $(U'^i, V'^i \times Y, \phi'^i)$ a subchart and $(U^1, V^1, \phi^1)$ be $G$-$\rho$-compatible with $(U^2, V^2, \phi^2)$. Then $(U'^1, V'^1 \times Y, \phi'^1)$ is $G$-$\rho$-compatible with $(U'^2, V'^2 \times Y, \phi'^2)$.*

*Proof.* If $U^1 \cap U^2 = \emptyset$ then $U'^1 \cap U'^2 = \emptyset$. If $U'^1 \cap U'^2 = \emptyset$ then $(U'^1, V'^1 \times Y, \phi'^1)$ is $G$-$\rho$-compatible with $(U'^2, V'^2 \times Y, \phi'^2)$. Otherwise, the transition function $t^1_2 \stackrel{\text{def}}{=} \phi^1 \circ \phi^{2-1}\restriction_{\phi^2[U^1 \cap U^2]}$ is generated by the group action and hence $t^1_2\restriction_{\phi^2[U'^1 \cap U'^2]} \colon \phi^2[U'^1 \cap U'^2] \overset{\tilde{=}}{\rightarrowtail\!\!\!\!\!\rightarrow} \phi^1[U'^1 \cap U'^2]$ is generated by the group action. $\square$

**Corollary 24.8** ($G$-$\rho$-compatibility with subcharts). *Let $(U, V \times Y, \phi)$, be a $Y$-$\pi$-bundle chart of $E$ in the coordinate space $X \times Y$ and $(U', V' \times Y, \phi')$ a subchart. Then $(U', V' \times Y, \phi')$ is $G$-$\rho$-compatible with $(U, V \times Y, \phi)$.*

*Proof.* $(U, V \times Y, \phi)$ is $G$-$\rho$-compatible with itself and is a subchart of itself. $\square$

**Definition 24.9** (Covering by $Y$-$\pi$-bundle charts). Let $A$ be a set of $Y$-$\pi$-bundle charts of $E$ in the coordinate space $X \times Y$. $A$ covers $E$ iff $E = \bigcup \pi_1[A]$.

**Lemma 24.10.** *Let $A$ be a set of $Y$-$\pi$-bundle charts of $E$ in the coordinate space $X \times Y$ that covers $E$ and $x \in X$. Then $\pi^{-1}[\{x\}]$ is homeomorphic to $Y$.*

*Proof.* Since $A$ covers $E$, there is a chart $(U, V, \phi)$ in $A$ containing $x$. Then $\pi^{-1}[\{x\}]$ is homeomorphic to $Y$ by lemma 24.2 (Properties of projection) on page 142. $\square$

## 25 Bundle atlases

A set of charts can be atlases for different fiber bundles even if it is for the same total model space, base space and fiber. In order to aggregate atlases into categories, there must be a way to distinguish them. Including the spaces[18], group and group action in the definitions of the categories serves the purpose.

**Definition 25.1** (Bundle atlases). Let $\boldsymbol{B} \stackrel{\text{def}}{=} (E, X, Y, \pi, G, \rho)$, be a protobundle. Then $\boldsymbol{A}$ is a bundle atlas of $B$, abbreviated isAtl$_{\text{Ob}}^{\text{Bun}}(\boldsymbol{A}, \boldsymbol{B})$ and $\boldsymbol{A}$ is a $\pi$-$G$-$\rho$-bundle atlas of $E$ in the coordinate space $X \times Y$, abbreviated isAtl$_{\text{Ob}}^{\text{Bun}}(\boldsymbol{A}, E, X, Y, \pi, G, \rho)$, iff it consists of a set of mutually $G$-$\rho$-compatible $Y$-$\pi$-bundle charts of $E$ in the coordinate space $X \times Y$ that covers $E$[19].

By abuse of language we write $U \in \boldsymbol{A}$ for $U \in \pi_1[\boldsymbol{A}]$.

*Remark* 25.2. The definition of a $\pi$-$G$-$\rho$-bundle atlas is by design similar to the definition of a coordinate bundle in [Steenrod, 1999, p. 7], but there are significant differences. This paper will use the term bundle atlas to avoid confusion.

Let $\boldsymbol{B}^\alpha \stackrel{\text{def}}{=} (E^\alpha, X^\alpha, Y^\alpha, \pi^\alpha, G^\alpha, \rho^\alpha), \alpha \prec A$, be a protobundle and $\boldsymbol{B} \stackrel{\text{def}}{=} \{\boldsymbol{B}^\alpha | \alpha \prec A\}$ be a set of protobundles. Then

$$\mathcal{A}t\ell_{\text{Ob}}^{\text{Bun}}(\boldsymbol{B}^\alpha) \stackrel{\text{def}}{=} \left\{(\boldsymbol{A}, \boldsymbol{B}^\alpha) \,\middle|\, \text{isAtl}_{\text{Ob}}^{\text{Bun}}(\boldsymbol{A}, E^\alpha, X^\alpha, Y^\alpha, G^\alpha, \pi^\alpha, \rho^\alpha)\right\} \quad (25.1)$$

$$\mathcal{A}t\ell_{\text{Ob}}^{\text{Bun}} \boldsymbol{B} \stackrel{\text{def}}{=} \bigcup\nolimits_{\alpha \prec A} \mathcal{A}t\ell_{\text{Ob}}^{\text{Bun}}(\boldsymbol{B}^\alpha) \quad (25.2)$$

**Lemma 25.3** (Bundle atlases). *Let $E, X, Y$ be topological spaces, $G$ a topological group, $\pi \colon E \twoheadrightarrow X$ surjective and $\rho \colon Y \times G \longrightarrow Y$ an effective right action of $G$ on $Y$. Then $\boldsymbol{A}$ is a $\pi$-$G$-$\rho$-bundle atlas of $E$ in the coordinate space $X \times Y$ iff it is an m-atlas of $\underset{\text{triv}}{E}$ in the coordinate space $\underset{G-\rho-\text{triv}}{X, Y}$ and every coordinate function preserves fibers.*

*Proof.* If $\boldsymbol{A}$ is a $\pi$-$G$-$\rho$-bundle atlas of $E$ in the coordinate space $X \times Y$ then

1. Every chart in $\boldsymbol{A}$ is a $Y$-$\pi$-bundle chart of $E$ in the coordinate space $X \times Y$, and hence its coordinate function preserves fibers.

2. The charts in $\boldsymbol{A}$ are mutually $G$-$\rho$-compatible; hence the transition functions are generated by the group action and are morphisms of $\underset{G-\rho-\text{triv}}{X, Y}$.

If $\boldsymbol{A}$ is an m-atlas of $\underset{\text{triv}}{E}$ in the coordinate space $\underset{G-\rho-\text{triv}}{X, Y}$ then the transition functions are generated by the group action and and thus the charts are mutually $G$-$\rho$-compatible.

---

[18]The spaces are redundant, but convenient.

[19]There is no need to introduce the concept of a full $\pi$-$G$-$\rho$-bundle atlas because a $\pi$-$G$-$\rho$-bundle atlas is automatically full.

If every coordinate function preserves fibers then the m-charts of $A$ are $Y$-$\pi$-bundle charts. $\square$

**Definition 25.4** (Compatibility of charts with bundle atlases). A $Y$-$\pi$-bundle chart $(U, V \times Y, \phi)$ of $E$ in the coordinate space $X \times Y$ is $G$-$\rho$-compatible with a $\pi$-$G$-$\rho$-bundle atlas $A$ iff it is $G$-$\rho$-compatible with every chart in the atlas.

**Lemma 25.5** (Compatibility of subcharts with bundle atlases). *Let $A$ be a $\pi$-$G$-$\rho$-bundle atlas of $E$ in the coordinate space $X \times Y$ and $C^1 = (U, V \times Y, \phi)$ a $Y$-$\pi$-bundle chart in $A$. Then any subchart of $C^1$ is $G$-$\rho$-compatible with $A$.*

*Proof.* Let $C'^1 = (U'^1, V'^1 \times Y, \phi'^1)$ be a subchart of $C^1$ and $C^2 = (U^2, V^2 \times Y, \phi^2)$ another chart in $A$.

1. If $U^1 \cap U^2 = \emptyset$, then $U' \cap U^2 = \emptyset$.

2. If $U' \cap U^2 = \emptyset$ then $C'^1$ is $G$-$\rho$-compatible with $C^2$.

3. Otherwise the transition function $t^1_2 \stackrel{\text{def}}{=} \phi^1 \circ \phi^{2-1}\restriction_{\phi^1[U^1 \cap U^2]}$ is generated by the group action and thus $t^1_2\restriction_{\phi^2[U'^1 \cap U^2]}$ is generated by the group action.

$\square$

**Lemma 25.6** (Extensions of bundle atlases). *Let $A$ be a $\pi$-$G$-$\rho$-atlas of $E$ in the coordinate space $X \times Y$ and $(U_i, V_i, \phi_i)$, $i = 1, 2$, be a $\pi$-$G$-$\rho$ chart of $E$ in the coordinate space $X \times Y$ $G$-$\rho$-compatible with $A$ in the coordinate space $X \times Y$. Then $(U_1, V_1, \phi_1)$ is $G$-$\rho$-compatible with $(U_2, V_2, \phi_2)$ in the coordinate space $X \times Y$.*

*Proof.* If $U_1 \cap U_2 = \emptyset$ then $(U_1, V_1, \phi_1)$ is $G$-$\rho$-compatible with $(U_2, V_2, \phi_2)$. Otherwise, $\phi_2 \circ \phi_1^{-1}\restriction_{\phi_1[U_1 \cap U_2]}: \phi_1[U_1 \cap U_2] \xrightarrow{\tilde{=}} \phi_2[U_1 \cap U_2]$ is a homeomorphism. It remains to show that $\phi_2 \circ \phi_1^{-1}\restriction_{\phi_1[U_1 \cap U_2]}$ is generated by the group action. Let $(U'_\alpha, V'_\alpha, \phi'_\alpha)$, $\alpha < A$, be charts in $A$ such that $U_1 \cap U_2 \subseteq \bigcup_{\alpha < A} U'_\alpha$ and $U_1 \cap U_2 \cap U'_\alpha \neq \emptyset$, $\alpha < A$. Since the charts are $G$-$\rho$-compatible with $(U'_\alpha, V'_\alpha, \phi'_\alpha)$, $\phi_2 \circ \phi'^{-1}_\alpha\restriction_{U_1 \cap U_2 \cap U'_\alpha}$ and $\phi'_\alpha \circ \phi_1^{-1}\restriction_{U^1 \cap U^2 \cap U'_\alpha}$ are generated by the group action and thus $\phi^2 \circ \phi_1^{-1} = \phi^2 \circ \phi'^{-1}_\alpha \circ \phi'_\alpha \circ \phi_1^{-1}$ is generated by the group action. $\square$

**Definition 25.7** (Maximal bundle atlases). Let $A$ be a $\pi$-$G$-$\rho$-bundle atlas of $E$ in the coordinate space $X \times Y$. $A$ is a maximal $\pi$-$G$-$\rho$-bundle atlas, abbreviated isAtl$^{\text{Bun}}_{\text{Ob}}$(A, E, X, Y, G, $\pi$, $\rho$), iff it cannot be extended by adding an additional $G$-$\rho$-compatible $Y$-$\pi$-bundle chart.

$A$ is a semi-maximal $\pi$-$G$-$\rho$-bundle atlas of $E$ in the coordinate space $C$, abbreviated isAtl$^{\text{Bun}}_{\text{Ob}}$(A, E, X, Y, G, $\pi$, $\rho$), iff whenever $(U, V \times Y, \phi) \in A$, $U' \subseteq U$, $V' \times Y \subseteq V \times Y$ and $V'' \times Y \subseteq X \times Y$ are open sets, $\phi[U'] = V' \times Y$, $\phi': V' \times Y \xrightarrow{\tilde{=}} V'' \times Y$ is a fiber preserving homeomorphism generated by the group action then $(U', V'' \times Y, \phi' \circ \phi) \in A$.

Let $\boldsymbol{B}^\alpha \stackrel{\text{def}}{=} (E^\alpha, X^\alpha, Y^\alpha, \pi^\alpha, G^\alpha, \rho^\alpha)$, $\alpha \prec A$, be a protobundle and $\boldsymbol{B} \stackrel{\text{def}}{=} \{\boldsymbol{B}^\alpha | \alpha \prec A\}$ be a set of protobundles. Then

$$\mathcal{A}t\ell^{\text{Bun}}_{\underset{\text{S-max (max)}}{\text{Ob}}}(\boldsymbol{B}^\alpha) \stackrel{\text{def}}{=} \left\{ (\boldsymbol{A}, \boldsymbol{B}^\alpha) \,\middle|\, \text{isAtl}^{\text{Bun}}_{\underset{\text{S-max (max)}}{\text{Ob}}}\left(\boldsymbol{A}, E^\alpha, X^\alpha, Y^\alpha, G^\alpha, \pi^\alpha, \rho^\alpha\right) \right\} \quad (25.3)$$

$$\mathcal{A}t\ell^{\text{Bun}}_{\underset{\text{S-max (max)}}{\text{Ob}}} \boldsymbol{B} \stackrel{\text{def}}{=} \bigcup\nolimits_{\alpha \prec A} \mathcal{A}t\ell^{\text{Bun}}_{\underset{\text{S-max (max)}}{\text{Ob}}} (\boldsymbol{B}^\alpha) \quad (25.4)$$

**Lemma 25.8** (Maximal $\pi$-$G$-$\rho$-bundle atlases are semi-maximal $\pi$-$G$-$\rho$-bundle atlases). *Let $E, X, Y$ be topological spaces, $G$ a topological group, $\pi\colon E \twoheadrightarrow X$ surjective, $\rho\colon Y \times G \longrightarrow Y$ an effective right action of $G$ on $Y$ and $\boldsymbol{A}$ a maximal $\pi$-$G$-$\rho$-bundle atlas of $E$ in the coordinate space $X \times Y$. Then $\boldsymbol{A}$ is a semi-maximal $\pi$-$G$-$\rho$-bundle atlas of $E$ in the coordinate space $X \times Y$.*

*Proof.* Let $(U, V \times Y, \phi) \in \boldsymbol{A}$, $U' \subseteq U$, $V' \times Y \subseteq V \times Y$ and $V'' \times Y \subseteq X \times Y$ be open sets and $\phi[U'] = V' \times Y$, $\phi'\colon V' \times Y \xrightarrow{\tilde{=}} V'' \times Y$ be a fiber preserving homeomorphism generated by the group action. $(U', V', \phi)$ is a subchart of $(U, V, \phi)$ and by lemma 25.5 (Compatibility of subcharts with bundle atlases) on page 145 is $G$-$\rho$-compatible with the charts of $\boldsymbol{A}$. Since $\phi'$ is a fiber preserving homeomorphism generated by the group action, $(U', V'' \times Y, \phi' \circ \phi)$ is $G$-$\rho$-compatible with the charts of $\boldsymbol{A}$. Since $\boldsymbol{A}$ is maximal, $(U', V'', \phi' \circ \phi)$ is a chart of $\boldsymbol{A}$. $\square$

**Theorem 25.9** (Existence and uniqueness of maximal $\pi$-$G$-$\rho$-bundle atlases). *Let $\boldsymbol{B} \stackrel{\text{def}}{=} (E, X, Y, \pi, G, \rho)$ be a protobundle and $\boldsymbol{A}$ $\pi$-$G$-$\rho$-bundle atlas of $E$ in the coordinate space $X \times Y$. Then there exists a unique maximal $\pi$-$G$-$\rho$-bundle atlas $\text{Atlas}^{\text{Bun}}_{\text{max}}(\boldsymbol{A}, \boldsymbol{B})$ of $E$ in the coordinate space $X \times Y$ $G$-$\rho$-compatible with $\boldsymbol{A}$.*

*Proof.* Let $\boldsymbol{P}$ be the set of all $\pi$-$G$-$\rho$-bundle atlases $E$ in the coordinate space $X \times Y$ containing $\boldsymbol{A}$ and $G$-$\rho$ compatible in the coordinate space $X \times Y$ with $\boldsymbol{A}$. Let $\boldsymbol{P}_{\text{max}}$ be a maximal chain of $\boldsymbol{P}$. Then $A' = \bigcup \boldsymbol{P}_{\text{max}}$ is a maximal $\pi$-$G$-$\rho$-bundle atlas of $\boldsymbol{E}$ in the coordinate space $X \times Y$ $G$-$\rho$ compatible with $\boldsymbol{A}$. Uniqueness follows from lemma 25.6 (Extensions of bundle atlases) on page 145. $\square$

**Lemma 25.10** (Existence and uniqueness of projection for atlases in $G$-$\rho$-model spaces). *Let $E$, $X$ and $Y$ be topological spaces, $G$ a topological group, $\rho$ an effective action of $G$ on $Y$, $\boldsymbol{C} \stackrel{\text{def}}{=} (X \times Y, \mathcal{XY})$ a $G$-$\rho$ model space of $X \times Y$ and $\boldsymbol{A}$ an m-atlas of $E$ in the coordinate model space $\boldsymbol{C}$. Then there exists a unique function $\pi\colon \boldsymbol{E} \longrightarrow X$ such that for any chart $(U, V, \phi)$ in $\boldsymbol{A}$, $\pi\restriction_U = \pi_1 \circ \phi$. If $\boldsymbol{A}$ is full then $\pi$ is surjective.*

*Proof.* Let $(U, V, \phi)$ be an arbitrary chart in $A$ and define $\pi(e \in U) = \pi_1 \circ \phi(e)$. $\pi(e)$ does not depend on the choice of chart because the morphisms of a $G$-$\rho$ model space preserve fibers. $\pi$ is continuous because it is continuous on each coordinate patch.

Let $x \in X$. If $A$ is full then there exists a chart $(U, V, \phi)$ in $A$ such that $x \in \pi_1[V]$. Let $u$ be an arbitrary point in $\phi^{-1}[\{x\} \times Y]$. Then $\pi(u) = \pi_1(\phi(u)) = x$. □

## 26 Bundle-atlas (near) morphisms and functors

This section introduces a taxonomy for morphisms between bundle atlases, defines categories of bundle atlases and constructs functors from them to categories of m-atlases. It only constructs reverse functors for subcategories of $\mathbf{B}_{\text{Bun-prod-}\mathcal{T}riv}$.

### 26.1 Bundle-atlas (near) morphisms

#### 26.1.1 Definitions of bundle-atlas (near) morphisms

**Definition 26.1** (Bundle-atlas near morphisms). Let $\mathbf{B}^i \stackrel{\text{def}}{=} (E^i, X^i, Y^i, \pi^i, G^i, \rho^i)$, $i = 1, 2$, be a protobundle and let $\mathbf{A}^i$ be a $\pi^i$-$G^i$-$\rho^i$-atlas of $E^i$ in the coordinate space $C^i = X^i \times Y^i$. Then $\mathbf{f} \stackrel{\text{def}}{=} (f_E \colon E^1 \longrightarrow E^2, f_X \colon X^1 \longrightarrow X^2, f_Y \colon Y^1 \longrightarrow Y^2, f_G \colon G^1 \longrightarrow G^2)$ is a $\mathbf{B}^1$-$\mathbf{B}^2$ bundle-atlas morphism from $\mathbf{A}^1$ to $\mathbf{A}^2$, abbreviated $\text{isAtl}_{\text{Ar}}^{\text{Bun,near}}(\mathbf{A}^1, \mathbf{B}^1, \mathbf{A}^2, \mathbf{B}^2, \mathbf{f})$, iff

1. all four functions are continuous

2. $f_G$ is a homomorphism

3. $\mathbf{f}$ commutes with $\pi^i$ and $\rho^i$, i.e.,

   (a) $\pi^2 \circ f_E = f_X \circ \pi^1$

   (b) $\left( \forall_{\substack{y \in Y^1 \\ g \in G^1}} \right) f_Y(y) \star^2 f_G(g) = f_Y(y \star^1 g)$

4. For any charts $(U^i, V^i, \phi^i \colon U^i \overset{\tilde{=}}{\rightarrowtail\!\!\!\twoheadrightarrow} V^i) \in \mathbf{A}^i$. $i = 1, 2$, with $I \stackrel{\text{def}}{=} U^1 \cap f_E^{-1}[U^2] \neq \emptyset$, diagram (26.1) below is locally nearly commutative in $X, Y, \pi, \rho$, i.e., for any $u^1 \in I$ there are open sets $U'^1 \subseteq I$, $V'^1 \subseteq V^1$, $U'^2 \subseteq U^2$, $V'^2 \subseteq V^2$, $\hat{V}'^2 \subseteq X^2 \times Y^2$ and a homeomorphism $\hat{f} \colon V'^2 \overset{\tilde{=}}{\rightarrowtail\!\!\!\twoheadrightarrow} \hat{V}'^2$ [20] such that eqs. (26.2) to (26.7) below hold, as shown in figs. 17 and 18.

---

[20]This reverses the direction of the arrow $\hat{f} \colon U_m \overset{\tilde{=}}{\rightarrowtail\!\!\!\twoheadrightarrow} V_n$ from section 23 (*G-rho*-nearly commutative diagrams) on page 141. This is permissible since $\hat{f}$ is a homeomorphism

$$D \stackrel{\text{def}}{=} \{$$

$$f_E \colon I \stackrel{\text{def}}{=} U^1 \cap f_E^{-1}[U^2] \longrightarrow U^2, \phi^2 \colon U^2 \stackrel{\tilde{=}}{\rightarrowtail} V^2, \quad (26.1)$$
$$\phi^1 \colon I \rightarrowtail V^1, f_X \times f_Y \colon V^1 \longrightarrow X^2 \times Y^2$$
$$\}$$

$$u^1 \in U'^1 \tag{26.2}$$

$$f_E[U'^1] \subseteq U'^2 \tag{26.3}$$

$$\phi^1[U'^1] = V'^1 \tag{26.4}$$

$$f_X \times f_Y[V'^1] \subseteq \hat{V}'^2 \tag{26.5}$$

$$\phi^2[U'^2] = V'^2 \tag{26.6}$$

$$\hat{f} \circ \phi^2 \circ f_E = f_X \times f_Y \circ \phi^1 \tag{26.7}$$

*Remark* 26.2. $(U'^i, V'^i, \phi^i \restriction_{U'^i, V'^i})$ need not be a chart of $\boldsymbol{A}^i$.

Figure 17: Uncompleted m-atlas near morphism

If $\boldsymbol{A}^1$ and $\boldsymbol{A}^2$ are semi-maximal (maximal) atlases then $\boldsymbol{f}$ is also a semi-maximal (maximal) $\boldsymbol{B}^1$-$\boldsymbol{B}^2$ bundle-atlas near morphism from $\boldsymbol{A}^1$ to $\boldsymbol{A}^2$, abbreviated $\text{isAtl}_{\text{Ar}}^{\text{Bun,near}}(\boldsymbol{A}^1, \boldsymbol{B}^1, \boldsymbol{A}^2, \boldsymbol{B}^2, \boldsymbol{f})$
S-max (max)

The identity morphism of $(\boldsymbol{A}^i, \boldsymbol{B}^i)$ is

$$\text{Id}_{(\boldsymbol{A}^i, \boldsymbol{B}^i)} \stackrel{\text{def}}{=} \Big( (\text{Id}_{E^i}, \text{Id}_{X^i}, \text{Id}_{Y^i}, \text{Id}_{G^i}), (\boldsymbol{A}^i, \boldsymbol{B}^i), (\boldsymbol{A}^i, \boldsymbol{B}^i) \Big) \tag{26.8}$$

This nomenclature will be justified later.

Figure 18: Completed m-atlas near morphism

**Definition 26.3** (Bundle-atlas morphisms). Let $\mathbf{B}^i \stackrel{\text{def}}{=} (E^i, X^i, Y^i, \pi^i, G^i, \rho^i), i = 1, 2$, be a protobundle and let $\mathbf{A}^i$ be a $\pi^i$-$G^i$-$\rho^i$-atlas of $E^i$ in the coordinate space $C^i = X^i \times Y^i$. Then $\mathbf{f} \stackrel{\text{def}}{=} (f_E \colon E^1 \longrightarrow E^2, f_X \colon X^1 \longrightarrow X^2, f_Y \colon Y^1 \longrightarrow Y^2, f_G \colon G^1 \longrightarrow G^2)$ is a $\mathbf{B}^1$-$\mathbf{B}^2$ bundle-atlas morphism from $\mathbf{A}^1$ to $\mathbf{A}^2$, abbreviated isAtl$_{\text{Ar}}^{\text{Bun}}(\mathbf{A}^1, \mathbf{B}^1, \mathbf{A}^2, \mathbf{B}^2, \mathbf{f})$, iff

1. all four functions are continuous

2. $f_G$ is a homomorphism

3. $\mathbf{f}$ commutes with $\pi^i$ and $\rho^i$, i.e.,

   (a) $\pi^2 \circ f_E = f_X \circ \pi^1$

   (b) $\left( \displaystyle\forall_{\substack{y \in Y^1 \\ g \in G^1}} \right) f_Y(y) \star^2 f_G(g) = f_Y(y \star^1 g)$

4. For any charts $(U^i, V^i, \phi^i \colon U^i \stackrel{\tilde{=}}{\rightarrowtail\!\!\!\twoheadrightarrow} V^i) \in \mathbf{A}^i$. $i = 1, 2$, and any $u^1 \in I \stackrel{\text{def}}{=} U^1 \cap f_E^{-1}[U^2]$, there exists a subchart $(U'^1, V'^1, \phi^1 \colon U'^1 \stackrel{\tilde{=}}{\rightarrowtail\!\!\!\twoheadrightarrow} V'^1) \in \mathbf{A}^1$ at $u^1$, an open set $U'^2 \subseteq U^2$ and a chart $(U'^2, \hat{V}'^2, \phi'^2 \colon U'^2 \stackrel{\tilde{=}}{\rightarrowtail\!\!\!\twoheadrightarrow} \hat{V}'^2) \in \mathbf{A}^2$ at $f_E(u^1)$ such that $f_E[U'^1] \subseteq U'^2$ and diagram (26.9) below is commutative, as shown in fig. 19.

$$D \stackrel{\text{def}}{=} \left\{ \begin{array}{c} i\colon \{u^1\} \rightarrowtail U'^1, \phi^1\colon U'^1 \stackrel{\cong}{\rightarrowtail\!\!\!\!\twoheadrightarrow} V'^1, f_X \times f_Y\colon V'^1 \longrightarrow \hat{V}'^2, \\ f_E\colon U'^1 \longrightarrow U'^2, \phi'^2\colon U'^2 \stackrel{\cong}{\rightarrowtail\!\!\!\!\twoheadrightarrow} \hat{V}'^2 \end{array} \right\} \qquad (26.9)$$

Figure 19: Completed Bundle-atlas morphism

If $A^1$ and $A^2$ are semi-maximal (maximal) atlases then $f$ is also a semi-maximal (maximal) $B^1$-$B^2$ bundle-atlas morphism from $A^1$ to $A^2$, abbreviated
$\text{isAtl}^{\text{Bun}}_{\text{Ar}}\, (A^1, B^1, A^2, B^2, f)$
S-max (max)

The identity morphism of $(A^i, B^i)$ is

$$\text{Id}_{(A^i, B^i)} \stackrel{\text{def}}{=} \Big( \big(\text{Id}_{E^i}, \text{Id}_{X^i}, \text{Id}_{Y^i}, \text{Id}_{G^i}\big), \big(A^i, B^i\big), \big(A^i, B^i\big) \Big) \qquad (26.10)$$

This nomenclature will be justified later.

**Definition 26.4** (Equivalence of bundle-atlas (near) morphisms). Let
$B^i \stackrel{\text{def}}{=} (E^i, X^i, Y^i, \pi^i, G^i, \rho^i), i = 1, 2$, be a protobundle, $A^i$ be a $\pi^i$-$G^i$-$\rho^i$-atlas of $E^i$ in the coordinate space $C^i = X^i \times Y^i$, and
$f \stackrel{\text{def}}{=} (f_E\colon E^1 \longrightarrow E^2, f_X\colon X^1 \longrightarrow X^2, f_Y\colon Y^1 \longrightarrow Y^2, f_G\colon G^1 \longrightarrow G^2)$ and
$g \stackrel{\text{def}}{=} (g_E\colon E^1 \longrightarrow E^2, g_X\colon X^1 \longrightarrow X^2, g_Y\colon Y^1 \longrightarrow Y^2, g_G\colon G^1 \longrightarrow G^2)$ be $B^1$-$B^2$ bundle-atlas morphisms from $A^1$ to $A^2$.

Then $f$ is bundle equivalent to $g$ as a $B^1$-$B^2$ bundle-atlas morphism from $A^1$ to $A^2$ iff $f_E = g_E, f_X = g_X$ and $f_G = g_G$,

By abuse of language we shall write $f$ is bundle equivalent to $g$ when the protobundles and bundle atlases are understood by context.

### 26.1.2 Proclamations on bundle-atlas (near) morphisms

**Lemma 26.5** (Bundle-atlas morphisms). *Let $B^i \overset{\text{def}}{=} (E^i, X^i, Y^i, \pi^i, G^i, \rho^i)$, $i = 1, 2$, be a protobundle and let $A^i$ be a $\pi^i$-$G^i$-$\rho^i$-atlas of $E^i$ in the coordinate space $C^i = X^i \times Y^i$. Then $\boldsymbol{f} \overset{\text{def}}{=} (f_E\colon E^1 \longrightarrow E^2, f_X\colon X^1 \longrightarrow X^2, f_Y\colon Y^1 \longrightarrow Y^2, f_G\colon G^1 \longrightarrow G^2)$ is a $\boldsymbol{B^1}$-$\boldsymbol{B^2}$ bundle-atlas (near) morphism from $A^1$ to $A^2$ iff $f_X \times f_Y$ is a $G^1$-$G^2$-$\rho^1$-$\rho^2$ morphism from $X^1, Y^1_{G^1-\rho^1-\text{triv}}$ to $X^2, Y^2_{G^2-\rho^2-\text{triv}}$ and $(f_E, f_X \times f_Y)$ is a $E^1_{\text{triv}}$-$E^2_{\text{triv}}$ m-atlas (near) morphism from $A^1$ to $A^2$ in the coordinate spaces $X^1, Y^1_{G^1-\rho^1-\text{triv}}$, $X^2, Y^2_{G^2-\rho^2-\text{triv}}$.*

*Proof.* If $\boldsymbol{f}$ is a $\boldsymbol{B^1}$-$\boldsymbol{B^2}$ bundle-atlas morphism then $f_X \times f_Y$ is a model function and $f_G$ is the function asserted to exist in eq. (22.3) (preservation of group action) of definition 22.3 (Morphisms of $G$-$\rho$-model spaces) on page 136, so $f_X \times f_Y$ is a $G^1$-$G^2$-$\pi^1$-$\pi^2$ morphism.

A diagram is locally nearly commutative in $X^2, Y^2, \pi^2, \rho^2$ iff it is m-locally nearly commutative in $\pi_2\left(X^2, Y^2_{G^2-\rho^2-\text{triv}}\right)$, thus $(f_E, f_X \times f_Y)$ is an m-atlas morphism in the coordinate spaces $X^1, Y^1_{G^1-\rho^1-\text{triv}}$, $X^2, Y^2_{G^2-\rho^2-\text{triv}}$.

If $f_X \times f_Y$ is a $G^1$-$G^2$-$\pi^1$-$\pi^2$ morphism from $X^1, Y^1_{G^1-\rho^1-\text{triv}}$ to $X^2, Y^2_{G^2-\rho^2-\text{triv}}$ then $\boldsymbol{f}$ commutes with $\rho^i$.

If $(f_E, f_X \times f_Y)$ is a $E^1_{\text{triv}}$-$E^2_{\text{triv}}$ m-atlas morphism from $A^1$ to $A^2$ in the coordinate spaces $X^1, Y^1_{G^1-\rho^1-\text{triv}}$, $X^2, Y^2_{G^2-\rho^2-\text{triv}}$ then $\boldsymbol{f}$ commutes with $\pi^i$.

A diagram is locally nearly commutative in $X^2, Y^2, \pi^2, \rho^2$ iff it is m-locally nearly commutative in $\pi_2\left(X^2, Y^2_{G^2-\rho^2-\text{triv}}\right)$ thus $(f_E, f_X \times f_Y)$ is a m-atlas morphism in the coordinate space $X^2, Y^2_{G^2-\rho^2-\text{triv}}$, so $\boldsymbol{f}$ is a $\boldsymbol{B^1}$-$\boldsymbol{B^2}$ bundle-atlas morphism from $A^1$ to $A^2$. □

**Corollary 26.6** (Bundle-atlas morphisms). *Let $B^i \overset{\text{def}}{=} (E^i, X^i, Y^i, \pi^i, G^i, \rho^i)$, $i = 1, 2, 3$, be a protobundle, $A^i$ be a $\pi^i$-$G^i$-$\rho^i$-atlas of $E^i$ in the coordinate space $C^i = X^i \times Y^i$ and $\boldsymbol{f}^i \overset{\text{def}}{=} (f_E^i\colon E^i \longrightarrow E^{i+1}, f_X^i\colon X^i \longrightarrow X^{i+1}, f_Y^i\colon Y^i \longrightarrow Y^{i+1}, f_G^i\colon G^i \longrightarrow G^{i+1})$ a $\boldsymbol{B^1}$-$\boldsymbol{B^2}$ bundle-atlas (near) morphism from $A^i$ to $A^{i+i}$.*

*If each $A^i$ is semi-maximal or if each $\boldsymbol{f}^i$ is a morphism then $\boldsymbol{f}^2 \overset{()}{\circ} \boldsymbol{f}^1$ is a bundle-atlas (near) morphism from $A^1$ to $A^3$.*

*Proof.* $f_X^1 \times f1_Y \circ f_X^2 \times f_Y^2$ is a $G^1$-$G^3$-$\pi^1$-$\pi^3$ morphism from $X^1, Y^1_{G^1-\rho^1-\text{triv}}$ to $X^3, Y^3_{G^3-\rho^3-\text{triv}}$ by item 4 of lemma 22.4 (Morphisms of $G$-$\rho$-model spaces) on page 137 and $(f_E^1, f_X^1 \times f_Y^1) \overset{()}{\circ} (f_E^2, f_X^2 \times f_Y^2)$ is a $E^1_{\text{triv}}$-$E^3_{\text{triv}}$ m-atlas morphism from $A^1$ to $A^3$ in the coordinate space $X^3, Y^3$ by item 2 of lemma 11.16 (Composition of m-atlas $G^3-\rho^3-\mathbf{triv}$ (near) morphisms) on page 72. □

**Lemma 26.7** (Composition of equivalent bundle-atlas (near) morphisms). *Let*

1. $\mathbf{B}^i \stackrel{\text{def}}{=} (E^i, X^i, Y^i, \pi^i, G^i, \rho^i)$, $i = 1, 2,$, *be a protobundle*

2. $\mathbf{A}^i$ *be a* $\pi^i$-$G^i$-$\rho^i$-*atlas of* $E^i$ *in the coordinate space* $C^i = X^i \times Y^i$

3. $\mathbf{f}^i \stackrel{\text{def}}{=} (f_E^i \colon E^i \longrightarrow E^{i+1}, f_X^i \colon X^i \longrightarrow X^{i+1}, f_Y^i \colon Y^i \longrightarrow Y^{i+1}, f_G^i \colon G^i \longrightarrow G^{i+1})$
   and $\mathbf{g}^i \stackrel{\text{def}}{=} (g_E^i \colon E^i \longrightarrow E^{i+1}, g_X^i \colon X^i \longrightarrow X^{i+1}, g_Y^i \colon Y^i \longrightarrow Y^{i+1}, g_G^i \colon G^i \longrightarrow G^{i+1})$
   *be bundle equivalent* $\mathbf{B}^i$-$\mathbf{B}^{i+1}$ *bundle-atlas (near) morphisms from* $\mathbf{A}^1$ *to* $\mathbf{A}^2$

*Then* $\mathbf{f}^{i+1} \stackrel{()}{\circ} \mathbf{f}^1$ *is* $C^k$-*equivalent to* $\mathbf{g}^2 \stackrel{()}{\circ} \mathbf{g}^1$.

*Proof.* $f_E^1 = g_E^1$, and $f_E^2 = g_E^2$, hence $f_E^2 \circ f_E^1 = g_E^2 \circ g_E^1$. $f_X^1 = g_X^1$, and $f_X^2 = g_X^2$, hence $f_X^2 \circ f_X^1 = g_X^2 \circ g_X^1$. $f_G^1 = g_G^1$, and $f_G^2 = g_G^2$, hence $f_G^2 \circ f_G^1 = g_G^2 \circ g_G^1$. □

## 26.2 Categories of bundle atlases and functors

This subsection only defines categories whose morphisms are bundle atlas morphisms; bundle atlas near morphisms do not satisfy the requiremens for forming a category.

**Definition 26.8** (Sets of bundle atlas morphisms). Let
$\mathbf{B}^\alpha \stackrel{\text{def}}{=} (E^\alpha, X^\alpha, Y^\alpha, \pi^\alpha, G^\alpha, \rho^\alpha)$, $\alpha \prec A$, be a protobundle and $\mathbf{B} \stackrel{\text{def}}{=} \{\mathbf{B}^\alpha | \alpha \prec A\}$ be a set of protobundles.

$$\mathcal{A}t\ell_{\substack{\text{Ar} \\ (\text{S-max,max})}}^{\text{Bun}}(\mathbf{B}^\alpha, \mathbf{B}^\beta) \stackrel{\text{def}}{=} \left\{ (\mathbf{f}, (\mathbf{A}^\alpha, \mathbf{B}^\alpha), (\mathbf{A}^\beta, \mathbf{B}^\beta)) \,\middle|\, \text{isAtl}_{\substack{\text{Ar} \\ (\text{S-max,max})}}^{\text{Bun}}(\mathbf{A}^\alpha, \mathbf{B}^\alpha, \mathbf{A}^\beta, \mathbf{B}^\beta, \mathbf{f}) \right\} \quad (26.11)$$

$$\mathcal{A}t\ell_{\substack{\text{Ar} \\ (\text{S-max,max})}}^{\text{Bun}}(\mathbf{B}^\alpha) \stackrel{\text{def}}{=} \left\{ (\mathbf{f}, (\mathbf{A}^\alpha, \mathbf{B}^\alpha), (\mathbf{A}^\alpha, \mathbf{B}^\alpha)) \,\middle|\, \text{isAtl}_{\substack{\text{Ar} \\ (\text{S-max,max})}}^{\text{Bun}}(\mathbf{A}^\alpha, \mathbf{B}^\alpha, \mathbf{A}^\alpha, \mathbf{B}^\alpha, \mathbf{f}) \right\} \quad (26.12)$$

$$\mathcal{A}t\ell_{\substack{\text{Ar} \\ (\text{S-max,max})}}^{\text{Bun}} \mathbf{B} \stackrel{\text{def}}{=} \bigcup_{\substack{\mathbf{B}^\mu \in \mathbf{B} \\ \mathbf{B}^\nu \in \mathbf{B}}} \mathcal{A}t\ell_{\substack{\text{Ar} \\ (\text{S-max,max})}}^{\text{Bun}} (\mathbf{B}_\mu, \mathbf{B}_\nu) \quad (26.13)$$

**Definition 26.9** (Categories of bundle atlases). Let
$\boldsymbol{B}^\alpha \stackrel{\text{def}}{=} (E^\alpha, X^\alpha, Y^\alpha, \pi^\alpha, G^\alpha, \rho^\alpha)$, $\alpha \prec A$, be a protobundle and $\boldsymbol{B} \stackrel{\text{def}}{=} \{\boldsymbol{B}^\alpha | \alpha \prec A\}$ be a set of protobundles. Then

$$\mathcal{A}t\ell^{\text{Bun}}_{\text{(S-max,max)}}(\boldsymbol{B}^\alpha) \stackrel{\text{def}}{=} \left( \mathcal{A}t\ell^{\text{Bun}}_{\substack{\text{Ob}\\\text{(S-max,max)}}}(\boldsymbol{B}^\alpha), \mathcal{A}t\ell^{\text{Bun}}_{\substack{\text{Ar}\\\text{(S-max,max)}}}(\boldsymbol{B}^\alpha, \boldsymbol{B}^\alpha), \stackrel{A}{\circ} \right) \qquad (26.14)$$

$$\mathcal{A}t\ell^{\text{Bun}}_{\text{(S-max,max)}} \boldsymbol{B} \stackrel{\text{def}}{=} \left( \mathcal{A}t\ell^{\text{Bun}}_{\substack{\text{Ob}\\\text{(S-max,max)}}}(\boldsymbol{B}), \mathcal{A}t\ell^{\text{Bun}}_{\substack{\text{Ar}\\\text{(S-max,max)}}}(\boldsymbol{B}), \stackrel{A}{\circ} \right) \qquad (26.15)$$

**Lemma 26.10** ($\mathcal{A}t\ell^{\text{Bun}} \boldsymbol{B}$ is a category). Let $\boldsymbol{B}^\alpha \stackrel{\text{def}}{=} (E^\alpha, X^\alpha, Y^\alpha, \pi^\alpha, G^\alpha, \rho^\alpha)$, $\alpha \prec A$, be a protobundle and $\boldsymbol{B} \stackrel{\text{def}}{=} \{\boldsymbol{B}^\alpha | \alpha \prec A\}$ be a set of protobundles. Then $\mathcal{A}t\ell^{\text{Bun}}_{\text{(S-max,max)}} \boldsymbol{B}$ is a category.

Let $(\boldsymbol{A}^\alpha, \boldsymbol{B}^\alpha) \in \mathcal{A}t\ell^{\text{Bun}}_{\text{Ob}} \boldsymbol{B}$. Then $\text{Id}_{(\boldsymbol{A}^\alpha, \boldsymbol{B}^\alpha)}$ is the identity morphism for $(\boldsymbol{A}^\alpha, \boldsymbol{B}^\alpha)$.

*Proof.* Let $(\boldsymbol{A}^i, \boldsymbol{B}^i)$, $i \in [1, 3]$, be an object of $\mathcal{A}t\ell^{\text{Bun}} \boldsymbol{B}$ and let $m^i \stackrel{\text{def}}{=} (\boldsymbol{f}^i, (\boldsymbol{A}^i, \boldsymbol{B}^i), (\boldsymbol{A}^{i+1}, \boldsymbol{B}^{i+1}))$ be a morphism of $\mathcal{A}t\ell^{\text{Bun}} \boldsymbol{B}$. Then

**Composition:** $(\boldsymbol{f}^2 \stackrel{()}{\circ} \boldsymbol{f}^1, (\boldsymbol{A}^1, \boldsymbol{E}^1, \boldsymbol{C}^1), (\boldsymbol{A}^3, \boldsymbol{E}^3, \boldsymbol{C}^3))$ is a morphism of $\mathcal{A}t\ell^{\text{Bun}} \boldsymbol{B}$ by corollary 26.6 (Bundle-atlas morphisms) on page 151.

**Associativity:** Composition is associative by lemma 1.19 (Tuple composition for labeled morphisms) on page 12.

**Identity:** $\text{Id}_{(\boldsymbol{A}^i, \boldsymbol{E}^i, \boldsymbol{C}^i)}$ is an identity morphism by lemma 1.19.

$\square$

**Definition 26.11** (Functor from Bundle atlases to m-atlases). Let $\boldsymbol{B}^i \stackrel{\text{def}}{=} (E^i, X^i, Y^i, \pi^i, G^i, \rho^i)$, $i = 1, 2$, be a protobundle, let $\boldsymbol{A}^i$ be a $\pi^i$-$G^i$-$\rho^i$-atlas of $E^i$ in the coordinate space $C^i = X^i \times Y^i$ and let
$\boldsymbol{f} \stackrel{\text{def}}{=} (f_E \colon E^1 \longrightarrow E^2, f_X \colon X^1 \longrightarrow X^2, f_Y \colon Y^1 \longrightarrow Y^2, f_G \colon G^1 \longrightarrow G^2)$ be a $\boldsymbol{B^1}$-$\boldsymbol{B^2}$ bundle-atlas morphism from $\boldsymbol{A}^1$ to $\boldsymbol{A}^2$. Then

$$\mathcal{F}^{\text{Bun}}_{\text{Bun,M}}(\boldsymbol{A}^i, \boldsymbol{B}^i) \stackrel{\text{def}}{=} \left( \boldsymbol{A}^i, \underset{\text{triv}}{E^i}, \underset{G^i - \rho^i - \text{triv}}{X^i, Y^i} \right) \qquad (26.16)$$

$$\mathcal{F}^{\text{Bun}}_{\text{Bun,M}}\big(\boldsymbol{f}, (\boldsymbol{A}^1, \boldsymbol{B}^1), (\boldsymbol{A}^2, \boldsymbol{B}^2)\big) \stackrel{\text{def}}{=}$$
$$\left( (f_E, f_X \times f_Y), \left(\boldsymbol{A}^1, \underset{\text{triv}}{E^1}, \underset{G^1-\rho^1-\text{triv}}{C^1}\right), \left(\boldsymbol{A}^2, \underset{\text{triv}}{E^2}, \underset{G^2-\rho^2-\text{triv}}{C^2}\right) \right) \qquad (26.17)$$

**Theorem 26.12** (Functor from Bundle atlases to m-atlases). *Let $E$ be a set of topological spaces, $B^\alpha \stackrel{\text{def}}{=} (E^\alpha \in E, X^\alpha, Y^\alpha, \pi^\alpha, G^\alpha, \rho^\alpha)$, $\alpha \prec A$, be a protobundle, $B \stackrel{\text{def}}{=} \{B^\alpha | \alpha \prec A\}$ be a set of protobundles and $C^\alpha \stackrel{\text{def}}{=} \underset{G^\alpha - \rho^\alpha - \text{triv}}{(X^\alpha, Y^\alpha)}$. Then $\mathscr{F}^{\text{Bun}}_{\text{Bun,M}}$ is a functor from $\underset{(\text{S-max,max})}{\mathscr{At}\ell^{\text{Bun}}} B$ to $\mathscr{At}\ell\left(\underset{\mathscr{T}riv}{E}, \underset{\underset{(\text{S-max,max})}{\text{Bun}-\mathscr{T}riv}}{B}\right)$ and a functor from $\underset{(\text{S-max,max})}{\mathscr{At}\ell^{\text{Bun}}} B$ to $\mathscr{At}\ell\left(\underset{\mathscr{T}riv}{E}, \underset{\underset{(\text{S-max,max})}{\text{Bun}-\text{prod}-\mathscr{T}riv}}{B}\right)$.*

*Proof.* Let $o^i \stackrel{\text{def}}{=} (A^i, B^i)$, $i \in [1, 3]$, be an object of $\mathscr{At}\ell^{\text{Bun}}(B)$, $o'^i \stackrel{\text{def}}{=} \mathscr{F}^{\text{Bun}}_{\text{Bun,M}} o^i = \left(A^i, \underset{\text{triv}}{E^i}, \underset{G^i - \rho^i - \text{triv}}{C^i}\right)$ be the corresponding object of $\mathscr{At}\ell\left(\underset{\mathscr{T}riv}{E}, \underset{\text{Bun}-\mathscr{T}riv}{B}\right)$, $m^i \stackrel{\text{def}}{=} (f^i, o^i, o^{i+1})$, $i = 1, 2$, be a morphism of $\mathscr{At}\ell^{\text{Bun}}(B)$ from $o^i$ to $o^{i+1}$ and

$$m'^i \stackrel{\text{def}}{=} \mathscr{F}^{\text{Bun}}_{\text{Bun,M}} m^i$$
$$= \left((f_E^i, f_X^i \times f_Y^i), \left(A^i, \underset{\text{triv}}{E^i}, \underset{G^i - \rho^i - \text{triv}}{C^i}\right), \left(A^{i+1}, \underset{\text{triv}}{E^{i+1}}, \underset{G^{i+1} - \rho^{i+1} - \text{triv}}{C^2}\right)\right)$$
$$= ((f_E^i, f_X^i \times f_Y^i), o'^i, o'^{i+1})$$

be the corresponding morphism of $\mathscr{At}\ell\left(\underset{\mathscr{T}riv}{E}, \underset{\text{Bun}-\mathscr{T}riv}{B}\right)$.

**Preservation of endpoints:** $\underset{\text{triv}}{E^i} \overset{\text{Ob}}{\in} \underset{\mathscr{T}riv}{E}$. $\underset{G^i - \rho^i - \text{triv}}{C^i} \overset{\text{Ar}}{\in} \underset{\text{Bun}-\text{prod}-\mathscr{T}riv}{B}$.

By theorem 26.12 (Functor from Bundle atlases to m-atlases) on page 154, $(f_E, f_X \times f_Y)$ is a $\underset{\text{triv}}{E^i} - \underset{\text{triv}}{E^{i+1}}$ m-atlas morphism from $A^i$ to $A^{i+1}$ in the coordinate spaces $\underset{G^i - \rho^i - \text{triv}}{X^i, Y^i}$, $\underset{G^{i+1} - \rho^{i+1} - \text{triv}}{X^{i+1}, Y^{i+1}}$.

$$\mathscr{F}^{\text{Bun}}_{\text{Bun,M}} m^i = \mathscr{F}^{\text{Bun}}_{\text{Bun,M}}(f^i, o^i, o^{i+1})$$
$$= ((f_E^i, f_X^i \times f_Y^i), o'^i, o'^{i+1})$$

**Identity:**

$$\mathscr{F}^{\text{Bun}}_{\text{Bun,M}} \text{Id}_{(A^i, B^i)} = \mathscr{F}^{\text{Bun}}_{\text{Bun,M}}((\text{Id}_{E^i}, \text{Id}_{X^i}, \text{Id}_{Y^i}, \text{Id}_{G^i}), o^i, o^i)$$
$$= ((\text{Id}_{\underset{\text{triv}}{E^i}}, \text{Id}_{C^i}), o'^i, o'^i)$$
$$= ((\text{Id}_{\underset{\text{triv}}{E^i}}, \text{Id}_{C^i}), \mathscr{F}^{\text{Bun}}_{\text{Bun,M}} o^i, \mathscr{F}^{\text{Bun}}_{\text{Bun,M}} o^i)$$
$$= \text{Id}_{\mathscr{F}^{\text{Bun}}_{\text{Bun,M}} o^i}$$

**Composition:**

$$m^2 \overset{A}{\circ} m^1 = ((f_0^2 \circ f_0^1, f_X^2 \circ f_X^1, f_Y^2 \circ f_Y^1, f_G^2 \circ f_G^1), o^1, o^3)$$

$$\mathcal{F}_{\text{Bun,M}}^{\text{Bun}}(m^2 \circ m^1) = ((f_E^2 \circ f_E^1, (f_X^2 \circ f_X^1) \times (f_Y^2 \circ f_Y^1)), o'^1, o'^3)$$

$$= ((f_E^2 \circ f_E^1, (f_X^2 \times f_Y^2) \circ (f_X^1 \times f_Y^1)), o'^1, o'^3)$$

$$= \mathcal{F}_{\text{Bun,M}}^{\text{Bun}} m^2 \overset{A}{\circ} \mathcal{F}_{\text{Bun,M}}^{\text{Bun}} m^1$$

□

**Lemma 26.13** (Base space functions derived from bundle-atlas morphisms). *Let $E^i$, $i = 1, 2$, be a model space, $X^i$, $Y^i$ topological spaces, $G^i$ a topological group, $\rho^i$ an effective action of $G^i$ on $Y^i$, $C^i \overset{\text{def}}{=} (X^i \times Y^i, \mathcal{XY}^i)$ a $G^i$-$\rho^i$ model space of $X^i \times Y^i$, $f_C \colon C^1 \longrightarrow C^2$ a $G^1$-$G^2$-$\rho^1$-$\rho^2$ morphism of $X^1 \times Y^1$ to $X^2 \times Y^2$, i.e., a model function that preserves group action, $A^i$ an m-atlas of $E^i$ in the coordinate model space $C^i$ and $f \overset{\text{def}}{=} (f_E \colon E^1 \longrightarrow E^2, f_C \colon C^1 \longrightarrow C^2)$ an $E^1$-$E^2$ m-atlas morphism of $A^1$ to $A^2$ in the coordinate spaces $C^1$, $C^2$.*

*Then there exists a unique function $f_X \colon X^1 \longrightarrow X^2$ such that $f_X \circ \pi_1 = \pi_1 \circ f_C$.*

*Proof.* Let $x$ in $X^1$, $y, y'$ in $Y^1$, Then $f_C$ preserves fibers, i.e., $\pi_1(f_C(x, y)) = \pi_1(f_C(x, y'))$. by lemma 22.4 (Morphisms of $G$-$\rho$-model spaces) on page 136. Define $f_X(x) \overset{\text{def}}{=} \pi_1(f_C(x, y))$ □

*Remark* 26.14. There need not exist $f_Y \colon Y^1 \longrightarrow Y^2$ such that $f_C = f_X \times f_Y$.

**Definition 26.15** (Functor from m-atlases to Bundle atlases). Let $\mathscr{E}$ be a trivial model category, $\mathscr{C}$ be a trivial product coordinate model category with objects $\left\{ C^\alpha \overset{\text{def}}{=} (X^\alpha \times Y^\alpha, \mathcal{XY}^\alpha) \,\middle|\, \alpha \prec A \right\}$, $G$ a group valued function on $\text{Ob}(\mathscr{C})$ and $\rho$ a function valued function on $\text{Ob}(\mathscr{C})$ such that for every $C^\alpha \overset{\text{def}}{=} (X^\alpha \times Y^\alpha, \mathcal{XY}^\alpha) \overset{\text{Ob}}{\in} \mathscr{C}$, $\rho(C^\alpha)$ is an effective action of $G(C^\alpha)$ on $Y^\alpha$ and $\underset{G(C^\alpha) - \rho(C^\alpha) - \text{triv}}{X^\alpha, Y^\alpha} = C^\alpha$.

Let $E^i \overset{\text{def}}{=} (E^i, \mathscr{E}^i) \overset{\text{Ob}}{\in} \mathscr{E}$, $i = 1, 2$, be a trivial model space, $C^i \overset{\text{def}}{=} (X^i \times Y^i, \mathcal{XY}^i) \overset{\text{Ob}}{\in} \mathscr{C}$ $G^i \overset{\text{def}}{=} G(C^i)$, $\rho^i \overset{\text{def}}{=} \rho(C^i)$, $A^i$ a full m-atlas of $E^i$ in the coordinate space $C^i$, $\pi^i \colon E^i \twoheadrightarrow X^i$ the unique function asserted in lemma 25.10 and $B^i \overset{\text{def}}{=} (E^i, X^i, Y^i, \pi^i, G^i, \rho^i)$. Then define

$$\mathcal{F}_{\text{M}-G-\rho,\text{Bun}}^{\text{Bun}}(A^i, E^i, C^i) \overset{\text{def}}{=} (A^i, B^i) \tag{26.18}$$

Let $f \overset{\text{def}}{=} (f_E \colon E^1 \longrightarrow E^2, f_C \overset{\text{def}}{=} f_X \times f_Y \colon C^1 \longrightarrow C^2)$ be an $E^1$-$E^2$ m-atlas morphism of $A^1$ to $A^2$ in the coordinate spaces $C^1$, $C^2$ that preserves the group action, $f_G \colon G^1 \longrightarrow G^2$ the unique function asserted to exist in eq. (22.3) (preservation of group action) of definition 22.3 (Morphisms of $G$-$\rho$-model spaces) on page 136.

Then

$$\mathscr{F}^{\mathrm{Bun}}_{\mathrm{M}-G-\rho,\mathrm{Bun}}\bigl(\boldsymbol{f},(\boldsymbol{A}^1,\mathbf{E^1},\boldsymbol{C}^1),(\boldsymbol{A}^2,\mathbf{E^2},\boldsymbol{C}^2)\bigr) \stackrel{\mathrm{def}}{=}$$
$$\bigl((f_E,f_X,f_Y,f_G),(\boldsymbol{A}^1,\boldsymbol{B}^1),(\boldsymbol{A}^2,\boldsymbol{B}^2)\bigr) \quad (26.19)$$

**Theorem 26.16** (Functor from m-atlases to bundle atlases). *Let*

1. $\mathscr{C}$ *be a trivial product coordinate model category,*

2. $\mathscr{E}$ *be a model category,*

3. $G$ *be a group valued function on* $\mathrm{Ob}(\mathscr{C})$

4. $\rho$ *be a function valued function on* $\mathrm{Ob}(\boldsymbol{C})$ *such that for every* $\boldsymbol{C}^\alpha \stackrel{\mathrm{def}}{=} (X^\alpha \times Y^\alpha, \mathcal{X}\mathcal{Y}^\alpha) \stackrel{\mathrm{Ob}}{\in} \mathscr{C}$, $\rho(\boldsymbol{C}^\alpha)$ *is an effective action of* $G(\boldsymbol{C}^\alpha)$ *on* $Y^\alpha$ *and* $\boldsymbol{C}^\alpha = \underset{G(\boldsymbol{C}^\alpha)-\rho(\boldsymbol{C}^\alpha)-\mathbf{triv}}{X^\alpha, Y^\alpha}$.

5. $\pi$ *be the unique function valued function on* $\mathcal{A}t\ell_{\mathrm{Ob}}^{\mathrm{full}}(\mathscr{E},\mathscr{C})$ *such that for every* $(\boldsymbol{A}^\alpha,(E^\alpha,\mathscr{E}^\alpha),(C^\alpha,\mathscr{C}^\alpha)) \in \mathcal{A}t\ell_{\mathrm{Ob}}^{\mathrm{full}}(\mathscr{E},\mathscr{C})$ *and every* $(U,V,\phi) \in \boldsymbol{A}^\alpha$, $\pi_1 \circ \phi = \pi(\boldsymbol{A}^\alpha)\!\upharpoonright_U$

6. $\boldsymbol{B} \stackrel{\mathrm{def}}{=} \left\{ \boldsymbol{B}^\alpha \stackrel{\mathrm{def}}{=} (E^\alpha, X^\alpha, Y^\alpha, G(\boldsymbol{C}^\alpha), \pi(\boldsymbol{A}^\alpha), \rho(\boldsymbol{C}^\alpha)) \,\Big|\, \right.$
$\left. (\boldsymbol{A}^\alpha,(E^\alpha,\mathscr{E}^\alpha),(C^\alpha \stackrel{\mathrm{def}}{=} X^\alpha \times Y^\alpha, \mathscr{C}^\alpha)) \in \mathcal{A}t\ell_{\mathrm{Ob}}^{\mathrm{full}}(\mathscr{E},\mathscr{C}) \right\}$

*Then* $\mathscr{F}^{\mathrm{Bun}}_{\mathrm{M}-G-\rho,\mathrm{Bun}}$ *is a functor from* $\mathcal{A}t\ell\left(\underset{\mathscr{T}\!\boldsymbol{riv}}{\boldsymbol{E}}, \underset{\mathrm{Bun}-\mathrm{prod}-\mathscr{T}\!\boldsymbol{riv}}{\boldsymbol{B}}\right)$ *to* $\mathcal{A}t\ell^{\mathrm{Bun}}(\boldsymbol{B})$.

*Proof.* Let $o^i \stackrel{\mathrm{def}}{=} (\boldsymbol{A}^i, \boldsymbol{E}^i, \boldsymbol{C}^i \stackrel{\mathrm{def}}{=} (X^i \times Y^i, \mathcal{X}\mathcal{Y}^i)) \stackrel{\mathrm{Ob}}{\in} \mathcal{M})$, $i \in [1,3]$, $G^i \stackrel{\mathrm{def}}{=} G(o^i)$, $\pi^i \stackrel{\mathrm{def}}{=} \pi(o^i)$, $\rho^i \stackrel{\mathrm{def}}{=} \rho(o^i)$, $\boldsymbol{A}^i$ an m-atlas of $E^i$ in the coordinate model space $\boldsymbol{C}^i$, $\boldsymbol{B}^i \stackrel{\mathrm{def}}{=} (E^i, X^i, Y^i, \pi^i, G^i, \rho^i)$, $o'^i \stackrel{\mathrm{def}}{=} (\boldsymbol{A}^i, \boldsymbol{B}^i)$, $m^i \stackrel{\mathrm{def}}{=} (\boldsymbol{f}^i, o^i, o^{i+1}) \stackrel{\mathrm{Ar}}{\in} \mathcal{M}$, $i = 1,2$, an $E^i$-$E^{i+1}$ m-atlas morphism of $\boldsymbol{A}^i$ to $\boldsymbol{A}^{i+1}$ in the coordinate spaces $\boldsymbol{C}^i, \boldsymbol{C}^{i+1}$ that preserves the group action, $f_G \colon G^i \longrightarrow G^{i+1}$ the unique function asserted to exist in eq. (22.3) (preservation of group action) of definition 22.3 (Morphisms of $G$-$\rho$-model spaces) on page 136 and $\pi^i \colon E^i \twoheadrightarrow X^i$ the unique function asserted in lemma 25.10.

Let $f_G^i \colon G^i \longrightarrow G^{i+1}$ be the function asserted by eq. (22.3) (Morphisms of $G$-$\rho$-model spaces) on page 136; let $f_1^i = f_X^i \times f_Y^i$ be the decomposition given by definition 22.7 (Trivial $G$-$\rho$-model spaces) on page 138, and let $m'^i \stackrel{\mathrm{def}}{=} \mathscr{F}^{\mathrm{Bun}}_{\mathrm{M}-G-\rho,\mathrm{Bun}} m^i = (\boldsymbol{f}^i, (f_0^i f_X^i f_Y^i f_G^i), o'^i, o'^{i+1})$.

**Preservation of endpoints:** $E^i, X^i$ and $Y^i$ are topological spaces.

$G^i$ is a topological group,

$\pi^i \colon E^i \twoheadrightarrow X^i$ is a continuous surjection.

$\rho^i \stackrel{\text{def}}{=} \rho(C^i)$ is an effective action of $G^i \stackrel{\text{def}}{=} G(C^i)$ on $Y^\alpha$.

Hence $B^i$ is a protobundle.

**Identity:**

$$\mathcal{F}^{\text{Bun}}_{\text{M}-G-\rho,\text{Bun}} \,\text{Id}_{(A^i,B^i)} = \text{Id}_{(A^i,B^i)}\big((\text{Id}_{E^i},\text{Id}_{X^i},\text{Id}_{Y^i},\text{Id}_{G^i}), o^i, o^i\big)$$

$$= \big((\text{Id}_{\underset{\text{triv}}{E^i}},\text{Id}_{C^i}), o'^i, o'^i\big)$$

$$= \Big((\text{Id}_{\underset{\text{triv}}{E^i}},\text{Id}_{C^i}), \mathcal{F}^{\text{Bun}}_{\text{M}-G-\rho,\text{Bun}} o^i, \mathcal{F}^{\text{Bun}}_{\text{M}-G-\rho,\text{Bun}} o^i\Big)$$

$$= \text{Id}_{\mathcal{F}^{\text{Bun}}_{\text{M}-G-\rho,\text{Bun}} o^i}$$

**Composition:**

$\mathcal{F}^{\text{Bun}}_{\text{M}-G-\rho,\text{Bun}} m^i$ is a morphism from $\mathcal{F}^{\text{Bun}}_{\text{M}-G-\rho,\text{Bun}} o^i$ to $\mathcal{F}^{\text{Bun}}_{\text{M}-G-\rho,\text{Bun}} o^{i+1}$:

$$\mathcal{F}^{\text{Bun}}_{\text{M}-G-\rho,\text{Bun}}(m^i) = \mathcal{F}^{\text{Bun}}_{\text{M}-G-\rho,\text{Bun}}(f^i, o^i, o^{i+1})$$

$$= \big((f^i_0, f^i_X, f^i_Y, f^i_G),(A^i,B^i),(A^{i+1},B^{i+1})\big)$$

$$= \Big((f^i_0, f^i_X, f^i_Y, f^i_G), \mathcal{F}^{\text{Bun}}_{\text{M}-G-\rho,\text{Bun}} o^i, \mathcal{F}^{\text{Bun}}_{\text{M}-G-\rho,\text{Bun}} o^{i+1}\Big)$$

$\mathcal{F}^{\text{Bun}}_{\text{M}-G-\rho,\text{Bun}}$ maps identity functions to identity functions:

$$\begin{aligned}
&\mathcal{F}^{\text{Bun}}_{\text{M}-G-\rho,\text{Bun}} \,\text{Id}_{o^i} = \\
&\mathcal{F}^{\text{Bun}}_{\text{M}-G-\rho,\text{Bun}}\big((\text{Id}_{E^i},\text{Id}_{C^i}), o^i, o^i\big) = \\
&\big((\text{Id}_{E^i},\text{Id}_{X^i},\text{Id}_{Y^i},\text{Id}_{G^i}),(A^i,B^i),(A^i,B^i)\big) = \\
&\big((\text{Id}_{E^i},\text{Id}_{X^i},\text{Id}_{Y^i},\text{Id}_{G^i}), \mathcal{F}^{\text{Bun}}_{\text{M}-G-\rho,\text{Bun}} o^i, \mathcal{F}^{\text{Bun}}_{\text{M}-G-\rho,\text{Bun}} o^i\big) = \\
&\text{Id}_{\mathcal{F}^{\text{Bun}}_{\text{M}-G-\rho,\text{Bun}} o^i}
\end{aligned} \qquad (26.20)$$

$\mathcal{F}^{\text{Bun}}_{\text{M}-G-\rho,\text{Bun}} m^2 \stackrel{A}{\circ} \mathcal{F}^{\text{Bun}}_{\text{M}-G-\rho,\text{Bun}} m^1 = \mathcal{F}^{\text{Bun}}_{\text{M}-G-\rho,\text{Bun}}(m^2 \stackrel{A}{\circ} m^1)$:

1. $m^2 \stackrel{A}{\circ} m^1 = \big((f^2_0 \circ f^1_0, f^2_1 \circ f^1_1), o^1, o^3\big)$

2. $\mathcal{F}^{\text{Bun}}_{\text{M}-G-\rho,\text{Bun}}(o^i) = (A^i, B^i)$

3. $\mathcal{F}^{\text{Bun}}_{\text{M}-G-\rho,\text{Bun}}(m^i) = $
   $\big((f^i_0, f^i_X, f^i_Y, f^i_G),(A^i,B^i),(A^{i+1},B^{i+1})\big)$

4. $\mathcal{F}^{\text{Bun}}_{\text{M}-G-\rho,\text{Bun}} m^2 \stackrel{A}{\circ} \mathcal{F}^{\text{Bun}}_{\text{M}-G-\rho,\text{Bun}} m^1 = $
   $\big((f^2_0 \circ f^1_0, f^2_X \circ f^1_X, f^2_Y \circ f^1_Y, f^2_G \circ f^1_G),(A^1,B^1),(A^3,B^3)\big)$

5. $\mathcal{F}^{\text{Bun}}_{\text{M}-G-\rho,\text{Bun}}(m^2 \circ m^1) =$
   $((f_0^2 \circ f_0^1, f_X^2 \circ f_X^1, f_Y^2 \circ f_Y^1, f_G^2 \circ f_G^1), (\boldsymbol{A}^1, \boldsymbol{B}^1), (\boldsymbol{A}^3, \boldsymbol{B}^3))$

$\square$

$\mathcal{F}^{\text{Bun}}_{\text{Bun,M}} \circ \mathcal{F}^{\text{Bun}}_{\text{M}-G-\rho,\text{Bun}} = \text{Id}.$

*Proof.* Expanding the definitions, we have

1. $\mathcal{F}^{\text{Bun}}_{\text{M}-G-\rho,\text{Bun}}(\boldsymbol{A}^i, \boldsymbol{E}^i, \boldsymbol{C}^i) = (\boldsymbol{A}^i, \boldsymbol{B}^i)$

2. $\mathcal{F}^{\text{Bun}}_{\text{Bun,M}}(\boldsymbol{A}^i, \boldsymbol{B}^i) = (\boldsymbol{A}^i, \underset{\textbf{triv}}{\boldsymbol{E}^i}, \underset{\textbf{G}^i-\rho^i-\textbf{triv}}{\boldsymbol{X}^i, \boldsymbol{Y}^i})(\boldsymbol{A}^i, \boldsymbol{E}^i, \boldsymbol{C}^i)$, since by definition 26.15, $\boldsymbol{E}^i$ and $\boldsymbol{C}^i$ are trivial.

3. $\mathcal{F}^{\text{Bun}}_{\text{M}-G-\rho,\text{Bun}}((\boldsymbol{f}_E^1, \boldsymbol{f}_X^1 \times \boldsymbol{f}_Y^1), (\boldsymbol{A}^1, \boldsymbol{E}^1, \boldsymbol{C}^1), (\boldsymbol{A}^2, \boldsymbol{E}^2, \boldsymbol{C}^2)) =$
   $((f_E^1, f_X^1, f_Y^1, f_G^1), (\boldsymbol{A}^1, \boldsymbol{B}^1), (\boldsymbol{A}^2, \boldsymbol{B}^2))$

4. $\mathcal{F}^{\text{Bun}}_{\text{Bun,M}}((f_E^1, f_X^1, f_Y^1, f_G^1), (\boldsymbol{A}^1, \boldsymbol{B}^1), (\boldsymbol{A}^2, \boldsymbol{B}^2)) =$
   $((\boldsymbol{f}_E^1, \boldsymbol{f}_X^1 \times \boldsymbol{f}_Y^1), (\boldsymbol{A}^1, \underset{\textbf{triv}}{\boldsymbol{E}^1}, \underset{\textbf{G}^1-\rho^1-\textbf{triv}}{\boldsymbol{C}^1}), (\boldsymbol{A}^2, \underset{\textbf{triv}}{\boldsymbol{E}^2}, \underset{\textbf{G}^2-\rho^2-\textbf{triv}}{\boldsymbol{C}^2})) =$
   $((\boldsymbol{f}_E^1, \boldsymbol{f}_X^1 \times \boldsymbol{f}_Y^1), (\boldsymbol{A}^1, \boldsymbol{E}^1, \boldsymbol{C}^1), (\boldsymbol{A}^2, \boldsymbol{E}^2, \boldsymbol{C}^2))$, since by definition 26.15, $\boldsymbol{E}^i$ and $\boldsymbol{C}^i$ are trivial.

$\square$

## 27 Associated model spaces and functors

**Definition 27.1** (Coordinate model spaces associated with bundle atlases)**.** Let $\boldsymbol{B}^i \overset{\text{def}}{=} (E^i, X^i, Y^i, \pi^i, G^i, \rho^i)$ be a protobundle $\boldsymbol{A}^i$ be a $\pi^i$-$G^i$-$\rho^i$-bundle atlas of $E^i$ in the coordinate space $C^i = X^i \times Y^i$ and $\boldsymbol{f}$
$\overset{\text{def}}{=} (f_E \colon E^1 \longrightarrow E^2, f_X \colon X^1 \longrightarrow X^2, f_Y \colon Y^1 \longrightarrow Y^2, f_G \colon G^1 \longrightarrow G^2)$ be a $\boldsymbol{B}^1$-$\boldsymbol{B}^2$ bundle-atlas (near) morphism from $\boldsymbol{A}^1$ to $\boldsymbol{A}^2$. Then:

The minimal $G$-$\rho$ coordinate model space with neighborhoods in $\boldsymbol{A}^i$ is $\mathcal{F}_2^{\text{Bun}}(\boldsymbol{A}^i, \boldsymbol{B}^i)$.

$$\mathcal{F}_2^{\text{Bun}}(\boldsymbol{A}^i, \boldsymbol{B}^i) \overset{\text{def}}{=}$$

$$\underset{\min}{\text{Mod}}\left(C^i, \pi_2[\boldsymbol{A}^i], \left\{\phi' \circ \phi^{-1} \middle| \left(\underset{(U',V',\phi') \in \boldsymbol{A}^i}{\exists_{(U,V,\phi) \in \boldsymbol{A}^i}}\right) U \cap U' \neq \emptyset\right\}\right) \quad (27.1)$$

The coordinate mapping associated with the $\boldsymbol{B}^1$-$\boldsymbol{B}^2$ bundle-atlas (near) morphism $\boldsymbol{f}$ from $\boldsymbol{A}^1$ to $\boldsymbol{A}^2$ is

$$\mathcal{F}_2^{\text{Bun}}\bigl(\boldsymbol{f},(\boldsymbol{A}^1,\boldsymbol{B}^1),(\boldsymbol{A}^2,\boldsymbol{B}^2)\bigr) \overset{\text{def}}{=} f_X \times f_Y \colon \mathcal{F}_2^{\text{Bun}}(\boldsymbol{A}^1,\boldsymbol{B}^1) \longrightarrow \mathcal{F}_2^{\text{Bun}}(\boldsymbol{A}^2,\boldsymbol{B}^2) \quad (27.2)$$

If it is a model function then it is also the coordinate $G^1$-$G^2$-$\rho^1$-$\rho^2$ (near) morphism associated with the $\boldsymbol{B}^1$-$\boldsymbol{B}^2$ bundle-atlas (near) morphism $\boldsymbol{f}$ from $\boldsymbol{A}^1$ to $\boldsymbol{A}^2$.

**Lemma 27.2** (Coordinate model spaces associated with bundle atlases). *Let $\boldsymbol{B} \overset{\text{def}}{=} (E,X,Y,G,\pi,\rho)$ and let $\boldsymbol{A}$ be a $\pi$-$G$-$\rho$-bundle atlas of $E$ in the coordinate space $C = X \times Y$. Then $\mathcal{F}_2^{\text{Bun}}(\boldsymbol{A},\boldsymbol{B})$ is a model space.*

*Proof.* $\mathcal{F}_2^{\text{Bun}}(\boldsymbol{A},\boldsymbol{B})$ satisfies the conditions for a model space. for a model space. Let $\mathscr{C} \overset{\text{def}}{=} \text{Cat}\bigl(\mathcal{F}_2^{\text{Bun}}(\boldsymbol{A},\boldsymbol{B})\bigr)$.

1. Since $\pi_2[\boldsymbol{A}]$ is an open cover of $\bigcup \pi_2[\boldsymbol{A}]$, the set of finite intersections is also an open cover.

2. Finite intersections of finite intersections are finite intersections

3. Restrictions of continuous functions are continuous

4. If $f\colon A \longrightarrow B$ is a morphism of $\mathcal{F}_2^{\text{Bun}}(\boldsymbol{A},\boldsymbol{B})$ $A, A', B, B'$ objects of $\mathscr{C}$ $A' \subseteq A$, $B' \subseteq B$ and $f[A'] \subseteq B'$ then since $f\colon A \longrightarrow B$ is a morphism it is a restriction of a transition function between its restrictions to sets in $\pi_2[\boldsymbol{A}]$ and its restrictions are also, hence morphisms, and thus $f\!\restriction_{A'}\, A' \longrightarrow B'$ is a morphism.

5. If $(U,V,\phi) \in \boldsymbol{A}$ then $\text{Id}_V = \phi \circ \phi^{-1}A$ is a transition function and hence a morphism of $\mathcal{F}_2^{\text{Bun}}(\boldsymbol{A},\boldsymbol{B})$. If $A, A'$ objects of $\mathscr{C}$ and $A' \subseteq A$ then the inclusion map $i\colon A' \hookrightarrow A$ is a restriction of an identity morphism of $\mathcal{F}_2^{\text{Bun}}(\boldsymbol{A},\boldsymbol{B})$ and hence a morphism.

6. Restricted sheaf condition: let

    (a) $U_\alpha, V_\alpha, \alpha \prec A$, be an object of $\mathscr{C}$
    
    (b) $f_\alpha\colon U_\alpha \longrightarrow V_\alpha$ be a morphism of $\mathscr{C}$
    
    (c) $U \overset{\text{def}}{=} \bigcup_{\alpha \prec A} U_\alpha$
    
    (d) $V \overset{\text{def}}{=} \bigcup_{\alpha \prec A} V_\alpha$
    
    (e) $f\colon U \longrightarrow V$ be continuous and $\left(\forall_{\substack{\alpha \prec A \\ x \in U_\alpha}}\right) f(x) = f_\alpha(x)$

    Then $f$ is generated by the group action and hence a morphism of $\mathscr{C}$

□

**Definition 27.3** (Model spaces associated with bundle atlases). Let
$\boldsymbol{B}^i \stackrel{\text{def}}{=} (E^i, X^i, Y^i, \pi^i, G^i, \rho^i)$, $i = 1, 2$, and $\boldsymbol{A}^i$ be a $\pi^i$-$G^i$-$\rho^i$-bundle atlas of $E^i$ in the coordinate space $C^i = X^i \times Y^i$. Then

The minimal model space with neighborhoods in the atlas $\boldsymbol{A}^i$ is

$$\mathcal{F}_1^{\text{Bun}}(\boldsymbol{B}^i, \boldsymbol{A}^i) \stackrel{\text{def}}{=}$$

$$\underset{\min}{\text{Mod}}\left(E^i, \pi_1[\boldsymbol{A}^i], \left\{\phi'^{-1} \circ \phi \,\middle|\, \left(\exists_{\substack{(U,V,\phi)\in \boldsymbol{A}^i \\ (U',V',\phi')\in \boldsymbol{A}^i}}\right) U \cap U' \neq \emptyset\right\}\right) \quad (27.3)$$

The mapping associated with the $\boldsymbol{B}^1$-$\boldsymbol{B}^2$ bundle-atlas (near) morphism $\boldsymbol{f}$ from $\boldsymbol{A}^1$ to $\boldsymbol{A}^2$ is

$$\mathcal{F}_1^{\text{Bun}}(\boldsymbol{f}, (\boldsymbol{A}^1, \boldsymbol{B}^1), (\boldsymbol{A}^2, \boldsymbol{B}^2)) \stackrel{\text{def}}{=} f_E \colon \mathcal{F}_1^{\text{Bun}}(\boldsymbol{A}^1, \boldsymbol{B}^1) \longrightarrow \mathcal{F}_1^{\text{Bun}}(\boldsymbol{A}^2, \boldsymbol{B}^2) \quad (27.4)$$

If it is a model function then it is also the m-atlas (near) morphism associated with the $\boldsymbol{B}^1$-$\boldsymbol{B}^2$ bundle-atlas (near) morphism $\boldsymbol{f}$ from $\boldsymbol{A}^1$ to $\boldsymbol{A}^2$.

**Lemma 27.4** (Model spaces associated with bundle atlases). *Let*
$\boldsymbol{B} \stackrel{\text{def}}{=} (E, X, Y, G, \pi, \rho)$ *and let $\boldsymbol{A}$ be a $\pi$-$G$-$\rho$-bundle atlas of $E$ in the coordinate space $C = X \times Y$. Then $\mathcal{F}_1^{\text{Bun}}(\boldsymbol{A}, \boldsymbol{B})$ is a model space.*

*Proof.* Lemma 5.4 (Minimal model spaces are model spaces) on page 27 □

**Theorem 27.5** (Functors from bundle atlases to model spaces). *Let*
$\boldsymbol{B}^\alpha \stackrel{\text{def}}{=} (E^\alpha, X^\alpha, Y^\alpha, \pi^\alpha, G^\alpha, \rho^\alpha)$, $\alpha \prec A$, *be a protobundle, and $\boldsymbol{B} \stackrel{\text{def}}{=} \{\boldsymbol{B}^\alpha | \alpha \prec A\}$ be a set of protobundles. Then $\mathcal{F}_1^{\text{Bun}}$ is a functor from $\mathcal{A}t\ell^{\text{Bun}}\boldsymbol{B}$ to $\underset{\text{triv}}{\boldsymbol{E}}$ and $\mathcal{F}_2^{\text{Bun}}$ is a functor from $\mathcal{A}t\ell^{\text{Bun}}\boldsymbol{B}$ to $\underset{\text{Bun-}\mathcal{T}riv}{\boldsymbol{B}}$.*

*Proof.* Let $o^i \stackrel{\text{def}}{=} (\boldsymbol{A}^i, E^i, C^i)$, $i \in [1, 3]$, be objects in $\text{Bun}(E, X, Y, \pi, G, \rho)$ and let $m^i \stackrel{\text{def}}{=} ((f_0^i, f_1^i), o^i, o^{i+1})$ be a morphism in $\mathcal{A}t\ell^{\text{Bun}}(\boldsymbol{E}, \boldsymbol{C})$.

$\mathcal{F}_1^{\text{Bun}}$ is a functor from $\mathcal{A}t\ell^{\text{Bun}}\boldsymbol{B}$ to $\underset{\text{triv}}{\boldsymbol{E}}$.

**Preservation of endpoints:**

$$\mathcal{F}_2^{\text{Bun}} m^i = \mathcal{F}_1^{\text{Bun}}(\boldsymbol{f}^i, o^i, o^{i+1})$$
$$= f_E^i \colon \mathcal{F}_1^{\text{Bun}} o^i \longrightarrow o^{i+1}$$

by definition 27.3 (Model spaces associated with bundle atlases) on page 160
.

**Composition:**

$$\mathcal{F}_1^{\text{Bun}}(m^2 \overset{A}{\circ} m^1) = \mathcal{F}_1^{\text{Bun}}(\boldsymbol{f}^2 \circ \boldsymbol{f}^1(\boldsymbol{A}^1, E^1, C^1), (\boldsymbol{A}^3, E^3, C^3))$$
$$= f_E^2 \circ f_E^1 \colon \mathcal{F}_1^{\text{Bun}} o^1 \longrightarrow \mathcal{F}_1^{\text{Bun}} o^3$$
$$= (f_E^2 \colon \mathcal{F}_1^{\text{Bun}} o^2 \longrightarrow \mathcal{F}_1^{\text{Bun}} o^3) \circ$$
$$(f_E^1 \colon \mathcal{F}_1^{\text{Bun}} o^1 \longrightarrow \mathcal{F}_1^{\text{Bun}} o^2)$$
$$= \mathcal{F}_1^{\text{Bun}}(\boldsymbol{f}^2, (\boldsymbol{A}^2, E^2, C^2), (\boldsymbol{A}^3, E^3, C^3)) \circ$$
$$\mathcal{F}_1^{\text{Bun}}(\boldsymbol{f}^1, (\boldsymbol{A}^1, E^1, C^1), (\boldsymbol{A}^2, E^2, C^2))$$
$$= \mathcal{F}_1^{\text{Bun}} m^2 \circ \mathcal{F}_1^{\text{Bun}} m^1$$

**Identity:**

1. $\mathcal{F}_1^{\text{Bun}}(\text{Id}_{o^i}) = \mathcal{F}_1^{\text{Bun}}((\text{Id}_{E^i}, \text{Id}_{C^i}), (\boldsymbol{A}^i, E^i, C^i), (\boldsymbol{A}^i, E^i, C^i)) = \text{Id}_{E^i} \colon \mathcal{F}_1^{\text{Bun}} o^i \longrightarrow \mathcal{F}_1^{\text{Bun}} o^i$
2. $\text{Id}_{\mathcal{F}_1^{\text{Bun}} o^i} = \text{Id}_{E^i} \colon \mathcal{F}_1^{\text{Bun}} o^i \longrightarrow \mathcal{F}_1^{\text{Bun}} o^i$

$\mathcal{F}_2^{\text{Bun}}$ is a functor from $\mathcal{A}t\ell^{\text{Bun}}\boldsymbol{B}$ to $\underset{\text{Bun}-\mathcal{T}riv}{\boldsymbol{B}}$.

**Preservation of endpoints:** $\mathcal{F}_2^{\text{Bun}}(m^i) = f_1^i \colon \mathcal{F}_2^{\text{Bun}} o^i \longrightarrow \mathcal{F}_2^{\text{Bun}} o^{i+1}$

**Composition:**

$$\mathcal{F}_2^{\text{Bun}}(m^2 \overset{A}{\circ} m^1) = \mathcal{F}_2^{\text{Bun}}((f_0^2 \circ f_0^1, f_1^2 \circ f_1^1)(\boldsymbol{A}^1, E^1, C^1), (\boldsymbol{A}^3, E^3, C^3))$$
$$= f_0^2 \circ f_0^1 \colon \mathcal{F}_2^{\text{Bun}} o^1 \longrightarrow \mathcal{F}_2^{\text{Bun}} o^3$$
$$= (f_1^2 \colon \mathcal{F}_2^{\text{Bun}} o^2 \longrightarrow \mathcal{F}_2^{\text{Bun}} o^3) \circ (f_1^1 \colon \mathcal{F}_2^{\text{Bun}} o^1 \longrightarrow \mathcal{F}_2^{\text{Bun}} o^2)$$
$$= \mathcal{F}_2^{\text{Bun}}((f_0^2, f_1^2)(\boldsymbol{A}^2, E^2, C^2), (\boldsymbol{A}^3, E^3, C^3)) \circ$$
$$\mathcal{F}_2^{\text{Bun}}((f_0^1, f_1^1)(\boldsymbol{A}^1, E^1, C^1), (\boldsymbol{A}^2, E^2, C^2))$$
$$= \mathcal{F}_2^{\text{Bun}}(m^2) \circ \mathcal{F}_2^{\text{Bun}}(m^1)$$

**Identity:**

1. $\mathcal{F}_2^{\text{Bun}} \text{Id}_{o^i} = \mathcal{F}_2^{\text{Bun}}((\text{Id}_{E^i}, \text{Id}_{C^i}), (\boldsymbol{A}^i, E^i, C^i), (\boldsymbol{A}^i, E^i, C^i)) = \text{Id}_{C^i} \colon \mathcal{F}_2^{\text{Bun}} o^i \longrightarrow \mathcal{F}_2^{\text{Bun}} o^i$
2. $\text{Id}_{\mathcal{F}_2^{\text{Bun}} o^i} = \text{Id}_{C^i} \colon \mathcal{F}_2^{\text{Bun}} o^i \longrightarrow \mathcal{F}_2^{\text{Bun}} o^i$

$\square$

# 28 Fiber bundles

Conventionally a fiber bundle is different from its atlases, but definition 26.9 (Categories of bundle atlases) on page 153 encourages treating them on an equal footing. All of the results for maximal bundle atlases carry directly over to results for fiber bundles.

## 28.1 Definitions

**Definition 28.1** (fiber bundles)**.** Let $E, X$ and $Y$ be topological spaces, $\pi\colon E \twoheadrightarrow X$ surjective, $G$ a topological group, $\rho\colon Y \times G \longrightarrow Y$ an effective right action of $G$ on $Y$ and $A$ a maximal $\pi$-$G$-$\rho$-bundle atlas of $E$ in the coordinate space $X \times Y$. Then $(E, X, Y, \pi, G, \rho, A)$ is a fiber bundle.

Let $\boldsymbol{B} \stackrel{\text{def}}{=} \left\{\boldsymbol{B}^\alpha \stackrel{\text{def}}{=} (E^\alpha, X^\alpha, Y^\alpha, \pi^\alpha, G^\alpha, \rho^\alpha) \,\middle|\, \alpha \prec A\right\}$, where $E^\alpha, X^\alpha, Y^\alpha$ are topological spaces, $G^\alpha$ a topological group, $\pi^\alpha\colon E^\alpha \twoheadrightarrow X^\alpha$ surjective and $\rho^\alpha\colon Y^\alpha \times G^\alpha \longrightarrow Y^\alpha$ an effective right action of $G^\alpha$ on $Y^\alpha$. Then

$$\mathrm{Bun}_{\mathrm{Ob}}\boldsymbol{B} \stackrel{\text{def}}{=} \mathcal{A}t\ell^{\mathrm{Bun}}_{\substack{\mathrm{Ob}\\ \max}}\boldsymbol{B} \quad (28.1)$$

**Definition 28.2** (Bundle maps)**.** Let
$\boldsymbol{B} \stackrel{\text{def}}{=} \left\{\boldsymbol{B}^\alpha \stackrel{\text{def}}{=} (E^\alpha, X^\alpha, Y^\alpha, \pi^\alpha, G^\alpha, \rho^\alpha) \,\middle|\, \alpha \prec A\right\}$, where $E^\alpha, X^\alpha, Y^\alpha$ are topological spaces, $G^\alpha$ a topological group, $\pi^\alpha\colon E^\alpha \twoheadrightarrow X^\alpha$ surjective and $\rho^\alpha\colon Y^\alpha \times G^\alpha \longrightarrow Y^\alpha$ an effective right action of $G^\alpha$ on $Y^\alpha$. Then

$$\mathrm{Bun}_{\mathrm{Ar}}\boldsymbol{B} \stackrel{\text{def}}{=} \mathcal{A}t\ell^{\mathrm{Bun}}_{\substack{\mathrm{Ar}\\ \max}}\boldsymbol{B} \quad (28.2)$$

$$\mathrm{Bun}\boldsymbol{B} \stackrel{\text{def}}{=} \left(\mathrm{Bun}_{\mathrm{Ob}}\boldsymbol{B}, \mathrm{Bun}_{\mathrm{Ar}}\boldsymbol{B}, \overset{A}{\circ}\right) \quad (28.3)$$

Let $(A^i, B^i) \in \mathrm{Bun}_{\mathrm{Ob}}\boldsymbol{B}$, $i = 1, 2$. Then
$\boldsymbol{f} \stackrel{\text{def}}{=} (f_E\colon E^1 \longrightarrow E^2, f_X\colon X^1 \longrightarrow X^2, f_Y\colon Y^1 \longrightarrow Y^2, f_G\colon G^1 \longrightarrow G^2)$ is a bundle map from $(A^1, B^1)$ to $(A^1, B^1)$ iff it is a bundle-atlas morphism from $A^1$ to $A^2$. The identity morphism for $(A^i, B^i)$ is

$$\mathrm{Id}_{(A^i, B^i)} \stackrel{\text{def}}{=} \left((\mathrm{Id}_{E^i}, \mathrm{Id}_{X^i}, \mathrm{Id}_{Y^i}, \mathrm{Id}_{G^i}), (A^i, B^i), (A^i, B^i)\right) \quad (28.4)$$

## 28.2 Categories of fiber bundles

**Theorem 28.3** (Categories of fiber bundles)**.** *Let*
$\boldsymbol{B} \stackrel{\text{def}}{=} \left\{\boldsymbol{B}^\alpha \stackrel{\text{def}}{=} (E^\alpha, X^\alpha, Y^\alpha, G^\alpha, \pi^\alpha, \rho^\alpha) \,\middle|\, \alpha \prec A\right\}$, *where* $E^\alpha, X^\alpha, Y^\alpha$ *are topological spaces,* $G^\alpha$ *a topological group,* $\pi^\alpha\colon E^\alpha \twoheadrightarrow X^\alpha$ *surjective and* $\rho^\alpha\colon Y^\alpha \times G^\alpha \longrightarrow Y^\alpha$ *an effective* *right action of $G^\alpha$ on $Y^\alpha$. Then* Bun**B** *is a category and* $Id_{B^\alpha}$ *is the identity morphism for $B^\alpha$.*

*Proof.* The result follows directly from definition 28.1 (fiber bundles), definition 28.2 (Bundle maps) and lemma 26.10 ($\mathcal{A}t\ell^{\text{Bun}}$**B** is a category) on page 153
. □

# Part IX
# Future directions

If this paradigm proves useful, it can be extended to include a set of admissible functions on the model neighborhoods of the charts, possibly using the language of sheaves. That might be desirable for coordinate spaces more general than Fréchet spaces.

Further work is needed to determine whether there is a productive way to define a notion of classic m-atlas near morphisms.

Further work is needed to determine whether it is productive to define (semi-)strict m-atlas morphisms in terms of concrete categories over a category of model spaces rather than directly in terms of a model category.

The extension of paracompactness to model spaces is intended to be useful for partitions of unity on fiber bundles. Further work is needed to determine whether that is actually the case.

The definitions given here include some fairly strong conditions, e.g., AOC. Further work is needed to determine whether they should be relaxed for applications beyond manifolds and fiber bundles.

Further work is needed to determine whether the concept of category-based atlases[21] of model spaces has general utility.

If the concept of nearly commutative diagrams proves useful, further work is needed to determine whether a more general definition has utility.

Further work is needed to determine where it is most fruitful to require full subcategories and where just subcategories.

Further work is needed to determine conditions for mappings associated with atlas morphisms to be model functions.

This paper uses the language of category theory as an organizing principle, but defines various notions concretely with sets. It may be desirable to abstract away some of the details, in the spirit of, e.g., topoi.

---

[21] As opposed to pseudogroup based